\documentclass[aos,preprint]{imsart}

\RequirePackage[OT1]{fontenc}
\RequirePackage{amsthm,amsmath,amsfonts}
\RequirePackage[numbers]{natbib}
\usepackage{amsfonts}
\usepackage{tikz}
\usepackage{bm}
\usepackage{graphicx}
\usepackage{mathtools}
\usepackage{enumitem}
\usepackage[titlenumbered,ruled]{algorithm2e}
\usepackage{bigints}
\usepackage{mathrsfs}
\usepackage{fancybox}
\usepackage{color}
\usetikzlibrary{arrows}
\usetikzlibrary{positioning}
\usepackage[colorlinks,citecolor=blue,urlcolor=blue]{hyperref}
\usepackage{appendix}
\usepackage{etoc}

\startlocaldefs
\usepackage{enumitem, hyperref}
\makeatletter
\def\namedlabel#1#2{\begingroup
    #2%
    \def\@currentlabel{#2}%
    \phantomsection\label{#1}\endgroup
}
\makeatother
\newcommand{\sone}{\ref{sone} }
\newcommand{\stwo}{\ref{stwo} }
\newcommand{\bDelta}{\bm{\Delta}}
\newcommand{\cvec}{\bc(\nat,\, \bx)}
\newcommand{\cvect}{\bc(\nat,\, \bx)^{\top}}
\newcommand{\bmat}{\bA(\nat,\, \bx)}
\newcommand{\vvar}{v(\nat,\, \bx)}
\newcommand{\schur}{\xi(\nat,\, \bx)}
\newcommand{\dyads}{\mathfrak{D}_N}

\usepackage[mathscr]{eucal}
\usepackage{tikz}
\usepackage{dlfltxbcodetips}
\usepackage{bbm} 
\usepackage{mathtools}
\usepackage{subfigure} 
\usepackage{amsfonts}
\usepackage{amsmath} 
\usepackage{amssymb}
\usepackage{fancybox} 
\usepackage{graphicx}
\usepackage{bm} 
\usepackage{latexsym} 
\usepackage{centernot}

\usetikzlibrary{fit,positioning,shapes,backgrounds,calc}

\tikzset{node/.style={circle, fill=white, draw, minimum size=.85cm, text = black, inner sep=1pt, thick}}

\definecolor{myyellow}{RGB}{250, 218, 94}
\definecolor{myred}{RGB}{255, 75, 70}
\definecolor{myblue}{RGB}{0, 0, 250}
\definecolor{mypink}{RGB}{252, 137, 172}
\definecolor{mygrey}{RGB}{211,211,211} 
\definecolor{myorange}{RGB}{255, 0, 0}

\tikzset{circle split part fill/.style  args={#1,#2}{%
 alias=tmp@name, 
  postaction={%
    insert path={
     \pgfextra{%
     \pgfpointdiff{\pgfpointanchor{\pgf@node@name}{center}}%
                  {\pgfpointanchor{\pgf@node@name}{east}}%
     \pgfmathsetmacro\insiderad{\pgf@x}
      \fill[#1] (\pgf@node@name.base) ([xshift=-\pgflinewidth]\pgf@node@name.east) arc
                          (0:180:\insiderad-\pgflinewidth)--cycle;
      \fill[#2] (\pgf@node@name.base) ([xshift=\pgflinewidth]\pgf@node@name.west)  arc
                           (180:360:\insiderad-\pgflinewidth)--cycle;            
         }}}}}

\newcommand{\dpath}{\centernot{\longleftrightarrow}}

\newcommand{\Mle}{\widehat\bTheta}
\newcommand{\mle}{\widehat\btheta}
\newcommand{\mple}{\widetilde\btheta}
\newcommand{\Mple}{\widetilde\bTheta(\gamma_N)}
\newcommand{\truth}{\nat^\star}
\newcommand{\mI}{\mathcal{I}}

\newcommand{\pvec}{\bm{\xi}}
\newcommand{\pp}{U}

\newcommand{\su}{\sup\limits}
\newcommand{\logit}{\mbox{logit}}
\newcommand{\tv}{{\tiny\mbox{TV}}}
\newcommand{\lip}{{\tiny\mbox{Lip}}}

\newcommand{\mA}{\mathscr{A}}
\newcommand{\mB}{\mathscr{B}}

\newcommand{\mD}{\mathscr{D}}
\newcommand{\mE}{\mathscr{E}}

\newcommand{\mG}{\mathscr{G}}
\newcommand{\mmG}{\mathfrak{N}}

\newcommand{\mM}{\mathscr{M}}
\newcommand{\mN}{\mathscr{N}}

\newcommand{\mR}{\mathbb{R}}
\newcommand{\mS}{\mathscr{S}}
\newcommand{\mT}{\mathscr{T}}

\newcommand{\mUTA}{\widetilde{\mathscr{C}}}
\newcommand{\mV}{\mathscr{V}}

\newcommand{\mbB}{\mathbb{B}}
\newcommand{\mbC}{\mathbb{C}}

\newcommand{\mbE}{\mathbb{E}}

\newcommand{\mbG}{\mathbb{G}(\gamma_N)}
\newcommand{\mbH}{\mathbb{H}}

\newcommand{\mbM}{\mathbb{M}}

\newcommand{\mbP}{\mathbb{P}}
\newcommand{\mbQ}{\mathbb{Q}}
\newcommand{\mbR}{\mathbb{R}}

\newcommand{\mbT}{\mathbb{T}}

\newcommand{\mbV}{\mathbb{V}}

\newcommand{\mbX}{\mathbb{X}}

\newcommand{\bc}{\bm{c}}

\newcommand{\bg}{\bm{g}}

\newcommand{\bv}{\bm{v}}
\newcommand{\bw}{\bm{w}}
\newcommand{\bx}{\bm{x}}

\newcommand{\bz}{\bm{z}}

\newcommand{\bA}{\bm{A}}
\newcommand{\bB}{\bm{B}}

\newcommand{\bI}{\bm{I}}

\newcommand{\bX}{\bm{X}}

\newcommand{\bmu}{\bm\mu}

\newcommand{\bTheta}{\bm\Theta}

\newcommand{\btheta}{\boldsymbol{\theta}}

\newcommand{\bXi}{\mbox{\boldmath$\Xi$}}

\newcommand{\msim}{\mathop{\rm \sim}}
\newcommand{\one}{\mathbbm{1}}

\DeclareMathOperator*{\argmax}{arg\,max}

\newcommand\orth{\protect\mathpalette{\protect\independenT}{\perp}}

\def\independenT#1#2{\mathrel{\rlap{$#1#2$}\mkern2mu{#1#2}}}

\newcommand{\iid}{\msim\limits^{\mbox{\tiny iid}}}

\newcommand{\dsum}{\displaystyle\sum\limits}
\newcommand{\dint}{\displaystyle\int\limits}
\newcommand{\dprod}{\displaystyle\prod\limits}
\newcommand{\dd}{\mathop{\mbox{d}}\nolimits}

\DeclareMathOperator*{\bd}{\mathrm{bd}}

\DeclareMathOperator*{\interior}{\mbox{int}}
\newcommand{\var}{\mathbb{V}}
\newcommand{\cov}{\mathbb{C}}
\usepackage{lipsum}
\usepackage{relsize}
\usepackage{mathrsfs}
\usepackage{bigints}
\newcommand{\norm}[1]{|\!|#1|\!|}
\newcommand{\mnorm}[1]{|\!|\!|#1|\!|\!|}

\definecolor{me}{RGB}{255,36,0}

\newcommand{\btt}{\begin{box}}
\newcommand{\ett}{\end{box}}

\newcommand{\btheorem}{\begin{bclogo}[couleur={rgb:orange,0;yellow,0;white,1},arrondi=0.1,logo=\bcplume,ombre=true]{Theorem}}
\newcommand{\ettheorem}{\end{bclogo}}

\newcommand{\bst}{\begin{bclogo}[couleur={rgb:orange,1;yellow,1;white,0.5},arrondi=0.1,logo=\bcpanchant]}
\newcommand{\est}{\end{bclogo}}

\newcommand{\ba}{\begin{array}{llllllllll}}
\newcommand{\ea}{\end{array}}

\newcommand{\bea}{\begin{equation}\begin{array}{llllllllll}}
\newcommand{\eea}{\end{array}\end{equation}}

\newcommand{\be}{\begin{equation}\begin{array}{lllllllllllllllll}}
\newcommand{\beno}{\begin{equation}\begin{array}{lllllllllllll}\nonumber}
\newcommand{\ee}{\end{array}\end{equation}}

\newcommand{\bel}{\begin{equation}\begin{array}{lllllllllllll}\nonumber}
\newcommand{\eel}{\Box\end{array}\end{equation}}

\newcommand{\bi}{\begin{itemize}}
\newcommand{\ei}{\end{itemize}}

\newcommand{\ben}{\begin{enumerate}}
\newcommand{\een}{\end{enumerate}}

\newcommand{\bbq}{\begin{quote}\bf\em}
\newcommand{\ebq}{\end{quote}}

\renewcommand{\=}{&=&}

\newcommand{\lte}{&\leq&}
\newcommand{\gte}{&\geq&}

\newcommand{\hide}[1]{}
\newcommand{\ghost}[1]{}
\newcommand{\blue}{\textcolor{black}}
\newcommand{\alert}{\textcolor{black}}
\newcommand{\s}{\vspace{0.25cm}}

\newcounter{ex}

\newcounter{counterexample}

\newcounter{definition}
\newenvironment{definition}[1][]{\refstepcounter{definition}\par\smallskip\indent%
\textbf{Definition~\thedefinition #1.} \rmfamily}{\smallskip}

\newcounter{theorem}
\newenvironment{theorem}[1][]{\refstepcounter{theorem}\par\smallskip\indent%
\textbf{Theorem~\thetheorem #1}.\em \rmfamily}{}

\newcounter{proposition}
\newenvironment{proposition}[1][]{\refstepcounter{proposition}\par\smallskip\indent%
\textbf{Proposition~\theproposition #1}.\em \rmfamily}{}

\newcounter{ttproof}
\newcommand{\ttproof}{%
\addtocounter{ttproof}{1}%
{\medskip}\textsc{Proof of Theorem}{}
}

\newcounter{pproof}
\newcommand{\pproof}{%
\addtocounter{pproof}{1}%
{\medskip}\textsc{Proof of Proposition}{}
}

\newcounter{corollary}
\newenvironment{corollary}[1][]{\refstepcounter{corollary}\par\smallskip\indent%
\textbf{Corollary~\thecorollary #1}.\em \rmfamily}{}

\newcounter{ccproof}
\newcounter{ccsproof}
\newcommand{\ccsproof}{%
\addtocounter{ccsproof}{1}%
{\medskip}\textsc{Proof of Corollaries}{}
}

\newcounter{llproof}
\newcommand{\llproof}{%
\addtocounter{llproof}{1}%
{}\textsc{Proof of Lemma}{}
}

\newcounter{lemma}
\newenvironment{lemma}[1][]{\refstepcounter{lemma}\par\medskip\noindent%
\textbf{Lemma~\thelemma #1}.\em \rmfamily}{}

\newcounter{example}

\newcounter{com}
\newenvironment{com}[1][]{\refstepcounter{com}\par\medskip\indent%
\textit{Remark}. \rmfamily}{\medskip}

\newcounter{assumption}

\newcommand{\nat}{\btheta}
\newcommand{\Nat}{\bTheta}

\newcommand{\llambda}{\llambdamin}

\newcommand{\llambdamin}{\Lambda_{N}(\truth)}

\newcommand{\tllambdamin}{\widetilde\Lambda_{N}(\truth)}

\newcommand{\ip}{\mathop{\to}\limits^{\mbox{\footnotesize p}}}

\newcounter{daggerfootnote}

\numberwithin{equation}{section}
\endlocaldefs

\usepackage{amsthm}

\newcommand{\longtitle}{Pseudo-likelihood-based $M$-estimation of random graphs with dependent edges and parameter vectors of increasing dimension}
\newcommand{\shorttitle}{Pseudo-likelihood-based $M$-estimation}

\begin{document}

\begin{frontmatter}
\thankstext{T2}{MSC2010 subject classifications. Primary 05C80; secondary 62B05, 62F10, 91D30.}

\title{\longtitle}
\runtitle{\shorttitle}

\begin{aug}
\author{\fnms{Jonathan R.} \snm{Stewart}\ead[label=e2]{jrstewart@fsu.edu}}
\and
\author{\fnms{Michael} \snm{Schweinberger}\ead[label=e1]{mus47@psu.edu}}
\affiliation{Florida State University \and Penn State University}
\address{Jonathan R.\ Stewart\\
Department of Statistics\\
Florida State University\\
117 N Woodward Ave\\
Tallahassee, FL 32306-4330 \\ 
E-mail:\ jrstewart@fsu.edu\\
\mbox{}\\
Michael Schweinberger\\
Department of Statistics\\
Penn State University\\
326 Thomas Building\\
University Park, PA 16802\\
E-mail:\ mus47@psu.edu}
\end{aug}

\runauthor{Jonathan R.\ Stewart and Michael Schweinberger}

\begin{abstract}
An important question in statistical network analysis is how to estimate models of discrete and dependent network data with intractable likelihood functions, without sacrificing computational scalability and statistical guarantees. We demonstrate that scalable estimation of random graph models with dependent edges is possible, by establishing convergence rates of pseudo-likelihood-based $M$-estimators for discrete undirected graphical models with exponential parameterizations and parameter vectors of increasing dimension in single-observation scenarios. We highlight the impact of two complex phenomena on the convergence rate: phase transitions and model near-degeneracy. The main results have possible applications to discrete and dependent network, spatial, and temporal data. To showcase convergence rates, we introduce a novel class of generalized $\beta$-models with dependent edges and parameter vectors of increasing dimension, which leverage additional structure in the form of overlapping subpopulations to control dependence. We establish convergence rates of pseudo-likelihood-based $M$-estimators for generalized $\beta$-models in dense- and sparse-graph settings.
\end{abstract}


\end{frontmatter}

\section{Introduction}
\label{sec:intro}

Network data have garnered considerable attention in recent years,
driven by the growth of the internet and online social networks that can serve as echo chambers and facilitate polarization,
and applications in science, technology, and public health (e.g., pandemics).

During the past two decades,
substantial progress has been made on models of network data,
including $\beta$- and $p_1$-models \citep[e.g.,][]{HpLs81,ChDiSl11,YaXu13,RiPeFi13,MuMuSe18,ChKaLe19};
exchangeable random graph models \citep[e.g.,][]{BiCh09}; 
stochastic block models \citep[e.g.,][]{
BiCh09,
RoChYu11,
AmChBiLe13,Gaetal18};
latent space models \citep[e.g.,][]{HpRaHm01};
and exponential-family models of random graphs \citep[e.g.,][]{Ha03p,
Sc09b,ChDi11,Mu13,ScSt19}.
Other models are small-world networks and scale-free networks with power law degree distributions. 
That said, despite strides in modeling and inference,
fundamental questions arising from the statistical analysis of non-standard and dependent network data have remained unanswered.

\subsection{Three questions}
\label{motivation}

Since the dawn of statistical network analysis in the 1980s \citep[][]{HpLs81,FoSd86},
three questions have loomed large:
\bi
\item[I.] How can one construct models that allow the propensities of nodes to form edges and other subgraphs to vary across nodes?
\item[II.] How can one construct models that do justice to the fact that network data are dependent data?
\item[III.] How can one learn models from a single observation of a random graph with dependent edges and parameter vectors of increasing dimension,
regardless of whether the likelihood function is tractable?
\ei
We take steps to answer these questions by building on the statistical exponen\-tial-family platform \citep{Br86},
which has long served as a convenient mathematical platform for obtaining first answers to statistical questions involving discrete and dependent data and hosts Bernoulli random graphs,
$\beta$- and $p_1$-models \citep{HpLs81,ChDiSl11},
generalized linear models of random graphs,
and undirected graphical models of random graphs \citep{FoSd86,LaRiSa17}.
An alternative route,
not considered here, 
is provided by the Hoover-Aldous representation theorem via exchangeable random graphs \citep{BiCh09},
which can likewise induce dependence (as demonstrated by stochastic block and latent space models).

On the statistical exponential-family platform,
research has focused on $\beta$- and $p_1$-models,
which provide answers to the first question but assume that edges are independent;
and exponential-family random graph models,
which allow edges to be dependent and can capture observed heterogeneity via covariates,
but are less suited to capturing unobserved heterogeneity and often give rise to intractable likelihood functions.
An additional issue is that theoretical properties of statistical procedures
-- well-established in the literature on $\beta$- and $p_1$-models \citep[e.g.,][]{ChDiSl11,YaXu13,RiPeFi13,
YaLeZh11,
MuMuSe18,ChKaLe19} --
are scarce in the literature on exponential-family random graph models,
with two recent exceptions.
\citet{Mu13} considered models with functions of degrees as sufficient statistics,
which allow edges to be dependent,
but have two parameters and do not capture network features other than degrees.
\citet{ScSt19} considered models with dependent edges,
but constrained dependence to non-overlapping subpopulations of nodes.
While both works provide statistical guarantees,
these works focus on the second question rather than the first question.

We aim to provide tentative answers to all three questions,
leveraging the statistical exponential-family platform.

\subsection{Probabilistic framework}

On the modeling side,
we consider a flexible approach to specifying random graph models with complex dependence from simple building blocks.
\alert{We demonstrate the probabilistic framework by extending the $\beta$-model of \citet{ChDiSl11}
-- studied by \citet{RiPeFi13},
\citet{YaXu13},
\citet{MuMuSe18},
\citet{ChKaLe19},
and others --
to generalized $\beta$-models capturing dependence among edges along with heterogeneity in the propensities of nodes to form edges.
To control the dependence among edges,
generalized $\beta$-models leverage additional structure in the form of overlapping subpopulations.
The $\beta$-model and generalized $\beta$-models have in common that the number of parameters increases with the number of nodes.
Having said that,
the closest relative of generalized $\beta$-models with dependent edges is not the $\beta$-model with independent edges,
but are statistical exponential-family models for discrete and dependent random variables:
e.g.,
Ising models,
Markov random fields,
and undirected graphical models for discrete and dependent network, spatial, and temporal data \citep[e.g.,][]{GhMu20}.
}

\subsection{Computational scalability and statistical guarantees}

On the statistical side,
we demonstrate that computational scalability and statistical guarantees need not be sacrificed in order to estimate random graph models with dependent edges and parameter vectors of increasing dimension.

We do so by focusing on pseudo-likelihood-based $M$-estimators,
which possess convenient factorization properties and are more scalable than estimators based on intractable likelihood functions.
Despite computational advantages,
the properties of pseudo-likelihood-based $M$-estimators for random graphs with dependent edges and parameter vectors of increasing dimension are unknown.
In the related literature on Ising models and discrete Markov random fields in single-observation scenarios,
consistency 
of maximum pseudo-likelihood estimators 
has been established \citep{
Co92,Ch07a,BhMu18,GhMu20},
but those results are limited to a fixed number of parameters.

We demonstrate that scalable estimation of random graph models with dependent edges is possible,
by establishing convergence rates of pseudo-likelihood-based $M$-estimators for discrete undirected graphical models with exponential parameterizations and parameter vectors of increasing dimension in single-observation scenarios.
In contrast to high-dimensional graphical models,
we do not assume that independent replications are available.
The main results have possible applications to discrete and dependent network,
spatial,
and temporal data.
We highlight the impact of two complex phenomena on the convergence rate: phase transitions and model near-degeneracy.
To showcase convergence rates,
we establish convergence rates for generalized $\beta$-models with dependent edges and parameter vectors of increasing dimension in dense- and sparse-graph settings.

\subsection{Structure}

Section \ref{sec:models} introduces the probabilistic framework.
Section \ref{sec:stat_inf} establishes convergence rates for pseudo-likelihood-based $M$-estimators.
\hide{ 
Simulation results can be found in the supplement \citep{StSc20}.
}

\subsection{Notation}
\label{sec:notation}

Let $\mN \coloneqq \{1, \dots, N\}$ ($N \geq 3$) be a finite set of nodes and $\bX$ be a random graph defined on $\mN$ with sample space $\mbX \coloneqq \{0, 1\}^{\binom{N}{2}}$,
where $X_{i,j} = 1$ if nodes $i \in \mN$ and $j \in \mN$ are connected by an edge and $X_{i,j} = 0$ otherwise.
We focus on random graphs with undirected edges and without self-edges.
and the vector $\bm{0} \in \mR^d$ denotes the $d$-dimensional null vector in $\mR^d$ ($d \geq 1$). 
We denote the $\ell_1$-, $\ell_2$-, and $\ell_{\infty}$-norm of vectors in $\mR^d$ by $|\!|\,\cdot\,|\!|_{1}$, $|\!|\,\cdot\,|\!|_2$, and $|\!|\,\cdot\,|\!|_{\infty}$,
respectively.
For any matrix $\bA \in \mR^{d \times d}$,
let $|\!|\!|\bA|\!|\!|_1 \coloneqq \max_{1 \leq j \leq d}\, \sum_{i=1}^d |A_{i,j}|$,
$|\!|\!|\bA|\!|\!|_{\infty} \coloneqq \max_{1 \leq i \leq d}\, \sum_{j=1}^d |A_{i,j}|$,
and $|\!|\!|\bA|\!|\!|_2 \coloneqq \sup_{\bm{u} \in \mR^d:\; |\!|\bm{u}|\!|_2=1}\, |\!|\bA\, \bm{u}|\!|_2$.
For any vector norm $|\!|\cdot|\!|$,
the open ball in $\mR^d$ centered at $\bm{c} \in \mR^d$ with radius $\rho > 0$ is denoted by $\mB(\bm{c},\, \rho) \coloneqq \{\bm{a} \in \mR^d : |\!|\bm{a} - \bm{c}|\!| < \rho\}$:
e.g.,
the open hypercube in $\mR^d$ is $\mB_\infty(\bm{c},\, \rho) \coloneqq \{\bm{a} \in \mR^d : |\!|\bm{a} - \bm{c}|\!|_\infty < \rho\}$.
For any subset $\mS \subset \mR^d$, 
$\interior \mS$ and $\bd \mS$ denote the interior and boundary of $\mS$ in $\mR^d$,
respectively.
The total variation distance between two probability measures $\mbP_1$ and $\mbP_2$ defined on a common measurable space is denoted by $|\!|\mbP_1 - \mbP_2|\!|_{\tv}$.
Expectations,
variances,
and covariances are denoted by $\mbE$,
$\mbV$, 
and $\mbC$,
respectively.
For any finite set $\mS$,
the number of elements of $\mS$ is denoted by $|\mS|$.
The function $\one(\,\cdot\,)$ is an indicator function,
which is $1$ if its argument is true and is $0$ otherwise.
Uppercase letters $A, B, C, \dots$ denote finite constants.
We write $a(n) = O(b(n))$ if there exists a finite constant $C > 0$ such that $|a(n)\, /\, b(n)| \leq C$ for all large enough $n$,
and write $a(n) = o(b(n))$ if,
for all $\epsilon > 0$,
$|a(n)\, /\, b(n)|\, \,<\, \epsilon$ for all large enough $n$.

\section{Probabilistic framework}
\label{sec:models}

We consider a simple and flexible approach to specifying random graph models with complex dependence from simple building blocks.
Let $\{\mbP_{\nat}, \, \nat\, \in\, \Nat\}$ be a family of probability measures dominated by a $\sigma$-finite measure $\nu$,
with densities of the form 
\be
\label{eq:factorization}
\hide{
\dfrac{\dd \mbP_{\nat}}{\dd \nu}(\bx)
\=
}
f_{\nat}\left(\bx\right)
&\propto&
\dprod_{i<j}^N \varphi_{i,j}(x_{i,j}, \, \bx_{\mS_{i,j}};\, \nat),
&& \bx \in \mbX,
\ee
where $\varphi_{i,j}: \{0, 1\}^{|\mS_{i,j}|+1} \times \Nat \mapsto [0,\, \infty)$ is a function that specifies how edge variable $X_{i,j}$ depends on a subset of edge variables $\bX_{\mS_{i,j}}$. \alert{Here,
$\mS_{i,j}$ denotes a subset of unordered pairs of nodes $\{a, b\} \subset \mN$,
and $\bX_{\mS_{i,j}}$ denotes a set of indicators of edges between the unordered pairs of nodes in $\mS_{i,j}$.
} 
We allow the dimension $p \geq 1$ of parameter vector $\nat\, \in\, \Nat \subseteq \mR^p$ to increase as a function of the number of nodes $N$,
i.e.,
$p \to \infty$ as $N \to \infty$. 
A natural choice of reference measure $\nu$ is the counting measure.
\hide{
although other reference measures can be chosen,
e.g.,
reference measures that place more mass on sparse than dense graphs.
}

It is worth noting that the factorization of \eqref{eq:factorization} does not imply that edges are independent,
because each $\varphi_{i,j}$ can be a function of multiple edges and can hence induce dependence among edges.
That said,
the factorization of \eqref{eq:factorization} implies conditional independence properties \citep{Da79},
and the resulting models can be viewed as undirected graphical models of random graphs \citep{FoSd86,LaRiSa17}.
In contrast to the undirected graphical models of random graphs by \citet{FoSd86},
which allow edges to depend on many other edges and can give rise to undesirable behavior \citep[e.g., model near-degeneracy,][]{Ha03p,
Sc09b,
ChDi11},
we leverage additional structure to control dependence among edges.
The additional structure consists of a population with overlapping subpopulations and comes with two benefits.
First,
it facilitates the construction of novel models with non-trivial dependence.
Second,
it helps control the dependence among edges.
To demonstrate,
we introduce a novel class of generalized $\beta$-models with dependent edges in Sections \ref{model1}--\ref{model3}.

\subsection{Parameterizations}
\label{general.parameterizations}

It is convenient to parameterize the functions of edges $\varphi_{i,j}$ by using exponential parameterizations.
Exponential parameterizations are widely used in the literature on undirected graphical models:
see,
e.g.,
\citet{LaRiSa17}.
We therefore assume that 
\be
\label{parameterization}
\varphi_{i,j}(x_{i,j},\, \bx_{\mS_{i,j}};\, \nat)
&\coloneqq& a_{i,j}(x_{i,j},\, \bx_{\mS_{i,j}}) \, \exp(\langle\nat,\, s_{i,j}(x_{i,j},\, \bx_{\mS_{i,j}})\rangle),
\ee
where $a_{i,j}: \{0, 1\}^{|\mS_{i,j}|+1} \mapsto [0,\, \infty)$ is a function of $x_{i,j}$ and $\bx_{\mS_{i,j}}$,
which can be used to induce sparsity by penalizing edges,
and $\langle\nat,\, s_{i,j}(x_{i,j},\, \bx_{\mS_{i,j}})\rangle$ is the inner product of a vector of parameters $\nat \in \Nat \subseteq \mR^p$ and a vector of statistics $s_{i,j}: \{0, 1\}^{|\mS_{i,j}|+1} \mapsto \mR^p$ ($\{i, j\} \subset \mN$).
The probability density function \eqref{eq:factorization} with parameterization \eqref{parameterization} can be written in exponential-family form:
\be
\label{exponential.family}
f_{\nat}(\bx)
&=& a(\bx)\, \exp\left(\langle\nat,\, s(\bx)\rangle - \psi(\nat)\right),
&& \bx \in \mbX,
\ee
where $a: \mbX \mapsto [0,\, \infty)$ is given by
$a(\bx)
\coloneqq \prod_{i<j}^N\, a_{i,j}(x_{i,j},\, \bx_{\mS_{i,j}})$
and the vector of sufficient statistics $s: \mbX \mapsto \mR^p$ is given by
\be
\label{eq:sufficient}
s(\bx)
&\coloneqq& \dsum_{i<j}^N s_{i,j}(x_{i,j},\, \bx_{\mS_{i,j}}).
\ee
The function $\psi: \Nat \mapsto (0,\infty)$ ensures that $\int_{\mbX}\, f_{\nat}(\bx) \dd \nu(\bx) = 1$:
\beno
\psi(\nat)
&\coloneqq& \log \dint_{\mbX}\, a(\bx)\, \exp\left(\langle\nat,\, s(\bx)\rangle\right) \dd \nu(\bx),
&& \nat\, \in\, \Nat.
\ee
The parameter space is $\Nat \coloneqq \{\nat \in \mR^p \,:\, \psi(\nat) < \infty\} = \mR^p$,
because the family of densities is an exponential family of densities with respect to a $\sigma$-finite measure with a finite support \citep{Br86}.
To ensure that $\nat \in \Nat$ is identifiable,
we assume that the exponential family is minimal in the sense of \citet[][p.\ 2]{Br86}.
The assumption of a minimal exponential family involves no loss of generality,
because all non-minimal exponential families can be reduced to minimal exponential families 
\citep[][Theorem 1.9, p.\ 13]{Br86}.

We demonstrate the probabilistic framework by developing a novel class of generalized $\beta$-models with dependent edges and $p \geq N \to \infty$ parameters.

\subsection{Model 1: $\beta$-model with independent edges}
\label{model1}

To introduce generalized $\beta$-models with dependent edges,
we first review the $\beta$-model with independent edges \citep{ChDiSl11}.
The $\beta$-model assumes that edges between nodes $i\in\mN$ and $j\in\mN$ are independent Bernoulli$(\mu_{i,j})$ ($\mu_{i,j} \in (0, 1)$) random variables,
where
\beno
\log\dfrac{\mu_{i,j}}{1-\mu_{i,j}}
\= \theta_i + \theta_j,
&& \theta_i\in\mR,
&& \theta_j\in\mR.
\ee
The parameters $\theta_i$ and $\theta_j$ can be interpreted as the propensities of nodes $i$ and $j$ to form edges.
The $\beta$-model is a special case of the probabilistic framework introduced above,
corresponding to
\beno
\varphi_{i,j}(x_{i,j};\, \nat)
\= a_{i,j}(x_{i,j})\, \exp((\theta_i + \theta_j) \, x_{i,j}),
& \nat = (\theta_1, \dots, \theta_N) \in \mR^N,
\ee
where $a_{i,j}(x_{i,j})$ is $1$ if $x_{i,j} \in \{0, 1\}$ and is $0$ otherwise.
The $\beta$-model captures heterogeneity in the propensities of nodes to form edges,
but assumes that edges are independent.

\subsection{Model 2: generalized $\beta$-model with dependent edges}
\label{model2}

We introduce a generalization of the $\beta$-model,
which captures dependence among edges induced by brokerage in networks,
in addition to heterogeneity in the propensities of nodes to form edges.
Brokerage can influence economic and political outcomes of interest and has therefore been studied by economists, political scientists, and other network scientists since at least the 1980s.
An example of brokerage is given by faculty members of universities with appointments in both computer science and statistics,
who can facilitate collaborations between faculty members in computer science and faculty members in statistics and can hence facilitate interdisciplinary research.

To capture dependence among edges induced by brokerage in networks,
consider a finite population of nodes $\mN$ consisting of $K \geq 2$ known subpopulations $\mA_1, \ldots, \mA_K$,
which may overlap in the sense that the intersections of subpopulations are non-empty.
As a consequence,
nodes may belong to multiple subpopulations:
e.g.,
faculty members of universities may have appointments in multiple departments,
which implies that the faculties of departments overlap.
Subpopulation structure is inherent to many real-world networks,
in part because people tend to build communities, 
and in part because organizations tend to divide large bodies of people into small bodies of people (e.g., divisions, subdivisions). 
It is worth noting that we focus on known subpopulations that can overlap,
in contrast to the literature on stochastic block models \citep{BiCh09}.
In applications,
it is often possible to observe subpopulation structure:
e.g.,
the appointments of faculty members can be determined by scraping the websites of universities.

\begin{figure}[t]
\centering
\includegraphics[width = .31 \linewidth, keepaspectratio]{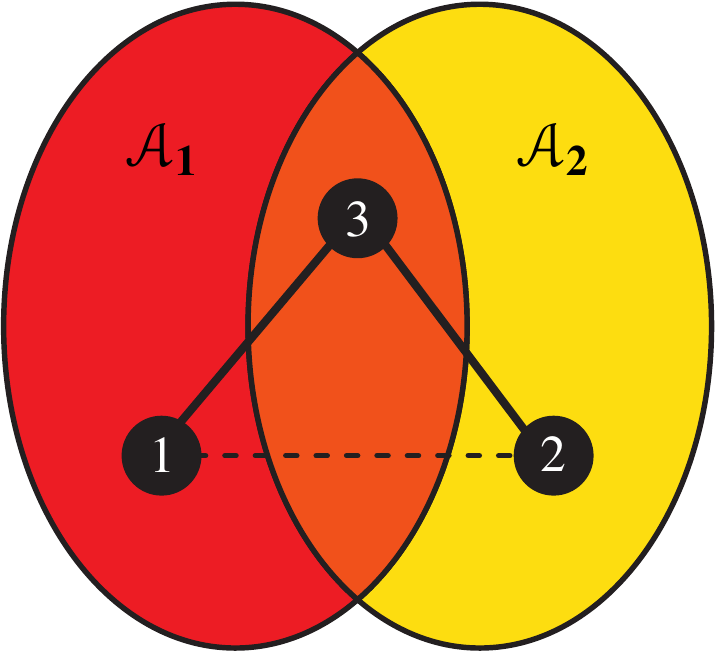}  
\caption{A graphical representation of the dependencies among edges induced by brokerage.
Consider two overlapping subpopulations $\mA_1$ and $\mA_2$.
The nodes $1 \in \mA_1 \setminus\, \mA_2$ and $2 \in \mA_2 \setminus\, \mA_1$ do not belong to the same subpopulation,
but the shared partner $3 \in \mA_1\, \cap\, \mA_2$ in the intersection of subpopulations $\mA_1$ and $\mA_2$ can facilitate an edge between nodes $1$ and $2$,
indicated by the dashed line between nodes $1$ and $2$.}
\label{fig:trans}
\end{figure}

Define,
for each node $i \in \mN$,
its neighborhood $\mN_i$ as the subset of all other nodes $j \in \mN \setminus \{i\}$ that share at least one subpopulation with node $i \in \mN$:
\beno
\mN_i
\;\coloneqq\; \left\{j \in \mN \setminus \{i\}:\, \exists\; k \in \{1, \dots, K\} \mbox{ such that } i \in \mA_k \mbox{ and } j \in \mA_k\right\}.
\ee
To capture dependence among edges induced by shared partners in the intersections of neighborhoods,
we consider functions of edges $\varphi_{i,j}$ of the form
\beno
\varphi_{i,j}(x_{i,j},\, \bx_{\mS_{i,j}};\, \nat)
\,\coloneqq\, a_{i,j}(x_{i,j})\, \exp\left((\theta_i + \theta_j)\, x_{i,j} + \theta_{N+1}\, b_{i,j}(x_{i,j},\, \bx_{\mS_{i,j}})\right),
\ee
where 
$\nat \coloneqq (\theta_1, \dots, \theta_{N+1}) \in \mR^{N+1}$, 
$\mS_{i,j} \subset \mN$ is the set of unordered pairs of nodes such that one node is an element of $\{i,\, j\}$ and the other node is an element of $\mN_i\, \cap\, \mN_j$,
$a_{i,j}(x_{i,j})$ is $1$ if $x_{i,j} \in \{0, 1\}$ and is $0$ otherwise,
and 
\be
\label{brokerage.statistic}
b_{i,j}(x_{i,j}, \bx_{\mS_{i,j}})
\coloneqq
\begin{cases}
0 & \mbox{if } \mN_i \cap \mN_j = \emptyset  
\\
x_{i,j}\, \one\left(\dsum_{h \in \mN_i\, \cap\, \mN_j} x_{i,h}\, x_{j,h} \geq 1\right)
& \mbox{if } \mN_i \cap \mN_j \neq \emptyset.
\end{cases}
\ee
Here,
$\one(\sum_{h\, \in\, \mN_i\, \cap\, \mN_j} x_{i,h}\, x_{j,h} \geq 1)$
is $1$ if nodes $i$ and $j$ have at least one shared partner in the 
intersection of neighborhoods $\mN_i$ and $\mN_j$, 
and is $0$ otherwise. 

\newpage 

\com {\em Generalized $\beta$-model captures brokerage in networks.}
The generalized $\beta$-model captures brokerage in networks,
along with heterogeneity in the propensities of nodes to form edges.
To demonstrate,
consider the two overlapping subpopulations $\mA_1$ and $\mA_2$ shown in Figure \ref{fig:trans}.
The nodes $1 \in \mA_1 \setminus\, \mA_2$ and $2 \in \mA_2 \setminus\, \mA_1$ do not belong to the same subpopulation,
but the shared partner $3 \in \mA_1\, \cap\, \mA_2$ in the intersection of subpopulations $\mA_1$ and $\mA_2$ can facilitate an edge between nodes $1$ and $2$,
provided $\theta_{N+1} > 0$.
In the language of network science,
nodes in the intersection $\mA_1\, \cap\, \mA_2$ of subpopulations $\mA_1$ and $\mA_2$ can act as brokers,
facilitating edges between nodes in $\mA_1 \setminus \mA_2$ and nodes in $\mA_2 \setminus \mA_1$.
In fact,
the generalized $\beta$-model can capture an excess in the expected number of brokered edges relative to the $\beta$-model,
in the sense that
\be
\label{inequality}
\begin{array}{ccccc}
\underbrace{\mbE_{\theta_1, \dots, \theta_{N}, \theta_{N+1}>0}\; b(\bX)}
&>&
\underbrace{\mbE_{\theta_1, \dots, \theta_{N}, \theta_{N+1}=0}\; b(\bX)},\s
\\
\mbox{\em generalized $\beta$-model}
&&
\mbox{\em $\beta$-model}
\end{array}
\ee
where $b(\bX) = \sum_{i<j}^N \, b_{i,j}(X_{i,j}, \, \bX_{\mS_{i,j}})$ and $\mbE_{\theta_1, \dots, \theta_{N}, \theta_{N+1}}\, b(\bX)$ is the expectation of $b(\bX)$ under $(\theta_1, \dots, \theta_{N}, \theta_{N+1}) \in \mR^{N+1}$.
In other words,
the generalized $\beta$-model with $\theta_{N+1}>0$ generates graphs that have,
on average,
more brokered edges than the $\beta$-model,
assuming that the propensities $\theta_1, \dots, \theta_{N}$ of nodes $1, \dots, N$ to form edges are the same under both models.
The inequality in \eqref{inequality} follows from the fact that the generalized $\beta$-model is an exponential-family model along with Corollary 2.5 of \citet[][p.\ 37]{Br86}.

\subsection{Model 3: sparse generalized $\beta$-models with dependent edges}
\label{model3}

\hide{
Many real-world networks are sparse in the sense that few of the possible edges are present,
which suggests that random graph models should place more mass on sparse graphs with few edges than dense graphs with many edges.
}
Sparse random graphs have been studied since the pioneering work of Erd\H{o}s and R\'enyi \citep[e.g.,][]{RiPeFi13, 
MuMuSe18, Mu13, ChKaLe19}.
To develop sparse versions of generalized $\beta$-models,
it makes sense to penalize edges between nodes $i \in \mN$ and $j \in \mN$ that are distant in the sense that $\mN_i\, \cap\, \mN_j = \emptyset$,
without penalizing edges between nodes that are close in the sense that $\mN_i\, \cap\, \mN_j \neq \emptyset$.
We therefore induce sparsity by considering Model 2 with
\beno
a_{i,j}(x_{i,j})
&\coloneqq&
\begin{cases}
N^{-\alpha\; x_{i,j}\, \one(\mN_i\, \cap\, \mN_j = \emptyset)} & \mbox{ if } x_{i,j} \in \{0, 1\}\s
\\
0 & \mbox{ otherwise,}
\end{cases}
\ee
where $\alpha \in (0,\, 1]$ is called the level of sparsity of the random graph.
\hide{
If $\alpha = 0$,
Model 3 reduces to Model 2,
whereas $\alpha > 0$ implies that Model 3 imposes a penalty on edges among nodes $i \in \mN$ and $j \in \mN$ with $\mN_i\, \cap\, \mN_j = \emptyset$,
which increases with the number of nodes $N$.
Thus,
Model 3 with $\alpha > 0$ places less mass on graphs with edges among nodes $i \in \mN$ and $j \in \mN$ satisfying $\mN_i\, \cap\, \mN_j = \emptyset$ than Model 3 with $\alpha = 0$.
}

To demonstrate that Model 3 encourages random graphs to be sparse,
we bound the expected degrees of nodes.
\begin{proposition}
\label{p:expected.degrees}
Consider Model 3 with $\nat \in \mR^{N+1}$ and $\alpha \in (0, 1]$. 
Then
\beno
\max\limits_{1 \leq i \leq N}\, \mbE_{\nat}\left(\dsum_{j \neq i}^N \, X_{i,j}\right)
\;\leq\; 2\; \exp(3\, \norm{\nat}_{\infty}) \; \left(\left(\max\limits_{1 \leq h \leq N} |\mN_h|\right)^2 + N^{1-\alpha}\right). 
\ee
\end{proposition}

Proposition \ref{p:expected.degrees} reveals that when the neighborhoods $\mN_h$ of nodes $h \in \mN$ are not too large,
the random graph is sparse in the sense that the expected degrees of all nodes are $o(N)$.
For example,
if $\max_{1 \leq h \leq N} |\mN_h|$ and $|\!|\nat|\!|_\infty$ are bounded above,
the expected degrees of nodes are $O(N^{1-\alpha})$.

\section{Statistical guarantees}
\label{sec:stat_inf}

We establish consistency results and convergence rates of maximum likelihood and pseudo-likelihood-based $M$-esti\-mators in Sections \ref{sec:consistency} and \ref{mple},
respectively.
We then present applications to $\beta$- and generalized $\beta$-models with dependent edges in Section \ref{sec:corollaries}.
To prepare the ground,
we first discuss how the dependence among edges and the smoothness of sufficient statistics can be quantified.
To ease the presentation,
we replace the double subscripts of edge variables by single subscripts and write $(X_m)_{1 \leq m \leq M}$ instead of $(X_{i,j})_{i<j:\, i\in\mN,\, j\in\mN}$,
where $M \coloneqq \binom{N}{2}$.
The data-generating parameter vector is denoted by $\truth \in \Nat = \mR^p$.

\subsection{Controlling dependence and smoothness} 
\label{sec:background}

To obtain consistency results and convergence rates based on a single observation of a random graph with dependent edges,
we need to control the dependence among edges along with the smoothness of the sufficient statistics of the model.

The dependence among edges can be controlled by bounding the total variation distance between conditional probability mass functions of edge variables,
quantifying how much the conditional probability mass functions of edge variables are affected by changes of other edge variables.
Define $\bX_{a:b} \coloneqq (X_a, \dots, X_b) \in \{0,\, 1\}^{b - a + 1}$,
where $a \leq b$ and $a,\, b \in \{1, \ldots, M\}$.
For each $i \in \{1, \dots, M\}$,
we denote the conditional probability mass function of subgraph $\bX_{i+1:M}$ given subgraph $(\bX_{1:i-1},\, X_i) = (\bx_{1:i-1},\, x_i)$ by $\mbP_{\truth,\bx_{1:i-1},x_i}$:
\beno
\mbP_{\truth,\bx_{1:i-1},x_i}(\bX_{i+1:M} = \bm{a})
\,\coloneqq\, \mbP_{\truth}(\bX_{i+1:M} = \bm{a} \mid (\bX_{1:i-1}, X_i) = (\bx_{1:i-1}, x_i)),
\ee
where $\bm{a} \in \{0, 1\}^{M-i}$. 
\alert{ 
We quantify the dependence among edges by bounding the total variation distance between the conditional probability mass functions $\mbP_{\truth,\bx_{1:i-1},0}$\, and $\mbP_{\truth,\bx_{1:i-1},1}$\, by using coupling methods \citep{Li02}:
\beno
\norm{\mbP_{\truth, \bx_{1:i-1}, 0} - \mbP_{\truth, \bx_{1:i-1}, 1}}_{\tv}
\;\leq\; \mbQ_{\truth,i,\bx_{1:i-1}}(\bX^\star_{i+1:M} \neq \bX^{\star\star}_{i+1:M}),
\hide{
\s
\\
=\; \mbQ_{\truth,i,\bx_{1:i-1}}\left(\bigcup\limits_{j = i+1}^M \{X_j^{\star} \neq X_j^{\star\star}\}\right)
\;\leq\; \dsum_{j = i + 1}^M \mbQ_{\truth,i,\bx_{1:i-1}}(X_j^{\star} \neq X_j^{\star\star}), 
}
\ee
where the pair of random vectors $(\bX^\star_{i+1:M}, \bX^{\star\star}_{i+1:M}) \in \{0, 1\}^{M-i} \times \{0, 1\}^{M-i}$ with joint probability mass function $\mbQ_{\truth,i,\bx_{1:i-1}}$ is a coupling of\, $\mbP_{\truth,\bx_{1:i-1},0}$\, and $\mbP_{\truth,\bx_{1:i-1},1}$ \citep{Li02}.
The coupling $\mbQ_{\truth,i,\bx_{1:i-1}}$ is constructed in Lemma \ref{prop:D_bound} in the supplement \citep{StSc20}.
Based on the coupling $\mbQ_{\truth,i,\bx_{1:i-1}}$,
we quantify the dependence among edges by the spectral norm $|\!|\!|\mD_N(\truth)|\!|\!|_2$ of the upper triangular $M \times M$ coupling matrix $\mD_N(\truth)$ with elements
\beno
\mD_{i,j}(\truth)
&\coloneqq&
\begin{cases}
0 & \mbox{ if } j < i
\\
1 & \mbox{ if } j = i
\\
\max\limits_{\bx_{1:i-1} \in \{0, 1\}^{i-1}} \mbQ_{\truth,i,\bx_{1:i-1}}(X_j^{\star} \neq X_j^{\star\star}) & \mbox{ if } j > i.
\end{cases}
\ee
While the definition of $\mD_N(\truth)$ depends on the ordering of edge variables,
it is possible to obtain bounds on the spectral norm $|\!|\!|\mD_N(\truth)|\!|\!|_2$ of $\mD_N(\truth)$ that hold for all orderings.
We describe in Section \ref{sec:coupling.matrix} how $|\!|\!|\mD_N(\truth)|\!|\!|_2$ can be bounded by using coupling methods from percolation theory \citep{BeMa94}.
}

\alert{To control the smoothness of the sufficient statistics of the model,
define
\beno
\Xi_{i,j}
&\coloneqq& \max\limits_{(\bx,\, \bx^\prime)\, \in\, \mbX \times \mbX:\;\;
x_k = x_k^\prime \mbox{\footnotesize\, for all } k \neq j}\;
|s_i(\bx) - s_i(\bx^\prime)|,
& j = 1, \dots, M,
\ee
where $s_1(\bx), \dots, s_p(\bx)$ 
are the coordinates of the sufficient statistic vector 
$s(\bx) \in \mR^p$ defined in \eqref{eq:sufficient}.
Let $\bm\Xi_i = (\Xi_{i,1}, \dots, \Xi_{i,M})$ and define
\beno
\Psi_N
&\coloneqq& \max\limits_{1 \leq i \leq p} \, |\!|\bm\Xi_i|\!|_2.
\ee
To exclude the trivial case where $\Psi_N = 0$,
we assume that there exists an integer $N_0 \geq 3$ such that $\Psi_N > 0$ for all $N > N_0$.}

\subsection{Maximum likelihood estimators} 
\label{sec:consistency}

Consider a single observation $\bx$ of a random graph $\bX$ with dependent edges.
Let $\ell(\nat;\, \bx) \coloneqq \log f_{\nat}(\bx)$ and
\beno
\Mle
&\coloneqq& \left\{\nat\, \in\, \mR^p:\; |\!|\nabla_{\nat}\; \ell(\nat;\, \bx)|\!|_\infty\, =\, 0\,\right\}.
\ee
We develop a novel approach to establishing consistency results and convergence rates of maximum likelihood estimators for discrete undirected graphical models with exponential parameterizations and parameter vectors of increasing dimension in single-observation scenarios.
These results serve as a stepping stone for establishing consistency results and convergence rates of pseudo-likelihood-based $M$-estimators in Section \ref{mple}.

Let $\mI(\nat) \coloneqq \nabla_{\nat}^2\, \psi(\nat) = \mbC_{\nat}\, s(\bX)
= -\mbE_{\nat}\, \nabla_{\nat}^2\,\, \ell(\nat;\, \bX)$ \citep[][Corollary 2.3, pp.\ 34--36]{Br86}.
Assume that there exists a constant $\epsilon^\star \in (0,\, \infty)$, 
independent of $N$ and $p$,
such that $\mI(\nat)$ is invertible for all $\nat \in \mB_{\infty}(\truth, \, \epsilon^\star)$.
Define
\be
\label{definition.phi}
\Lambda_N(\truth) &\coloneqq& \sup\limits_{\nat\, \in\, \mB_{\infty}(\truth, \, \epsilon^\star)} \, \mnorm{\mI(\nat)^{-1}}_{\infty}\s
\\
\Phi_N(\truth) &\coloneqq& \Lambda_N(\truth)\, \mnorm{\mD_N(\truth)}_2\, \Psi_N\, \sqrt{\log \, \max\{N, p\}}.
\ee

\vspace{-.75cm}

\alert{
\begin{theorem}
\label{theorem.mle}
Consider a single observation of a random graph with $N$ nodes and dependent edges.
Assume that $\truth \in \Nat = \mR^p$,
where $p \to \infty$ as $N \to \infty$ is allowed.
If\, $\Phi_N(\truth) \to 0$ as $N \to \infty$,
there exists an integer $N_0 \geq 3$ such that,
for all $N > N_0$,
the random set $\Mle$ is non-empty and its unique element $\mle$ satisfies
\beno
\norm{\mle - \truth}_{\infty}
&\leq& \sqrt{3 / 2}\;\, \Phi_N(\truth)
\ee
with probability at least $1 - 2\, / \max\{N,\, p\}^2$.
\end{theorem}
}

\s

\hide{ 
While Theorem \ref{theorem.mle} is stated in terms of random graphs,
Theorem \ref{theorem.mle} covers discrete undirected graphical models with exponential parameterizations and parameter vectors of increasing dimension in single-observation scenarios.
Theorem \ref{theorem.mle} suggests that the convergence rate of maximum likelihood estimators depends on the dimension $p$ of the parameter space $\Nat = \mR^p$ and
\bi
\item the inverse Fisher information matrix in a neighborhood of the data-generating parameter vector $\truth \in \mR^p$,
quantified by $\Lambda_N(\truth)$;\vspace{.1cm} 
\item the dependence induced by the model,
quantified by $|\!|\!|\mD_N(\truth)|\!|\!|_2$;\vspace{.1cm}
\item the sensitivity of sufficient statistics,
quantified by $\Psi_N$.
\ei}
We highlight the impact of two complex phenomena on the convergence rate:
phase transitions and model near-degeneracy \citep{Ha03p,
Sc09b,
ChDi11}.
It is known that some random graph models with dependent edges \citep[e.g., the ill-posed edge-and-triangle model,][]{Ha03p,
Sc09b,
ChDi11} exhibit phase transitions and model near-degeneracy.
To examine the impact of phase transitions and model near-degeneracy on the convergence rate,
consider a model with a parameter space $\Nat = \mR^p$ divided into two or more subsets (regimes) inducing very different distributions,
some of which may place almost all mass on a small subset of graphs (e.g., 
near-empty or near-complete graphs).
\vspace{.2cm}

{\em Phase transitions.}
On subsets of $\Nat$ where transitions between such regimes occur,
small changes of natural parameters $\nat$ can lead to large changes of mean-value parameters $\bmu(\nat) \coloneqq \nabla_{\nat}\, \psi(\nat) = \mbE_{\nat}\, s(\bX)$.
In such cases,
$\mI(\nat) \coloneqq \nabla_{\nat}^2\,\, \psi(\nat)$ can become ill-posed and non-invertible,
in which case Theorem \ref{theorem.mle} does not establish consistency. 
\vspace{.2cm}

{\em Model near-degeneracy.}
On subsets of $\Nat$ inducing near-degenerate distributions,
the variances of sufficient statistics (e.g., the number of edges) can be small,
so that the elements on the main diagonal of\, $\mI(\nat) = \mbC_{\nat}\, s(\bX)$ can be small for some or all $\nat\, \in\, \mB_{\infty}(\truth,\, \epsilon^\star)$. 
In such cases,
the convergence rate is reduced via $\llambdamin$.
In addition,
model near-degeneracy is sometimes associated with strong dependence and high sensitivity of sufficient statistics \citep{Sc09b},
depressing the convergence rate via $|\!|\!|\mD_N(\truth)|\!|\!|_2$ and $\Psi_N$.
An example is the ill-posed edge-and-triangle model \citep{Ha03p,
Sc09b,
ChDi11}.
We are interested in well-posed models that are amenable to scalable estimation with statistical guarantees.
Therefore,
the applications in Section \ref{sec:corollaries} focus on models that leverage additional structure to control all relevant quantities.
\vspace{.2cm}

To prove Theorem \ref{theorem.mle} along with Theorem \ref{thm:mple_consistency} in Section \ref{mple},
we first prove two lemmas.
Both lemmas are applicable to any homeomorphism\break 
$\bg: \mR^p \mapsto \mR^p$ (i.e., $\bg$ is bijective, 
continuous,
and its inverse $\bg^{-1}$ exists and is continuous),
any vector norm $|\!|\cdot|\!|$ with induced matrix norm $\mnorm{\cdot}$,\break
and any ball $\mB(\truth,\, \epsilon)\, \coloneqq\, \{\nat \in \mR^p:\, |\!|\nat - \truth|\!| < \epsilon\}$ in $\mR^p$.
\begin{lemma}
\label{lemma.inequality}
Let $\bg: \mR^p \mapsto \mR^p$ be any homeomorphism,
and let $|\!|\cdot|\!|$ be any vector norm with induced matrix norm $\mnorm{\cdot}$.
Consider any $\truth \in \mR^p$ and any $\epsilon \in (0,\, \infty)$, 
and define
\be
\label{delta.epsilon.definition}
\delta(\epsilon)
&\coloneqq& \inf\limits_{\nat\, \in\, {\footnotesize \mathrm{bd}}\, \mB(\truth,\, \epsilon)}\, |\!|\bg(\nat) - \bg(\truth)|\!|.
\ee
If $\bg(\nat)$ is continuously differentiable for all $\nat \in \mB(\truth,\, \epsilon)$ and the Jacobian matrix $\bm{J}(\nat) \coloneqq \nabla_{\nat}\; \bg(\nat)$ is invertible for all $\nat \in \mB(\truth,\, \epsilon)$,
then
\beno
\dfrac{\epsilon}{\sup_{\nat\, \in\, \mB(\truth,\, \epsilon)}\, \mnorm{\bm{J}(\nat)^{-1}}}
&\leq& \delta(\epsilon).
\ee
\end{lemma}

{\bf Proof of Lemma \ref{lemma.inequality}.}
Pick any $\truth \in \mR^p$ and any $\epsilon \in (0,\, \infty)$.
Then $\delta(\epsilon)$ can be expressed as
\beno
\delta(\epsilon)
&\coloneqq& \inf\limits_{\nat\, \in\, \bd \mB(\truth,\, \epsilon)} \, \norm{\bg(\nat) - \bg(\truth)}
\= \inf\limits_{\bg^\prime \in\, \bd \bg(\mB(\truth,\, \epsilon))} \, \norm{\bg^\prime - \bg(\truth)},
\ee
because $\bg(\bd \mB(\truth,\, \epsilon)) = \bd \bg(\mB(\truth,\, \epsilon))$ by Lemma \ref{lem:hom_bound_to_bound} in the supplement \citep{StSc20}.
By the invariance of domain theorem, 
$\epsilon > 0$ implies $\delta(\epsilon) > 0$.
Next,
pick any element $\nat^\prime\, \in\, \bd \mB(\truth,\, \epsilon)$ of the boundary $\bd \mB(\truth,\, \epsilon)$.
By the mean-value theorem \citep[Theorem 3.2.3, p.~69]{OrRh20} applied to $\bg^{-1}$ at $\bg(\nat^\prime)$ and $\bg(\truth)$,
\be
\label{key4}
\norm{\nat^\prime - \truth}
&\leq& \sup\limits_{0\, \leq\, \lambda\, \leq\, 1}\, \mnorm{\bm{J}(\lambda\, \nat^\prime + (1-\lambda)\, \truth)^{-1}}\; \norm{\bg(\nat^\prime) - \bg(\truth)}\s
\\
&\leq& \sup\limits_{\nat\, \in\, \mB(\truth,\, \epsilon)}\, \mnorm{\bm{J}(\nat)^{-1}}\; \norm{\bg(\nat^\prime) - \bg(\truth)},
\ee
where the second inequality follows from the fact that $\lambda\, \nat^\prime + (1-\lambda)\, \truth$ is an element of the convex set $\mB(\truth,\, \epsilon)$ and 
$\bm{J}(\nat)$
is invertible for all $\nat \in \mB(\truth,\, \epsilon)$ by assumption.
Since $\nat^\prime$ is an element of the boundary $\bd \mB(\truth,\, \epsilon)$,
the left-hand side of Equation \eqref{key4} is equal to $\norm{\nat^\prime - \truth} = \epsilon$,
which implies that
\be
\label{key5}
\epsilon
&\leq& \sup\limits_{\nat\, \in\, \mB(\truth,\, \epsilon)}\, 
\mnorm{\bm{J}(\nat)^{-1}}\;\, \norm{\bg(\nat^\prime) - \bg(\truth)}.
\ee
Taking the infimum over $\nat^\prime\, \in\, \bd \mB(\truth,\, \epsilon)$ on both sides of Equation \eqref{key5} and invoking the definition of $\delta(\epsilon)$ in Equation \eqref{delta.epsilon.definition} shows that
\beno
\epsilon
&\leq& \sup\limits_{\nat\, \in\, \mB(\truth,\, \epsilon)}\, \mnorm{\bm{J}(\nat)^{-1}}\;\, \delta(\epsilon).
\qed
\ee

\vspace{-.3cm}

\begin{lemma}
\label{lemma.boundary}
Let $\bg: \mR^p \mapsto \mR^p$ be any homeomorphism, 
and let $|\!|\cdot|\!|$ be any vector norm.
Consider any $\truth \in \mR^p$ and any $\epsilon \in (0,\, \infty)$,
and let $\delta(\epsilon)$ be defined by Equation \eqref{delta.epsilon.definition}.
Then $\mB(\bg(\truth),\, \delta(\epsilon))\, \subseteq\, \bg(\mB(\truth,\, \epsilon))$.
\end{lemma}

\vspace{.34cm}

{\bf Proof of Lemma \ref{lemma.boundary}.}
Lemma \ref{lem:hom_bound_to_bound} in the supplement \citep{StSc20} proves that
\be
\label{subset.of}
\bg(\bd \mB(\truth,\, \epsilon))
&=& \bd \bg(\mB(\truth,\, \epsilon)).
\ee
By the definition of $\delta(\epsilon)$ in Equation \eqref{delta.epsilon.definition} along with Equation \eqref{subset.of},
\beno
\delta(\epsilon) 
&\coloneqq& \inf\limits_{\nat\, \in\, {\footnotesize \bd}\, \mB(\truth,\, \epsilon)}\, |\!|\bg(\nat) - \bg(\truth)|\!|
&=& \inf\limits_{\bg^\prime \in\, {\footnotesize \bd}\, \bg(\mB(\truth,\, \epsilon))}\, |\!|\bg^\prime - \bg(\truth)|\!|.
\ee
As a result,
$\delta(\epsilon)$ is the shortest distance from $\bg(\truth)$ to the boundary of $\bg(\mB(\truth,\, \epsilon))$.
Thus,
all elements $\bg^{\prime\prime} \in \mR^p$ for which 
$|\!|\bg^{\prime\prime} - \bg(\truth)|\!| < \delta(\epsilon)$ are elements of the interior of $\bg(\mB(\truth,\, \epsilon))$,
which implies that
\beno
\mB(\bg(\truth),\, \delta(\epsilon))
&\subseteq& \bg(\mB(\truth,\, \epsilon)).
\qed
\ee

\vspace{.1cm}

{\bf Proof of Theorem \ref{theorem.mle}.}
Let $\nat \in \Nat$ be the natural parameter vector and $\bmu(\nat) \coloneqq \mbE_{\nat}\, s(\bX)$ be the mean-value parameter vector of an exponential family of the form \eqref{exponential.family}.
To ensure that $\nat \in \Nat$ is identifiable,
we assume that the exponential family is minimal in the sense of \citet[][p.\ 2]{Br86}.
The assumption of a minimal exponential family involves no loss of generality,
because all non-minimal exponential families can be reduced to minimal ones \citep[][Theorem 1.9, p.\ 13]{Br86}.
The map $\bmu: \Nat \mapsto \mbM$ from $\Nat \coloneqq \{\nat \in \mR^p:\, \psi(\nat) < \infty\} = \mR^p$ to $\mbM\, \coloneqq\, \bmu(\Nat)$ is a homeomorphism \citep[][Theorem 3.6, p.~74]{Br86} and is continuously differentiable \citep[Theorem 2.2, pp.~34--35]{Br86}.
By assumption,
there exists a constant $\epsilon^\star \in (0,\, \infty)$,
independent of $N$ and $p$,
such that $\mI(\nat) = \nabla_{\nat}^2\, \psi(\nat) = \nabla_{\nat}\, \bmu(\nat)$ is invertible for all $\nat\, \in\, \mB_{\infty}(\truth, \, \epsilon^\star)$.
We will show that we can focus on the subset $\mB_\infty(\truth,\, \epsilon^\star) \subset \Nat$ by proving in Equation \eqref{completion} that
\beno
\label{probprob}
\mbP(\mle \in \mB_{\infty}(\truth,\, \epsilon^\star))
&\geq& 1 - \dfrac{2}{\max\{N,\, p\}^2}\;\;
\mbox{ for all large enough $N$.}
\ee
Consider any $\epsilon \in (0, \, \epsilon^\star)$ and define
\beno
\delta(\epsilon)
&\coloneqq& \inf\limits_{\nat\, \in\, {\footnotesize \mathrm{bd}}\, \mB_\infty(\truth,\, \epsilon)}\, |\!|\bmu(\nat) - \bmu(\truth)|\!|_\infty.
\ee
We bound the probability of event $\mle \in \mB_{\infty}(\truth,\, \epsilon)$ in four steps.
\vspace{.1cm}

{\bf Step 1:}
The fact that $\bmu: \Nat \mapsto \mbM$ is a homeomorphism \citep[][Theorem 3.6, p.~74]{Br86} implies that
\beno
\label{fundamental.identity1}
\mbP(\mle \in \mB_{\infty}(\truth,\, \epsilon))
&=& \mbP(s(\bX) \in \bmu(\mB_{\infty}(\truth,\, \epsilon))),
\ee
noting that $\mle$ exists, is unique, and solves $\bmu(\mle) = s(\bX)$ in the event $s(\bX) \in \bmu(\mB_{\infty}(\truth,\, \epsilon)) \subset \mbM$ \citep[][Theorem 5.5, p.\ 148]{Br86}.
\vspace{.1cm}

{\bf Step 2:} 
Since $\bmu: \Nat \mapsto \mbM$ is a homeomorphism,
Lemma \ref{lemma.boundary} establishes
\beno
\mB_{\infty}(\bmu(\truth),\, \delta(\epsilon))
&\subseteq& \bmu(\mB_{\infty}(\truth,\, \epsilon)),
\ee
which implies that
\beno
\label{fundamental.identity2}
\mbP(s(\bX) \in \bmu(\mB_{\infty}(\truth,\, \epsilon)))
&\geq& \mbP(s(\bX) \in \mB_{\infty}(\bmu(\truth),\, \delta(\epsilon))).
\ee

{\bf Step 3:} 
Lemma \ref{theorem.mle.lemma} in the supplement \citep{StSc20} shows that
\beno
\label{fundamental.identity3}
\mbP(s(\bX) \in \mB_{\infty}(\bmu(\truth),\, \delta(\epsilon)))
&\geq& 1 - 2\, \exp\left( - \dfrac{2\; \delta(\epsilon)^2}{|\!|\!|\mD_N(\truth)|\!|\!|_2^2 \; \Psi_N^2} + \log p\right).
\ee

{\bf Step 4:} 
The homeomorphism $\bmu: \Nat \mapsto \mbM$ is continuously differentiable for all $\nat \in \mB(\truth,\, \epsilon)$ \citep[Theorem 2.2, pp.~34--35]{Br86} and $\mI(\nat) = \nabla_{\nat}^2\, \psi(\nat) = \nabla_{\nat}\, \bmu(\nat)$ is invertible for all $\nat\, \in\, \mB_\infty(\truth,\, \epsilon) \subset \mB_\infty(\truth,\, \epsilon^\star)$ by assumption.
Therefore,
Lemma \ref{lemma.inequality} can be invoked to establish 
\beno
\label{fundamental.identity4}
\delta(\epsilon)
&\geq& \dfrac{\epsilon}{\sup_{\nat\, \in\, \mB_\infty(\truth,\, \epsilon)}\, \mnorm{\mI(\nat)^{-1}}_\infty}.
\ee

{\bf Combining Steps 1--4} proves that,
for all $\epsilon \in (0,\, \epsilon^\star)$,
\beno
\label{fundamental.identity1}
\mbP(\mle \in \mB_{\infty}(\truth,\, \epsilon))
&=& \mbP(s(\bX) \in \bmu(\mB_{\infty}(\truth,\, \epsilon)))\s
\\
&\geq& \mbP(s(\bX) \in \mB_{\infty}(\bmu(\truth),\, \delta(\epsilon)))\s
\\
&\geq& 1 - 2\, \exp\left( - \dfrac{2\; \delta(\epsilon)^2}{|\!|\!|\mD_N(\truth)|\!|\!|_2^2 \; \Psi_N^2} + \log p\right)\s
\\
&\geq& 1 - 2 \, \exp\left( - \dfrac{2 \, \epsilon^2}{\Lambda_N(\truth)^2 \, \mnorm{\mD_N(\truth)}_2^2 \, \Psi_N^2} + \log \, p \right),
\ee
where $\Lambda_N(\truth) \coloneqq \sup_{\nat\, \in\, \mB_{\infty}(\truth, \, \epsilon^\star)} \, \mnorm{\mI(\nat)^{-1}}_{\infty}$.
Set
\beno
\epsilon
&\coloneqq& \sqrt{3 / 2}\; \Phi_N(\truth)
&=& \sqrt{3 / 2}\; \Lambda_N(\truth)\, \mnorm{\mD_N(\truth)}_2\, \Psi_N\, \sqrt{\log \, \max\{N, p\}}.
\ee
The chosen $\epsilon \coloneqq \sqrt{3 / 2}\; \Phi_N(\truth)$ satisfies $\epsilon < \epsilon^\star$ for all large enough $N$,
because the assumption $\Phi_N(\truth) \to 0$ as $N \to \infty$ implies that there exists an integer $N_0 \geq 3$ such that $\epsilon \coloneqq \sqrt{3 / 2}\;\, \Phi_N(\truth) < \epsilon^\star$ for all $N > N_0$.
Thus:
\be
\label{mle.convergence}
\mbP(\mle \in \mB_{\infty}(\truth, \sqrt{3 / 2}\, \Phi_N(\truth)))
\,\geq\, 1 - \dfrac{2}{\max\{N, p\}^2}
\mbox{ for all } N > N_0.
\ee
To complete the proof, 
we show that the focus on $\mB_{\infty}(\truth,\, \epsilon^\star) \subset \Nat$ is legitimate.
First,
the fact that $\sqrt{3 / 2}\; \Phi_N(\truth) < \epsilon^\star$ for all $N > N_0$ implies
\be
\label{subset.relation}
\mB_{\infty}(\truth, \sqrt{3 / 2}\; \Phi_N(\truth))
&\subset& \mB_{\infty}(\truth,\, \epsilon^\star)
&\mbox{for all}& N > N_0.
\ee
Second,
combining Equation \eqref{subset.relation} with Equation \eqref{mle.convergence} establishes
\be
\label{completion}
\mbP(\mle \in \mB_{\infty}(\truth,\, \epsilon^\star))
&\geq& \mbP(\mle\, \in\, \mB_{\infty}(\truth,\, \sqrt{3 / 2}\; \Phi_N(\truth)))\s
\\
&\geq& 1 - \dfrac{2}{\max\{N,\, p\}^2}\; \mbox{ for all }\; N > N_0.
\ee
To conclude,
for all $N > N_0$,
$\mle$ exists,
is unique,
and satisfies\break
$\norm{\mle - \truth}_{\infty} \leq \sqrt{3 / 2}\; \Phi_N(\truth)$ with probability at least $1 - 2 \,/ \max\{N, \, p\}^2$.$\qed$

\s

\vspace{-.25cm}

\subsection{Pseudo-likelihood-based $M$-estimators}
\label{mple}

Maximum likelihood estimators are unappealing on computational grounds,
because evaluating $\ell(\nat;\, \bx)$ requires evaluating the normalizing constant of $f_{\nat}(\bx)$.
The normalizing constant of $f_{\nat}(\bx)$ is a sum over $\exp(M \log 2)$ possible graphs and cannot be computed unless $M \coloneqq \binom{N}{2}$ is small or the model makes restrictive independence assumptions.
As a scalable alternative,
consider $M$-estimators 
\beno
\Mple
&\coloneqq& \left\{\nat\, \in\, \mR^p:\; |\!|\nabla_{\nat}\; \widetilde\ell(\nat;\, \bx)|\!|_\infty\; \leq\; \gamma_N\right\},
&& \gamma_N \in [0,\, \infty)
\ee
based on the pseudo-loglikelihood function 
\beno
\widetilde\ell(\nat;\, \bx)
&\coloneqq& \log \dprod_{i=1}^M f_{\nat}(x_i \mid \bx_{-i}),
\ee
where $f_{\nat}(x_i \mid \bx_{-i})$ is the conditional probability of $X_i = x_i$ given all other edge variables $\bX_{-i} = \bx_{-i}$ ($i = 1, \dots, M$).

To bound the statistical error of pseudo-likelihood-based $M$-estimators in single-observation scenarios with $p \to \infty$ parameters,
let $i \in \{1, \dots, M\}$ and $\mmG_{i}\; \subseteq\; \{1, \dots, M\}\, \setminus\, \{i\}$ be the smallest subset of $\{1, \dots, M\}\, \setminus\, \{i\}$ such that
\beno
X_i &\orth& \bX_{\{1, \dots, M\}\, \setminus\, (\{i\} \,\cup\, \mmG_i)} \mid \bX_{\mmG_i}.
\ee
Let $\epsilon^\star \in (0,\, \infty)$ be a constant,
independent of $N$ and $p$,
and 
\beno
\mbH
&\subseteq& 
\left\{\bx \in \mbX:\; -\nabla_{\nat}^2\; \widetilde\ell(\nat;\, \bx) \mbox{ is invertible for all $\nat \in \mB_{\infty}(\truth,\, \epsilon^\star)$}\right\},
\ee
and define 
\beno
\widetilde{\Lambda}_{N,\bx}(\truth)
&\coloneqq& \sup\limits_{\nat\, \in\, \mB_{\infty}(\truth,\, \epsilon^\star)}\, \mnorm{(-\nabla_{\nat}^2\; \widetilde\ell(\nat;\, \bx))^{-1}}_{\infty},\;\;\; \bx \in \mbH\s
\\
\widetilde{\Lambda}_{N}(\truth)
&\coloneqq& \max\limits_{\bx\, \in\, \mbH}\, \widetilde\Lambda_{N,\bx}(\truth)\s
\\
\widetilde{\Phi}_N(\truth)
&\coloneqq& \tllambdamin\,\, (1 + D_N) \,\, |\!|\!|\mD_N(\truth)|\!|\!|_2\,\, \Psi_N\; \sqrt{\log \max\{N,\, p\}},
\ee
where $D_N \coloneqq \max\{|\mathfrak{N}_1|, \dots, |\mathfrak{N}_M|\} \in \{0, \dots, M-1\}$.
It is worth noting that the set $\mbH$ can be a proper subset of the set of all $\bx \in \mbX$ for which $-\nabla_{\nat}^2\; \widetilde\ell(\nat;\, \bx)$ is invertible on $\mB_{\infty}(\truth,\, \epsilon^\star)$,
provided $\mbH$ is a high probability subset of $\mbX$ (see Theorem \ref{thm:mple_consistency} below).
The fact that $\mbH$ can be a proper subset of all $\bx \in \mbX$ for which the Hessian is invertible on $\mB_{\infty}(\truth,\, \epsilon^\star)$ has two advantages.
First,
characterizing the set of all $\bx \in \mbX$\, for which the Hessian is invertible may be hard,
but finding a sufficient condition for invertibility may well be possible.
Second, 
subsets of $\bx \in \mbX$ for which $-\nabla_{\nat}^2\; \widetilde\ell(\nat;\, \bx)$ is invertible but $\sup_{\nat\, \in\, \mB_{\infty}(\truth,\, \epsilon^\star)}\, \mnorm{(-\nabla_{\nat}^2\; \widetilde\ell(\nat;\, \bx))^{-1}}_{\infty}$ is large can be excluded from $\mbH$,
as long as those subsets are low probability subsets of\, $\mbX$.

\setcounter{theorem}{1}
\begin{theorem}
\label{thm:mple_consistency}
Consider a single observation of a random graph with $N$ nodes and dependent edges.
Let $\truth \in \Nat = \mR^p$,
where $p \to \infty$ as $N \to \infty$ is allowed.
Assume that there exists an integer $N_0 \geq 3$,
independent of $N$ and $p$,
such that, 
for all $N > N_0$,
$\sqrt{96}\; \widetilde\Phi_N(\truth) < \epsilon^\star$ and 
\be
\label{cond1}
\mbP\left(\bX \in \mbH\right)
&\geq& 1 - \dfrac{2}{\max\{N,\, p\}^2}.
\ee
Then,
for all $N > N_0$,
the random set $\Mple$ is non-empty and satisfies
\beno
\Mple 
&\subseteq& \mB_{\infty}(\truth,\, \sqrt{96}\;\, \widetilde\Phi_N(\truth))
&\subset& \mB_{\infty}(\truth,\, \epsilon^\star)
\ee 
with probability at least $1 - 4\, / \max\{N,\, p\}^2$,
provided
\beno
\gamma_N
&=& \sqrt{24}\;\, (1 + D_N) \; |\!|\!|\mD_N(\truth)|\!|\!|_2\; \Psi_N\; \sqrt{\log \max\{N,\, p\}}
&>& 0.
\ee
\end{theorem}

\hide{

{\bf Remark.}
The subset $\mbH \subseteq \mbX$ depends on the model and is,
in general,
a proper subset of $\mbX$.
As a case in point,
consider generalized $\beta$-models:
If $\bx \in \mbX$ is the empty graph without edges,
then changing an edge cannot change the number of brokered edges.
As a result,
$-\nabla_{\nat}^2\; \widetilde\ell(\nat;\, \bx)$ is not invertible.
That said,
under reasonable models (e.g., $\beta$- and generalized $\beta$-models),
the subset $\mbH$ is a high probability subset of $\mbX$.

}

\hide{ 

A proof of Theorem \ref{thm:mple_consistency} is provided in the supplement \citep{StSc20}.
While stated in terms of random graphs,
Theorems \ref{theorem.mle} and \ref{thm:mple_consistency} cover discrete undirected graphical models with exponential parameterizations and parameter vectors of increasing dimension in single-observation scenarios.
As a result,
Theorems \ref{theorem.mle} and \ref{thm:mple_consistency} have possible applications to discrete and dependent network,
spatial,
and temporal data. 

}

\alert{We first provide a simple application of Theorems \ref{theorem.mle} and \ref{thm:mple_consistency} in Section \ref{sec:dimension} and explore how fast the dimension $p$ of the parameter space $\Nat = \mR^p$ can grow as a function of the number of nodes $N$.
We then explain in Sections \ref{sec:coupling.matrix}, 
\ref{sec:determinant},
and \ref{sec:psi} how $\widetilde\Phi_N(\truth)$ can be bounded.
Applications to generalized $\beta$-models with dependent edges and $p \geq N \to \infty$ parameters are presented in Section \ref{sec:corollaries}.
These applications demonstrate that $\widetilde\Phi_N(\truth) \to 0$ as $N \to \infty$ provided $D_N$ does not grow too fast.
We conclude Section \ref{sec:corollaries} with a comparison with related statistical exponential-family models for discrete and dependent random variables in single-observation scenarios.}

\subsubsection{Example: growth of $p$ as a function of $N$}
\label{sec:dimension}

To showcase Theorems \ref{theorem.mle} and \ref{thm:mple_consistency} in one of the simplest possible scenarios and explore how fast the dimension $p$ of the parameter space $\Nat = \mR^p$ can grow as a function of $N$,
we consider inhomogeneous Bernoulli random graphs in the dense-graph regime.
Inhomogeneous Bernoulli random graphs assume that edge variables $X_i$ are independent Bernoulli$(\mu_i)$ random variables,
with edge probabilities $\mu_i \coloneqq \mbE\, X_i$ satisfying $0 < C_1 < \mu_i < C_2 < 1$ for finite constants $C_1$ and $C_2$,
independent of $N$.
Suppose that each edge variable $X_i$ belongs to one of $p \leq M$ distinct categories $k \in \{1, \ldots, p\}$ with edge probabilities $\pi_k \in (0,\, 1)$,
and that $\mu_i = \pi_k$ if edge variable $X_i$ is assigned to category $k$.
Inhomogeneous Bernoulli random graphs are statistical exponential families with natural parameters $\theta_k \coloneqq \logit(\pi_k)$ and sufficient statistics $s_k(\bx) \coloneqq \sum_{i=1}^M\, \one_k(i)\, x_i$ ($k = 1, \dots, p$),
where $\nat \coloneqq (\theta_1, \dots, \theta_p) \in \mR^p$,\,
$s(\bx) \coloneqq (s_1(\bx), \dots, s_p(\bx)) \in \mR^p$,\,
and $\one_k(i)$ is $1$ if edge variable $X_i$ is assigned to category $k$ and is $0$ otherwise.
Since edges are independent,
$\widetilde\ell(\nat;\, \bx) = \ell(\nat;\, \bx)$ for all $(\nat,\, \bx) \in \Nat \times \mbX$ and hence $\nabla_{\nat}^2\, \widetilde\ell(\nat;\, \bx) = \nabla_{\nat}^2\, \ell(\nat;\, \bx) = \mI(\nat) = \mbC_{\nat}\, s(\bX)$ for all $(\nat,\, \bx) \in \Nat \times \mbX$.
By the independence of edges,
$\mbC_{\nat}\, s(\bX)$ is a diagonal matrix,
so the variances $\var_{\nat}\, s_1(\bX)$, $\dots$, $\var_{\nat}\, s_p(\bX)$ are the eigenvalues of $\mbC_{\nat}\, s(\bX)$.
To bound them,
assume that there exist finite constants $0 < C_3 < C_4$ such that
\beno
\dfrac{C_3\, N^2}{p} 
\lte \dsum_{i=1}^M \one_k(i) 
\lte \dfrac{C_4\, N^2}{p},
&& k = 1, \dots, p,
\ee
that is,
the $p$ categories are balanced,
in the sense that the sizes of the $p$ categories are of the same order of magnitude.
Then there exists a finite constant $C_5 > 0$,
independent of $N$ and $p$,
such that 
\beno
\Lambda_N(\truth) 
&=& \widetilde\Lambda_N(\truth) 
\;\;=\;\; \sup\limits_{\nat\, \in\, \mB_{\infty}(\truth, \, \epsilon^\star)} \, \mnorm{\mI(\nat)^{-1}}_{\infty}\s
\\
&=& \sup\limits_{\nat\, \in\, \mB_{\infty}(\truth, \, \epsilon^\star)} \; \max\limits_{1 \leq k \leq p} \, \dfrac{1}{\var_{\nat} \, s_k(\bX)}
\;\;\leq\;\; \dfrac{C_5\,\, p}{N^2}.
\ee
By the independence of edges,
$D_N = 0$ and the coupling matrix $\mD_N(\truth)$ is the $M \times M$ identity matrix with spectral norm $|\!|\!|\mD_N(\truth)|\!|\!|_2 = 1$.
The quantity $\Psi_N \coloneqq \max_{1 \leq k \leq p} \, |\!|\bm\Xi_k|\!|_2$ can be bounded as follows.
First,
adding or deleting an edge in any category $k$ can change the number of edges $s_k(\bx)$ in category $k$ by $-1$ or $+1$,
while changes of edges in other categories leave $s_k(\bx)$ unchanged.
Second,
each category $k$ contains at most $C_4\, N^2 /\, p$\, edges,
so $|\!|\bm\Xi_k|\!|_2 \leq \sqrt{C_4\, N^2 /\, p}$ for all $k$ and hence $\Psi_N \leq \sqrt{C_4\, N^2 /\, p}$.
Thus,
there exists a finite constant $C > 0$, 
independent of $N$ and $p$,
such that 
\beno
\Phi_N(\truth)
&=& \widetilde\Phi_N(\truth)
&\leq& \dfrac{\sqrt{p \, \log \max\{N,\, p\}}}{C\, N}.
\ee
If $p = o(N^2 / \log N)$,
then $\Phi_N(\truth) = \widetilde\Phi_N(\truth) \to 0$ and the maximum likelihood and pseudo-likelihood estimators $\mle$ and $\mple$ are consistent estimators of $\truth \in \Nat = \mR^p$ by Theorem \ref{theorem.mle};
note that $\mle$ and $\mple$ are equal with probability $1$ when edges are independent.
Thus,
Theorems \ref{theorem.mle} and \ref{thm:mple_consistency} confirm the intuition that the number of parameters $p$ we can estimate (without assuming $\truth$ to be sparse) is less than $N^2$ (ignoring logarithmic terms).
These results dovetail with the 
results of \citet[][Theorem 2.1]{Po88} based on $n \to \infty$ independent 
observations from a statistical exponential family with $p \to \infty$ parameters,
which suggest that consistency results can be obtained as long as $p = o(n)$;
note that the number of independent observations under inhomogeneous Bernoulli random graphs is $n = \binom{N}{2}$.
While the example is limited to inhomogeneous Bernoulli random graphs,
we conjecture that $p$ can grow at most as fast when edges are dependent and the random graph is sparse,
because dependence increases $|\!|\!|\mD_N(\truth)|\!|\!|_2$ while sparsity decreases information $\mI(\nat)$ and hence increases $\llambda$.
\hide{
Remark:
Consider any epsilon in (0, 1), however small, and p < epsilon^2 * N^2 / log(N).
Case max(N,p) = N: trivial.
Case max(N,p) = p:
p * log(p)
< (epsilon^2 * N^2 / log(N)) * log(epsilon^2 * N^2 / log(N))
< (epsilon^2 * N^2 / log(N)) * log(N^2)
= 2 * (epsilon^2 * N^2 / log(N)) * log(N)
= 2 * epsilon^2 * N^2
so
sqrt(p * log(p)) / (C * N)
< sqrt(2 * epsilon^2 * N^2) / (C * N)
= (sqrt(2) * epsilon * N) / (C * N)
= sqrt(2) * epsilon / C
}

\subsubsection{Bounding the spectral norm of the coupling matrix}
\label{sec:coupling.matrix}

\alert{If edges are independent,
the spectral norm $|\!|\!|\mD_N(\truth)|\!|\!|_2$ of the coupling matrix $\mD_N(\truth)$ is $1$,
otherwise $|\!|\!|\mD_N(\truth)|\!|\!|_2$ needs to be bounded from above.
We transform the hard problem of bounding $|\!|\!|\mD_N(\truth)|\!|\!|_2$ into the more convenient problem of studying paths in a conditional independence graph $\mG$ that represents the conditional independence structure of a random graph \citep{FoSd86,LaRiSa17}.
A conditional independence graph $\mG$ consists of a set of vertices $\mV \coloneqq \{X_1, \dots, X_M\}$ and a set of undirected edges $\mE \subset \mV \times \mV$ indicating the absence of conditional independencies among edge variables $X_1, \dots, X_M$ \citep[see, e.g.,][]{FoSd86,LaRiSa17}.

We begin with the observation that the concentration results of \citet{Chetal07} leveraged in Theorems \ref{theorem.mle} and \ref{thm:mple_consistency} hold for all possible couplings $\mbQ_{\truth,i,\bx_{1:i-1}}$ of\, $\mbP_{\truth,\bx_{1:i-1},0}$\, and $\mbP_{\truth,\bx_{1:i-1},1}$,
and all possible couplings bound the total variation distance between $\mbP_{\truth,\bx_{1:i-1},0}$\, and $\mbP_{\truth,\bx_{1:i-1},1}$: 
\beno
\norm{\mbP_{\truth, \bx_{1:i-1}, 0} - \mbP_{\truth, \bx_{1:i-1}, 1}}_{\tv}
\;\leq\; \mbQ_{\truth,i,\bx_{1:i-1}}(\bX^\star_{i+1:M} \neq \bX^{\star\star}_{i+1:M})\s
\\
=\; \mbQ_{\truth,i,\bx_{1:i-1}}\left(\bigcup\limits_{j = i+1}^M \{X_j^{\star} \neq X_j^{\star\star}\}\right) 
\;\leq\; \dsum_{j = i + 1}^M \mD_{i,j}(\truth).
\ee
We can therefore replace optimal couplings (which provide the tightest bounds on the total variation distance) by suboptimal but more convenient couplings that facilitate bounds on the spectral norm $|\!|\!|\mD_N(\truth)|\!|\!|_2$ of $\mD_N(\truth)$.
To do so,
we adapt the coupling approach of \citet[][pp.~759--760]{BeMa94} from Markov random fields to random graphs.
The resulting coupling $\mbQ_{\truth,i,\bx_{1:i-1}}$ is described in Lemma \ref{prop:D_bound} in the supplement \citep{StSc20} and may not be optimal,
but it helps translate the hard problem of bounding the spectral norm $|\!|\!|\mD_N(\truth)|\!|\!|_2$ of $\mD_N(\truth)$ into the more convenient problem of studying paths in the conditional independence graph $\mG$.

We start with the inequality
\beno
\label{basic.inequality}
|\!|\!|\mD_N(\truth)|\!|\!|_2
&\leq& \sqrt{|\!|\!|\mD_N(\truth)|\!|\!|_1 \; |\!|\!|\mD_N(\truth)|\!|\!|_{\infty}}.
\ee
We then bound the quantities $|\!|\!|\mD_N(\truth)|\!|\!|_1$ and $|\!|\!|\mD_N(\truth)|\!|\!|_{\infty}$ by bounding the above-diagonal elements $\mD_{i,j}(\truth)$ of $\mD_N(\truth)$,
using paths of disagreement $i \centernot{\longleftrightarrow} j$ between vertices $X_i$ and $X_j$ in the conditional independence graph $\mG$;
note that the below-diagonal and diagonal elements of $\mD_N(\truth)$ are $0$ and $1$.
A path of disagreement $i \centernot{\longleftrightarrow} j$ between vertices $X_i$ and $X_j$ is a sequence of two or more distinct vertices $(X_i, \dots, X_j)$ in the conditional independence graph $\mG$ starting at vertex $X_i$ and ending at vertex $X_j$,
such that 
\bi
\item each subsequent pair of vertices $(X_{v}, X_{w})$ in the sequence is connected by an edge in the conditional independence graph $\mG$,
which indicates the absence of conditional independence of vertices $X_{v}$ and $X_{w}$;\vspace{.1cm}
\item the coupling $(\bX^\star_{i+1:M}, \bX^{\star\star}_{i+1:M}) \in \{0, 1\}^{M-i} \times \{0, 1\}^{M-i}$ with joint probability mass function $\mbQ_{\truth,i,\bx_{1:i-1}}$ 
disagrees at each vertex $X_v$ in the sequence,
in the sense that $X_{v}^{\star} \neq X_{v}^{\star\star}$.
\ei

\noindent
Theorem 1 of \citet[][]{BeMa94} implies that the coupling $\mbQ_{\truth,i,\bx_{1:i-1}}$ constructed in the supplement \citep{StSc20} satisfies
\be
\label{vandenbergandmaes}
\mbQ_{\truth,i,\bx_{1:i-1}}^{}(X_{j}^{\star} \neq X_{j}^{\star\star})
= \mbQ_{\truth,i,\bx_{1:i-1}}^{}(i \centernot{\longleftrightarrow} j)
\leq \mbB_{\bm\pi(\truth)}(i \centernot{\longleftrightarrow} j),
\ee
where $\mbB_{\bm\pi(\truth)}$ is a Bernoulli product measure on $\{0, 1\}^{M}$ with probability vector $\bm\pi(\truth) \in [0, 1]^{M}$.
The coordinates $\pi_v(\truth)$ of $\bm\pi(\truth)$ are given by
\beno
\pi_v(\truth)
\coloneqq
\begin{cases}
0 & \mbox{ if } v \leq i-1
\\
1 & \mbox{ if } v = i
\\
\max\limits_{(\bx_{-v},\, \bx_{-v}^\prime)\, \in\, \{0, 1\}^{M-1} \times \{0, 1\}^{M-1}}\, \pi_{v,\, \bx_{-v},\, \bx_{-v}^\prime}(\truth)
& \mbox{ if } v \geq i+1,
\end{cases}
\ee
where
\beno
\pi_{v,\, \bx_{-v},\, \bx_{-v}^\prime}(\truth)
&\coloneqq&
|\!|\mbP_{\truth}(\,\cdot \mid \bX_{-v} = \bx_{-v}) - \mbP_{\truth}(\,\cdot \mid \bX_{-v} = \bx_{-v}^\prime)|\!|_{\tv}
\ee
is the total variation distance between the conditional probability mass functions of vertex $X_v$ given $\bX_{-v} = \bx_{-v}$ and $\bX_{-v} = \bx_{-v}^\prime$.
Leveraging \eqref{vandenbergandmaes},
we can bound the above-diagonal elements $\mD_{i,j}(\truth)$ of $\mD_N(\truth)$ as follows:
\beno
\label{boundingdij}
\mD_{i,j}(\truth)
&\coloneqq& \max\limits_{\bx_{1:i-1}\, \in\, \{0,1\}^{i-1}}\, \mbQ_{\truth,i,\bx_{1:i-1}}^{}(X_{j}^{\star} \neq X_{j}^{\star\star})
&\leq& \mbB_{\bm\pi(\truth)}(i \centernot{\longleftrightarrow} j).
\ee
In other words,
the spectral norm $|\!|\!|\mD_N(\truth)|\!|\!|_2$ of $\mD_N(\truth)$ can be bounded by using paths of disagreement $i \centernot{\longleftrightarrow} j$ in the conditional independence graph $\mG$,
and by bounding the probabilities of those paths by Bernoulli product measures.
Specific bounds depend on the data-generating model with parameter vector $\truth \in \Nat = \mR^p$.
Applications to generalized $\beta$-models with dependent edges can be found in the supplement \citep{StSc20}.
}

\subsubsection{Bounding the $\ell_{\infty}$-norm of inverse negative (expected) Hessians} 
\label{sec:determinant}

\alert{\blue{
To establish convergence rates, 
$\llambdamin$ and $\tllambdamin$ need to be bounded,
which amounts to bounding the suprema of\, $\mnorm{(-\mbE\, \nabla_{\nat}^2\, \ell(\nat; \bX))^{-1}}_{\infty}$ and\break 
$\mnorm{(-\nabla_{\nat}^2\; \widetilde\ell(\nat;\, \bx))^{-1}}_{\infty}$ for all $\bx \in \mbH$ on $\mB_{\infty}(\truth,\, \epsilon^\star) \subset \Nat = \mR^p$. 
}

\blue{
As a case in point,
consider the $\ell_{\infty}$-induced matrix norm of $(-\nabla_{\nat}^2\; \widetilde\ell(\nat;\, \bx))^{-1}$,
where $-\nabla_{\nat}^2\; \widetilde\ell(\nat;\, \bx)$ is invertible for all $(\nat,\, \bx) \in \mB_\infty(\truth,\, \epsilon^\star) \times \mbH$.
Then
\be
\label{standard.inequalities}
\mnorm{(-\nabla_{\nat}^2\; \widetilde\ell(\nat;\, \bx))^{-1}}_{\infty} 
\lte \sqrt{p}\; \mnorm{(-\nabla_{\nat}^2\; \widetilde\ell(\nat;\, \bx))^{-1}}_{2}\s
\\
\= \dfrac{\sqrt{p}}{\raisebox{-4pt}{$\lambda_{\min}(-\nabla_{\nat}^2\; \widetilde\ell(\nat;\, \bx))$}},
\ee
where $\lambda_{\min}(-\nabla_{\nat}^2\; \widetilde\ell(\nat;\, \bx)) > 0$ is the smallest eigenvalue of\, $-\nabla_{\nat}^2\; \widetilde\ell(\nat;\, \bx)$.
}
That said,
bounds of $\tllambdamin$ based on \eqref{standard.inequalities} may be loose when $p \to \infty$ as $N \to \infty$,
as is the case with generalized $\beta$-models with dependent edges.

To establish bounds on $\tllambdamin$ in scenarios with $p \to \infty$ parameters,
we leverage the fact that generalized $\beta$-models with dependent edges and $p = N + 1 \to \infty$ parameters include the $\beta$-model with independent edges and $p = N \to \infty$ parameters as a special case,
along with the fact that the negative expected Hessian of the $\beta$-model is diagonally dominant in the sense of \citet{HiWi15}. 
By leveraging these properties,
Lemma \ref{lemma.max} in the supplement \citep{StSc20} establishes the bound $\tllambdamin \leq C\, D_N^9\, /\, N^{1-(\alpha+\vartheta)}$,
where the constants $C \in (0,\, \infty)$,\,
$\alpha \in [0,\, 1/2)$, 
and $\vartheta \in [0,\, 1/2 - \alpha)$ are independent of $N$ and $p$,
while $D_N$ satisfies $D_N = O(\log N)$.
\hide{
The resulting bound $\tllambdamin \leq C\, D_N^9\, /\, N^{1-(\alpha+\vartheta)}$ reveals a trade-off between the sparsity of the random graph controlled by $\alpha \in [0,\, 1/2)$ and the growth of $\norm{\truth}_\infty$ controlled by $\vartheta \in [0,\, 1/2 - \alpha)$,
as discussed in Section \ref{sec:corollaries}.
}
} 

\subsubsection{Bounding the smoothness of the sufficient statistics}
\label{sec:psi}

\alert{The quantity $\Psi_N \coloneqq \max_{1 \leq i \leq p} \, |\!|\bm\Xi_i|\!|_2$ can be bounded by bounding the coordinates $\Xi_{i,j}$ of $\bm\Xi_i$.
Bounding $\Xi_{i,j}$ amounts to bounding changes of sufficient statistics.
}

\subsection{Applications}
\label{sec:corollaries}

\alert{We present applications of pseudo-likelihood-based $M$-estimators to $\beta$- and generalized $\beta$-models with dependent edges and $p \geq N \to \infty$ parameters,
in dense- and sparse-graph settings.
Throughout,
we assume that 
the data-generating parameter vector $\truth \in \Nat = \mR^p$ satisfies
\be
\label{condition.theta0}
\norm{\truth}_{\infty} 
&\leq& \dfrac{L + \vartheta \, \log N}{12\, (3 + D_N)} - \epsilon^\star,
\ee
where $L \in [0,\, \infty)$,\,
$\vartheta \in [0,\, \infty)$,\,
and $\epsilon^\star \in (0,\, \infty)$ are constants,
independent of $N$ and $p$.
The constant $\epsilon^\star \in (0,\, \infty)$ is identical to the constant $\epsilon^\star$ in the definition of $\tllambdamin$ and Theorem \ref{thm:mple_consistency}.
The quantity $D_N \coloneqq \max\{|\mathfrak{N}_1|, \dots, |\mathfrak{N}_M|\}$ is identical to the quantity $D_N$ in the definition of $\widetilde\Phi_N(\truth)$ and satisfies $D_N = 0$ under Model 1, 
but can increase as a function of $N$ under Models 2 and 3.
To ensure that $\norm{\truth}_{\infty} > 0$, 
we assume that $D_N$ satisfies
\beno
1  
\lte D_N
&<& \dfrac{L + \vartheta \, \log N}{12\; \epsilon^\star} - 3
\ee 
under Models 2 and 3.
}
\indent

\alert{
We start with the $\beta$-model \citep{ChDiSl11},
because its theoretical properties have been studied and it is therefore a convenient benchmark.
}

\vspace{-.25cm}

\alert{\begin{corollary}{\bf {} $\beta$-model.}
\label{cor:model_1}
Consider Model 1 with $\truth \in \mbR^{N}$ satisfying \eqref{condition.theta0} with $\vartheta \in [0,\, 3)$.
Then there exist finite constants $C > 0$ and $N_0 \geq 3$,
independent of $N$ and $p$, 
such that,
for all $N > N_0$,\,
\beno
\Phi_N(\truth)
&=& \widetilde\Phi_N(\truth)
&\leq& C\, \sqrt{\dfrac{\log N}{N^{1 - \vartheta / 3}}}.
\ee
\end{corollary}
}

\vspace{-.25cm}

\alert{
Corollary \ref{cor:model_1} shows that the convergence rate is highest when $|\!|\truth|\!|_\infty$ is bounded above ($\vartheta=0$).
Condition \eqref{condition.theta0} is the weakest known condition on $|\!|\truth|\!|_\infty$:
\citet[][Theorem 1.3]{ChDiSl11} report a non-asymptotic error bound of the form $|\!|\mle - \truth|\!|_\infty \leq C\, \sqrt{\log N\, /\, N}$ assuming that $|\!|\truth|\!|_\infty$ is bounded above ($\vartheta=0$),
while \citet[][Theorem 1]{YaXu13} report asymptotic consistency and normality results assuming that $|\!|\truth|\!|_\infty = o(\log\log N)$.
By contrast,
condition \eqref{condition.theta0} assumes that $|\!|\truth|\!|_\infty\, <\, (1 / 12)\, \log N$ ($\truth \in \mR^N$,\, $L = 0$, $\vartheta < 3$, $D_N = 0$),
which dovetails with the condition 
$|\!|\truth|\!|_\infty < (1 / 24)\, \log N$ ($\truth \in \mR^{2\, N - 1}$) of \citet[][Theorem 1]{YaLeZh11} based on the $p_1$-model for directed random graphs;
note that the $\beta$-model for undirected random graphs can be viewed as a relative of the $p_1$-model for directed random graphs,
because both models are statistical exponential-family models of degree sequences.
These results, 
along with the results on the dimension $p$ of the parameter space $\Nat = \mR^p$ in Section \ref{sec:dimension}, 
demonstrate that Theorems \ref{theorem.mle} and \ref{thm:mple_consistency} recover the sharpest known results for random graphs with independent edges and $p \to \infty$ parameters,
suggesting that the generality of Theorems \ref{theorem.mle} and \ref{thm:mple_consistency} comes at a low cost. It is worth noting that it is unknown whether the constants mentioned above are sharp.
While it would be of interest to investigate whether these constants are sharp,
constants do not affect convergence rates and the question of whether these constants are sharp is therefore not pertinent to the main results of the paper.

}

\alert{To demonstrate that Theorem \ref{thm:mple_consistency} covers random graph models with non-trivial dependence,
we turn to generalized $\beta$-models with dependent edges.
Throughout,
we assume that the size $|\mA_k|$ of each subpopulation $\mA_k$ satisfies $|\mA_k| \geq 3$ ($k = 1, \dots, K$).
We start with non-overlapping subpopulations in Corollary \ref{cor:model_2} and deal with overlapping subpopulations in Corollary \ref{cor:model_3}.}
\alert{\begin{corollary}{\bf {} Generalized $\beta$-models with dependent edges.}
\label{cor:model_2}
Consider Models 2 and 3 with non-overlapping subpopulations,
level of sparsity $\alpha \in [0,\, 1/2)$, 
and $\truth \in \mbR^{N+1}$ satisfying \eqref{condition.theta0} with $\vartheta \in [0, \, 1/2 - \alpha)$. 
Then 
\beno
4 \, D_N \, \sqrt{N \log N}
\lte \gamma_N
\lte 28 \, D_N^{5} \, \sqrt{N \log N},
\ee
and there exist finite constants $C > 0$ and $N_0 \geq 3$,
independent of $N$ and $p$,
such that,
for all $N > N_0$,
\textcolor{black}{\eqref{cond1} holds provided $D_N = O(\log N)$,
and} 
\beno
\widetilde\Phi_N(\truth)
&\leq& C\, D_N^{14}\; \sqrt{\dfrac{\log N}{N^{1 - 2\, (\alpha + \vartheta)}}}.
\ee
\end{corollary}
} 

\vspace{-.25cm}

\alert{
\noindent
Corollary \ref{cor:model_2} shows that the convergence rate of pseudo-likelihood-based $M$-estimators under generalized $\beta$-models with dependent edges and non-overlapping subpopulations resembles the convergence rate under the $\beta$-model with independent edges when the random graph is dense ($\alpha = 0$) and $\norm{\truth}_\infty$ is bounded above ($\vartheta=0$),
ignoring logarithmic terms;
note that $D_N$ needs to satisfy $D_N = O(\log N)$ to ensure $\norm{\truth}_\infty > 0$.
In addition,
Corollary \ref{cor:model_2} reveals a trade-off between the sparsity of the random graph controlled by $\alpha \in [0,\, 1/2)$ and the growth of $\norm{\truth}_\infty$ controlled by $\vartheta \in [0,\, 1/2 - \alpha)$.
}

\alert{We turn to overlapping subpopulations.
To bound $|\!|\!|\mD_N(\truth)|\!|\!|_2$ in scenarios with overlapping subpopulations,
we need to control the amount of overlap of subpopulations,
because the dependence among edges can propagate through overlapping subpopulations.
To do so,
we introduce a subpopulation graph $\mG_{\mA}$ with a set of vertices $\mV_{\mA} \coloneqq \{\mA_1, \ldots, \mA_K\}$,
where a pair of distinct subpopulations $\mA_k$ and $\mA_l$ is connected by an edge if 
$\mA_k\, \cap\, \mA_l\, \neq\, \emptyset$. 
}
Denote by $d_{\mG_{\mA}}: \mV_{\mA} \times \mV_{\mA} \mapsto \{0, 1, \ldots\} \, \cup\, \{\infty\}$ the length of the shortest path between pairs of subpopulations in $\mG_{\mA}$,
called the graph distance;
note that $d_{\mG_{\mA}}(\mA_k,\, \mA_k) \coloneqq 0$ and $d_{\mG_{\mA}}(\mA_k,\, \mA_l) \coloneqq \infty$ if there is no path of finite length between two distinct subpopulations $\mA_k$ and $\mA_l$.
Let $\mV_{\mA_k, l}$ be the subset of subpopulations at graph distance $l$ from a given subpopulation $\mA_k$:
\beno
\mV_{\mA_k, l} 
&\coloneqq& \left\{\mA^\star \in \{\mA_1, \ldots, \mA_K\} \setminus \{\mA_k\}:\; d_{\mG_{\mA}}(\mA_k,\, \mA^\star) = l \right\}.
\ee

\alert{
{\bf Assumption A.}
{\em
Define $U \coloneqq  1 \,/\, (1 + \exp(-L))$, 
where $L \in [0,\, \infty)$ is identical to the constant $L$ in \eqref{condition.theta0} and is independent of $N$ and $p$. 
Assume that $D_N \in [1,\, \infty)$ and that there exist finite constants $\omega_1 \in [0,\, \infty)$ and
\beno
0
&\leq& \omega_2 
&\leq& \min\left\{ \omega_1, \; \dfrac{1}{(\omega_1 + 1) \, |\log(1 - U)|} \right\},
\ee
independent of $N$ and $p$, 
such that  
\beno
\max\limits_{1 \leq k \leq K}\, |\mV_{\mA_k,l}|
\lte \omega_1 + \dfrac{\omega_{2}}{2 \, D_N^3} \, \log l, 
&& l \in \{1, \dots, K-1\}. 
\ee
}
\noindent
}
Assumption A covers tree- and non-tree subpopulation graphs in which, 
for each subpopulation, 
the number of subpopulations at graph distance $l$ is either constant or grows slowly as a function of $l$ (depending on $D_N$).

\vspace{-.4cm}

\alert{\begin{corollary}{\bf {} Generalized $\beta$-models with dependent edges.}
\label{cor:model_3}
Consider Models 2 and 3 with overlapping subpopulations and level of sparsity $\alpha \in [0,\, 1/2)$.
Assume that $\truth \in \mbR^{N+1}$ satisfies \eqref{condition.theta0} with $\vartheta = 0$ and that Assumption A is satisfied. 
Then there exist finite constants $A > 0$, $B > 0$, $C > 0$, and $N_0 \geq 3$, 
independent of $N$ and $p$, 
such that,  
for all $N > N_0$, 
\beno
4 \, D_N \, \sqrt{N \log N}
\lte \gamma_N
\lte B\, \exp(A\, D_N^3) \, \sqrt{N \log N},
\ee
\textcolor{black}{\eqref{cond1} holds provided $D_N = o((\log(N\, /\, \log N))^{1/3})$,
and} 
\beno
\widetilde\Phi_N(\truth)
&\leq& C\, \exp(A\, D_N^3)\; \sqrt{\dfrac{\log N}{N^{1-2\, \alpha}}}.
\ee
\end{corollary} 
} 

\vspace{-.25cm}

A comparison of Corollaries \ref{cor:model_2} and \ref{cor:model_3} reveals that overlap comes at a cost.
First,
the convergence rate is lower due to the factor $\exp(A\, D_N^3)$ in the overlapping subpopulation scenario, 
compared with the factor $D_N^{14}$ in the non-overlapping subpopulation scenario.
Second,
overlap requires stronger restrictions on $D_N$.
For example,
consider the best-case scenario when the random graph is dense ($\alpha = 0$) and $\norm{\truth}_\infty$ is bounded above ($\vartheta=0$).
Then,
to ensure $\norm{\truth}_\infty > 0$ and $\norm{\mple - \truth}_\infty \ip 0$,\,
$D_N$ needs to satisfy 
\begin{itemize}[topsep=2pt,itemsep=2pt]
\item $D_N = O(\log N)$ when the subpopulations do not overlap;
\item $D_N = o((\log(N\, /\, \log N))^{1/3})$ when the subpopulations do overlap.
\ei
These results dovetail with results on other statistical exponential-family models for discrete and dependent random variables in single-observation scenarios.
For example,
\citet{ChDi11} considered the edge-and-triangle model with $p=2$ parameters and unbounded $D_N$ of order $O(N)$,
but concluded that the edge-and-triangle model possesses undesirable properties and did not report consistency results.
Likewise,
the recent results of \citet{GhMu20} on Ising models with $p=2$ parameters suggest that consistency results may not be obtainable unless $D_N$ is bounded or other restrictions are imposed.
By contrast,
we 
\begin{itemize}[topsep=3pt,itemsep=2pt]
\item allow $D_N \to \infty$ as $N \to \infty$ provided $D_N = O(\log N)$ (non-overlapping subpopulations) or $D_N = o((\log(N\, /\, \log N))^{1/3})$ (overlapping subpopulations),
as discussed above;
\item allow $p \to \infty$ as $N \to \infty$ provided $p = o(N^2 / \log N)$,
as discussed in Section \ref{sec:dimension};
\item cover a wide range of model specifications, 
beyond the pairwise interaction terms of discrete graphical models (e.g., Ising models).
\ei

\vspace{-.25cm}

\section*{Supplementary materials}

Proofs of Proposition \ref{p:expected.degrees},
Theorem \ref{thm:mple_consistency},
and Corollaries \ref{cor:model_1}--\ref{cor:model_3} along with simulations can be found in the supplement \citep[][]{StSc20}.
 
\section*{Acknowledgements}

We acknowledge support via NSF awards\break 
SES-2345043, 
DMS-1812119,   
DMS-1513644, 
and DoD award ARO\break 
W911NF-21-1-0335.
We are indebted to an anonymous referee who suggested the proof of Lemma \ref{lemma.inequality} and helped sharpen Theorem \ref{theorem.mle}.

\bibliographystyle{asa1}

\bibliography{base}

\pagebreak

\appendix

\makeatletter

\setcounter{page}{1}

\setcounter{section}{0}

\setcounter{com}{0}

\begin{center}
\Large\bf\textsc\bf
\vspace{-.75in}\\
Supplement to:\\
\longtitle\s\s

\normalsize
{\normalfont\textsc{By Jonathan R.\ Stewart and Michael Schweinberger}}\s
\\
{\normalfont\em Florida State University and Penn State University}\s\s
\end{center}

\thispagestyle{empty}

\noindent
{Appendix \ref{sup-sec:mle}: Auxiliary results for Theorem 1}\dotfill\pageref{sup-sec:mle}\s
\\
{Appendix \ref{sec:mple}: Proof of Theorem 2}\dotfill\pageref{sec:mple}\s
\\
\mbox{\hspace{0.75cm}}Appendix \ref{mple.auxiliary}: Auxiliary results for Theorem 2\dotfill\pageref{mple.auxiliary}\s
\\
{Appendix \ref{sec:conc}: Proofs of Corollaries 1--3}\dotfill\pageref{sec:conc}\s
\\
\mbox{\hspace{0.75cm}}Appendix \ref{sup-sec:scale_pl}: Bounding \texorpdfstring{$\widetilde\Lambda_N(\truth)$}{Lambda}\dotfill\pageref{sup-sec:scale_pl}\s
\\
\mbox{\hspace{0.75cm}}Appendix \ref{sup-sec:mrf_conc}: Bounding \texorpdfstring{$\mnorm{\mD_N(\truth)}_2$}{D}\dotfill\pageref{sup-sec:mrf_conc}\s
\\
{Appendix \ref{proof.degrees}: Proof of Proposition 1}\dotfill\pageref{proof.degrees}\s
\\
{Appendix \ref{sec:simulationsstudy}: Simulation results}\dotfill\pageref{sec:simulationsstudy}\s

\s

In Appendices \ref{sup-sec:mle}, 
\ref{sec:mple}, 
\ref{D} and \ref{sup-sec:D_aux_proofs},  
we adopt the notation used in Section \ref{sec:stat_inf} of the manuscript,
by denoting the number of  edge variables by $M \coloneqq \binom{N}{2}$ and edge variables by $X_1, \dots, X_M$. 
In addition,
we denote the data-generating parameter vector by $\truth \in \Nat = \mR^p$ and the data-generating probability measure and expectation by $\mbP \equiv \mbP_{\truth}$ and $\mbE \equiv \mbE_{\truth}$,
respectively.
Throughout,
we assume that $\min_{1 \leq k \leq K} |\mA_k| \geq 3$.
\hide{
To simplify the proofs,
we suppress the argument $\truth$ of all other quantities that depend on the data-generating parameter vector $\truth$,
and suppress the subscript $N$ of all quantities that depend on the number of nodes $N$.
}

\hide{

\s

{\bf References.}
The references in the supplement are identical to the references in the main manuscript.
A list of all references in the main manuscript and supplement can be found on pages 22--25 of the main manuscript.

\s
}

\setcounter{section}{0}

\hide{

\section{Additional results}
\label{model4}

We present an additional version of generalized $\beta$-models with size-dependent parameterizations,
called Model 4.

Model 2,
introduced in Section \ref{model2} of the manuscript, 
assumes that the weight of the brokerage statistic $b_{i,j}$ does not depend on the sizes of neighborhood intersections.
While convenient on mathematical grounds,
the assumption that small and large neighborhood intersections have the same brokerage parameter may be unwarranted.
To allow small and large neighborhood intersections to have distinct brokerage parameters,
we consider functions of edges $\varphi_{i,j}$ of the form
\beno
\varphi_{i,j}(x_{i,j}, \bx_{\mS_{i,j}}; \nat)
= a_{i,j}(x_{i,j}) \exp\left((\theta_i + \theta_j)\, x_{i,j} + \eta_{i,j}(\theta_{N+1}) \, b_{i,j}(x_{i,j}, \bx_{\mS_{i,j}})\right),
\ee
where $a_{i,j}$ is $1$ if $x_{i,j} \in \{0, 1\}$ and is $0$ otherwise, 
$b_{i,j}$ is defined by \eqref{brokerage.statistic},
and
\beno
\label{eta.parameterization}
\eta_{i,j}(\theta_{N+1})
&=& \theta_{N+1} \, \log\left(1 + \dfrac{\log |\mN_i\, \cap\, \mN_j|}{|\mN_i\, \cap\, \mN_j|}\right),
&& \theta_{N+1}  \in \mbR.
\ee

\vspace{-.25cm}

\begin{corollary}
\label{cor:model_4}
Consider Model 4,
the generalized $\beta$-model with dependent edges and a size-dependent parameterization,
with $\truth \in \mbR^{N+1}$ and $\gamma$ satisfying the same conditions as in Corollary \ref{cor:model_3}.
Let $C > 0$,\,
$C_1 > 0$,\,
$C_2 > 0$,\,
$C_3 > 0$,\,
and $N_0 \geq 3$ be finite constants,
independent of $N$ and $p$.
Then,
for all $N > N_0$,
\beno
\widetilde\Phi
&\geq& \dfrac{C}{\raisebox{-4pt}{$\mUTA\; D^3\; \mnorm{\mD}_2$}}\;\, \sqrt{\dfrac{N^{2\, \vartheta - 1}}{\log N}},
\ee
where edges are dependent provided $\theta_{N+1} \neq 0$, 
$D_N \coloneqq \max_{1 \leq i \leq M} |\mmG_i| \geq 1$ can increase as a function of $N$ provided $\theta_{N+1} \neq 0$,
and $\min_{1 \leq k \leq K} |\mA_k| \geq 3$. 
If Assumption A is satisfied,
then:
\bi
\item Scenario 1: If the subpopulations do not intersect ($\omega_1 = \omega_{2} = 0$) and $\truth \in \mbR^{N+1}$ satisfies \eqref{condition.theta0} with $\vartheta \in (4/5,\, 1]$,
then
\beno
\widetilde\Phi
&\geq& \dfrac{C_1}{\raisebox{-4pt}{$\mUTA\; D^5$}}\;\, \sqrt{\dfrac{N^{2\, \vartheta - 1}}{\log N}}.
\ee
\item Scenario 2: If the subpopulations do intersect and $\truth \in \mbR^{N+1}$ satisfies \eqref{condition.theta0} with $\vartheta = 1$,
then
\beno
\widetilde\Phi
&\geq& \dfrac{C_2}{\raisebox{-4pt}{$\mUTA\; D^5\, \exp(C_3\, D^2)$}}\;\, \sqrt{\dfrac{N}{\log N}}.
\ee
\ei
\end{corollary}
Corollary \ref{cor:model_4} demonstrates that theoretical guarantees can be obtained for models with size-dependent parameterizations,
which allow parameters to depend on the sizes of neighborhood intersections.

\s

}

\section{Auxiliary results for Theorem 1} 
\label{sup-sec:mle}

\textcolor{black}{Theorem \ref{theorem.mle} along with Lemmas \ref{lemma.inequality} and \ref{lemma.boundary} are stated and proved in Section \ref{sec:consistency} of the manuscript.
Here,
we state and prove Lemmas \ref{lem:hom_bound_to_bound} and \ref{theorem.mle.lemma},
which are used to prove Theorem \ref{theorem.mle} and Lemmas \ref{lemma.inequality} and \ref{lemma.boundary}.
}

\vspace{-.65cm}

\textcolor{black}{ 
\begin{lemma}
\label{lem:hom_bound_to_bound}
Let $\bg: \mD \mapsto \mathscr{R}$ be a homeomorphism between $\mD \subseteq \mbR^p$ and $\mathscr{R} \subseteq \mbR^p$.
Consider any $\truth \in \mD$, 
any $\rho \in (0,\, \infty)$,
and any vector norm $\norm{\cdot}$,
and define $\mB(\bm{c},\, \rho) \coloneqq \left\{\bm{v} \in \mbR^p:\, \norm{\bm{v} - \bm{c}} < \rho\right\}$ ($\bm{c} \in \mR^p$). 
Then 
\beno
\bg(\bd \mB(\truth,\, \rho)) 
\= \bd \bg(\mB(\truth,\, \rho)).
\ee
\end{lemma}
\llproof \ref{lem:hom_bound_to_bound}.
\hide{
We consider $\mbR^p$ equipped with the topology induced by the metric $d(\bv,\, \bw) \coloneqq \norm{\bv - \bw}$ ($\bv,\, \bw \in \mbR^p$).
The balls $\mB(\bv,\, \rho)$ form a basis for the topology, 
so any open set in the topology can be written as a union of balls in $\mbR^p$.  
As a result, 
the boundary of a proper subset $\mS \subset \mbR^p$ can be written as
\beno
\bd \mS 
\;=\; \left\{\bv \in \mbR^p:\, 
\mB(\bv,\, \rho)\, \cap\, \mS \neq \emptyset  
\;\mbox{ and }\; 
\mB(\bv,\, \rho)\, \cap\, \mS^c \neq \emptyset 
\,\mbox{ for all }\, \rho > 0\right\}.  
\ee 
Since $\bg: \mD \mapsto \mathscr{R}$ is a homeomorphism,
it is bijective,
hence
\be
\label{eq:disjoint}
\bg(\mB(\truth,\, \rho)) 
\;\cap\; \bg(\mD \setminus \mB(\truth,\, \rho))
\;=\; \emptyset
\;\;\mbox{ for all }\;\; \mB(\truth,\, \rho) 
\;\subset\; \mD. 
\ee 
}
We prove Lemma \ref{lem:hom_bound_to_bound} by proving that
\beno
1.\;\;\; \bg(\bd \mB(\truth,\, \rho)) 
&\subseteq& \bd \bg(\mB(\truth,\, \rho))\vspace{.3cm}
\\
2.\;\;\; \bd \bg(\mB(\truth,\, \rho)) 
&\subseteq& \bg(\bd \mB(\truth,\, \rho)),
\ee
which establishes the desired result:
\beno
\bg(\bd \mB(\truth,\, \rho)) 
&=& \bd \bg(\mB(\truth,\, \rho)).
\ee 

\s

{\bf 1. Proving $\bg(\bd \mB(\truth,\, \rho)) \subseteq \bd \bg(\mB(\truth,\, \rho))$.}
Consider any
\beno
\nat_{\mbox{\footnotesize $\bd$}} &\in& \bd \mB(\truth,\, \rho)
\= \left\{\nat \in \mbR^p:\, \norm{\nat - \truth} = \rho\right\}. 
\ee
Since $\bg: \mD \mapsto \mathscr{R}$ is a homeomorphism,
it is continuous and one-to-one.
Thus,
for any $\epsilon > 0$,
there exists a real number $\delta(\epsilon) > 0$ such that 
\be
\label{ggggg}
\bg(\mB(\nat_{\mbox{\footnotesize $\bd$}}, \, \delta(\epsilon)))
&\subseteq& \mB(\bg(\nat_{\mbox{\footnotesize $\bd$}}), \, \epsilon).
\ee
In light of the fact that $\nat_{\mbox{\footnotesize $\bd$}}$ is an element of the boundary $\bd \mB(\truth,\, \rho)$ of $\mB(\truth,\, \rho)$, 
there exist points 
\begin{equation}
\label{ffff}
\begin{array}{ccccccc}
\nat_1 
&\in& \mB(\nat_{\mbox{\footnotesize $\bd$}},\, \delta(\epsilon)) 
&\cap& \mB(\truth,\, \rho)\s
\\
\nat_2 
&\in& \mB(\nat_{\mbox{\footnotesize $\bd$}},\, \delta(\epsilon)) 
&\cap& (\mathscr{D} \setminus \mB(\truth,\, \rho)).
\end{array}
\end{equation}
The fact that $\bg(\nat_1) \in \mB(\bg(\nat_{\mbox{\footnotesize $\bd$}}),\, \epsilon)$ and $\bg(\nat_2) \in \mB(\bg(\nat_{\mbox{\footnotesize $\bd$}}),\, \epsilon)$ by Equations \eqref{ggggg} and \eqref{ffff} implies that
\begin{equation}
\begin{array}{cccccc}
\label{eq:show_bound}
{ \bg(\nat_1) \;\;\in\;\; } 
\mB(\bg(\nat_{\mbox{\footnotesize $\bd$}}),\, \epsilon)
&\cap& \bg(\mB(\truth,\, \rho))
&\neq& \emptyset \s \\
{ \bg(\nat_2) \;\;\in\;\;} 
\mB(\bg(\nat_{\mbox{\footnotesize $\bd$}}),\, \epsilon)
&\cap& (\mathscr{R} \setminus \bg(\mB(\truth,\, \rho)))
&\neq& \emptyset.
\end{array}
\end{equation}
Equation \eqref{eq:show_bound} holds for all $\epsilon > 0$,
so $\bg(\nat_{\mbox{\footnotesize $\bd$}}) \in \bd  \bg(\mB(\truth,\, \rho))$.
Since\break 
$\nat_{\mbox{\footnotesize $\bd$}} \in \bd \mB(\truth,\, \rho)$ was arbitrary,
$\bg(\bd  \mB(\truth,\, \rho)) \subseteq \bd  \bg(\mB(\truth,\, \rho))$.

\s

{\bf 2. Proving $\bd  \bg(\mB(\truth,\, \rho)) \subseteq \bg(\bd  \mB(\truth,\, \rho))$.} 
Consider any
\beno
\bg_{\mbox{\footnotesize $\bd$}} 
&\in& \bd \bg(\mB(\truth,\, \rho)).
\ee
By assumption,
$\bg: \mD \mapsto \mathscr{R}$ is a homeomorphism, 
so $\bg^{-1} : \mathscr{R} \mapsto \mD$ exists and is continuous and one-to-one.
As a consequence,
for any $\delta > 0$,
there exists a real number $\epsilon(\delta) > 0$ such that 
\be
\label{eeee0}
\bg^{-1}\left(\mB(\bg_{\mbox{\footnotesize $\bd$}}, \, \epsilon(\delta)) \right)
&\subseteq& \mB(\bg^{-1}(\bg_{\mbox{\footnotesize $\bd$}}), \, \delta). 
\ee
Since $\bg_{\mbox{\footnotesize $\bd$}} \in \bd  \bg(\mB(\truth,\, \rho))$, 
there exist points 
\be
\label{eeee1}
\begin{array}{ccccccc}
\bg_1
&\in& \mB(\bg_{\mbox{\footnotesize $\bd$}}, \, \epsilon(\delta)) 
&\cap& \bg(\mB(\truth,\, \rho))\s
\\
\bg_2
&\in& \mB(\bg_{\mbox{\footnotesize $\bd$}}, \, \epsilon(\delta)) 
&\cap& (\mathscr{R} \setminus \bg(\mB(\truth,\, \rho))).
\end{array}
\ee
Since $\bg^{-1}(\bg_1) \in \mB(\bg^{-1}(\bg_{\mbox{\footnotesize $\bd$}}),\, \delta)$ and $\bg^{-1}(\bg_2) \in \mB(\bg^{-1}(\bg_{\mbox{\footnotesize $\bd$}}),\, \delta)$ by Equations \eqref{eeee0} and \eqref{eeee1}, 
\hide{
we conclude that
\beno
\bg^{-1}(\bg(\mB(\truth,\, \rho)))
&\cap& 
\bg^{-1}(\mathscr{R} \setminus \bg(\mB(\truth,\, \rho)))
&=& \emptyset,  
\ee
which implies that
}
\begin{equation}
\begin{array}{cccccc}
\label{eq:show_bound_2}
{ \bg^{-1}(\bg_1) \;\;\in\;\; }
\mB(\bg^{-1}(\bg_{\mbox{\footnotesize $\bd$}}),\, \delta)
&\cap& \mB(\truth,\, \rho) 
&\neq& \emptyset\s
\\
{ \bg^{-1}(\bg_2) \;\;\in\;\; }
\mB(\bg^{-1}(\bg_{\mbox{\footnotesize $\bd$}}),\, \delta)
&\cap& (\mD \setminus \mB(\nat,\, \rho))
&\neq& \emptyset. 
\end{array}
\end{equation}
Equation \eqref{eq:show_bound_2} holds for all $\delta > 0$,
hence $\bg^{-1}(\bg_{\mbox{\footnotesize $\bd$}}) \in \bd  \mB(\truth,\, \rho)$.
Since $\bg_{\mbox{\footnotesize $\bd$}} \in \bd  \bg(\mB(\truth,\, \rho))$ was arbitrary,
$\bg^{-1}(\bd  \bg(\mB(\truth,\, \rho))) \subseteq \bd \mB(\truth,\, \rho)$.\break
We have therefore established that $\bd \bg(\mB(\truth,\, \rho)) \subseteq\, \bg(\bd \mB(\truth,\, \rho))$,\break
because $\bg^{-1}: \mathscr{R} \mapsto \mD$ is the inverse map corresponding to $\bg: \mD \mapsto \mathscr{R}$.
\qed

}

\vspace{-.5cm}

\textcolor{black}{ 
\begin{lemma}
\label{theorem.mle.lemma}
Under the assumptions of Theorem \ref{theorem.mle},
for all $t > 0$,
\vspace{-.05cm}
\beno
\mbP\left(s(\bX)\; \in\; \mB_{\infty}(\bmu(\truth),\, t)\right)
&\geq& 1 - 2\; \exp\left( - \dfrac{2\; t^2}{\mnorm{\mD_N(\truth)}_2^2 \; \Psi_N^2} + \log p\right),
\ee
where $\mnorm{\mD_N(\truth)}_2 \geq 1$ and $\Psi_N > 0$ provided $N$ is large enough.
\end{lemma}

}

\textcolor{black}{ 
\hide{ 

\vspace{-.1cm}

\begin{lemma}
\label{theorem.mle.lemma}
Under the assumptions of Theorem \ref{theorem.mle}, 
for all $t > 0$,
\vspace{-.05cm}
\beno
\mbP\left(s(\bX)\; \in\; \mB_{\infty}(\bmu(\truth),\, t)\right)
&\geq& 1 - 2\; \exp\left( - \dfrac{2\; t^2}{\mnorm{\mD_N(\truth)}_2^2 \; \Psi_N^2} + \log p\right),
\ee
where $\mnorm{\mD_N(\truth)}_2 \geq 1$ and $\Psi_N > 0$ provided $N$ is large enough.
\end{lemma}

\s

{\em Remark.}
The concentration results of \citet{Chetal07} used in the proof of Lemma \ref{theorem.mle.lemma} hold for all possible couplings $\mbQ_{\truth,i,\bx_{1:i-1}}$ of the conditional probability mass functions $\mbP_{\truth,\bx_{1:i-1},0}$\, and $\mbP_{\truth,\bx_{1:i-1},1}$ defined in Section \ref{sec:background} of the manuscript.
Having said that,
non-optimal couplings can lead to less sharp concentration results compared with optimal couplings.

\s

}

\llproof \ref{theorem.mle.lemma}.
By Theorem 1 of \citet[][p.\ 207]{Chetal07},
\beno
\mbP\left(|s_i(\bX) - \mbE \, s_i(\bX)|\,\geq\, t\right)
\,\leq\, 2 \, \exp\left( - \dfrac{2 \; t^2}{\mnorm{\mD_N(\truth)}_2^2 \; \norm{\bm\Xi_i}_2^2} \right),\;
 i = 1, \dots, p.
\ee
A union bound over the $p$ coordinates of $s(\bX)$ shows that
\beno
\mbP\left(\norm{s(\bX) - \mbE \, s(\bX)}_\infty\,\geq\, t\right)
\lte 2 \, \exp\left( - \dfrac{2 \; t^2}{\mnorm{\mD_N(\truth)}_2^2\,\; \Psi_N^2} + \log p\right),
\ee
where $\Psi_N \coloneqq \max_{1 \leq i \leq p}\; \norm{\bm\Xi_i}_2$.
As a result,
we obtain
\beno
\mbP\left(s(\bX)\; \in\; \mB_{\infty}(\bmu(\truth),\, t)\right)
\gte 1 -  2\; \exp\left( - \dfrac{2 \; t^2}{\mnorm{\mD_N(\truth)}_2^2\,\; \Psi_N^2} + \log p\right),
\ee
using $\bmu(\truth) \coloneqq \mbE_{\truth}\, s(\bX)$.
\qed
 
}

\com {\em Extensions to dependent random variables with countable and uncountable sample spaces.}
\label{exten}
Theorem \ref{theorem.mle} is not restricted to random graphs with dependent edges.
It covers models of dependent random variables with finite sample spaces,
and can be extended to countable sample spaces:
e.g.,
the concentration result of \citet{Chetal07} used in Theorem \ref{theorem.mle} assumes that the sample spaces are finite---motivated by applications to Ising models---but could be extended to countable sample spaces.
Uncountable sample spaces could be accommodated by replacing the concentration result of \citet{Chetal07} by other suitable concentration results,
e.g.,
Subgaussian concentration results.
Likewise,
the exponential-family properties used in Theorem \ref{theorem.mle} are neither restricted to finite nor countable sample spaces \citep{Br86}.

\vspace{.125cm}

\section{Proof of Theorem 2}
\label{sec:mple}

We prove Theorem \ref{thm:mple_consistency} stated in Section \ref{mple} of the manuscript.
Auxiliary results are proved in Appendix \ref{mple.auxiliary}.

\s

\blue{ 
\ttproof \ref{thm:mple_consistency}. 
Define
\beno
\bg(\nat;\, \bx) 
&\coloneqq& \nabla_{\nat}\; \widetilde\ell(\nat;\, \bx),
&& \nat \in \Nat, 
& \bx \in \mbX, 
\ee
where $\widetilde\ell(\nat;\, \bx)$ is the pseudo-loglikelihood of $\nat \in \Nat$ based on $\bx \in \mbX$,
and 
\beno
\mbG
&\coloneqq& \left\{\bx \in \mbX:\, \norm{\bg(\truth;\, \bx)}_{\infty}\; \leq\; \gamma_N\right\},
& \gamma_N \in [0,\, \infty).
\ee
Recall that $\mbH$ is a subset of $\mbX$ satisfying 
\beno
\mbH
&\subseteq&
\left\{\bx \in \mbX:\; -\nabla_{\nat}^2\; \widetilde\ell(\nat;\, \bx) \mbox{ is invertible for all $\nat \in \mB_{\infty}(\truth,\, \epsilon^\star)$}\right\}
\ee 
and
\beno
\label{r1}
\mbP\left(\bX \in \mbH\right)
\gte 1 - \dfrac{2}{\max\{N,\, p\}^2},
\ee
where $\epsilon^\star \in (0,\, \infty)$ is a constant independent of $N$ and $p$. 
The conditional distributions of edge variables are exponential families. 
Thus,
$\widetilde\ell(\nat;\, \bx)$ is twice continuously differentiable in $\nat \in \Nat$,
because $\widetilde\ell(\nat;\, \bx)$ is a sum of exponential-family loglikelihoods,
each of which is twice continuously  differentiable. 

\s

{\bf I.\ In the event $\bX \in \mbG$,
the set $\widetilde{\bm\Theta}(\gamma_N)$ is non-empty.}
By the construction of the sets  $\mbG$ and $\Mple$,
the set $\Mple$ is non-empty for all $\bx \in \mbG$, 
because 
\beno
\truth 
\,\in\, \Mple
\,\coloneqq\, \left\{\nat \in \Nat:\; \norm{\bg(\nat;\, \bx)}_{\infty}\, \leq\, \gamma_N\right\}.
\ee

{\bf II.\ In the event $\bX \in \mbG\, \cap\,\, \mbH$,\,
the set $\Mple$ satisfies $\Mple \subseteq \mB_{\infty}(\truth,\, \sqrt{96}\;\, \widetilde\Phi_N(\truth))$ provided $N > N_0$.}
Let $\epsilon = \sqrt{96}\;\, \widetilde\Phi_N(\truth) > 0$,
where 
\beno
\widetilde{\Phi}_N(\truth)
&\coloneqq& \tllambdamin\,\, (1 + D_N) \,\, |\!|\!|\mD_N(\truth)|\!|\!|_2\,\, \Psi_N\; \sqrt{\log \max\{N,\, p\}},
\ee 
recalling the definitions 
\beno
\widetilde{\Lambda}_{N,\bx}(\truth)
&\coloneqq& \sup\limits_{\nat\, \in\, \mB_{\infty}(\truth,\, \epsilon^\star)}\, \mnorm{(-\nabla_{\nat}^2\; \widetilde\ell(\nat;\, \bx))^{-1}}_{\infty}\;\;\; \mbox{for all } \bx \in \mbH\s
\\
\widetilde{\Lambda}_{N}(\truth)
&\coloneqq& \max\limits_{\bx\, \in\, \mbH}\, \widetilde\Lambda_{N,\bx}(\truth).
\ee
By assumption,
$\widetilde\Phi_N(\truth) \to 0$ as $N \to \infty$,
which implies that there exists a constant $N_0 \geq 3$,
independent of $N$ and $p$, 
such that 
\beno
\epsilon
&=& \sqrt{96}\;\, \widetilde\Phi_N(\truth) 
&<& \epsilon^\star
&\mbox{for all}& N > N_0.
\ee
By the definition of $\mbH$,
the matrix $-\nabla_{\nat}^2\; \widetilde\ell(\nat;\, \bx)$ is invertible
on $\mB_{\infty}(\truth, \epsilon^\star)$ 
for all $\bx \in \mbH$. 
By Lemma \ref{lem:support_mple},
$\widetilde{\ell}(\, \cdot\, ;\, \bx)$ is strictly concave on $\Nat$ for any given $\bx \in \mbH$,
which implies that $\bg(\nat;\, \bx) = - \nabla_{\nat} \, \widetilde\ell(\nat;\, \bx)$,
considered as a function of $\nat \in\Nat$ for fixed $\bx \in \mbG\, \cap\, \mbH$,
is continuous and injective, 
and is thus a homeomorphism 
by the invariance of domain theorem. 
Since the inverse $\bg^{-1}(\, \cdot\, ;\, \bx)$ of\, $\bg(\, \cdot\, ;\, \bx)$ exists and is continuous on $\Nat$,
there exists, 
for each $\bx \in \mbH$, 
a real number $\delta_{\bx}(\epsilon) \in (0,\, \infty)$ such that 
\be
\label{continuity}
\norm{\bg(\nat;\, \bx) - \bg(\truth;\, \bx)}_{\infty}
&\leq& \delta_{\bx}(\epsilon)
&\mbox{implies}&
\norm{\nat - \truth}_{\infty} 
&\leq& \epsilon. 
\ee
Since the pseudo-loglikelihood function is a sum of conditional Bernoulli loglikelihood functions (i.e., conditional exponential-family loglikelihood functions) and each of them is twice continuously differentiable \citep[Theorem 2.2, pp.\ 34--35,][]{Br86},
we can invoke Lemma \ref{lemma.inequality} in the manuscript to conclude that $\delta_{\bx}(\epsilon)$ is related to $\epsilon$ by the following inequality:
\be
\label{eq:lower_relation}
\dfrac{\epsilon}{\raisebox{-3pt}{$\widetilde\Lambda_{N,\bx}(\truth)$}}
&\leq& \delta_{\bx}(\epsilon).
\ee
To leverage \eqref{eq:lower_relation},
note that,
for all $\nat \in \Mple$ and all $\bx \in \mbG \cap \mbH$,
\begin{equation}
\begin{array}{cccccccc}
\label{bound00}
\norm{\bg(\nat;\, \bx) - \bg(\truth;\, \bx)}_{\infty}
&\leq& \norm{\bg(\nat;\, \bx)}_{\infty} &+& \norm{\bg(\truth;\, \bx)}_{\infty}\s
\\
&\leq& \underbrace{\gamma_N} &+& \underbrace{\gamma_N}\s
\\
\mbox{\em by virtue of} && \mbox{\em $\nat \in \Mple$} && \mbox{\em $\bx \in \mbG$}\vspace{.25cm}
\\
&=& \dfrac{\epsilon}{\raisebox{-3pt}{$\widetilde\Lambda_{N}(\truth)$}},
\end{array}
\end{equation}
because $\epsilon = \sqrt{96} \; \widetilde\Phi_N(\truth) > 0$ and $\gamma_N$ is assumed to be of the form 
\beno
\gamma_N
&=& \dfrac{\epsilon}{\raisebox{-3pt}{$2\; \widetilde\Lambda_{N}(\truth)$}}\vspace{.4cm}
\\
&=& \sqrt{24}\;\, (1 + D_N) \; |\!|\!|\mD_N(\truth)|\!|\!|_2\; \Psi_N\; \sqrt{\log \max\{N,\, p\}}.
\ee
So,
by the definition of 
\beno
\widetilde\Lambda_{N,\bx}(\truth) 
\;>\; 0 
&&\mbox{and}&&
\widetilde\Lambda_{N}(\truth) 
\;\coloneqq\; \max\limits_{\bx \in \mbH}\, \widetilde\Lambda_{N,\bx}(\truth) 
\;>\; 0,
\ee
we have 
\be
\label{bound00}
\norm{\bg(\nat;\, \bx) - \bg(\truth;\, \bx)}_{\infty}
&\leq& \dfrac{\epsilon}{\raisebox{-3pt}{$\widetilde\Lambda_{N}(\truth)$}}
&\leq& \dfrac{\epsilon}{\raisebox{-3pt}{$\widetilde\Lambda_{N,\bx}(\truth)$}}
\ee
and, 
according to \eqref{eq:lower_relation},
\be
\label{delta.x}
\norm{\bg(\nat;\, \bx) - \bg(\truth;\, \bx)}_{\infty}
&\leq& \dfrac{\epsilon}{\raisebox{-3pt}{$\widetilde\Lambda_{N,\bx}(\truth)$}}
&\leq& \delta_{\bx}(\epsilon).
\ee
Since
\beno
\norm{\bg(\nat;\, \bx) - \bg(\truth;\, \bx)}_{\infty}
&\leq& \delta_{\bx}(\epsilon)
&\mbox{implies}&
|\!|\nat - \truth|\!|_\infty 
&\leq& \epsilon
\ee
and
\beno
\epsilon
&=& \sqrt{96}\;\, \widetilde\Phi_N(\truth)
&<& \epsilon^\star
&\mbox{for all}& N > N_0,
\ee
the random set $\Mple$ is non-empty and satisfies
\be
\label{keyresult}
\Mple 
\,\subseteq\, \mB_{\infty}(\truth,\, \sqrt{96}\;\, \widetilde\Phi_N(\truth))
\,\subset\, \mB_{\infty}(\truth,\, \epsilon^\star)
\mbox{ for all } N > N_0
\ee
in the event $\bX \in \mbG\, \cap\, \mbH$.

\s

{\bf III. The event $\bX\, \in\, \mbG\, \cap\, \mbH$ occurs with probability at least $1 - 4\, / \max\{N,\, p\}^2$ provided $N > N_0$.}
A union bound shows that
\be
\label{union.bound}
\mbP\left(\bX \in \mbG \cap\, \mbH\right)
\geq 1 - \mbP\left(\bX \in \mbX \setminus \mbG\right) - \mbP\left(\bX \in \mbX \setminus \mbH\right).
\ee
The first probability on the right-hand side of \eqref{union.bound} can be bounded above by invoking Lemma \ref{uc.lemma} along with
\beno
\epsilon 
&=& \sqrt{96}\;\, \widetilde\Phi_N(\truth)
&>& 0, 
\ee
which leads to the following upper bound: 
\be
\label{r2}
&& \mbP\left(\bX \in \mbX\, \setminus\, \mbG\right)\s
\\
&\leq& \mbP\left(\norm{\bg(\truth;\, \bX)}_{\infty}\; \geq\; \dfrac{\epsilon}{\raisebox{-3pt}{$2\; \widetilde\Lambda_{N}(\truth)$}}\right)\s
\\
&=& \mbP\left(\norm{\bg(\truth; \bX) - \mbE\,\, \bg(\truth; \bX)}_{\infty}\; \geq\; \dfrac{\epsilon}{\raisebox{-3pt}{$2\; \widetilde\Lambda_{N}(\truth)$}}
\right)\s
\\
&\leq& 2\, \exp\left(-\dfrac{\epsilon^2}{\raisebox{-3pt}{32\, $\widetilde\Lambda_N(\truth)^2\; (1 + D_N)^2\; \mnorm{\mD_N(\truth)}_2^2\; \Psi_N^2$}} + \log \, p \right)\s
\\
&=& \dfrac{2}{\max\{N,\, p\}^2},
\ee 
using the fact that $\mbE\,\, \bg(\truth; \bX) = \bm{0}$ by Lemma \ref{lem:concave_pl}.
The second probability on the right-hand side of \eqref{union.bound} is bounded above by assumption \eqref{cond1}:
\be
\label{r1}
\mbP\left(\bX \in \mbX\, \setminus\, \mbH\right)
\lte \dfrac{2}{\max\{N,\, p\}^2}.
\ee

\s

{\bf IV. Conclusion.}
Combining \eqref{keyresult} with \eqref{r2} and \eqref{r1} establishes that,
for all $N > N_0$,
the random set $\Mple$ is non-empty and satisfies
\beno
\Mple 
&\subseteq& \mB_{\infty}(\truth,\, \sqrt{96}\;\, \widetilde\Phi_N(\truth))
&\subset& \mB_{\infty}(\truth,\, \epsilon^\star)
\ee 
with probability at least $1 - 4\, / \max\{N,\, p\}^2$,
provided
\beno
\gamma_N
\hide{
&=& \dfrac{\epsilon}{2\; \Lambda_{N}(\truth)} 
}
&=& \sqrt{24}\;\, (1 + D_N) \; |\!|\!|\mD_N(\truth)|\!|\!|_2\; \Psi_N\; \sqrt{\log \max\{N,\, p\}}
&>& 0.\qed
\ee

\hide{

{\bf Remark.}
The assumptions of Theorem \ref{thm:mple_consistency} imply a form of stochastic equicontinuity in the sense of \citet{Ne91},
because \eqref{delta.x} implies that,
for all $\nat \in \Mple$ and all $\bx\, \in\, \mbG\, \cap\, \mbH$,
\beno
\norm{\bg(\nat;\, \bx) - \bg(\truth;\, \bx)}_{\infty}
&\leq& \delta_{\bx}(\epsilon),
\ee
which in turn implies that
\beno
\sup\limits_{\nat\, \in\, \mB_{\infty}(\truth,\, \epsilon^\star)}\, \norm{\bg(\nat;\, \bx) - \bg(\truth;\, \bx)}_{\infty} 
&\leq& \delta_{\bx}(\epsilon),
\ee
with probability at least $1 - 4\, / \max\{N,\, p\}^2$. 
That said,
we deal with $p \geq N \to \infty$ parameters,
in contrast to \citet{Ne91}.

}

}

\hide{

\vspace{-.15cm}

\com {\em Extensions to dependent random variables with countable and uncountable sample spaces.}
Neither Theorem \ref{theorem.mle} nor Theorem \ref{thm:mple_consistency} are restricted to random graphs with dependent edges.
Both theorems can be extended to dependent random variables with countable and uncountable sample spaces:
see the remark in Appendix \ref{sup-sec:mle}.

}

\subsection{Auxiliary results for Theorem 2}
\label{mple.auxiliary}

\mbox{}

\vspace{-.5cm}

\blue{ 
\begin{lemma}
\label{lem:support_mple} 
$\widetilde{\ell}(\, \cdot\, ;\, \bx)$ is strictly concave on $\Nat$ provided $\bx \in \mbH$,
where
\beno
\mbH
&\subseteq& \left\{\bx \in \mbX:\; -\nabla_{\nat}^2\; \widetilde\ell(\nat;\, \bx) \mbox{ is invertible for all $\nat \in \mB_{\infty}(\truth,\, \epsilon^\star)$}\right\}.
\ee
\end{lemma}
\llproof \ref{lem:support_mple}. 
The pseudo-loglikelihood function $\widetilde{\ell}(\, \cdot\, ;\, \bx)$ is a sum of loglikelihoods based on conditional Bernoulli distributions,
each of which is concave on $\Nat$ \citep[][Lemma 5.3, p.~146]{Br86}.
As a result,
$\widetilde{\ell}(\, \cdot\, ;\, \bx)$ is concave on $\Nat$.
We first establish strict concavity of $\widetilde{\ell}(\, \cdot\, ;\, \bx)$ on $\mB_{\infty}(\truth,\, \epsilon^\star)$ and then extend the result from $\mB_{\infty}(\truth,\, \epsilon^\star)$ to $\Nat$.
All of the following results focus on the subset $\mbH \subseteq \mbX$.
By construction of $\mbH$,
$-\nabla_{\nat}^2\; \widetilde\ell(\nat;\, \bx)$ is invertible for all $(\nat,\, \bx) \in \mB_{\infty}(\truth,\, \epsilon^\star) \times \mbH$.

\s

{\bf Characterizing the strict concavity of $\widetilde{\ell}(\nat;\, \bx)$ provided $\bx \in \mbH$.} 
Observe that
\beno
\widetilde{\ell}(\nat;\, \bx)
\= \dsum_{i=1}^{M} \, \left(\langle \nat, s(x_i,\, \bx_{-i}) \rangle - \psi_i(\nat;\, \bx_{-i}) \right), 
\ee 
where 
\be
\label{eq:psi_cond}
\psi_i(\nat;\, \bx_{-i})
\= \log\, \dsum_{x_i=0}^{1}\, \exp\left(\langle\nat,\, s(x_i,\, \bx_{-i})\rangle\right),
& i \in \{1, \ldots, M\}.
\ee
We establish strict concavity of $\widetilde{\ell}(\, \cdot\, ;\, \bx)$ by demonstrating that at least one of the $\psi_i(\, \cdot\, ;\, \bx_{-i})$ is strictly convex on $\Nat$.  
H\"older's inequality shows that,
for any $\bx_{-i} \in \{0,1\}^{M-1}$,
the function $\psi_i(\nat;\, \bx_{-i})$ is a convex function of $\nat \in \Nat$ \citep[][Theorem 1.13, p.\ 19]{Br86}:
\be
\label{psi.convex}
\psi_i(\lambda \, \nat_1 + (1- \lambda) \, \nat_2;\, \bx_{-i})
\,\leq\, \lambda\, \psi_i(\nat_1;\, \bx_{-i}) + (1 - \lambda)\, \psi_i(\nat_2;\, \bx_{-i}),
\ee
where $\lambda \in (0,\, 1)$ and $(\nat_1,\, \nat_2) \in \Nat \times \Nat$.
The inequality \eqref{psi.convex} is an equality if and only if
\be
\label{cond0}
\exp\left(\langle\nat_1,\, s(x_{i},\, \bx_{-i})\rangle\right)
\propto \exp\left(\langle\nat_2,\, s(x_i,\, \bx_{-i})\rangle\right)
\mbox{ for all } x_i \in \{0,\, 1\}.  
\ee 
Upon taking the natural logarithm on both sides, 
condition \eqref{cond0} can be written as follows: 
\be 
\label{eq:prop_cond}
\langle\nat_1 - \nat_2,\, s(x_i,\, \bx_{-i})\rangle 
\= \log\, C
& \mbox{for all}& x_i \in \{0,\, 1\},
\ee 
where $C > 0$ is a constant.
For $\widetilde{\ell}(\, \cdot\, ;\, \bx)$ to be concave---but not strictly concave---condition \eqref{eq:prop_cond} 
must be satisfied for all $i \in \{1, \ldots, M\}$.  
By construction of $\mbH$,
the negative Hessian $-\nabla_{\nat}^2\; \widetilde{\ell}(\nat;\, \bx)$ 
is invertible for all $(\nat,\, \bx) \in \mB_{\infty}(\truth,\, \epsilon^\star) \times \mbH$,
implying that $\widetilde{\ell}(\, \cdot\, ;\, \bx)$ ($\bx \in \mbH$) 
is strictly concave on 
$\mB_{\infty}(\truth,\, \epsilon^\star)$ by virtue of the concavity of $\widetilde{\ell}(\, \cdot\, ;\, \bx)$ ($\bx \in \mbX$)
on $\Nat$,
discussed above.  
As a result, 
there exists an integer $i \in \{1, \ldots, M\}$ such that condition \eqref{eq:prop_cond} is not satisfied
for each $\bx \in \mbH$. 
We leverage the above characterization to extend the strict concavity of 
$\widetilde{\ell}(\,\cdot\,;\, \bx)$ on $\mB_{\infty}(\truth,\, \epsilon^\star)$ to $\Nat$.

\s

{\bf Extending the strict concavity of $\widetilde{\ell}(\nat;\, \bx)$ on $\mB_{\infty}(\truth,\, \epsilon^\star)$ to $\Nat$ provided $\bx \in \mbH$.}
We extend the above result from $\mB_{\infty}(\truth,\, \epsilon^\star)$ to $\Nat$ by using a proof by contradiction.
Consider any parameter vector $\nat \in \Nat \setminus \mB_{\infty}(\truth,\, \epsilon^\star)$.
Since $\Nat$ and $\mB_{\infty}(\truth,\, \epsilon^\star)$ are convex sets \citep[][Theorem 5.8, p.\ 154]{Br86},
there exists a real number $\lambda \in (0,\, 1)$ and a parameter vector $\dot{\nat} \in \mB_{\infty}(\truth,\, \epsilon^\star)$ such that 
\begin{equation}
\begin{array}{cccccc}
\label{lambda.cond}
\dot{\nat} 
\=  \lambda \, \nat + (1 - \lambda) \, \truth
&\in& \mB_{\infty}(\truth,\, \epsilon^\star).
\ee
Therefore,
$\nat$ can be represented as
\be
\nat 
\= \dfrac{\dot{\nat}}{\lambda} - \dfrac{1 - \lambda}{\lambda}\, \truth
&\in& \Nat \setminus \mB_{\infty}(\truth,\, \epsilon^\star).
\end{array} 
\end{equation}
Next,
consider any pair of parameter vectors $(\nat_1,\, \nat_2) \in \Nat \times (\Nat \setminus \mB_{\infty}(\truth,\, \epsilon^\star))$ and observe that both $\nat_1$ and $\nat_2$ can be represented in the form \eqref{lambda.cond}, 
and that the same $\lambda$ can be used to represent them (by choosing $\lambda$ small enough).
In other words,
there exists a real number $\lambda \in (0,\, 1)$ and a pair of parameter vectors $(\dot{\nat}_1,\, \dot{\nat}_2)\, \in\, \mB_{\infty}(\truth,\, \epsilon^\star) \times \mB_{\infty}(\truth,\, \epsilon^\star)$ such that 
\beno
\dot{\nat}_1
\= \lambda \, \nat_1 + (1 - \lambda) \, \truth
&\in& \mB_{\infty}(\truth,\, \epsilon^\star)\s\s
\\
\dot{\nat}_2 
\= \lambda \, \nat_2 + (1 - \lambda) \, \truth
&\in& \mB_{\infty}(\truth,\, \epsilon^\star),
\ee
so $\nat_1$ and $\nat_2$ can be represented as
\beno
\nat_1
\= \dfrac{\dot{\nat}_1}{\lambda} - \dfrac{1 - \lambda}{\lambda}\, \truth
&\in& \Nat\s\\
\nat_2
\= \dfrac{\dot{\nat}_2}{\lambda} - \dfrac{1 - \lambda}{\lambda}\, \truth
&\in& \Nat \setminus \mB_{\infty}(\truth,\, \epsilon^\star).
\end{array} 
\end{equation}
Assume that
\be
\label{condd1}
\langle \nat_1 - \nat_2,\, s(x_i,\, \bx_{-i}) \rangle
\= \log\, C
&\mbox{for all}& x_i \in \{0,\, 1\}\s
\\
&&&\mbox{and all}& i \in \{1, \dots, M\},
\ee 
which in turn implies that
\be
\label{impll1}
\langle\dot{\nat}_1 - \dot{\nat}_2,\, s(x_i,\, \bx_{-i})\rangle
\= \lambda\, \log\, C
&\mbox{for all}& x_i \in \{0,\, 1\}\s
\\
&&&\mbox{and all}& i \in \{1, \dots, M\},
\ee
because
\beno
\nat_1 - \nat_2 
\= \left(\dfrac{\dot{\nat}_1}{\lambda} - \dfrac{1 - \lambda}{\lambda}\, \truth\right) - \left(\dfrac{\dot{\nat}_2}{\lambda} - \dfrac{1 - \lambda}{\lambda} \, \truth\right)
\= \dfrac{1}{\lambda} \, (\dot{\nat}_1 - \dot{\nat}_2).
\ee
The conclusion \eqref{impll1} contradicts the strict concavity of $\widetilde\ell(\, \cdot\, ;\, \bx)$ on $\mB_{\infty}(\truth,\, \epsilon^\star)$,
because both $\dot{\nat}_1$ and $\dot{\nat}_2$ are elements of $\mB_{\infty}(\truth,\, \epsilon^\star)$
and $\lambda > 0$ is a constant independent of $x_i \in \{0, 1\}$. 
Therefore,
the assumption \eqref{condd1} cannot be satisfied,
hence $\widetilde\ell(\, \cdot\, ;\, \bx)$ is strictly concave on $\Nat$ provided $\bx \in \mbH$.
\qed

}

\vspace{-.25cm}

\blue{
\begin{lemma}
\label{uc.lemma}
Under the assumptions of Theorem \ref{thm:mple_consistency},
for any $\nat \in \Nat$ and $t > 0$,
\beno
\mbP\left(\norm{\bg(\nat;\, \bX) - \bg(\nat)}_{\infty} \,\geq\, t\right)
\,\leq\, 2\, \exp\left( - \dfrac{t^2}{8 \, (1+D_N)^2 \, \mnorm{\mD_N(\truth)}_2^2 \, \Psi_N^2} + \log\, p \right), 
\ee
where $\bg(\nat; \bX)$ and $\bg(\nat)$ are defined by 
\[
\begin{array}{cll}
\bg(\nat; \bX) 
&\coloneqq& \nabla_{\nat}\; \widetilde\ell(\nat; \bX)\s
\\
\bg(\nat) 
&\coloneqq& \mbE\; \nabla_{\nat}\; \widetilde\ell(\nat; \bX),
\end{array}
\]
while $D_N \geq 0$,\,
$\mnorm{\mD_N(\truth)}_2 \geq 1$,\,
and $\Psi_N > 0$ provided $N$ is large enough.
\end{lemma}
}

\s\s

\llproof \ref{uc.lemma}.
We prove Lemma \ref{uc.lemma} by leveraging concentration results of \citet{Chetal07} along with conditional independence properties of models with factorization properties of the form \eqref{eq:factorization}.
 
Consider any $\nat \in \Nat$ and any $\bx \in \mbX$.
By definition,
\beno
\widetilde\ell(\nat;\, \bx)
&\coloneqq& \dsum_{i=1}^{M} \log \, \mbP_{\nat}(X_i = x_i \,|\, \bX_{-i} = \bx_{-i}),
\ee
which implies that
\beno
\bg(\nat;\, \bx)
&\coloneqq& \nabla_{\nat}\; \widetilde\ell(\nat;\, \bx) 
\= \dsum_{i=1}^M \nabla_{\nat} \log \, \mbP_{\nat}(X_i = x_i \,|\, \bX_{-i} = \bx_{-i}).
\ee
Observe that
\be
\label{eq:cond_prob_score}
\nabla_{\nat} \log \, \mbP_{\nat}(X_i = x_i \,|\, \bX_{-i} = \bx_{-i})
\= s(\bx) - \mbE_{\nat, \bx_{-i}} \, s(\bX),
\ee
where $\mbE_{\nat, \bx_{-i}}$ denotes the expectation with respect to the conditional probability distribution of $X_{i}$ given $\bX_{-i} = \bx_{-i}$.
\hide{
The result in \eqref{eq:cond_prob_score} follows from exponential-family properties \citep{Br86},
because the conditional distributions of $X_i$ given $\bX_{-i} = \bx_{-i}$ are Bernoulli distributions ($i = 1, \ldots, M$),
and Bernoulli distributions are exponential-family distributions.  
}
The result in \eqref{eq:cond_prob_score} follows from exponential-family properties \citep{Br86},
because  the conditional distribution of $X_i$ given $\bX_{-i} = \bx_{-i}$ 
is an exponential-family distribution 
with sufficient statistic vector $s(\bx)$ and natural parameter vector $\nat$. 

We are interested in events of the form
\be
\label{event.interest}
\left\{
\norm{\bg(\nat;\, \bX) - \bg(\nat)}_{\infty} \,\geq\, t\right\},
&& t > 0,
& \nat \in \Nat. 
\ee
To bound the probabilities of events of the form \eqref{event.interest},
we leverage concentration results of \citet{Chetal07}.
Theorem 1 of \citet{Chetal07} states that,
for each $k \in \{1, \ldots, p\}$ and $t > 0$,
\beno
\mbP\left(\left| g_k(\nat, \bX) - \mbE \, g_k(\nat, \bX) \right| \,\geq\, t \right)
&\leq& 2 \, \exp\left(- \dfrac{2 \, t^2}{\norm{\bDelta_k}_2^2\; \mnorm{\mD_N(\truth)}_2^2} \right),
\ee
where $\mD_N(\truth)$ is defined in Section \ref{sec:background} and $\bDelta_k \in [0,\, \infty)^{M}$ is defined by  
\beno
\Delta_{k,i} 
\;\coloneqq\; \max\limits_{(\bx, \bx^\prime)\, \in\, \mbX\times\mbX:\;\, x_l = x_l^\prime \mbox{\footnotesize\, for all } l \neq i}\; |g_k(\nat; \bx) - g_k(\nat; \bx^\prime)|,  
& i \in \{1, \ldots, M\}. 
\ee
We bound the probability of event \eqref{event.interest} by bounding
\beno
\norm{\bDelta_k}_2^2
&=& \dsum_{i=1}^{M} \left(\, \max\limits_{(\bx, \bx^\prime)\, \in\, \mbX\times\mbX:\; x_l = x_l^\prime \mbox{\footnotesize\, for all } l \neq i}\; |g_k(\nat; \bx) - g_k(\nat; \bx^\prime)| \right)^2.
\ee
Consider any $i \in \{1, \dots, M\}$ and any $(\bx, \bx^\prime)\, \in\, \mbX \times \mbX$ such that $x_i = 0$ and $x_i^\prime = 1$ while $x_l = x_l^\prime$ for all $l \neq i$.
Write
\beno
\bg(\nat;\, \bx) - \bg(\nat;\, \bx^\prime)
&=& \dsum_{j=1}^{M} \nabla_{\nat} \, \lambda_j(\nat; \bx, \bx^\prime),  
\ee
where 
\beno
\lambda_j(\nat; \bx, \bx^\prime)
&\coloneqq&
\log \, 
\dfrac{\mbP_{\nat}(X_j = x_j \,|\, \bX_{-j} = \bx_{-j})}
{\mbP_{\nat}(X_j = x_j^\prime \,|\, \bX_{-j} = \bx_{-j}^\prime)},
&& j \in \{1, \ldots, M\}. 
\ee
By definition,
for any given $j \in \{1, \dots, M\}$,\,
the set $\mmG_{j}\; \subseteq\; \{1, \dots, M\}\, \setminus\, \{j\}$ is the smallest subset of indices such that
\be
\label{eq:cond_prop_conc}
X_j &\orth& \bX_{\{1,\ldots,M\} \setminus (\{j\} \,\cup\, \mmG_j)} \mid \bX_{\mmG_j}. 
\ee
Therefore,
for all $j \in \{1, \ldots, M\} \setminus (\{i\}\, \cup\, \mmG_i)$,
the conditional probability mass function of $X_j$ is unaffected by $X_i$,
so \eqref{eq:cond_prop_conc} implies that 
\beno
\mbP_{\nat}(X_j = x_j \,|\, \bX_{-j} = \bx_{-j})
\= \mbP_{\nat}(X_j = x_j^\prime \,|\, \bX_{-j} = \bx_{-j}^\prime),
\ee
which in turn implies,
for all $j \in \{1, \ldots, M\} \setminus (\{i\}\, \cup\, \mmG_i)$,
that 
\beno
\lambda_j(\nat; \bx, \bx^\prime) 
&\coloneqq& \log\, \dfrac{\mbP_{\nat}(X_j = x_j \,|\, \bX_{-j} = \bx_{-j})}
{\mbP_{\nat}(X_j = x_j^\prime \,|\, \bX_{-j} = \bx_{-j}^\prime)}
&=& 0,
\ee
noting that $x_l = x_l^\prime$ for all $l \neq i$. 
As a result,
\beno
\bg(\nat;\, \bx) - \bg(\nat;\, \bx^\prime)
\= \dsum_{j=1}^{M} \nabla_{\nat} \, \lambda_j(\nat; \bx, \bx^\prime) 
\= \dsum_{j\, \in\, \{i\}\, \cup\, \mmG_i} \nabla_{\nat} \, \lambda_j(\nat; \bx, \bx^\prime). 
\ee
The triangle inequality 
and \eqref{eq:cond_prob_score} 
imply, 
for each $k \in \{1, \ldots, p\}$,
that  
\beno
&& |g_k(\nat;\, \bx) - g_k(\nat;\, \bx^\prime)|\s\s
\\
\lte \dsum_{j\, \in\, \{i\}\, \cup\, \mmG_i}\, \left( \left| s_k(\bx) - s_k(\bx^\prime) \right| + \left| \mbE_{\nat, \bx_{-j}} \, s_k(\bX) - \mbE_{\nat, \bx_{-j}^\prime} \, s_k(\bX) \right| \right).
\ee 
We bound the terms of the above sum one by one.\s

{\bf Bounding $|s_k(\bx) - s_k(\bx^\prime)|$.}
Consider any $i \in \{1, \dots, M\}$ and any $(\bx, \bx^\prime)\, \in\, \mbX \times \mbX$ such that $x_i = 0$ and $x_i^\prime = 1$ while $x_l = x_l^\prime$ for all $l \neq i$.

By definition, 
\beno
\Xi_{k,i}
&\coloneqq& \max\limits_{(\bx,\, \bx^\prime)\, \in\, \mbX \times \mbX:\;\, x_l = x_l^\prime \mbox{\footnotesize\, for all } l \neq i}\;
|s_k(\bx) - s_k(\bx^\prime)|,
\ee
providing the following bound:
\beno
\left|s_k(\bx) - s_k(\bx^\prime) \right| 
&\leq& \Xi_{k,i},
\ee 
provided $(\bx, \bx^\prime) \in \mbX \times \mbX$ 
satisfies 
$x_l = x_l^\prime$ for all $l \neq i$ with $x_i = 0$ and $x_i^\prime = 1$.

\s

{\bf Bounding $|\mbE_{\nat, \bx_{-j}} \, s_k(\bX) - \mbE_{\nat, \bx_{-j}^\prime} \, s_k(\bX)|$.}
We take advantage of the coupling argument in Section 2.1 of \citet{Chetal07} to bound deviations of conditional expectations.

Consider any $i \in \{1, \dots, M\}$ and any $(\bx, \bx^\prime)\, \in\, \mbX \times \mbX$ such that $x_i = 0$ and $x_i^\prime = 1$ while $x_l = x_l^\prime$ for all $l \neq i$.
Define
\beno
\mbP_{j,\nat,\bx_{-j}}(a) 
&\coloneqq& \mbP_{\nat}(X_j = a \,|\, \bX_{-j} = \bx_{-j}),
&& a \in \{0, 1\}.
\ee
Let $(\bX^\star, \bX^{\star\star}) \in \{0,1\}^M \times \{0,1\}^M$ be an optimal coupling of the conditional probability mass functions $\mbP_{j,\nat,\bx_{-j}}$ and $\mbP_{j,\nat,\bx^\prime_{-j}}$
such that
\bi
\item the marginal probability mass function of $X_j^\star$ is $\mbP_{j,\nat,\bx_{-j}}$;\s
\item the marginal probability mass function of $X_j^{\star\star}$ is $\mbP_{j,\nat,\bx^\prime_{-j}}$;\s
\item the coupling ensures the  following events occur with probability 1: \vspace{.1cm} 
\bi
\item $\{X_{i}^\star = x_i = 0\}$,\vspace{.1cm}
\item $\{X_{i}^{\star\star} = x_i^\prime = 1\}$,\vspace{.1cm}
\item $\{X^\star_l = X_l^{\star\star} = x_l\}$ for all $l \in \{1, \ldots, M\}\, \setminus\, \{i,\, j\}$.\vspace{.1cm}
\ei
\ei
An optimal coupling is guaranteed to exist,
but it may not be unique \citep[][pp.\ 99--107]{Li02}. 
That said,
any optimal coupling will suffice.
We denote the joint probability mass function of $(\bX^\star, \bX^{\star\star})$ by $\mbT_{j,\nat,\bx_{-j},\bx_{-j}^\prime}$.

An important property of the optimal coupling $\mbT_{j,\nat,\bx_{-j},\bx_{-j}^\prime}$ is that,
for all\break 
$j \in \{1, \ldots, M\}\, \setminus\, (\{i\}\, \cup\, \mathfrak{N}_i)$,
\be
\label{tv.coupling}
\mbT_{j,\nat,\bx_{-j},\bx_{-j}^\prime}(X_j^\star \neq X_j^{\star\star})
\= \norm{\mbP_{j,\nat,\bx_{-j}} - \mbP_{j,\nat,\bx^\prime_{-j}}}_{\tv}
\= 0,
\ee
because
\bi
\item $\mbT_{j,\nat,\bx_{-j},\bx_{-j}^\prime}(X_j^\star \neq X_j^{\star\star})
= \norm{\mbP_{j,\nat,\bx_{-j}} - \mbP_{j,\nat,\bx^\prime_{-j}}}_{\tv}$, 
by virtue of 
$\mbT_{j,\nat,\bx_{-j},\bx_{-j}^\prime}$ being an optimal coupling of 
$\mbP_{j,\nat,\bx_{-j}}$ and $\mbP_{j,\nat,\bx^\prime_{-j}}$; \s 
\item the conditional independence of $X_i$ and $X_j$ implies 
$\mbP_{j,\nat,\bx_{-j}} = \mbP_{j,\nat,\bx^\prime_{-j}}$ 
for all $(\bx, \bx^\prime)\, \in\, \mbX \times \mbX$ such that $x_l = x_l^\prime$ for all $l \neq i$,
further implying that $\norm{\mbP_{j,\nat,\bx_{-j}} - \mbP_{j,\nat,\bx^\prime_{-j}}}_{\tv} = 0$.
\ei
By construction of the coupling $\mbT_{j,\nat,\bx_{-j},\bx_{-j}^\prime}$,
we can write
\beno
\mbE_{\nat, \bx_{-j}} \, s_k(\bX) - \mbE_{\nat, \bx_{-j}^\prime} \, s_k(\bX)
\,=\, \mbE_{\mbT_{j,\nat,\bx_{-j},\bx_{-j}^\prime}} \, s_k(\bX^\star) - \mbE_{\mbT_{j,\nat,\bx_{-j},\bx_{-j}^\prime}}\, s_k(\bX^{\star\star})\s
\\
\hspace{5.1cm}=\; \mbE_{\mbT_{j,\nat,\bx_{-j},\bx_{-j}^\prime}}(s_k(\bX^\star) - s_k(\bX^{\star\star})),
\ee 
where $\mbE_{\mbT_{j,\nat,\bx_{-j},\bx_{-j}^\prime}}$ 
is the expectation operator under the coupling probability distribution $\mbT_{j,\nat,\bx_{-j},\bx_{-j}^\prime}$
of $(\bX^\star, \bX^{\star\star})$.  
Taking advantage of the telescoping identity on page 205 and the bounding argument on page 206 of \citet{Chetal07} gives rise to the bound
\beno
\left|\mbE_{\nat, \bx_{-j}} \, s_k(\bX) - \mbE_{\nat, \bx_{-j}^\prime} \, s_k(\bX) \right|
\= \left|\mbE_{\mbT_{j,\nat,\bx_{-j},\bx_{-j}^\prime}}(s_k(\bX^\star) - s_k(\bX^{\star\star}))\right|\s
\\
\lte \dsum_{l=1}^{M}\, \Xi_{k,l}\; \mbT_{j,\nat,\bx_{-j},\bx_{-j}^\prime}(X_l^\star \neq X_l^{\star\star}).
\ee
The construction of the coupling $\mbT_{j,\nat,\bx_{-j},\bx_{-j}^\prime}$ implies that
\beno
\left|\mbE_{\nat, \bx_{-j}} \, s_k(\bX) - \mbE_{\nat, \bx_{-j}^\prime} \, s_k(\bX) \right|
\lte \dsum_{l=1}^{M}\, \Xi_{k,l}\; \mbT_{j,\nat,\bx_{-j},\bx_{-j}^\prime}(X_l^\star \neq X_l^{\star\star})\s
\\
\lte 
\begin{cases}
\Xi_{k,i}
& \mbox{if } i = j\s
\\
\Xi_{k,i} + \Xi_{k,j} 
& \mbox{if } i \neq j \mbox{ and } j\, \in\, \mmG_i\s
\\
0 
& \mbox{if } i \neq j \mbox{ and } j\, \not\in\, \mmG_i.
\end{cases}
\ee

\s

{\bf Collecting terms.}
Upon collecting terms,
we obtain the bounds
\beno
&& |g_k(\nat;\, \bx) - g_k(\nat;\, \bx^\prime)|\s\s
\\
\lte \dsum_{j\, \in\, \{i\}\, \cup\, \mmG_i}\, \left( \left| s_k(\bx) - s_k(\bx^\prime) \right| + \left| \mbE_{\nat, \bx_{-j}} \, s_k(\bX) - \mbE_{\nat, \bx_{-j}^\prime} \, s_k(\bX) \right| \right)\s\s
\\
\hide{
&\leq& 2\; \Xi_{k,i} + \dsum_{j\, \in\, \mmG_i}\, \left(\left|s_k(\bx) - s_k(\bx^\prime)\right| + \left|\mbE_{\nat, \bx_{-j}} \, s_k(\bX) - \mbE_{\nat, \bx_{-j}^\prime}\, s_k(\bX) \right|\right)\s\s
\\
}
&\leq& 2\; \Xi_{k,i} + \dsum_{j\, \in\, \mmG_i}\, \left(\Xi_{k,i} + (\Xi_{k,i} + \Xi_{k,j})\right)\s
\\
&\leq& 2\, \left(\Xi_{k,i} + \dsum_{j\, \in\, \mmG_i} (\Xi_{k,i} + \Xi_{k,j})\right)
\ee
and
\[
\begin{array}{ccc}
&& \max\limits_{(\bx, \bx^\prime)\, \in\, \mbX\times\mbX:\; x_l = x_l^\prime\mbox{\footnotesize\, for all } l \neq i}\, |g_k(\nat;\, \bx) - g_k(\nat;\, \bx^\prime)|\s\s
\\
\lte 2 \left(\Xi_{k,i} + \dsum_{j\, \in\, \mmG_i} (\Xi_{k,i} + \Xi_{k,j})\right).  
\end{array}
\]
The Cauchy–Schwarz inequality implies that
\beno
\left(\max\limits_{(\bx, \bx^\prime)\, \in\, \mbX\times\mbX:\;\, x_l = x_l^\prime \mbox{\footnotesize\, for all } l \neq i}\; |g_k(\nat;\, \bx) - g_k(\nat;\, \bx^\prime)|\right)^2\s\s
\\
\leq\; 4 \, \left(\Xi_{k,i} + \dsum_{j\, \in\, \mmG_i} (\Xi_{k,i} + \Xi_{k,j})\right)^2\s\s
\\
\leq\; 4\, (1 + 2\, D_N) \left(\Xi_{k,i}^2 + \dsum_{j\, \in\, \mmG_i} (\Xi_{k,i}^2 + \Xi_{k,j}^2)\right)\s\s
\\
\leq\; 8\, (1 + D_N) \left(\Xi_{k,i}^2 + \dsum_{j\, \in\, \mmG_i} (\Xi_{k,i}^2 + \Xi_{k,j}^2)\right)\s\s
\\
\leq\; 8\, (1 + D_N) \left((1 + D_N)\, \Xi_{k,i}^2 + \dsum_{j\, \in\, \mmG_i}\, \Xi_{k,j}^2 \right), 
\ee
using $D_N \coloneqq \max\{|\mmG_1|, \dots, |\mmG_M|\}$. 
We hence obtain
\beno
\norm{\bDelta_k}_2^2 
&\coloneqq& \dsum_{i=1}^{M} \, 
\left(\, \max\limits_{(\bx, \bx^\prime)\, \in\, \mbX\times\mbX:\;\, x_l = x_l^\prime \mbox{\footnotesize\, for all } l \neq i}\; |g_k(\nat; \bx) - g_k(\nat; \bx^\prime)| \right)^2\s\s\\
&\leq& \dsum_{i=1}^{M} \, 8 \, (1 + D_N) \left( (1 + D_N) \, \Xi_{k,i}^2 + \dsum_{j\, \in\, \mmG_i} \, \Xi_{k,j}^2 \right) \s\s\\
\= 8\, (1 + D_N)^2 \, \norm{\bXi_k}_2^2 + 8 \, (1 + D_N) \, \dsum_{j=1}^{M} \, \Xi_{k,j}^2 \, \dsum_{i=1}^{M} \, \one(j\, \in\, \mmG_i). 
\ee
To bound the second term on the right-hand side,
note that edge variable $X_j$ can be in the dependence neighborhoods 
$\mmG_i$ ($i = 1, \ldots, M$) of at most $D_N \coloneqq \max\{|\mmG_1|, \dots, |\mmG_M|\}$ other edge variables $X_i$,
which implies that
\beno
8 \, (1 + D_N) \, \dsum_{j=1}^{M} \, \Xi_{k,j}^2 \, \dsum_{i=1}^{M} \, \one(j\, \in\, \mmG_i) 
&\leq& 8 \, (1 + D_N) \, \dsum_{j=1}^{M} \, \Xi_{k,j}^2\, D_N\s\s\\ 
&\leq& 8 \, (1 + D_N)^2 \; \norm{\bXi_k}_2^2.  
\ee
We hence arrive at the following bound on $\norm{\bDelta_k}_2^2$: 
\beno
\norm{\bDelta_k}_2^2 
\lte 8 \, (1 + D_N)^2 \, \norm{\bXi_k}_2^2 + 8 \, (1 + D_N) \, \dsum_{j=1}^{M} \, \Xi_{k,j}^2 \, \dsum_{i=1}^{M} \, \one(j\, \in\, \mmG_i)\s
\\
\lte 16 \, (1 + D_N)^2 \, \norm{\bXi_k}_2^2 \s\s
\\
\lte 16 \, (1 + D_N)^2 \, \Psi_N^2,
\ee
where $\Psi_N \coloneqq \max_{1 \leq k \leq p} \, \norm{\bXi_k}_2$. \s

\s

{\bf Concentration result.} By applying Theorem 1 of \citet{Chetal07} to each coordinate $g_k(\nat;\, \bX) - g_k(\nat)$ of\, $\bg(\nat;\, \bX) - \bg(\nat)$ ($k = 1, \dots, p$) using the above bound on $\norm{\bDelta_k}_2$,
we have,
for all $t > 0$,
\beno
\mbP\left( \left|g_k(\nat;\, \bX) - g_k(\nat) \right| \geq  t\right)
\lte 2 \, \exp\left( - \dfrac{2 \, t^2}{16 \, (1+D_N)^2 \, \mnorm{\mD_N(\truth)}_2^2 \, \Psi_N^2} \right)\s
\\
&=& 2 \, \exp\left( - \dfrac{t^2}{8\, (1+D_N)^2 \, \mnorm{\mD_N(\truth)}_2^2 \, \Psi_N^2} \right),
\ee
where $D_N \geq 0$,\,
$\mnorm{\mD_N(\truth)}_2 \geq 1$,\, 
and $\Psi_N > 0$ provided $N$ is large enough.
A union bound over the $p$ coordinates $g_k(\nat;\, \bX) - g_k(\nat)$ of\, $\bg(\nat;\, \bX) - \bg(\nat)$ gives rise to the bound
\beno
\mbP\left(\norm{\bg(\nat;\, \bX) - \bg(\nat)}_{\infty} \geq t\right) 
\,\leq\, 2 \, \exp\left( - \dfrac{t^2}{8 \, (1+D_N)^2 \, \mnorm{\mD_N(\truth)}_2^2 \,  \Psi_N^2} + \log \, p \right). 
\ee
\hide{
Since the above bound does not depend on $\nat \in \Nat$,
we conclude that
\beno
&& \mbP\left(\su_{\nat\, \in\, \Nat}\; \norm{\bg(\nat;\, \bX) - \bg(\nat)}_{\infty} \,\geq\, t\right)\s
\\
&\leq& 2\, \exp\left( - \dfrac{t^2}{8 \, (1+D_N)^2 \, \mnorm{\mD_N(\truth)}_2^2 \, \Psi_N^2} + \log\, p \right).
\ee
}
\qed

\s

\begin{lemma}
\label{lem:concave_pl}
The function $\mbE\; \widetilde\ell(\nat;\, \bX)$ is a strictly concave function on the convex set $\Nat = \mR^p$.
In addition,
the data-generating parameter vector $\truth \in \Nat$ maximizes the expected loglikelihood and pseudo-loglikelihood functions:
\beno
\truth
\= \argmax\limits_{\nat\, \in\, \Nat}\; \mbE\; \ell(\nat;\, \bX)  
\= \argmax\limits_{\nat\, \in\, \Nat}\; \mbE\; \widetilde\ell(\nat;\, \bX).
\ee
\end{lemma}

\llproof \ref{lem:concave_pl}.
Section \ref{sec:models} of the manuscript shows that the family of densities $\{f_{\nat},\, \nat\, \in\, \Nat\}$ parameterized by \eqref{eq:factorization} and \eqref{parameterization} is an exponential family of densities.
We take advantage of the properties of exponential families \citep{Br86} to prove Lemma \ref{lem:concave_pl},
and divide the proof into three parts:
\bi
\item[I.] $\mbE \; \widetilde\ell(\nat;\, \bX)$ is a strictly concave function on the convex set $\Nat$.\s 
\item[II.] $\truth$ is the unique maximizer of $\mbE\; \ell(\nat;\, \bX)$.\s
\item[III.] $\truth$ is the unique maximizer of $\mbE\; \widetilde\ell(\nat;\, \bX)$.
\ei

\s

\hide{
{\bf I.\ $\Nat$ is a convex set.}
By definition of $\Nat = \{\nat \in \mR^p:\, \psi(\nat) < \infty\}$ and H\"older's inequality,
$\Nat$ is a convex set \citep[][Theorem 1.13, p.\ 19]{Br86}.

\s

}

{\bf I.\ $\mbE \; \widetilde\ell(\nat;\, \bX)$ is a strictly concave function on the convex set $\Nat$.}
Let $\bx$ be an observation of a random graph $\bX$ with dependent edges.
Then,
by definition,
\beno
\widetilde\ell(\nat;\, \bx)
\= \dsum_{i=1}^M \widetilde\ell_i(\nat;\, \bx),
\ee
where
\beno
\widetilde\ell_i(\nat;\, \bx)
\= \langle \nat, \, s(\bx) \rangle - \psi_{i}(\nat;\; \bx_{-i})
\ee
and
\beno
\psi_{i}(\nat;\, \bx_{-i}) 
&=& \log \, \dsum_{x_i = 0}^{1} \, \exp\left( \langle \nat, \, s(\bx_{-i}, \, x_i) \rangle \right). 
\ee 
The set $\Nat$ is a convex set \citep[][Theorem 1.13, p.\ 19]{Br86}. 
We first show that $\mbE \; \widetilde\ell(\nat;\, \bX) = \sum_{i=1}^M \mbE \; \widetilde\ell_i(\nat;\, \bX)$ is a concave function on the convex set $\Nat$ 
by proving that the functions $\mbE \; \widetilde\ell_i(\nat;\, \bX)$ are concave on $\Nat$.
Observe that the functions $\mbE \; \widetilde\ell_i(\nat;\, \bX)$ are concave provided the functions $\widetilde\ell_i(\nat;\, \bx)$ are concave for all $\bx \in \mbX$.
To show that the functions $\widetilde\ell_i(\nat;\, \bx)$ are concave for all $\bx \in \mbX$,
consider any $i \in \{1, \dots, M\}$,
any $x_i \in \{0, 1\}$,
and any $\bx_{-i} \in \{0, 1\}^{M-1}$.
Each $\widetilde\ell_i(\nat;\, \bx)$ consists of two terms.
The first term,
$\langle \nat, \, s(\bx) \rangle$,
is a linear function of $\nat$,
so $\widetilde\ell_i(\nat;\, \bx)$ is a concave function of $\nat$ if the second term,
$\psi_{i}(\nat;\; \bx_{-i})$, 
is a convex function of $\nat$.
Consider any $(\nat^{(1)}, \nat^{(2)}) \in \Nat \times \Nat$ 
and any $\lambda \in (0, 1)$.
Then,
by H\"older's inequality,
\beno
\hide{
&& 
}
\psi_{i}\left(\lambda\; \nat^{(1)} + (1 - \lambda)\; \nat^{(2)};\; \bx_{-i}\right)
\hide{
\s\s
\\
&=& \log\left(\dsum_{x_{i}^\prime = 0}^{1} \, \exp\left(\left\langle\lambda\; \nat^{(1)} + (1 - \lambda)\; \nat^{(2)},\; s(\bx_{-i},\, x_{i}^\prime)\right\rangle\right)\right)\s\s
\\
\= \log\left( \dsum_{x_{i}^\prime = 0}^{1} \, \exp\left(\left\langle\nat^{(1)},\; s(\bx_{-i},\, x_{i}^\prime)\right\rangle\right)^{\lambda}\; \exp\left(\left\langle\nat^{(2)},\; s(\bx_{-i},\, x_{i}^\prime)\right\rangle\right)^{1 - \lambda}\right)\s\s
\\
\lte \log\left(\dsum_{x_{i}^\prime = 0}^{1} \exp\left(\left\langle\nat^{(1)},\; s(\bx_{-i},\, x_{i}^\prime)\right\rangle\right)\right)^{\lambda}\s\s
\\
&+& \log \left(\dsum_{x_{i}^\prime = 0}^{1} \exp\left(\left\langle\nat^{(2)},\; s(\bx_{-i},\, x_{i}^\prime)\right\rangle\right)\right)^{1-\lambda}\s\s
\\
&=& 
}
\leq \lambda\, \psi_{i}\left(\nat^{(1)};\, \bx_{-i}\right) + (1 - \lambda)\, \psi_{i}\left(\nat^{(2)};\, \bx_{-i}\right).
\ee
As a consequence,
for any $\bx_{-i} \in \{0, 1\}^{M-1}$,\;
$\psi_{i}(\nat;\, \bx_{-i})$ is a convex function on $\Nat$.
Hence,
for all $\bx \in \mbX$,\;
$\widetilde\ell(\nat;\, \bx)$ is a concave function on $\Nat$,
and so is $\mbE \; \widetilde\ell(\nat;\, \bX)$ as a finite sum of concave functions on $\Nat$. 

Second, 
we prove by contradiction that 
\beno
\mbE \; \widetilde\ell(\nat;\, \bX) 
\= \dsum_{i=1}^{M} \, \left[ \langle \nat, \, \mbE \, s(\bX) \rangle - \mbE \, \psi_{i}(\nat;\; \bX_{-i}) \right]
\ee
is a strictly concave function on $\Nat$. 
As discussed above, 
the first term $\langle \nat, \, \mbE \, s(\bX) \rangle$ is concave in $\nat$. 
Therefore, 
the strict concavity of $\mbE \; \widetilde\ell(\nat;\, \bX)$ 
is determined by the terms $\mbE \, \psi_{i}(\nat;\; \bX_{-i})$ ($i \in \{1, \ldots, M\}$). 
Suppose that there does not exist any $i^\star \in \{1, \dots, M\}$ such that $\mbE \; \psi_{i^\star}(\nat;\, \bX_{-i^\star})$ is strictly convex on $\Nat$.
Then,
there exists $(\nat^{(1)}, \nat^{(2)}) \in \Nat \times \Nat$ 
such that, 
for all $i \in \{1, \dots, M\}$,
all $\bx_{-i} \in \{0, 1\}^{M-1}$,
and all $x_{i} \in \{0, 1\}$,
\be
\label{eq:linearly_dependent} 
\exp\left(\left\langle \nat^{(1)},\; s(\bx_{-i},\, x_{i})\right\rangle\right)
\,\propto\, \exp\left(\left\langle \nat^{(2)},\; s(\bx_{-i},\, x_{i})\right\rangle\right),
\ee
as 
H\"older's inequality reduces to an equality if and only if  
\eqref{eq:linearly_dependent} holds, 
i.e.,
\beno
\psi_{i}\left(\lambda\; \nat^{(1)} + (1 - \lambda)\; \nat^{(2)};\; \bx_{-i}\right)
= \lambda\, \psi_{i}\left(\nat^{(1)};\, \bx_{-i}\right) + (1 - \lambda)\, \psi_{i}\left(\nat^{(2)};\, \bx_{-i}\right) 
\ee
if and only if \eqref{eq:linearly_dependent} holds. 
In other words,
for all $\bx \in \mbX$,
\be
\label{exf}
\exp\left(\left\langle \nat^{(1)}, \, s(\bx)\right\rangle\right)
\,\propto\, \exp\left(\left\langle \nat^{(2)}, \, s(\bx)\right\rangle\right).
\ee
The conclusion \eqref{exf} contradicts the assumption that the exponential family is minimal.
Therefore,
there exists $i^\star \in \{1, \dots, M\}$ such that $\mbE \; \psi_{i^\star}(\nat;\, \bX_{-i^\star})$ is strictly convex on $\Nat$,
which implies that $\mbE \; \widetilde\ell_{i^\star}(\nat;\, \bX)$ is strictly concave on $\Nat$,
and so is $\mbE \; \widetilde\ell(\nat;\, \bX) = \sum_{i=1}^M \mbE \; \widetilde\ell_i(\nat;\, \bX)$.

\s

{\bf II.\ $\truth$ is the unique maximizer of $\mbE\; \ell(\nat; \bX)$.}
Maximizing $\mbE\, \ell(\nat; \bX)$ is equivalent,
by Lemma \ref{lem:interchange},
to solving
\be
\label{solsol}
\nabla_{\nat}\; \mbE\; \ell(\nat;\, \bX)
\hide{
\= \nabla_{\nat}\, (\langle\nat,\, \mbE\, s(\bX)\rangle - \psi(\nat))\s
\\
}
\= \mbE\, s(\bX) - \mbE_{\nat}\, s(\bX)
\= \bm{0}.
\ee
The unique solution of \eqref{solsol} is $\truth \in \Nat = \mR^p$,
because $\mbE \equiv \mbE_{\truth}$. 
The fact that the solution is unique follows from the fact the map $\bmu: \Nat \mapsto \mbM$ defined by $\bmu(\nat) \coloneqq \mbE_{\nat}\, s(\bX)$ is one-to-one \citep[][Theorem 3.6, p.\ 74]{Br86}.
As a result,
$\truth \in \Nat = \mR^p$ is the unique maximizer of $\mbE\; \ell(\nat; \bX)$.

\s

{\bf III.\ $\truth$ is the unique maximizer of $\mbE\; \widetilde\ell(\nat;\, \bX)$.}
Observe that,
for any $\bx \in \mbX$,\;
$\widetilde\ell(\nat; \bx)$ is a sum of exponential-family loglikelihood functions,
because the conditional distributions of edge variables $X_i$ given $\bX_{-i} = \bx_{-i}$ 
($i = 1, \ldots, M$)
are exponential-family distributions with sufficient statistic vector $s(\bx)$ 
and natural parameter vector $\nat$. 
As a result,
$\widetilde\ell(\nat; \bx)$ is continuously differentiable on $\Nat$ for all $\bx \in \mbX$ \citep{Br86},
and so is $\mbE \; \widetilde\ell(\nat; \bX)$.
We then have
\beno
\bg(\nat)
&\coloneqq& \mbE \; \nabla_{\nat} \; \widetilde\ell(\nat; \bX) 
\= \mbE \, \dsum_{i=1}^M \, \left( s(\bX) - \mbE_{\nat, \bX_{-i}} \, s(\bX) \right),
\ee
where $\mbE_{\nat, \bx_{-i}}$ denotes the conditional expectation with respect to the conditional distribution of $X_i$ given $\bX_{-i} = \bx_{-i}$.
By the law of total expectation and the fact that $\mbE \equiv \mbE_{\truth}$, 
we have $\mbE \, \mbE_{\truth, \bX_{-i}}\, s(\bX) = \mbE \, s(\bX)$, 
which implies 
\be
\label{lte}
\bg(\truth)
\= \mbE \, \dsum_{i=1}^M \, \left( s(\bX) - \mbE_{\truth, \bX_{-i}} \, s(\bX) \right)
\= \bm{0}.
\ee
Thus,
a root of $\bg(\nat)$ exists,
and $\truth$ is a root of $\bg(\nat)$.
In addition,
$\mbE \; \widetilde\ell(\nat; \bX)$ is strictly concave on $\Nat$,
so $\truth$ is the unique root of $\bg(\nat)$.
As a consequence,
the maximizer of $\mbE \; \widetilde\ell(\nat; \bX)$ as a function of $\nat \in \Nat = \mR^p$ exists and is unique, 
and is given by $\truth \in \Nat = \mR^p$.
\qed

\s\s

\hide{ 

\begin{lemma}
\label{lem:one-to-one}
Let $\bg: \Nat \mapsto \mbR$ be any continuously differentiable function on the open and convex set $\Nat$.
If $\bg(\btheta)$ is strictly concave on $\Nat$, 
then its gradient $\nabla_{\btheta}\; \bg(\btheta)$ exists, 
is continuous, 
and is one-to-one.
\end{lemma}

\s 

\llproof \ref{lem:one-to-one}. 
The existence and continuity of $\nabla_{\btheta}\, \bg(\btheta)$ on $\Nat$ follow from the assumption that $\bg(\btheta)$ is continuously differentiable on the open and convex set $\Nat$.
We prove by contradiction that $\nabla_{\btheta}\; \bg(\btheta)$ is one-to-one on $\Nat$. 
Suppose that $\nabla_{\btheta}\; \bg(\btheta)$ is not one-to-one on $\Nat$,
that is,
there exists $(\btheta_1,\, \btheta_2) \in \Nat \times \Nat$ such that $\btheta_1 \neq \btheta_2$ and $\nabla_{\btheta}\; \bg(\btheta)\, |_{\btheta=\btheta_1} = \nabla_{\btheta}\; \bg(\btheta)\, |_{\btheta=\btheta_2}$. 
By the strict concavity of $\bg(\btheta)$ on $\Nat$,
\be
\label{eq1}
\langle \nabla_{\btheta}\; \bg(\btheta)\, |_{\btheta=\btheta_1}, \; \btheta_2 - \btheta_1 \rangle 
&>& \bg(\btheta_2) - \bg(\btheta_1)
\ee
and 
\be
\label{eq2}
\langle \nabla_{\btheta}\; \bg(\btheta)\, |_{\btheta=\btheta_2}, \; \btheta_1 - \btheta_2 \rangle 
&>& \bg(\btheta_1) - \bg(\btheta_2).
\ee
By multiplying both sides of \eqref{eq2} by $-1$,
we obtain 
\beno
\langle \nabla_{\btheta}\; \bg(\btheta)\, |_{\btheta=\btheta_2}, \; \btheta_2 - \btheta_1 \rangle 
&<& \bg(\btheta_2) - \bg(\btheta_1).
\ee
If $\nabla_{\btheta}\; g(\btheta)\, |_{\btheta=\btheta_1} = \nabla_{\btheta}\; g(\btheta)|_{\btheta=\btheta_2}$, 
then
\be
\label{eq3}
\langle \nabla_{\btheta}\; \bg(\btheta)\, |_{\btheta=\btheta_1}, \; \btheta_2 - \btheta_1 \rangle 
&<& \bg(\btheta_2) - \bg(\btheta_1).
\ee
The conclusion \eqref{eq3} contradicts \eqref{eq1},
so $\nabla_{\btheta}\, \bg(\btheta)$ is one-to-one on $\Nat$.
\qed

}

\vspace{.125cm}

\begin{lemma}
\label{lem:interchange} 
Under the assumptions of Theorem \ref{thm:mple_consistency},
\beno
\nabla_{\nat} \; \mbE\; \widetilde\ell(\nat; \bX)
\= \mbE\; \nabla_{\nat} \; \widetilde\ell(\nat; \bX)
&\mbox{for all}& \nat \in \Nat = \mR^p.
\ee
\end{lemma}

\vspace{-.25cm}

\llproof \ref{lem:interchange}.
We start with two observations.
First,
the exponential family introduced in Section \ref{general.parameterizations} of the manuscript is regular in the sense of \citet[][p.~2]{Br86},
because
\beno
\Nat \coloneqq \{\nat \in \mR^p:\, \psi(\nat) < \infty\} =  \mR^p
&\mbox{and}& \Nat = \mR^p \mbox{ is open}.
\ee
Second,
for any $\bx \in \mbX$,\;
$\widetilde\ell(\nat;\, \bx)$ is a sum of exponential-family loglikelihood functions,
because the conditional distribution of $X_i$ given $\bX_{-i} = \bx_{-i}$ 
is an exponential-family distribution  
with sufficient statistic vector $s(\bx)$ and natural parameter vector $\nat$.   
Thus,
for any $\bx \in \mbX$,\,
$\widetilde\ell(\,\cdot\,;\, \bx)$ is continuously differentiable on $\Nat$.

Consider $\widetilde\ell(\nat;\, \bx)$ as a function of $\bx \in \mbX$ for fixed $\nat \in \Nat$ and define
\beno
\bg(\nat;\, \bx)
&\coloneqq& \nabla_{\nat}\; \widetilde\ell(\nat;\, \bx)
\= \dsum_{i=1}^M\, \left(s(\bx) - \mbE_{\nat, \bx_{-i}}\, s(\bX)\right),
\ee
where $\mbE_{\nat, \bx_{-i}}$ denotes the expectation with respect to the conditional distribution of $X_i$ given $\bX_{-i} = \bx_{-i}$.
Here, 
$\bg(\nat;\, \bx)$ is considered as a function of $\bx \in \mbX$ for fixed $\nat \in \Nat$.
By the triangle inequality,
for each $k \in \{1, \ldots, p\}$, 
\beno
|g_k(\nat;\, \bx)|
&\leq& M\, |s_k(\bx)| + \dsum_{i=1}^M\, |\mbE_{\nat,\, \bx_{-i}}\, s_k(\bX)| 
&\eqqcolon& h_k(\bx),
\ee
where the dependence of $h_k(\bx)$ on $\nat$ is supressed.
Since the exponential family is regular in the sense of \citet[][p.~2]{Br86},
all moments of $s(\bX)$ exist \citep[Theorem 2.2, p.~34,][]{Br86},
implying that,
for all $k \in \{1, \ldots, p\}$,
$\mbE |s_k(\bX)| < \infty$ and $\mbE |\mbE_{\nat, \bX_{-i}} \, s_k(\bX)| < \infty$. 
As a result,
\beno
\mbE\, h_k(\bX) 
&=& M\, \mbE |s_k(\bX)| + \dsum_{i=1}^M\, \mbE |\mbE_{\nat, \bX_{-i}} \, s_k(\bX)|
&<& \infty.
\ee
Since $|g_k(\nat;\, \bx)| \leq h_k(\bx)$ for all $\bx \in \mbX$ and all $\nat \in \bTheta$ 
and $\mbE\, h_k(\bX) < \infty$,
Lebesgue's dominated convergence theorem implies that
\beno
\nabla_{\nat} \; \mbE\; \widetilde\ell(\nat; \bX)
\= \mbE\; \nabla_{\nat} \; \widetilde\ell(\nat; \bX)
&\mbox{for all}& \nat \in \Nat = \mR^p.
\ee
\qed

\s 

\section{Proofs of Corollaries 1--3}
\label{sec:conc}

We prove Corollaries \ref{cor:model_1}--\ref{cor:model_3} stated in Section \ref{sec:corollaries} of the manuscript,
using auxiliary results proved in Appendices \ref{sup-sec:scale_pl} and \ref{sup-sec:mrf_conc}.
To prove them,
it is convenient to return to the notation used in Section \ref{sec:models} of the manuscript,
denoting edge variables by $X_{i,j}$ ($\{i, j\} \subset \mN$).
Throughout Appendix \ref{sec:conc},
we assume that Models 2 and 3 satisfy one of the two following conditions.

\begin{description}[style=multiline, labelwidth=.75cm]
\item[\textcolor{black}{\namedlabel{sone}{S.1}}]
The subpopulations do not intersect ($\omega_1 = \omega_2 = 0$),\,
$\truth \in \mR^{N+1}$ satisfies condition \eqref{condition.theta0} with $\vartheta \in [0,\, 1/2 - \alpha)$ and $\alpha \in [0,\, 1/2)$,
and
\beno
D_N 
&=& O(\log N).
\ee
\item[\textcolor{black}{\namedlabel{stwo}{S.2}}]
The subpopulations intersect ($\omega_1 > 0$, $\omega_2 \geq 0$),\,
$\truth \in \mR^{N+1}$ satisfies condition \eqref{condition.theta0} with $\vartheta = 0$ and $\alpha \in [0,\, 1/2)$,
and
\beno
D_N
&=& o((\log(N / \log N)^{1/3}).
\ee
\end{description}
Condition \ref{sone} and \ref{stwo} ensure that the assumptions of Corollaries \ref{cor:model_2} and \ref{cor:model_3} are satisfied,
respectively.

\ccsproof \ref{cor:model_1}--\ref{cor:model_3}.
To prove Corollaries \ref{cor:model_1}--\ref{cor:model_3},
we bound
\beno
\widetilde\Phi_N(\truth)
&\coloneqq& \tllambdamin\,\, (1 + D_N)\,\, |\!|\!|\mD_N(\truth)|\!|\!|_2\,\, \Psi_N\; \sqrt{\log \max\{N,\, p\}}.
\ee
We first bound $\Psi_N$ and then prove Corollaries \ref{cor:model_1}--\ref{cor:model_3}.

\s\s 

{\bf Bounding $\Psi_N$.}
Recall the definition of $\Psi_N$:
For each $a \in \{1, \ldots, p\}$ and each pair of nodes $\{i, j\} \subset \mN$,
\beno
\Xi_{a,\{i,j\}}
&\coloneqq& \max\limits_{(\bx, \bx^\prime)\, \in\, \mbX \times \mbX:\,\;
x_{k,l} = x_{k,l}^\prime \mbox{\footnotesize\, for all } 
\{k,l\} \neq \{i,j\}}\;
|s_a(\bx) - s_a(\bx^\prime)|
\ee
and
\beno
\Psi_N
&\coloneqq& \max\limits_{1 \leq a \leq p} \, \norm{\bm\Xi_a}_2.
\ee
We show that 
$\Psi_N \leq \sqrt{N}$ under Model 1 and 
$\Psi_N \leq \norm{s_{N+1}}_{\lip}\; \sqrt{N}$ under Models 2 and 3 and bound $\norm{s_{N+1}}_{\lip}$,
where $\norm{s_{N+1}}_{\lip}$ is the Lipschitz coefficient of $s_{N+1}(\bX)$ with respect to the Hamming metric on $\mbX \times \mbX$:\s
\bi
\item Models 1, 2, and 3 have sufficient statistics $s_1(\bX), \dots, s_N(\bX)$,
the degrees of nodes $1, \dots, N$,
respectively.
Since the degrees of nodes are sums of $N-1$ edge variables $X_{i,j} \in \{0, 1\}$,
we have 
\beno
\norm{\bXi_a}_2 
&=& \sqrt{N-1} 
&\leq& \sqrt{N},
&& a = 1, \dots, N.
\ee
\item Models 2 and 3 include the additional sufficient statistic for brokerage $s_{N+1}(\bx) \coloneqq \sum_{i<j}^{N} \, X_{i,j} \; I_{i,j}(\bX)$,
where 
\beno
I_{i,j}(\bx) \;\coloneqq\; \one\left(\dsum_{h\, \in\, \mN_i\, \cap\, \mN_j} \, x_{i,h} \, x_{j,h} \,\geq\, 1 \right),
&& \{i,j\} \subset \mN. 
\ee
By the definition of $s_{N+1}(\bx)$,
we have $\Xi_{N+1,\{i,j\}} = 0$ for all pairs of nodes $\{i,j\} \subset \mN$ satisfying $\mN_i\, \cap\, \mN_j = \emptyset$.
The number of pairs of nodes $\{i,j\} \subset \mN$ satisfying $\mN_i\, \cap\, \mN_j \neq \emptyset$
is bounded above by $N D_N^2$:
For each of the $N$ nodes $i \in \mN$,
there are at most $D_N^2$ distinct nodes $j \in \mN_i$ such that $\mN_i\, \cap\, \mN_j \neq \emptyset$,
a fact established by Lemma \ref{lem:DN}.
In addition,
Lemma \ref{lem:brokerage_ham_bound} shows,
for each $\{i,j\} \subset \mN$, that  
$\Xi_{N+1,\{i,j\}} \leq 1+D_N$.
Thus,
\beno
\norm{\bm\Xi_{N+1}}_2 
\lte \sqrt{N\, D_N^2 \, (1+D_N)^2} 
\lte \sqrt{4\, N\, D_N^4}
&=& 2\; D_N^{2}\, \sqrt{N}.
\ee
\ei
As a result,
under Model 1,
\beno
\sqrt{N / 2}
\lte \Psi_N 
&\coloneqq& \max\limits_{1 \leq a \leq p} \, \norm{\bm\Xi_a}_2
&=& \sqrt{N-1}
\lte \sqrt{N},
\ee
whereas under Models 2 and 3, 
\beno
\sqrt{N / 2}
\lte \Psi_N
&\coloneqq& \max\limits_{1 \leq a \leq p} \, \norm{\bm\Xi_a}_2
\lte 2\; D_N^{2} \; \sqrt{N},
\ee
noting that $D_N \geq 1$ under Models 2 and 3. 

\s

{\bf Convergence rates.}
We obtain the following convergence rates using the auxiliary results in Appendices \ref{sup-sec:scale_pl} and \ref{sup-sec:mrf_conc}.
The following results hold for all large enough $N$.
The constants vary from model to model.

\bi
\item {\bf Corollary \ref{cor:model_1}:} 
The independence of edges under Model 1 implies that $D_N = 0$,\,
$\mnorm{\mD_N(\truth)}_2 = 1$,\,
and $\Lambda_N(\truth) = \widetilde\Lambda_N(\truth)$,
which in turn implies that $\Phi_N(\truth) = \widetilde\Phi_N(\truth)$.
We have $p = N$ and $\Psi_N \leq \sqrt{N}$.
By Lemma \ref{lemma.max} with $\alpha = 0$ and 
$\vartheta \in [0,\, 6\, (1/2 - \alpha)) = [0,\, 3)$, 
there exist constants $B > 0$ and $N_0 \geq 3$,
independent of $N$ and $p$, 
such that
\beno
\widetilde\Lambda_N(\truth)
&\leq& \dfrac{B}{N^{1 - \vartheta / 6}}
& \mbox{for all } N > N_0. 
\ee
As a result,
there exists a constant $C > 0$,
independent $N$ and $p$, 
such that,
for all $N > N_0$, 
\beno
\Phi_N(\truth)
\;=\; \widetilde\Phi_N(\truth)
\;\leq\; \dfrac{C \, \sqrt{N \log  N}}{N^{1-\vartheta/6}} 
\;=\; C \, \sqrt{\dfrac{\log \, N}{N^{1 - 2\, \vartheta / 6}}}
\;=\;  C \, \sqrt{\dfrac{\log \, N}{N^{1 - \vartheta / 3}}}.  
\ee

\s

\item {\bf Corollaries \ref{cor:model_2} and \ref{cor:model_3}:}
By assumption,
$D_N \geq 1$ under Models 2 and 3,
and $\Psi_N$ is bounded as follows: 
\beno
\sqrt{N/2} 
\lte \Psi_N 
\lte 2\; D_N^{2} \; \sqrt{N}. 
\ee
To bound $\gamma_N$,
recall that $\gamma_N$ is given by
\beno
\gamma_N
&=& \sqrt{24}\;\, (1 + D_N)\; \mnorm{\mD_N(\truth)}_2\; \Psi_N \; \sqrt{\log \max\{N,\, p\}}. 
\ee
Using $\mnorm{\mD_{N}(\truth)}_2 \geq 1$ along with $\Psi_N \geq \sqrt{N / 2}$,
we obtain the lower bound
\beno
\gamma_N
&\geq& \sqrt{12}\;\, D_N\; \sqrt{N \log N}
\ee
and the upper bound
\beno
\gamma_N
&\leq& \sqrt{768}\; D_N^{3}\; \mnorm{\mD_N(\truth)}_2\; \sqrt{N \log N},
\ee
using $1 + D_N \leq 2\, D_N$ ($D_N \geq 1$) along with $\Psi_N\, \leq\, 2\, D_N^{2} \sqrt{N}$ and
\beno 
\log \max\{N,\, p\} 
\= \log (N + 1)
&\leq& 2\, \log N.
\ee
Thus,
$\gamma_N$ satisfies
\beno
\sqrt{24}\;\, D_N\; \sqrt{N \log N}
\;\leq\; \gamma_N
\;\leq\; \sqrt{768}\;\, D_N^{3}\; \mnorm{\mD_N(\truth)}_2\; \sqrt{N \log N}.
\ee
We turn to bounding 
\beno
\widetilde\Phi_N(\truth)
&\coloneqq& \tllambdamin\,\, (1 + D_N)\,\, |\!|\!|\mD_N(\truth)|\!|\!|_2\,\, \Psi_N\; \sqrt{\log \max\{N,\, p\}}.
\ee
\blue{ 
By Lemma \ref{lem:Iscale},
there exists an integer $N_0 > 0$ and a subset $\mbH \subseteq \mbX$ such that,
for all $N > N_0$,
\bi
\item[$\bullet$] $-\nabla_{\nat}^2\; \widetilde\ell(\nat;\, \bx)$ is invertible for all $(\nat,\, \bx) \in \mB_\infty(\truth,\, \epsilon^\star) \times \mbH$,
\item[$\bullet$] the event $\bX \in \mbH$ occurs with probability at least $1 - 2\, / \max\{N,\, p\}^2$,
\ei
provided that either condition \sone or condition \stwo is satisfied.
}
By Lemma \ref{lemma.max},
there exist constants $C_1 > 0$ and $N_0 \geq 3$,
independent of $N$ and $p$,  
such that,
for all $N > N_0$, 
\beno
\widetilde\Lambda_N(\truth) 
\;\leq\; \dfrac{C_1 \, D_N^9}{N^{1 - (\alpha + \vartheta)}},
\ee
assuming that $\alpha \in [0,\, 1/2)$ and $\vartheta \in [0,\, 1/2-\alpha)$.
As a result,
there exists a constant $C_2 \coloneqq 4 \, C_1 > 0$,
independent of $N$ and $p$,
such that,
for all $N > N_0$,
\be
\label{eq:cor_big_bound}
\widetilde\Phi_N(\truth)
&\leq& \dfrac{C_2\; D_N^{9+1+2}\; \mnorm{\mD_N(\truth)}_2\; \sqrt{N \log N}}{N^{1-(\alpha + \vartheta)}}\s
\\
&=& C_2\; D_N^{12}\; \mnorm{\mD_N(\truth)}_2\; \sqrt{\dfrac{\log N}{N^{1-2\, (\alpha + \vartheta)}}},
\ee
using the inequalities $1 + D_N \leq 2 \, D_N$ and $\Psi_N\, \leq\, 2\, D_N^{2} \sqrt{N}$,
noting that $D_N \geq 1$ under Models 2 and 3.  
By Lemma \ref{prop:D_bound},  
there exist constants $C_3 > 0$ and $C_4 > 0$,
independent of $N$ and $p$, 
such that:
\bi
\item Corollary \ref{cor:model_2} with $\vartheta \in [0,\, 1/2-\alpha)$: 
\be
\label{eq:cor_D_1}
\mnorm{\mD_N(\truth)}_2
&\leq& 1 + 4 \, D_N^2
&\leq& 5\; D_N^2,
\ee
using the fact that $D_N$ satisfies $D_N \geq 1$ under Models 2 and 3.\vspace{.1cm}
\item Corollary \ref{cor:model_3} with $\vartheta = 0$:
If Assumption A is satisfied,
\be
\label{eq:cor_D_2}
\mnorm{\mD_N(\truth)}_2 
&\leq& 1 + 4 \, D_N^2 + \omega_1 \, C_3 \, \exp(C_4\, D_N^3)\s
\\
&\leq& 3\, \max\{4,\, \omega_1 \, C_3\} \, D_N^2 \,  \exp(C_4\, D_N^3)\s
\\
&\leq& B_1 \, \exp(2 \, \log D_N + C_4 \, D_N^3) \s\\ 
&\leq& B_1 \, \exp(A_1\, D_N^3),
\ee
using $D_N \geq 1$ along with $\log D_N \leq D_N \leq D_N^3$,
and defining $A_1 \coloneqq 2 + C_4 > 0$ and $B_1 \coloneqq 3\, \max\{4,\, \omega_1 \, C_3\} > 0$. 
Note that the constants $C_3 > 0$, $C_4 > 0$, and $\omega_1 \geq 0$ are independent of $N$ and $p$,
implying that $A_1$ and $B_1$ are likewise independent of $N$ and $p$.  \vspace{.1cm}
\ei
Upon collecting terms,
we conclude that there exist constants $C_5 > 0$ and $C_6 > 0$,
independent of $N$ and $p$, 
such that,
for all $N > N_0$:\vspace{.1cm}
\bi
\item Corollary \ref{cor:model_2} with $\alpha \in [0,\, 1/2)$ and $\vartheta \in [0,\, 1/2-\alpha)$: 
Equations \eqref{eq:cor_big_bound} and \eqref{eq:cor_D_1} provide the bound 
\beno
\widetilde\Phi_N(\truth)
&\leq& C_5\; D_N^{14}\, \sqrt{\dfrac{\log N}{N^{1-2\, (\alpha + \vartheta)}}},
\ee
where $C_5 \coloneqq 4\; C_2 > 0$ and 
$\gamma_N$ satisfies (noting $\sqrt{768} \leq 28$)
\beno
4 \, D_N \,\sqrt{N \log N}
\lte \gamma_N
\lte 28 \, D_N^{5}\, \sqrt{N \log N}.
\ee
\item Corollary \ref{cor:model_3} with $\alpha \in [0,\, 1/2)$ and $\vartheta = 0$:
If Assumption A is satisfied,
Equations \eqref{eq:cor_big_bound} and \eqref{eq:cor_D_2} provide the bound 
\beno
\widetilde\Phi_N(\truth)
&\leq& C_6\; D_N^{12} \exp(A_1\, D_N^3)\; \sqrt{\dfrac{\log N}{N^{1 - 2\, \alpha}}} \s\\
&\leq& C_6\, \exp(12 \log D_N + A_1 \, D_N^3) \;  \sqrt{\dfrac{\log N}{N^{1 - 2\, \alpha}}} \s\\
&\leq& C_6\, \exp(A_2 \, D_N^3) \;  \sqrt{\dfrac{\log N}{N^{1 - 2\, \alpha}}},
\ee
where $\gamma_N$ satisfies
\beno
4 \, D_N \, \sqrt{N \log N}
\lte \gamma_N\s
\\
\lte B_2\, \exp(A_3\, D_N^3) \, \sqrt{N \log N},
\ee
where $A_2 \coloneqq 12 + A_1 > 0$,
$C_6 \coloneqq B_1 \, C_2 > 0$,
and 
$B_2 \coloneqq 28 \, B_1 > 0$;
note that the upper bound on $\gamma_N$ leverages the fact that $D_N \geq 1$ under Models 2 and 3 
along with $\log D_N \leq D_N \leq D_N^3$,
so that
\beno 
D_N^{3} \exp(A_1 \, D_N^3) 
= \exp\left(3 \log D_N + A_1 \, D_N^3\right)
\leq \exp\left(A_3\, D_N^3 \right),
\ee
where $A_3 \coloneqq 3 + A_1 > 0$ is a constant,
independent of $N$ and $p$.
\ei
\ei
\qed

\s

\subsection{Bounding \texorpdfstring{$\widetilde{\Lambda}_N(\truth)$}{Lambda}}
\label{sup-sec:scale_pl}

\blue{ 
We bound 
\[
\begin{array}{cllll}
\widetilde\Lambda_N(\truth)
&\coloneqq& \max\limits_{\bx\, \in\, \mbH}\; \sup\limits_{\nat\, \in\, \mB_{\infty}(\truth,\, \epsilon^\star)}\, \mnorm{(-\nabla_{\nat}^2\; \widetilde\ell(\nat;\, \bx))^{-1}}_{\infty},
\end{array}
\]
using Lemma \ref{lemma.max},
which leverages auxiliary results supplied by Lemmas \ref{lem:beta}--\ref{lem:Iscale}. 
}
To do so,
we first introduce additional notation.
The negative Hessian $-\nabla_{\nat}^2\; \widetilde\ell(\nat;\, \bx)$ corresponding to Models 2 and 3 is of the form
\be
\label{expected.hessian}
-\nabla_{\nat}^2 \; \widetilde\ell(\nat;\, \bx) 
\= 
\left( 
\begin{matrix} 
\bmat        & \cvec\s
\\
\cvec^{\top} & \vvar \\ 
\end{matrix} \right),
\ee
where
\bi
\item the entries $A_{i,j}(\nat,\, \bx)$ of the matrix $\bmat \in \mR^{N \times N}$ are given by  
\beno
A_{i,j}(\nat,\, \bx) 
\= \dsum_{a<b}^{N} \, \cov_{\nat, \bx_{-\{a,b\}}}(s_i(\bX), \, s_j(\bX)), 
&& i,\, j = 1, \ldots, N;
\ee
\item the entries $c_i(\nat,\, \bx)$ of the vector $\cvec \in \mR^{N}$ are given by 
\beno
c_{i}(\nat,\, \bx) 
\= \dsum_{a<b}^{N} \, \cov_{\nat, \bx_{-\{a,b\}}}(s_{i}(\bX), \, s_{N+1}(\bX)),
&& i = 1, \ldots, N;
\ee
\item $\vvar  \in (0,\, \infty)$ is given by 
\beno
\vvar 
\= \dsum_{a<b}^{N} \,\var_{\nat, \bx_{-\{a,b\}}} \, s_{N+1}(\bX), 
\ee
\ei
where $\cov_{\nat, \bx_{-\{a,b\}}}$ and $\var_{\nat, \bx_{-\{a,b\}}}$ 
are the conditional covariance and variance operators with respect to the conditional probability distribution of edge variable 
$X_{a,b}$ given all other edge variables $\bX_{-\{a,b\}} = \bx_{-\{a,b\}}$. 
The negative Hessian under Model 1 is $-\nabla_{\nat}^2\; \widetilde\ell(\nat;\, \bx) = \bmat$.

\blue{ 
Throughout, 
the indicator function
\beno
I_{i,j}(\bx) 
&\coloneqq& \one\left(\dsum_{h\, \in\, \mN_i\, \cap\, \mN_j} \, x_{i,h} \, x_{j,h} \,\geq\, 1 \right),
&& \{i,j\} \subset \mN
\ee
indicates whether there exists a node $h \in \mN_i\, \cap\, \mN_j$ in the intersection of $i$'s neighborhood $\mN_i$ and $j$'s neighborhood $\mN_j$ such that both $i$ and $j$ are connected to $h$.
In addition,
define the vector
\beno
\bI_i(\bx)
&\coloneqq& (I_{i,1}(\bx), \ldots, I_{i,i-1}(\bx),\, I_{i,i+1}(\bx), \ldots, I_{i,N}(\bx)),
&& i \in \mN. 
\ee
As a result,
$\norm{\bI_i(\bx)}_{\infty}$ indicates whether,
for a given node $i \in \mN$, 
there exists a distinct node $j \in \mN \setminus\, \{i\}$ such that $\mN_i\, \cap\, \mN_j \neq \emptyset$ and $i$ and $j$ are both connected to a third node $h \in \mN \setminus \{i,\, j\}$ contained in $\mN_i\, \cap\, \mN_j$. 
}

{ 
\lemma
\label{lemma.max}
Assume that $\truth \in \mR^p$ satisfies condition \eqref{condition.theta0} in Section \ref{sec:corollaries}.
Then there exist constants $B > 0$,\,
$C > 0$,\,
and $N_0 \geq 3$, 
independent of $N$ and $p$,
such that $\widetilde\Lambda_N(\truth)$ is bounded above as follows:
\bi
\item Model 1: 
For all $N > N_0$, 
$-\nabla_{\nat}^2 \; \widetilde\ell(\nat;\, \bx)$ is invertible for all $\bx \in \mbX$ and
\beno
\widetilde\Lambda_N(\truth)
\;\coloneqq\; \max\limits_{\bx\, \in\, \mbX}\;\, \sup\limits_{\nat\, \in\, \mB_{\infty}(\truth,\, \epsilon^\star)}\, \mnorm{(-\nabla_{\nat}^2 \; \widetilde\ell(\nat;\, \bx))^{-1}}_{\infty}
\;\leq\; \dfrac{B}{N^{1 - (\alpha + \vartheta/6)}},
\ee
assuming $\alpha \in [0, \, 1/2)$ and $\vartheta \in [0,\, 6\, (1/2 - \alpha))$.\s
\item Models 2 and 3:
For all $N > N_0$,  
the event
\beno
\mbH 
\;\coloneqq\; \left\{\bx \in \mbX:\, \dsum_{i=1}^{N} \norm{\bI_i(\bx)}_{\infty} \geq \dfrac{N}{2\, (1 + \exp((3 + D_N) \, (\norm{\truth}_{\infty} + \epsilon^\star))^2}
\right\}
\ee
occurs with at least probability $1 - 2\, / \max\{N,\, p\}^2$.
In the event $\bx \in \mbH$,
$-\nabla_{\nat}^2 \; \widetilde\ell(\nat;\, \bx)$ is invertible and
\beno
\widetilde\Lambda_N(\truth)
\;\coloneqq\; \max\limits_{\bx\, \in\, \mbH}\;\, \sup\limits_{\nat\, \in\, \mB_{\infty}(\truth,\, \epsilon^\star)}\, \mnorm{(-\nabla_{\nat}^2 \; \widetilde\ell(\nat;\, \bx))^{-1}}_{\infty}
\;\leq\; \dfrac{C\, D_N^{9}}{N^{1 - (\alpha + \vartheta)}},
\ee
assuming that $\alpha \in [0,\, 1/2)$ and $\vartheta \in [0, \, 1/2 - \alpha)$,
$D_N$ satisfies
\beno
1 
&\leq& D_N
&<& \dfrac{L + \vartheta \, \log N}{12\; \epsilon^\star} - 3,
\ee
and either condition \sone or condition \stwo is satisfied.
\ei
}

{\em Remark.}
Under Model 1,
edges are independent and $D_N = 0$,
whereas under Models 2 and 3,
edges are dependent and $D_N \geq 1$.
The upper bound on $D_N$ under Models 2 and 3 ensures that $|\!|\truth|\!|_\infty > 0$.

\s

\s

\llproof \ref{lemma.max}.
Using \eqref{expected.hessian},  
we can write the negative Hessian $- \nabla_{\nat}^2 \; \widetilde\ell(\nat;\, \bx)$ 
corresponding to Models 2 and 3 as
\beno
- \nabla_{\nat}^2 \; \widetilde\ell(\nat;\, \bx)
&=& \left(
\begin{matrix}
\bmat & \cvec \s\\
\cvect & \vvar
\end{matrix} \right),
\ee
where $\bmat \in \mR^{N \times N}$,\,
$\cvec \in \mR^{N}$,\,
and $\vvar  \in \mR^{+}$ are defined above.

\paragraph{Bounding $\mnorm{\bmat^{-1}}_{\infty}$}

Lemma \ref{lem:beta} proves that the smallest eigenvalue of $\bmat$ is strictly positive on 
$\mB_{\infty}(\truth,\, \epsilon^\star)$ for all $\bx \in \mbX$,
so $\bmat$ is invertible on $\mB_{\infty}(\truth,\, \epsilon^\star)$ for all $\bx \in \mbX$.
Theorem 1.2 of \citet{HiWi15},
along with the bounds on the entries of $\bmat$ given in 
\eqref{eq:bmat_entry_bounds} of Lemma \ref{lem:beta},
reveal that
\beno
\mnorm{\bmat^{-1}}_{\infty}
&\leq& \dfrac{(3 \, N - 4) \, (1 + \exp((3 + D_N) \, \norm{\nat}_{\infty})))^2}
{2 \, N^{-\alpha} \, (N-2) \, (N-1)} \s\\
&\leq&  (1 + \exp((3 + D_N) \, \norm{\nat}_{\infty})))^2 \, \dfrac{9}{2 \, N^{1 - \alpha}}\s
\\
&\leq& (1 + \exp((3 + D_N) \, (\norm{\truth}_{\infty} + \epsilon^\star)))^2 \, \dfrac{9}{2 \, N^{1 - \alpha}}\s
\\
&\leq& (2\, \exp((3 + D_N) \, (\norm{\truth}_{\infty} + \epsilon^\star)))^2 \, \dfrac{9}{2 \, N^{1 - \alpha}}\s
\\
&=& \exp((3 + D_N) \, (\norm{\truth}_{\infty} + \epsilon^\star))^2 \, \dfrac{18}{N^{1 - \alpha}}\s
\\
&=& \dfrac{18\; \tau(\truth)^2}{N^{1 - \alpha}},
\ee
where
\beno
\tau(\truth)
&\coloneqq& \exp((3 + D_N) \, (\norm{\truth}_{\infty} + \epsilon^\star)).
\ee
The above exploits the fact that $\norm{\nat}_{\infty} \leq \norm{\truth}_{\infty} + \epsilon^\star$ for all $\nat \in \mB_{\infty}(\truth,\, \epsilon^\star)$,
along with the inequality
\beno
\dfrac{3 \, N - 4}{2\, (N-2) \, (N-1)}
\,=\, \dfrac{3 \, (N - 1) - 1}{2\, (N-2) \, (N-1)}
\,\leq\, \dfrac{3 \, (N - 1)}{2\, (N-2) \, (N-1)}
\,=\, \dfrac{3}{2\, (N-2)}, 
\ee
which is bounded above by 
\beno
\dfrac{3}{2\, (N-2)} 
&\leq& \dfrac{9}{2 \, N}, 
\ee
provided $N \geq 3$.
Using the assumption
\beno
\norm{\truth}_{\infty} 
&\leq& \dfrac{L + \vartheta \, \log N}{12\, (3 + D_N)} - \epsilon^\star,
\ee
we obtain
\beno 
\tau(\truth)^2
&=& \exp(2\, (3 + D_N) \, (\norm{\truth}_{\infty} + \epsilon^\star))
\;\leq\; \exp\left(2\, (3 + D_N)\, \dfrac{L + \vartheta \, \log N}{12\, (3 + D_N)}\right)\s
\\
&=& \exp\left(\dfrac{L + \vartheta \, \log N}{6}\right)\s
\;=\; \exp(L\, /\, 6)\; N^{\vartheta\, /\, 6}.
\ee 
As a consequence,
we find that,
for all $\bx \in \mbX$,  
\beno
\mnorm{\bmat^{-1}}_{\infty}
\lte \dfrac{18\; \tau(\truth)^2}{N^{1 - \alpha}}
\lte \dfrac{18\; \exp(L\, /\, 6)\; N^{\vartheta\, /\, 6}}{N^{1 - \alpha}}
&=& \dfrac{B}{N^{1 - (\alpha + \vartheta/6)}},
\ee
where $B \coloneqq 18 \exp(L / 6) > 0$.
Under Model 1,
we therefore obtain
\beno
\widetilde\Lambda_N(\truth)
&\coloneqq& \max\limits_{\bx\, \in\, \mbX}\;\, \sup\limits_{\nat\, \in\, \mB_{\infty}(\truth,\, \epsilon^\star)} \, \mnorm{(-\nabla_{\nat}^2 \; \widetilde\ell(\nat;\, \bx))^{-1}}_{\infty}\s\s
\\
&=& \max\limits_{\bx\, \in\, \mbX}\;\, \sup\limits_{\nat\, \in\, \mB_{\infty}(\truth,\, \epsilon^\star)} \, \mnorm{\bmat^{-1}}_{\infty}\s
\\
&\leq& \dfrac{B}{N^{1 - (\alpha + \vartheta/6)}},
\ee
assuming $\alpha \in [0, \, 1/2)$ and $\vartheta \in [0,\, 6\, (1/2 - \alpha))$.\s

\paragraph{Bounding $\mnorm{(- \nabla_{\nat}^2\; \widetilde\ell(\nat;\, \bx))^{-1}}_{\infty}$ under Models 2 and 3}

Let
\beno
\schur 
&\coloneqq& \vvar - \cvect \, \bmat^{-1} \, \cvec.
\ee 
Theorem 8.5.11 of \citet[][p.~99]{Ha97} implies that,
if the inverse of $- \nabla_{\nat}^2 \; \widetilde\ell(\nat;\, \bx)$ exists, 
then it can be written as
{\small 
\beno
\label{inverse.block.matrix}
(- \nabla_{\nat}^2 \; \widetilde\ell(\nat;\, \bx))^{-1} 
\;=\; \left(
\begin{matrix}
\bmat & \cvec \s\\
\cvect & \vvar
\end{matrix} \right)^{-1}
\;=\; \left(
\begin{matrix}
\bB_{1,1}(\nat,\, \bx) & \bB_{1,2}(\nat,\, \bx)\s
\\
\bB_{1,2}(\nat,\, \bx)^\top & \bB_{2,2}(\nat,\, \bx)
\end{matrix} \right),
\ee
}
where
\beno
\bB_{1,1}(\nat,\, \bx)
\coloneqq\; \bmat^{-1} + \schur^{-1} \, (\bmat^{-1} \, \cvec) \, (\bmat^{-1} \, \cvec)^{\top}\s
\\
\bB_{1,2}(\nat,\, \bx)
\coloneqq - \schur^{-1} \, \bmat^{-1} \, \cvec\s
\\
\bB_{2,2}(\nat,\, \bx)
\coloneqq\; \schur^{-1}.
\ee
To establish that $- \nabla_{\nat}^2 \; \widetilde\ell(\nat;\, \bx)$ is invertible,
note that $\bmat$ is invertible on $\mB_{\infty}(\truth,\, \epsilon^\star)$ for all $\bx \in \mbX$,
because its smallest eigenvalue is strictly positive by Lemma \ref{lem:beta}.
In addition,
we demonstrate below that there exists an integer $N_0 \geq 3$,
independent of $N$ and $p$, 
such that
\beno
\mbP(\bX \in \mbH)
&\geq& 1 - \dfrac{2}{\max\{N,\, p\}^2}
&\mbox{for all}& N > N_0
\ee
and 
\beno
\schur 
&>& 0 
& \mbox{for all}& \bx \in \mbH
\ee
provided that either condition \sone or condition \stwo is satisfied.
As a result,
$- \nabla_{\nat}^2 \; \widetilde\ell(\nat;\, \bx)$ is invertible on $\mB_{\infty}(\truth,\, \epsilon^\star)$ for all $\bx \in \mbH$ by virtue of Theorem 8.5.11 of \citet[][p.~99]{Ha97}. 
We proceed to bound $\mnorm{(- \nabla_{\nat}^2 \; \widetilde\ell(\nat;\, \bx))^{-1}}_{\infty}$ under Models 2 and 3.
Observe that
\beno
\mnorm{(-\nabla_{\nat}^2 \; \widetilde\ell(\nat;\, \bx))^{-1}}_{\infty}
&\leq& \max\left\{ T_1, \, T_2 \right\},  
\ee
where 
\beno
T_1 &\coloneqq& 
\mnorm{\bmat^{-1} + \schur^{-1} \, (\bmat^{-1} \, \cvec) \, (\bmat^{-1} \, \cvec)^{\top}}_{\infty}\s
\\
&+& \norm{\schur^{-1} \, \bmat^{-1} \, \cvec}_{\infty}\s
\\
T_2 &\coloneqq& 
\norm{\schur^{-1} \, \bmat^{-1} \, \cvec}_1 + |\schur|^{-1}. 
\ee
We bound the terms $T_1$ and $T_2$ one by one.

\paragraph*{Bounding $T_1$} 

The term $T_1$ is defined as
\beno
T_1 &\coloneqq& 
\mnorm{\bmat^{-1} + \schur^{-1} \, (\bmat^{-1} \, \cvec) \, (\bmat^{-1} \, \cvec)^{\top}}_{\infty}\s
\\
&+& \norm{\schur^{-1} \, \bmat^{-1} \, \cvec}_{\infty}.
\ee
We bound the first term of $T_1$ by using the triangle inequality:
\beno
&& \mnorm{\bmat^{-1} + \schur^{-1} \, (\bmat^{-1} \, \cvec) \, (\bmat^{-1} \, \cvec)^{\top}}_{\infty}\vspace{.35cm}
\\
&\leq& 
\hide{
\mnorm{\bmat^{-1}}_{\infty} + \mnorm{\schur^{-1} \, (\bmat^{-1} \, \cvec) \, (\bmat^{-1} \, \cvec)^{\top}}_{\infty}\vspace{.35cm}
\\
&=& 
}
\mnorm{\bmat^{-1}}_{\infty} + |\schur|^{-1} \, \mnorm{(\bmat^{-1} \, \cvec) \, (\bmat^{-1} \, \cvec)^{\top}}_{\infty}\vspace{.35cm}
\\
&=& \mnorm{\bmat^{-1}}_{\infty} + |\schur|^{-1} \, \norm{\bmat^{-1} \, \cvec}_{\infty} \, 
\norm{\bmat^{-1} \, \cvec}_{1}\vspace{.35cm}
\\
&\leq& \mnorm{\bmat^{-1}}_{\infty} + N\; |\schur|^{-1} \, \norm{\bmat^{-1} \, \cvec}_{\infty}^2\vspace{.35cm}
\\
&\leq& \mnorm{\bmat^{-1}}_{\infty} + N\, |\schur|^{-1}\, \norm{\cvec}_{\infty}^2\, \mnorm{\bmat^{-1}}_{\infty}^2\vspace{.35cm}
\\
\= \mnorm{\bmat^{-1}}_{\infty}\, (1 + N \, |\schur|^{-1} \, \norm{\cvec}_{\infty}^2 \, \mnorm{\bmat^{-1}}_{\infty}),
\ee
taking advantage of the identity
\beno
\mnorm{\bz\, \bz^{\top}}_{\infty}
\= \max\limits_{1 \leq i \leq N}\, \dsum_{j = 1}^N |z_i\, z_j|
\= \max\limits_{1 \leq i \leq N}\, |z_i|\, \dsum_{j = 1}^N |z_j|
\= \norm{\bz}_\infty\, \norm{\bz}_1
\ee
applied to the vector 
\beno
\bz
&\coloneqq& \bmat^{-1}\; \cvec,
\ee
along with the fact that $\mnorm{\bmat^{-1}}_1 = \mnorm{\bmat^{-1}}_{\infty}$,
thanks to the symmetry of $\bmat$. 
The second term of $T_1$ can be bounded as follows:
\beno
\norm{\schur^{-1} \, \bmat^{-1} \, \cvec}_{\infty}
&\leq& |\schur|^{-1} \, \norm{\cvec}_{\infty} \, \mnorm{\bmat^{-1}}_{\infty}. 
\ee
Combining these results gives the following bound on $T_1$:
\beno
T_1 
&\leq& \mnorm{\bmat^{-1}}_{\infty}\s
\\
&\times& \left(1 + N\, |\schur|^{-1}\, \norm{\cvec}_{\infty}^2\, \mnorm{\bmat^{-1}}_{\infty} + |\schur|^{-1}\, \norm{\cvec}_{\infty}\right).  
\ee

\paragraph*{Bounding $T_2$}

The term $T_2$ is defined as follows:
\beno
T_2 
&\coloneqq& \norm{\schur^{-1} \, \bmat^{-1} \, \cvec}_1 + |\schur|^{-1}.
\ee 
We bound $T_2$ by  
\beno 
T_2 
&=& \norm{\schur^{-1} \, \bmat^{-1} \, \cvec}_1 + |\schur|^{-1}\s
\\
&\leq& |\schur|^{-1} \, (1 + \norm{\cvec}_1 \, \mnorm{\bmat^{-1}}_{1})\s
\\
\= |\schur|^{-1} \, (1 + \norm{\cvec}_1 \, \mnorm{\bmat^{-1}}_{\infty}) \s
\\
&\leq& |\schur|^{-1} \, (1 + N \, \norm{\cvec}_{\infty} \, \mnorm{\bmat^{-1}}_{\infty}),
\ee
using the inequality $\norm{\bv}_1 \leq N \, \norm{\bv}_{\infty}$,
where the step from $\mnorm{\bmat^{-1}}_{1}$ to $\mnorm{\bmat^{-1}}_{\infty}$ follows from the symmetry of $\bmat$.
We bound the terms in $T_1$ and $T_2$ one by one.
The resulting bounds hold for all $\nat \in \mB_\infty(\truth,\, \epsilon^\star)$.

\paragraph{Bounding $\norm{\bmat^{-1}}_{\infty}$}

We have shown above that
\beno
\mnorm{\bmat^{-1}}_{\infty}
&\leq& \dfrac{18\; \tau(\truth)^2}{N^{1 - \alpha}}
& \mbox{for all } \bx \in \mbX. 
\ee




\paragraph{Bounding $\norm{\cvec}_{\infty}$}

Lemma \ref{lem:vctc} establishes that
\beno
\norm{\cvec}_{\infty} 
&\leq& 3\, D_N^3
&\mbox{for all } \bx \in \mbH.
\ee

\paragraph{Bounding $|\schur|^{-1}$}

The term $\schur$ is defined as
\beno
\schur 
&\coloneqq& \vvar - \cvect \, \bmat^{-1} \, \cvec\s
\\
&=& \cvect \, \cvec\s 
\\
&\times& \left(\dfrac{\vvar}{\cvect \, \cvec} - \dfrac{ \cvect \, \bmat^{-1} \, \cvec }{\cvect \cvec} \right). 
\ee
To bound $\schur$,
we leverage Lemmas \ref{lem:ral}--\ref{lem:vctc},
which show that the following bounds are satisfied for all $\bx \in \mbH$:
\beno
\cvect \cvec 
\;\geq\; \dfrac{N}{128 \; \tau(\truth)^6}\;\; \mbox{ by Lemma \ref{lem:ctc},}\s
\\
\dfrac{\vvar}{\cvect \cvec}
\;\geq\; \dfrac{1}{576\; D_N^6\; \tau(\truth)^4}\;\; \mbox{ by Lemma \ref{lem:vctc},}\s\s
\\
\dfrac{ \cvect \, \bmat^{-1} \, \cvec }{\cvect \cvec} \;\leq\; \dfrac{12\; \tau(\truth)^2}{N^{1-\alpha}}\;\; \mbox{ by Lemma 
\ref{lem:ral}.}
\ee
All of the above quantities are well-defined because $N \geq 3$,\,
$D_N \geq 1$,\,
and $\tau(\truth) > 0$ under Models 2 and 3. 
These results help bound $\schur$ as follows:
\beno
\schur 
&=& \cvect\, \cvec\s 
\\
&\times& \left( \dfrac{\vvar}{\cvect \, \cvec} - \dfrac{ \cvect \, \bmat^{-1} \, \cvec }{\cvect \cvec} \right)\s\s
\\
&\geq& \dfrac{N}{128 \, \tau(\truth)^6}\, \left(\dfrac{1}{576\; D_N^6\; \tau(\truth)^4} - \dfrac{12\; \tau(\truth)^2}{N^{1-\alpha}}\right)\s\s
\\
\= \dfrac{N}{(128) \, (576) \; D_N^6\; \tau(\truth)^{10}} \left(1 - \dfrac{(12) \, (576)\; D_N^6\; \tau(\truth)^6}{N^{1-\alpha}} \right).
\ee
To bound the term $1 - (12) \, (576) \; D_N^6\, \tau(\truth)^6\, /\, N^{1-\alpha}$,
observe that the assumption
\beno
\label{eq:theta-star-bound}
\norm{\truth}_{\infty} &\leq& \dfrac{L + \vartheta \, \log N}{12\, (3 + D_N)} - \epsilon^\star
\ee
implies that the term $\tau(\truth)^6$ is bounded above by
\beno
\tau(\truth)^6
\= \exp(6\, (3 + D_N)\, (\norm{\truth}_{\infty} + \epsilon^\star))\s
\\
\lte \exp\left(\dfrac{6}{12}\, (L + \vartheta \, \log  N)\right)\s
\\
\= \exp((1/2)\, L)\; N^{\vartheta / 2},
\ee
which in turn implies that
\beno
\dfrac{D_N^6\; \tau(\truth)^6}{N^{1-\alpha}}
\;\leq\; \dfrac{\exp((1/2)\, L)\; D_N^6\; N^{\vartheta / 2}}{N^{1-\alpha}}
\;\leq\; \dfrac{\exp((1/2)\, L)\; D_N^6}{N^{1/2}},
\ee
using the assumption that 
$\alpha \in [0,\, 1/2)$ and $\vartheta \in [0,\, 1/2-\alpha)$. 
Since conditions \sone and \stwo ensure that $D_N = O(\log N)$,
there exist constants $C_1 \in (0,\, 1)$ and $N_1 \geq N_0$,
independent of $N$ and $p$, 
such that,
for all $N > N_1$,
\beno
1 - \dfrac{(12) \, (576) \; D_N^6\; \tau(\truth)^6}{N^{1-\alpha}}
\;\geq\; 1 - \dfrac{(12) \, (576) \, \exp((1/2)\, L)\; D_N^6}{N^{1/2}}
\;\geq\; C_1 
\;>\; 0.
\ee
We then obtain,
for all $N > N_1$ and all $\bx \in \mbH$,
\beno
\schur
&\geq& \dfrac{N}{(128) \, (576) \; D_N^6\; \tau(\truth)^{10}}\, \left(1 - \dfrac{(12) \, (576) \; D_N^6\; \tau(\truth)^6}{N^{1-\alpha}} \right) \s\\
&\geq& \dfrac{C_1\; N}{(128) \, (576)  \; D_N^6\; \tau(\truth)^{10}},
\ee
which shows that $\schur > 0$ and hence 
\beno
\left| \schur \right|^{-1}
&=& \dfrac{1}{\schur}
&\leq& \dfrac{C_2\; D_N^6 \, \tau(\truth)^{10}}{N},
\ee
defining $C_2 \coloneqq (128) \, (576) \,/\, C_1 > 0$. 

\paragraph*{Bounding $\max\{T_1,\, T_2\}$}

We have shown that
{\small 
\beno
T_1
&\leq& \mnorm{\bmat^{-1}}_{\infty} \, 
\left(1 + N \, |\schur|^{-1}\norm{\cvec}_{\infty}^2 \, \mnorm{\bmat^{-1}}_{\infty} + |\schur|^{-1} \norm{\cvec}_{\infty}\right)\s
\\
T_2
&\leq& |\schur|^{-1} \, (1 + N \, \norm{\cvec}_{\infty} \, \mnorm{\bmat^{-1}}_{\infty}). \s\\
\ee
}
Using the bounds derived above,
we obtain,
for all $N \geq N_1$ and all $\bx \in \mbH$,
{\small
\beno
T_1
&\leq& \mnorm{\bmat^{-1}}_{\infty} \, \left(1 + N \, |\schur|^{-1} \, \norm{\cvec}_{\infty}^2 \, \mnorm{\bmat^{-1}}_{\infty} + |\schur|^{-1} \, \norm{\cvec}_{\infty}\right)\s\s
\\
&\leq& \dfrac{18 \, \tau(\truth)^2}{N^{1-\alpha}} \, 
\left(1 + N \left( \dfrac{C_2 \, D_N^6 \, \tau(\truth)^{10}}{N} \right) (3 \, D_N^3)^2 \left(\dfrac{18 \, \tau(\truth)^2}{N^{1-\alpha}}\right) + \left(\dfrac{C_2 \, D_N^6 \, \tau(\truth)^{10}}{N} \right) \, (3 \, D_N^3)\right)\s\s
\\
&\leq& \dfrac{18\; \tau(\truth)^2}{N^{1 - \alpha}} \, 
\left(1 + \dfrac{(18) \, (3)^2 \; C_2\; D_N^{12}\; \tau(\truth)^{12}}{N^{1-\alpha}}\; + \dfrac{3 \, C_2\; D_N^{9} \, \tau(\truth)^{10}}{N} \right)\s\s\s
\\
&=& \dfrac{(18)^2\, (3)^2 \, C_2\, D_N^{12}\, \tau(\truth)^{14}}{N^{1 - \alpha}} \, 
\left(\dfrac{1}{(18)\, (3)^2 \, C_2\, D_N^{12}\, \tau(\truth)^{12}} + \dfrac{1}{N^{1-\alpha}} + \dfrac{1}{(18)\, (3) \, D_N^3\,  \tau(\truth)^{2}\; N} \right)\s\s\s
\\ 
&\leq& \dfrac{(18)^2\, (3)^2 \; C_2\; D_N^{12}\; \tau(\truth)^{14}}{N^{1 - \alpha}}\, 
\left(\dfrac{1}{C_2\; D_N^{12}\; \tau(\truth)^{12}} + \dfrac{1}{N^{1-\alpha}} + \dfrac{1}{D_N^3\;  \tau(\truth)^{2}\; N} \right)\s\s
\\
\lte \dfrac{(3)^3 \, (18)^2 \; C_2\; D_N^{12}\; \tau(\truth)^{14}}{N^{1 - \alpha}}\; 
\dfrac{1}{\min\{C_2\; D_N^{12}\; \tau(\truth)^{12},\; N^{1-\alpha},\;D_N^3\; \tau(\truth)^{2}\; N\}}
\ee
}\hide{
\beno
T_1
&\leq& \mnorm{\bmat^{-1}}_{\infty}\s\s
\\
&\times& \left(1 + N \, |\schur|^{-1} \, \norm{\cvec}_{\infty}^2 \, \mnorm{\bmat^{-1}}_{\infty} + |\schur|^{-1} \, \norm{\cvec}_{\infty}\right)\s\s
\\
&\leq& \dfrac{18 \, \tau(\truth)^2}{N^{1-\alpha}}\s\s
\\
&\times& \left(1 + N \left( \dfrac{C_2 \, D_N^6 \, \tau(\truth)^{12}}{N} \right) (3 \, D_N^3)^2 \left(\dfrac{18 \, \tau(\truth)^2}{N^{1-\alpha}}\right) + \left(\dfrac{C_2 \, D_N^6 \, \tau(\truth)^{12}}{N} \right) \, (3 \, D_N^3)\right)\s\s
\\
&\leq& \dfrac{18\; \tau(\truth)^2}{N^{1 - \alpha}}\s\s
\\
&\times& \left(1 + \dfrac{(18) \, (3)^2 \; C_2\; D_N^{12}\; \tau(\truth)^{14}}{N^{1-\alpha}}\; + \dfrac{3 \, C_2\; D_N^{9} \, \tau(\truth)^{12}}{N} \right)\s\s\s
\\
&=& \dfrac{(18)^2\, (3)^2 \, C_2\, D_N^{12}\, \tau(\truth)^{16}}{N^{1 - \alpha}}\s\s
\\
&\times& \left(\dfrac{1}{(18)\, (3)^2 \, C_2\, D_N^{12}\, \tau(\truth)^{14}} + \dfrac{1}{N^{1-\alpha}} + \dfrac{1}{(18)\, (3) \, D_N^3\,  \tau(\truth)^{2}\; N} \right)\s\s\s
\\ 
&\leq& \dfrac{(18)^2\, (3)^2 \; C_2\; D_N^{12}\; \tau(\truth)^{16}}{N^{1 - \alpha}}\, 
\left(\dfrac{1}{C_2\; D_N^{12}\; \tau(\truth)^{14}} + \dfrac{1}{N^{1-\alpha}} + \dfrac{1}{D_N^3\;  \tau(\truth)^{2}\; N} \right)\s\s
\\
\lte \dfrac{(3)^3 \, (18)^2 \; C_2\; D_N^{12}\; \tau(\truth)^{16}}{N^{1 - \alpha}}\; 
\dfrac{1}{\min\{C_2\; D_N^{12}\; \tau(\truth)^{14},\; N^{1-\alpha},\;D_N^3\; \tau(\truth)^{2}\; N\}}
\ee
}
and
\beno
T_2
&\leq& |\schur|^{-1} \, (1 + N \, \norm{\cvec}_{\infty} \, \mnorm{\bmat^{-1}}_{\infty})\s\s
\\
\lte \dfrac{C_2 \, D_N^6 \, \tau(\truth)^{10}}{N}\, \left(1 + N\, D_N^3\, \dfrac{(3) \,  (18) \,\tau(\truth)^2}{N^{1 - \alpha}} \right)\s\s
\\
&=& \dfrac{C_2 \, D_N^6 \, \tau(\truth)^{10}}{N}\, \left(1 + (3) \, (18)\, D_N^3\, N^\alpha\, \tau(\truth)^2\right)\s\s
\\
\lte
\dfrac{C_3 \, D_N^{9}\, \tau(\truth)^{12}}{N^{1 - \alpha}}, \s
\ee
using the fact that $D_N^3\, N^\alpha\, \tau(\truth)^2 \geq 1$ under Models 2 and 3 and defining $C_3 \coloneqq (2) \, (3) \, (18) \, C_2 > 0$.
Define
\beno
U
&\coloneqq& \dfrac{C_3\;  D_N^{9} \, \tau(\truth)^{12}}{N^{1 - \alpha}}\s
\\
V
&\coloneqq& \dfrac{C_4\; D_N^3\; \tau(\truth)^2}{\min\{C_2\; D_N^{12}\; \tau(\truth)^{12},\; N^{1-\alpha},\;D_N^3\; \tau(\truth)^{2}\; N\}},
\ee
so that the bounds on $T_1$ and $T_2$ can be stated in terms of $U$ and $V$:
\beno
T_1 
\lte \dfrac{C_3 \, C_4 \; D_N^{12}\; \tau(\truth)^{14}}{N^{1 - \alpha}\; \min\{C_2\; D_N^{12}\; \tau(\truth)^{12},\; N^{1-\alpha},\; D_N^3\; \tau(\truth)^{2}\; N\}}
&=& U\,V\s
\\
T_2 \lte U,
\ee
where $C_4 \coloneqq (3)^3 \, (18)^2\; C_2\, /\, C_3 > 0$.
Thus, 
for all $N \geq N_1$ and all $\bx \in \mbH$, 
\beno
\label{eq:inf_mnorm_bound}
\sup\limits_{\nat\, \in\, \mB_{\infty}(\truth,\, \epsilon^\star)} \, 
\mnorm{(-\nabla_{\nat}^2 \; \widetilde\ell(\nat;\, \bx))^{-1}}_{\infty}
\;\leq\; \max\left\{T_1,\, T_2\right\}
\;\leq\; \max\left\{U\,V,\; U\right\}. 
\ee
To make the bound on $\max\left\{T_1,\, T_2\right\}$ as tight as possible,
we need constants $C_5 \geq 1$ and $N_2 \geq N_1$,
independent of $N$ and $p$, 
such that,
for all $N > N_2$,
\beno
V
&\coloneqq& \dfrac{C_4\; D_N^3\; \tau(\truth)^2}{\min\{C_2\; D_N^{12}\; \tau(\truth)^{12},\; N^{1-\alpha},\; D_N^3\; \tau(\truth)^{2}\; N\}}
&\leq& C_5.
\ee
Upon inspecting the denominator of $V$,
\beno
\min\{C_2\; D_N^{12}\; \tau(\truth)^{12},\; N^{1-\alpha},\; D_N^3\; \tau(\truth)^{2}\; N\},
\ee
and observing that $\alpha \in [0, 1/2)$,
it is evident that 
\bi
\item the first term $C_2\, D_N^{12}\, \tau(\truth)^{12}$ grows either slower or faster than $N^{1/2}$,
depending on the growth of $D_N$ and $\tau(\truth)$;\s
\item the second term $N^{1-\alpha}$ grows faster than $N^{1/2}$,
because $1 - \alpha > 1/2$ for all $\alpha \in [0, 1/2)$;\s  
\item the third term $D_N^3\, \tau(\truth)^{2}\, N$ grows faster than $N^{1/2}$,
because $D_N \geq 1$ and $\tau(\truth) \geq 1$ under Models 2 and 3.  
\ei
To bound the first term,
observe that the assumption
\beno
\label{eq:theta-star-bound}
\norm{\truth}_{\infty} &\leq& \dfrac{L + \vartheta \, \log N}{12\, (3 + D_N)} - \epsilon^\star
\ee
implies that the term $\tau(\truth)^{12}$ is bounded above by
\beno
\tau(\truth)^{12}
= \exp(12\, (3 + D_N)\, (\norm{\truth}_{\infty} + \epsilon^\star))
\leq \exp(L + \vartheta \, \log  N)
= \exp(L)\;  N^{\vartheta},
\ee
which in turn implies that the first term is bounded above by
\beno
C_2\; D_N^{12}\; \tau(\truth)^{12}
\lte C_2\, \exp(L)\; D_N^{12}\; N^{\vartheta}.
\ee
Since $\alpha \in [0, \, 1/2)$ and $\vartheta \in [0,\, 1/2 - \alpha)$ under Models 2 and 3,
the constant $\vartheta$ satisfies $\vartheta < 1/2$ while $D_N$ satisfies $D_N = O(\log N)$ by conditions \sone and \stwo \hspace{-.1cm}.
As a result,
the first term grows slower than $N^{1/2}$,
while the second and third term grow at least as fast as $N^{1/2}$.
Thus,
there exist constants $C_5 > 0$ and $N_2 \geq N_1$,
independent of $N$ and $p$,
such that,
for all $N > N_2$,
\beno
V
&\coloneqq& \dfrac{C_4\; D_N^3\; \tau(\truth)^2}{\min\{C_2\; D_N^{12}\; \tau(\truth)^{12},\; N^{1-\alpha},\; D_N^3\; \tau(\truth)^{2}\; N\}}\s\s
\\
\= \dfrac{C_4\; D_N^3\; \tau(\truth)^2}{C_2\; D_N^{12}\; \tau(\truth)^{12}}
\;\;=\;\; \dfrac{C_4}{C_2\; D_N^9\; \tau(\truth)^{10}} 
\;\;\leq\;\; C_5.
\ee
It is worth noting that $D_N$ and $\tau(\truth)$ may or may not increase as a function of $N$,
but both quantities are bounded below by $1$ under Models 2 and 3.

The above results show that,
for all $N > N_2 \geq 3$,
\beno
\mbP(\bX \in \mbH)
&\geq& 1 - \dfrac{2}{\max\{N,\, p\}^2},  
\ee
and,
for all $\bx \in \mbH$,  
$-\nabla_{\nat}^2 \; \widetilde\ell(\nat;\, \bx))^{-1}$ is invertible and
\beno
\label{eq:inf_mnorm_bound}
\widetilde\Lambda_N(\truth)
&\coloneqq& \max\limits_{\bx\, \in\, \mbH}\;\, \sup\limits_{\nat\, \in\, \mB_{\infty}(\truth,\, \epsilon^\star)}\, \mnorm{(- \nabla_{\nat}^2 \; \widetilde\ell(\nat;\, \bx))^{-1}}_{\infty}\s
\\
\hide{
\lte \max\left\{T_1,\, T_2\right\}
\\
}
\lte \max\left\{U\, V,\; U\right\}\s
\\
\lte \max\{1,\, C_5\}\; U\s
\\
\= \dfrac{C_3\; \max\{1,\, C_5\}\; D_N^{9}\; \tau(\truth)^{12}}{N^{1 - \alpha}}\s
\\
\lte \dfrac{C_3\; \max\{1,\, C_5\}\; \exp(L)\; D_N^{9}\; N^{\vartheta}}{N^{1 - \alpha}}\s
\\
\= \dfrac{C\, D_N^{9}}{N^{1 - (\alpha + \vartheta)}},
\ee
assuming $\alpha \in [0, \, 1/2)$ and $\vartheta \in [0,\, 1/2 - \alpha)$,
where the constant\break 
$C \coloneqq C_3\, \max\{1,\, C_5\}\, \exp(L) > 0$ is independent of $N$ and $p$. 
\qed
\hide{
Last,
but not least,
observe that
\beno
\hspace{.5in}
1 &\leq& D_N 
&\leq& \dfrac{L + \vartheta \, \log \, N}{12 \, \epsilon^\star} - 3
\ee
guarantees that $\norm{\truth}_{\infty} > 0$.
}

\s

\begin{lemma}
\label{lem:beta}
Consider Models 1, 2, and 3 with $\alpha \in [0,\, 1/2)$. 
Then,
for all $\bx \in \mbX$, 
\beno
\inf\limits_{\nat\, \in\, \mB_{\infty}(\truth,\, \epsilon^\star)}\; \lambda_{\min}(\bmat)
\;\geq\; \dfrac{N^{1-\alpha}}{12 \, \exp((3 + D_N) \, (\norm{\truth}_{\infty} + \epsilon^\star))^2} 
\;>\; 0,
\ee 
where $\lambda_{\min}(\bmat)$ is the smallest eigenvalue of $\bmat$.
\end{lemma} 

\s\s

\llproof \ref{lem:beta}.  
By definition, 
\beno
\widetilde\ell(\nat;\, \bx)
\= \dsum_{i < j}^{N} \,  \log \, \mbP_{\nat}(X_{i,j} = x_{i,j} \mid \bX_{-\{i,j\}} = \bx_{-\{i,j\}}),
&& \bx \in \mbX.
\ee
Note that the conditional distribution of edge variable $X_{i,j}$ 
conditional on the event $ \bX_{-\{i,j\}} = \bx_{-\{i,j\}}$ 
is an exponential-family distribution with sufficient statistic vector $s(\bx)$ 
and natural parameter vector $\nat$. 
Using standard properties of exponential families, 
it is straightforward to calculate,
for each pair of nodes $\{i, j\} \subset \mN$ and coordinates $(t,l) \in \{1, \ldots, N\}^2$:
\beno
&& -\dsum_{i < j}^{N}\, \dfrac{\partial}{\partial\theta_t\; \partial\theta_l}\, \log \mbP_{\nat}(X_{i,j} = x_{i,j} \mid \bX_{-\{i,j\}} = \bx_{-\{i,j\}})\s
\\
\= \dsum_{i < j}^{N} \cov_{\nat, \bx_{-\{i,j\}}}(s_t(\bX), s_l(\bX)), 
\ee
where $\cov_{\nat, \bx_{-\{i,j\}}}(s_t(\bX), s_l(\bX))$ denotes the conditional covariance of $s_t(\bX)$ and $s_l(\bX)$, 
computed with respect to the conditional distribution of $X_{i,j}$ given $\bX_{-\{i,j\}} = \bx_{-\{i,j\}}$. 
We have,
for all $\{i, j\} \subset \mN$ and $\bx_{-\{i,j\}} \in \{0, 1\}^{\binom{N}{2}-1}$,
\beno
\cov_{\nat, \bx_{-\{i,j\}}}(s_t(\bX), s_l(\bX))
\= \dsum_{h_1\, \in\, \mN\, \setminus\, \{t\}} \,
 \dsum_{h_2\, \in\, \mN\, \setminus\, \{l\}} \, 
\cov_{\nat, \bx_{-\{i, j\}}}(X_{t,h_1}, \, X_{l,h_2}),
\ee
as 
\beno
s_{t}(\bX) 
\= \dsum_{h \in \mN \setminus \{t\}} \, X_{t,h},
&& t \in \{1, \ldots, N\}. 
\ee
For each pair of nodes $\{i, j\} \subset \mN$,
we distinguish two cases:
\ben
\item If either $t \not\in \{i,j\}$ or $l \not\in \{i,j\}$,
then $s_t(\bX)$ and $s_l(\bX)$ cannot both be a function of $X_{i,j}$.
It then follows that,
conditional on $\bX_{-\{i,j\}} = \bx_{-\{i,j\}}$,
\beno
\cov_{\nat, \bx_{-\{i,j\}}}(s_t(\bX), s_l(\bX)) 
\= 0, 
\ee
as in this case either $s_t(\bX)$ or $s_l(\bX)$ will be almost surely constant. \s
\item If either $\{t,l\} = \{i,j\}$ or $t = l \in \{i,j\}$, 
then both $s_t(\bX)$ and $s_l(\bX)$ are functions of $X_{i,j}$.
Conditional on $\bX_{-\{i,j\}} = \bx_{-\{i,j\}}$,
edge variables $X_{a,b}$ corresponding to pairs of nodes $\{a,b\} \neq \{i,j\}$ are almost surely constant,
implying
\beno
\cov_{\nat, \bx_{-\{i,j\}}}(X_{t,h_1}, \, X_{l,h_2})
\= 0,
\ee
for all $\{t,h_1\} \neq \{i,j\}$ and all $\{l,h_2\} \neq \{i,j\}$.
We then have,
in the case $\{t,l\} = \{i,j\}$ ($t \neq l$),
that   
\beno
\cov_{\nat, \bx_{-\{i,j\}}}(s_t(\bX), s_l(\bX))
\= \var_{\nat, \bx_{-\{i,j\}}} \, X_{i,j},
\ee 
and in the case when $t = l \in \{i,j\}$, 
\beno
\cov_{\nat, \bx_{-\{i,j\}}}(s_t(\bX), s_l(\bX))
\= \var_{\nat, \bx_{-\{i,j\}}} \, s_t(\bX)
\= \var_{\nat, \bx_{-\{i,j\}}} \, X_{i,j}. 
\ee

\een
As a result,
for all $t \neq l \in \{1, \ldots, N\}$,
\be
\label{eq:off_diag}
\dsum_{i<j}^{N} \cov_{\nat, \bx_{-\{i,j\}}}(s_t(\bX), s_l(\bX))
\= \var_{\nat, \bx_{-\{t,l\}}} \, X_{t,l}
\ee
and all $t \in \{1, \ldots, N\}, $
\be
\label{eq:diag}
\dsum_{i<j}^{N} \var_{\nat, \bx_{-\{i,j\}}} \, s_t(\bX)
\= \dsum_{l\, \in\, \mN\, \setminus\, \{t\}} \, \var_{\nat, \bx_{-\{t,l\}}} \, X_{t,l}. 
\ee
An important consequence of \eqref{eq:off_diag} and \eqref{eq:diag} is that the matrix 
$\bmat$ given in \eqref{expected.hessian}
is diagonally balanced in the sense of \citet{HiWi15}.  
Observe that
\beno
\var_{\nat, \bx_{-\{i,j\}}} \, X_{i,j}
\= \mbP_{\nat}(X_{i,j} = 1 \,|\, \bX_{-\{i,j\}} = \bx_{-\{i,j\}}) \s \\  
&\times& (1 - \mbP_{\nat}(X_{i,j} = 1 \,|\, \bX_{-\{i,j\}} = \bx_{-\{i,j\}})). 
\ee 
Applying Lemma \ref{lem:bound_P}, 
for all $\bx_{-\{i,j\}} \in \{0,1\}^{\binom{N}{1}-1}$, 
\beno
\mbP_{\nat}(X_{i,j} = 1 \,|\, \bX_{-\{i,j\}} = \bx_{-\{i,j\}})
\gte \dfrac{N^{-\alpha}}{1 + \exp((3 + D_N) \norm{\nat}_\infty)}, \s\s \\
\mbP_{\nat}(X_{i,j} = 1 \,|\, \bX_{-\{i,j\}} = \bx_{-\{i,j\}}) 
&\leq& \dfrac{1}{1 + \exp(-(3 + D_N) \norm{\nat}_\infty)},
\ee
noting that $D_N = 0$ under Model 1, 
which implies that 
\beno
1 - \mbP_{\nat}(X_{i,j} = 1 \,|\, \bX_{-\{i,j\}} = \bx_{-\{i,j\}})
\gte \dfrac{1}{1 + \exp((3 + D_N) \, \norm{\nat}_{\infty})}. 
\ee
Thus, 
for all pairs of nodes $\{i,j\} \subset \mN$ and all $\bx_{-\{i,j\}} \in \{0, 1\}^{\binom{N}{2}-1}$,
\beno
\var_{\nat, \bx_{-\{i,j\}}} \, X_{i,j} 
&\geq& \dfrac{N^{-\alpha}}{(1 + \exp((3 + D_N) \, \norm{\nat}_{\infty}))^2}.  
\ee
As a result, 
each element $A_{t,l}(\nat, \bx)$ of $\bmat$ 
is bounded from below by 
\be
\label{eq:bmat_entry_bounds} 
A_{t,l}(\nat, \bx)
&\geq& \dfrac{N^{-\alpha}}{(1 + \exp((3 + D_N) \, \norm{\nat}_{\infty}))^2}
&>& 0.  
\ee
By invoking Lemma 2.1 of \citet{HiWi15} using the above bounds,
the smallest eigenvalue $\lambda_{\min}(\bmat)$ of the matrix $\bmat$ satisfies  
\beno
\lambda_{\min}(\bmat)
\gte \dfrac{N^{-\alpha} \, (N - 2)}{(1 + \exp((3 + D_N) \, \norm{\nat}_{\infty}))^2}
\gte \dfrac{N^{-\alpha} \, (N - 2)}{4 \, \exp((3 + D_N) \, \norm{\nat}_{\infty}))^2}.  
\ee 
Using the inequality $N - 2\, \geq\, N \,/\, 3$ (for $N \geq 3$),
we obtain,
for all $\bx \in \mbX$,  
\beno
\lambda_{\min}(\bmat)
\;\geq\; \dfrac{N^{-\alpha} \, (N - 2)}{4 \, \exp((3 + D_N) \, \norm{\nat}_{\infty}))^2}
\;\geq\; \dfrac{N^{1-\alpha}}{12 \, \exp((3 + D_N) \, \norm{\nat}_{\infty}))^2}. 
\ee
Since $\norm{\nat}_{\infty} \leq \norm{\truth}_{\infty} + \epsilon^\star$ for all $\nat \in \mB_{\infty}(\truth,\, \epsilon^\star)$,
we can conclude that,
for all $\bx \in \mbX$,
\beno
\inf\limits_{\nat\, \in\, \mB_{\infty}(\truth,\, \epsilon^\star)}\, \lambda_{\min}(\bmat)
&\geq& \dfrac{N^{1-\alpha}}{12 \, \exp((3 + D_N) \, (\norm{\truth}_{\infty} + \epsilon^\star))^2}.
\ee
\qed

\s

\begin{lemma}
\label{lem:ral}
Consider Models 2 and 3 with $\alpha \in [0,\, 1/2)$ and assume that $\truth \in \mR^p$ satisfies condition \eqref{condition.theta0} in Section \ref{sec:corollaries}.
Then there exists an integer $N_0 \geq 3$,
independent of $N$ and $p$, 
such that,
for all $N > N_0$ and all $\bx \in \mbH$,
\beno
\sup\limits_{\nat\, \in\, \mB_{\infty}(\truth,\, \epsilon^\star)}  
\dfrac{\cvec^{\top} \, \bmat^{-1} \, \cvec}{\cvec^{\top} \cvec}
\leq \dfrac{12 \, \exp((3 + D_N) \, (\norm{\truth}_{\infty} + \epsilon^\star))^2}{N^{1-\alpha}}.
\ee
\end{lemma} 

\llproof \ref{lem:ral}.
Based on \eqref{expected.hessian}, 
define 
\beno
R(\bmat^{-1},\, \cvec) 
&\coloneqq& \dfrac{\cvec^{\top}\,  \bmat^{-1} \, \cvec}{\cvec^{\top}\, \cvec},
\ee 
recognizing $R(\bmat^{-1},\, \cvec)$ to be 
the Rayleigh quotient of the matrix $\bmat^{-1} \in \mbR^{N \times N}$,  
assuming $\cvec \in \mR^N \setminus\, \bm{0}$,
where 
$\bm{0} \in \mR^N$ denotes the $N$-dimensional zero vector.
To bound 
$R(\bmat^{-1},\, \cvec)$,
note that if $\lambda_{1}, \ldots, \lambda_{N}$ are the eigenvalues of $\bmat$,
then $1 / \lambda_{1}$, $\ldots$, $1 / \lambda_{N}$ are the eigenvalues of $\bmat^{-1}$. 
Let $\lambda_{\min}(\bmat)$ be the smallest eigenvalue of $\bmat$, 
so that $1 / \lambda_{\min}(\bmat)$ is the largest eigenvalue of $\bmat^{-1}$.
Since the Rayleigh quotient of a matrix is bounded above by the largest eigenvalue of the matrix, 
Lemma \ref{lem:beta} shows that,
for all $\nat \in \mB_{\infty}(\truth,\, \epsilon^\star)$ and all $\bx \in \mbX$,  
\beno
R(\bmat^{-1},\, \cvec) 
&\leq& \dfrac{1}{\lambda_{\min}(\bmat)}\s\s
\\
&\leq& \dfrac{12 \, \exp((3+D_N) \, (\norm{\truth}_{\infty} + \epsilon^\star))^2}{N^{1-\alpha}},
\ee
which implies the bound 
\beno
\sup\limits_{\nat\, \in\, \mB_{\infty}(\truth,\, \epsilon^\star)}\, \dfrac{\cvec^{\top} \, \bmat^{-1} \, \cvec}{\cvec^{\top} \cvec}
\leq \dfrac{12 \, \exp((3 + D_N) \, (\norm{\truth}_{\infty} + \epsilon^\star))^2}{N^{1-\alpha}}.  
\ee
Last,
but not least,
we show that $\cvec \in \mbR^N \setminus \bm{0}$ for all $N > N_0$ and all $\bx \in \mbH$. 
By Lemma \ref{lem:ctc},
there exists an integer $N_0 \geq 3$,
independent of $N$ and $p$,
such that,
for all $N > N_0$ and all $\bx \in \mbH$,
\beno
\inf\limits_{\nat\, \in\, \mB_{\infty}(\truth,\, \epsilon^\star)}\, \cvec^{\top}\, \cvec
&\geq& \dfrac{N}{128\, (1 + \exp((3 + D_N) \, (\norm{\truth}_{\infty} + \epsilon^\star)))^6}\s\s
\\
&\geq& \dfrac{N}{128\, (2\, \exp((3 + D_N) \, (\norm{\truth}_{\infty} + \epsilon^\star)))^6}\s\s
\\
&=& \dfrac{N}{(64)\, (128)\, (\exp((3 + D_N) \, (\norm{\truth}_{\infty} + \epsilon^\star)))^6}.
\ee
The condition \eqref{condition.theta0} implies that
\beno
\norm{\truth}_{\infty} &\leq& \dfrac{L + \vartheta \, \log N}{12\, (3 + D_N)} - \epsilon^\star,
\ee 
which in turn implies that 
\beno
(\exp((3 + D_N) \, (\norm{\truth}_{\infty} + \epsilon^\star)))^6
&\leq& \exp(L / 2 + (\vartheta / 2)\, \log  N)\s\s
\\
&=& \exp(L / 2)\;  N^{\vartheta / 2}\s\s
\\
&<& \exp(L / 2)\;  N^{1 / 4}
\ee
using $\vartheta \in [0,\, 1/2)$.
As a result,
\beno
\inf\limits_{\nat\, \in\, \mB_{\infty}(\truth,\, \epsilon^\star)}\, \cvec^{\top}\, \cvec
&\geq& \dfrac{N}{(64)\, (128)\, (\exp((3 + D_N) \, (\norm{\truth}_{\infty} + \epsilon^\star)))^6}\s\s
\\
&\geq& \dfrac{N}{(64)\, (128)\, \exp(L / 2)\;  N^{1 / 4}}\s\s
\\
&=& \dfrac{N^{3/4}}{(64)\, (128)\, \exp(L / 2)},
\hide{
\\
&>& 0,
}
\ee
which implies that $\cvec \in \mbR^N \setminus\, \bm{0}$ for all $\nat \in \mB_{\infty}(\truth,\, \epsilon^\star)$ and all $\bx \in \mbH$.
\qed

\s

\begin{lemma}
\label{lem:ctc}
Consider Models 2 and 3.
Then there exists an integer $N_0 \geq 3$,
independent of $N$ and $p$, 
such that,
for all $N \geq N_0$ and all $\bx \in \mbH$,  
\beno
\inf\limits_{\nat\, \in\, \mB_{\infty}(\truth,\, \epsilon^\star)}\, \cvec^{\top}\, \cvec
&\geq& \dfrac{N}{128\, \exp((3 + D_N)\, (\norm{\truth}_{\infty} + \epsilon^\star))^6}.
\ee
\end{lemma}

\llproof \ref{lem:ctc}.
By \eqref{expected.hessian}, 
the coordinates of $\cvec$ are given by
\beno
c_t(\nat,\, \bx) \= \dsum_{i<j}^{N} \, 
\cov_{\nat, \bx_{-\{i,j\}}}(s_t(\bX),\, s_{N+1}(\bX)),
&& t \in \{1, \ldots, N\}.
\ee
Recall that 
\beno
s_t(\bX)
&\coloneqq& \dsum_{a\, \in\, \mN\, \setminus\, \{t\}} \, X_{t,a},
&& t \in \{1, \ldots, N\}. 
\ee
Then 
\beno
c_t(\nat,\, \bx) \= \dsum_{i<j}^{N}\,  \cov_{\nat, \bx_{-\{i,j\}}}(s_t(\bX),\, s_{N+1}(\bX))\s\s 
\\ 
\= \dsum_{i<j}^{N} \cov_{\nat, \bx_{-\{i,j\}}}\left( \dsum_{a\, \in\, \mN\, \setminus\, \{t\}} \, X_{t,a}, \, s_{N+1}(\bX) \right)\s\s 
\\ 
\= \dsum_{i<j}^{N}\; \dsum_{a\, \in\, \mN\, \setminus\, \{t\}} \,  \cov_{\nat, \bx_{-\{i,j\}}}(X_{t,a}, \, s_{N+1}(\bX)).  
\ee
Since $\cov_{\nat, \bx_{-\{i,j\}}}(X_{t,a}, \, s_{N+1}(\bX)) = 0$ almost surely
for all $\{i,j\} \neq \{t,a\}$,
\beno
\dsum_{i<j}^{N}\; \dsum_{a\, \in\, \mN\, \setminus\, \{t\}} \,  
\cov_{\nat, \bx_{-\{i,j\}}}(X_{t,a}, \, s_{N+1}(\bX))
\,=\, \dsum_{a\, \in\, \mN\, \setminus\, \{t\}} \,  \cov_{\nat, \bx_{-\{t,a\}}}(X_{t,a}, \, s_{N+1}(\bX)), 
\ee
as $\cov_{\nat, \bx_{-\{i,j\}}}$ is the conditional covariance operator with respect to the conditional distribution 
of $X_{i,j}$ given $\bX_{-\{i,j\}} = \bx_{-\{i,j\}}$,
implying $X_{t,a}$ is almost surely constant whenever $\{i,j\} \neq \{t,a\}$.  
Recall that 
\beno
s_{N+1}(\bX) 
&\coloneqq& \dsum_{i<j}^{N} \, X_{i,j} \, I_{i,j}(\bX),
\ee
where 
\beno
I_{i,j}(\bX) 
&\coloneqq& \one\left( \, \dsum_{h\, \in\, \mN_i\, \cap\, \mN_j} X_{i,h} \, X_{j,h} \,\geq\, 1 \right),
&& \{i,j\} \subset \mN.
\ee
Thus, 
we have 
\beno
\cov_{\nat, \bx_{-\{t,a\}}}(X_{t,a}, \, s_{N+1}(\bX))
\= \dsum_{i<j}^{N} \, \cov_{\nat, \bx_{-\{t,a\}}}(X_{t,a}, \, X_{i,j} \, I_{i,j}(\bX)), 
\ee
implying 
\beno
&& \dsum_{a\, \in\, \mN\, \setminus\, \{t\}} \, \cov_{\nat, \bx_{-\{t,a\}}}(X_{t,a}, \, s_{N+1}(\bX)) \s\\
\= \dsum_{a\, \in\, \mN\, \setminus\, \{t\}}\; \dsum_{i<j}^{N} \, \cov_{\nat, \bx_{-\{t,a\}}}(X_{t,a}, \, X_{i,j} \, I_{i,j}(\bX)). 
\ee
The FKG inequality
implies that
\beno
\cov_{\nat, \bx_{-\{t,a\}}}(X_{t,a}, \, X_{i,j} \, I_{i,j}(\bX)) \gte 0
&\mbox{for all}& \bx_{-\{t,a\}} \in \{0, 1\}^{\binom{N}{2}-1},
\ee
because the conditional covariance is computed with respect to the conditional distribution of $X_{t,a}$ and both $X_{t,a}$ and 
$X_{i,j} \, I_{i,j}(\bX)$ ($\{i,j\} \subset \mN$) are monotone non-decreasing functions of $X_{t,a}$.
As a result, 
\beno
\dsum_{a\, \in\, \mN\, \setminus\, \{t\}}\, \dsum_{i<j}^{N}  
\cov_{\nat, \bx_{-\{t,a\}}}(X_{t,a},  X_{i,j}  I_{i,j}(\bX))
\;\geq\; 
\cov_{\nat, \bx_{-\{t,a\}}}(X_{t,a}, X_{t,a} I_{t,a}(\bX)), 
\ee
for some $a\in \mN$ satisfying $\mN_a \cap \mN_t \neq \emptyset$.
Such a node $a \in \mN$ exists because each node $t \in \mN$ belongs to one or more subpopulations $\mA_k$ 
($k \in \{1, \dots, K\}$) and $|\mA_k| \geq 3$ for all $k \in \{1, \ldots, K\}$.  
We then obtain
\beno
\cvec^{\top}\, \cvec 
&\geq& \dsum_{t = 1}^{N} \, 
\left( \cov_{\nat, \bx_{-\{t,a\}}}(X_{t,a}, \, X_{t,a} \, I_{t,a}(\bX)) \right)^2. 
\ee
We can partition the sample space of $\bX_{-\{t,a\}} \in \{0, 1\}^{\binom{N}{2}-1}$ based on whether 
$ I_{t,a}(\bX) = 0$ or $I_{t,a}(\bX) = 1$.
When $I_{t,a}(\bX) = 0$, 
\beno
\cov_{\nat,\bx_{-\{t,a\}}}(X_{t,a}, \, X_{t,a} \, I_{t,a}(\bX))
\= \cov_{\nat,\bx_{-\{t,a\}}}(X_{t,a},\, 0) 
\= 0
\ee 
and when $I_{t,a}(\bX) = 1$,
\beno
\cov_{\nat,\bx_{-\{t,a\}}}(X_{t,a}, \, X_{t,a} \, I_{t,a}(\bX)) 
\= \var_{\nat, \bx_{-\{t,a\}}} \, X_{t,a}.
\ee 
By Lemma \ref{lem:bound_P}, 
for pairs $\{i,j\} \subset \mN$ with $\mN_i \cap \mN_j \neq \emptyset$,
we have the bounds   
\beno
\mbP_{\nat}(X_{i,j} = 1 \mid \bX_{-\{i,j\}} = \bx_{-\{i,j\}})
&\geq& \dfrac{1}{1 + \exp((3 + D_N) \, \norm{\nat}_{\infty})} 
\ee
and 
\beno
\mbP_{\nat}(X_{i,j} = 0 \mid \bX_{-\{i,j\}} = \bx_{-\{i,j\}})
&\geq& \dfrac{1}{1 + \exp((3 + D_N) \, \norm{\nat}_{\infty})}, 
\ee
for all $\bx_{-\{i,j\}} \in \{0, 1\}^{\binom{N}{2}-1}$. 
Using these bounds,
we obtain
\beno
\var_{\nat, \bx_{-\{t,a\}}} \, X_{t,a} 
&\geq& \left( \dfrac{1}{1 + \exp((3+ D_N) \, \norm{\nat}_{\infty})} \right)^2, 
\ee 
which shows that 
\beno
\cov_{\nat,\bx_{-\{t,a\}}}(X_{t,a}, \, X_{t,a} \, I_{t,a}(\bX)) 
&\geq& \dfrac{I_{t,a}(\bX)}{(1 + \exp((3+ D_N) \, \norm{\nat}_{\infty}))^2}. 
\ee
For all $\nat \in \mB_\infty(\truth,\, \epsilon^\star)$,
$\norm{\nat}_{\infty} \leq \norm{\truth}_{\infty} + \epsilon^\star$, 
which implies that
\beno
\dfrac{1}{1 + \exp((3 + D_N) \, \norm{\nat}_{\infty})} 
&\geq& \dfrac{1}{1 + \exp((3 + D_N) \, (\norm{\truth}_{\infty} + \epsilon^\star))}. 
\ee
We then obtain the following lower bound: 
\beno
\cvec^{\top}\, \cvec
&\geq& \dsum_{t = 1}^{N}\, \dfrac{\max_{a \in \mN \setminus \{t\}}\, I_{t,a}(\bx)}{(1 + \exp((3 + D_N)\, (\norm{\truth}_{\infty} + \epsilon^\star))^4}\s
\\
\= \dfrac{\sum_{t = 1}^{N} \norm{\bI_t(\bx)}_{\infty}}{(1 + \exp((3 + D_N)\, (\norm{\truth}_{\infty} + \epsilon^\star))^4},
\ee
where $\bI_t(\bx) \coloneqq (I_{t,1}(\bx),\ldots, I_{t,t-1}(\bx),\, I_{t,t+1}(\bx), \ldots,  I_{t,N}(\bx))$ ($t \in \mN$).  
By definition of $\mbH$,
\beno 
\mbH
&\coloneqq& \left\{\bx \in \mbX \;:\; \dsum_{i=1}^{N}\, \norm{\bI_i(\bx)}_{\infty}\, \geq\, \dfrac{N}{2\, (1 + \exp((3 + D_N) \, (\norm{\truth}_{\infty} + \epsilon^\star)))^2}
\right\},
\ee
all $\bx \in \mbH$ satisfy
\beno 
\cvec^{\top}\, \cvec
&\geq& \dfrac{N}{2 \, (1 + \exp((3+ D_N) \, (\norm{\truth}_{\infty} + \epsilon^\star))^6} \s\\
&\geq& \dfrac{N}{2 \, (2 \, \exp((3+ D_N) \, (\norm{\truth}_{\infty} + \epsilon^\star))^6} \s\\
&\geq& \dfrac{N}{128 \, \exp((3+ D_N) \, (\norm{\truth}_{\infty} + \epsilon^\star))^6}.  
\ee
\qed

\s

\begin{lemma}
\label{lem:vctc}
Consider Models 2 and 3.
Then there exists an integer $N_0 \geq 3$,
independent of $N$ and $p$, 
such that,
for all $N \geq N_0$ and all $\bx \in \mbH$,
\beno
\sup\limits_{\nat\, \in\, \mB_{\infty}(\truth, \, \epsilon^\star)} \, \norm{\cvec}_{\infty} 
\lte 3 \, D_N^3
\ee 
\vspace{-.2cm}
and 
\vspace{-.1cm}
\beno
\inf\limits_{\nat\, \in\, \mB_{\infty}(\truth,\, \epsilon^\star)}\; \dfrac{\vvar}{\cvec^{\top}\, \cvec}
&\geq& \dfrac{1}{576\; D^6_N\, \exp((3 + D_N)\, (\norm{\truth}_{\infty} + \epsilon^\star))^4},
\ee
noting that $D_N \geq 1$ under Models 2 and 3. 
\end{lemma}

\s\s

\llproof \ref{lem:vctc}.
Recall that $s_{N+1}(\bX)$ is defined by 
\beno
s_{N+1}(\bX) 
&\coloneqq& \dsum_{i<j}^{N} \, X_{i,j} \, I_{i,j}(\bX),
\ee
where 
\beno
I_{i,j}(\bX) \;=\; \one\left(\dsum_{h\, \in\, \mN_i\, \cap\, \mN_j} \, X_{i,h} \, X_{j,h} \,\geq\, 1 \right),
&& \{i,j\} \subset \mN.
\ee
According to \eqref{expected.hessian}, 
$\vvar$ is given by 
\beno
\vvar 
&\coloneqq& \dsum_{i<j}^{N} \, \var_{\nat, \bx_{-\{i,j\}}} \, s_{N+1}(\bX)\s\s
\\
\= \dsum_{i<j}^{N} \, \var_{\nat, \bx_{-\{i,j\}}} \left(\dsum_{a<b}^{N} \, X_{a,b} \, I_{a,b}(\bX)\right). 
\ee
Given any pair of nodes $\{i,\, j\} \subset \mN$, 
the FKG inequality implies,
for all pairs of nodes $\{a,\, b\} \subset \mN$ and $\{r,\, t\} \subset \mN$, 
that  
\beno
\cov_{\nat, \bx_{-\{i,j\}}}(X_{a,b}\, I_{a,b}(\bX),\; X_{r,t}\, I_{r,t}(\bX))
&\geq& 0,
\ee
because the conditional covariance is computed with respect to the conditional distribution of $X_{i,j}$ 
and each $X_{a,b} \, I_{a,b}(\bX)$ ($\{a,b\} \subset \mN$)
is a monotone non-decreasing function of $X_{i,j}$. 
Thus, 
\beno
\vvar
&=& \dsum_{i<j}^{N} \, \var_{\nat, \bx_{-\{i,j\}}}\left(\dsum_{a<b}^{N} \, X_{a,b} \, I_{a,b}(\bX)\right)\s
\\
\gte \dsum_{i<j}^{N} \, \dsum_{a<b}^{N} \, \var_{\nat, \bx_{-\{i,j\}}}\left(X_{a,b} \, I_{a,b}(\bX)\right) \s 
\\ 
\gte \dsum_{i<j}^{N} \,\var_{\nat, \bx_{-\{i,j\}}}\left(X_{i,j}\, I_{i,j}(\bX)\right) \s\\
\= \dsum_{i<j \,:\, \mN_i \,\cap\, \mN_j \neq \emptyset} \, 
\var_{\nat, \bx_{-\{i,j\}}}\left(X_{i,j} \, I_{i,j}(\bX)\right),  
\ee
noting that $I_{i,j}(\bX) = 0$ almost surely when $\mN_i \,\cap\, \mN_j = \emptyset$. 
Observe that $\var_{\nat, \bx_{-\{i,j\}}}(X_{i,j}\, I_{i,j}(\bX))$ can be simplified as follows:
\beno
\var_{\nat, \bx_{-\{i,j\}}}(X_{i,j}\, I_{i,j}(\bX))
&=& I_{i,j}(\bx)\; \var_{\nat, \bx_{-\{i,j\}}}\, X_{i,j},
\hide{
\= \begin{cases}
0 & I_{i,j}(\bx) = 0 \s\\
\var_{\nat, \bx_{-\{i,j\}}} \, X_{i,j} & I_{i,j}(\bX) = 1
\end{cases},
}
\ee 
because $I_{i,j}(\bX)$ is a function of $\bX_{-\{i,j\}}$ but is not a function of $X_{i,j}$,
and therefore $I_{i,j}(\bX)$ is almost surely constant conditional on $\bX_{-\{i,j\}} = \bx_{-\{i,j\}}$.
Hence,
\beno
\vvar
&\geq& \dsum_{i<j \,:\, \mN_i \,\cap\, \mN_j \,\neq\, \emptyset}^{N} \, 
I_{i,j}(\bx) \, \var_{\nat, \bx_{-\{i,j\}}} \, X_{i,j}. 
\ee
Using Lemma \ref{lem:bound_P} shows that,
for all $\bx_{-\{i,j\}} \in \{0, 1\}^{\binom{N}{2}-1}$
and $\{i,j\} \subset \mN$ satisfying $\mN_i \cap \mN_j \neq \emptyset$,
\beno
\mbP_{\nat}(X_{i,j} = 1 \mid \bX_{-\{i,j\}} = \bx_{-\{i,j\}})
&\geq& \dfrac{1}{1 + \exp((3 + D_N) \, \norm{\nat}_{\infty})}
\ee
and
\beno
\mbP_{\nat}(X_{i,j} = 0 \mid \bX_{-\{i,j\}} = \bx_{-\{i,j\}})
&\geq& \dfrac{1}{1 + \exp((3 + D_N) \, \norm{\nat}_{\infty})},
\ee
implying
\be
\label{eq:var_bound_}
\var_{\nat, \bx_{-\{i,j\}}} \, X_{i,j} 
&\geq& \dfrac{1}{(1 + \exp((3 + D_N) \, \norm{\nat}_{\infty}))^2}.
\ee
Therefore,
\beno
\vvar
&\geq& \dsum_{i<j:\; \mN_i\, \cap\, \mN_j\, \neq\, \emptyset}^{N} \, 
\dfrac{I_{i,j}(\bx)}{(1 + \exp((3 + D_N) \, \norm{\nat}_{\infty}))^2}\s\s
\\
&\geq& \dsum_{i \in \mN} \, \dsum_{j \in \mN \setminus \{i\}} \, 
\dfrac{I_{i,j}(\bx)}{2 \, (1 + \exp((3 + D_N) \, \norm{\nat}_{\infty}))^2} \s\s\\
&\geq& \dfrac{\sum_{i \in \mN} \max_{j \in \mN \setminus \{i\}} \, I_{i,j}(\bx)}{2 \, (1 + \exp((3 + D_N) \, \norm{\nat}_{\infty}))^2}\s\s
\\
\= \dfrac{\sum_{i \in \mN} \norm{\bI_i(\bx)}_{\infty}}{2 \, (1 + \exp((3 + D_N) \, \norm{\nat}_{\infty}))^2}\s\s
\\
&\geq& \dfrac{\sum_{i \in \mN} \norm{\bI_i(\bx)}_{\infty}}{2 \, (1 + \exp((3 + D_N)\, (\norm{\truth}_{\infty} + \epsilon^\star)))^2},
\ee 
defining $\bI_i(\bx) \coloneqq (I_{i,1}(\bx), \ldots,  I_{i,i-1}(\bx),\, I_{i,i+1}(\bx), \ldots,  I_{i,N}(\bx))$ for each $i \in \mN$ and using $\norm{\nat}_{\infty} \leq \norm{\truth}_{\infty} + \epsilon^\star$ for all $\nat \in \mB_{\infty}(\truth,\, \epsilon^\star)$.
The second inequality follows because $I_{i,j}(\bx) = 0$ ($\bx \in \mbX$)
for all $\{i,j\} \subset \mN$ satisfying $\mN_i \cap \mN_j = \emptyset$ 
and by noting that $I_{i,j}(\bx) = I_{j,i}(\bx)$ ($\{i,j\} \subset \mN$). 
By definition of $\mbH$,
\beno
\mbH
&\coloneqq& \left\{\bx \in \mbX:\; \dsum_{i=1}^{N}\, \norm{\bI_i(\bx)}_{\infty}\; \geq\; \dfrac{N}{2\, (1 + \exp((3 + D_N) \, (\norm{\truth}_{\infty} + \epsilon^\star)))^2}
\right\}.
\ee
Thus,
for all $N \geq N_0$,
all $\nat \in \mB_{\infty}(\truth,\, \epsilon^\star)$ and all $\bx \in \mbH$,  
\beno
\vvar
&\geq& \dfrac{N}{4\; (1 + \exp((3 + D_N)\, (\norm{\truth}_{\infty} + \epsilon^\star)))^4}.  
\ee

We proceed with $\cvec^{\top}\cvec$. 
Lemma \ref{lem:ctc} establishes that 
\beno
c_t(\nat,\, \bx)  
\= \dsum_{a\, \in\, \mN\, \setminus\, \{t\}} 
\, \cov_{\nat, \bx_{-\{t,a\}}}(X_{t,a}, \, s_{N+1}(\bX)) \s\\
\= \dsum_{a\, \in\, \mN\, \setminus\, \{t\}} \, \dsum_{i<j}^{N} \,
\cov_{\nat, \bx_{-\{t,a\}}}(X_{t,a}, \, X_{i,j} \, I_{i,j}(\bX)) \s\\
\= \dsum_{a\, \in\, \mN\, \setminus\, \{t\} \,:\, \mN_a\, \cap\, \mN_t\, \neq\, \emptyset}\; \dsum_{i<j}^{N} \,
\cov_{\nat, \bx_{-\{t,a\}}}(X_{t,a}, \, X_{i,j} \, I_{i,j}(\bX)), 
\ee
noting that, 
by Proposition \ref{prop:cond_ind}, 
$X_{t,a}$ is independent of all other edge variables
when $\mN_t \cap \mN_a = \emptyset$, 
in which case $\cov_{\nat, \bx_{-\{t,a\}}}(X_{t,a}, \, X_{i,j} \, I_{i,j}(\bX)) = 0$.
Hence,
\beno
&& \dsum_{a\, \in\, \mN\, \setminus\, \{t\} \,:\, \mN_a\, \cap\, \mN_t\, \neq\, \emptyset} \; \dsum_{i<j}^{N} \,
\cov_{\nat, \bx_{-\{t,a\}}}(X_{t,a}, \, X_{i,j} \, I_{i,j}(\bX)) \s\\ 
&\leq& D_N^2 \, \left(\max\limits_{a\, \in\, \mN\, \setminus\, \{t\}} \; \dsum_{i<j}^{N} \,
\cov_{\nat, \bx_{-\{t,a\}}}(X_{t,a}, \, X_{i,j} \, I_{i,j}(\bX)) \right).
\ee
This bound follows from Lemma \ref{lem:DN}, 
which shows that, 
for all $t \in \mN$,
\beno
\left|\left\{a \in \mN \setminus \{t\} \,:\, \mN_a \cap \mN_t \neq \emptyset \right\} \right|
&\leq& D_N^2. 
\ee
If $\mN_t \cap \mN_a \neq \emptyset$, 
then $\cov_{\nat, \bx_{-\{t,a\}}}(X_{t,a}, \, X_{i,j} \, I_{i,j}(\bX)) = 0$ if
\ben 
\item $\{i,j\} \neq \{t,a\}$, 
in which case $X_{i,j}$ is constant almost surely, and \s 
\item 
$I_{i,j}(\bX)$ is constant in $X_{t,a}$,
implying $I_{i,j}(\bX)$ is constant almost surely. 
\een 
The justification for the above statements regarding constancy is due to the fact that  
the conditional covariance $\cov_{\nat, \bx_{-\{t,a\}}}(X_{t,a}, \, X_{i,j} \, I_{i,j}(\bX))$ 
is computed
with respect to the conditional distribution of $X_{t,a}$ conditional on $\bX_{-\{t,a\}} = \bx_{-\{t,a\}}$. 
It is therefore enough to bound
\bi
\item the number of pairs $\{i,j\} \subset \mN$ which do not satisfy either point 1. or 2. above, 
for a  given $\{t,a\} \subset \mN$, and \s 
\item the quantity $\mbE \, \cov_{\nat, \bx_{-\{t,a\}}}(X_{t,a}, \, X_{i,j} \, I_{i,j}(\bX))$. 
\ei
Since $\cov_{\nat, \bx_{-\{t,a\}}}(X_{t,a}, \, X_{i,j} \, I_{i,j}(\bX)) \leq 1$,
we focus on the bounding the number of pairs $\{i,j\} \subset \mN$ for which
$X_{i,j} \, I_{i,j}(\bX)$ is a function of $X_{t,a}$: 
\bi 
\item First, 
$\{i,j\} = \{t,a\}$ for only one pair $\{i,j\} \subset \mN$. \s  
\item Second, 
\beno
I_{i,j}(\bX)
\= \one\left( \, \dsum_{h \in \mN_i \,\cap\, \mN_j} \, X_{i,h} \, X_{j,h} \,\geq\, 1 \right),
&& \{i,j\} \neq \{t,a\},  
\ee
is a function of $X_{t,a}$ 
if and only if one of the following holds: 
\ben
\item $\{i,j\} \cap \{t,a\} = a$ and $t \in \mN_i \cap \mN_j$, or  
\item $\{i,j\} \cap \{t,a\} = t$ and $a \in \mN_i \cap \mN_j$. 
\een
In either case, 
the number of possible pairs $\{i,j\} \neq \{t,a\}$ is bounded above by $\max\{|\mN_i|, \, |\mN_j|\} \leq D_N$
by Lemma \ref{lem:DN},
and hence is bounded above by $2 \, D_N$ in either case.  
\ei
Recalling $D_N \geq 1$ under Models 2 and 3,
we have the bound  
\beno
\dsum_{a\, \in\, \mN\, \setminus\, \{t\}} \, \dsum_{i<j}^{N} \,
\cov_{\nat, \bx_{-\{t,a\}}}(X_{t,a}, \, X_{i,j} \, I_{i,j}(\bX))
&\leq& D_N^2 \, (1 + 2 \, D_N)
&\leq& 3 \, D_N^3, 
\ee
which shows,
for all $\nat \in \mB_{\infty}(\truth,\, \epsilon^\star)$ and all $\bx \in \mbX$,
that  
\beno
\norm{\cvec}_{\infty} 
\lte 3\, D_N^3
\ee
and
\beno
\cvec^{\top} \cvec 
\lte N\, (3 \, D_N^3)^2 
&=& 9 \, D_N^6\, N.
\ee 

Collecting terms reveals that,
for all $\nat \in \mB_{\infty}(\truth,\, \epsilon^\star)$ and all $\bx \in \mbH$,  
\beno
\dfrac{\vvar}{\cvec^{\top}\, \cvec}
&\geq& \dfrac{N}{4\; (1 + \exp((3 + D_N) \, (\norm{\truth}_{\infty} + \epsilon^\star))^4}
 \left(\dfrac{1}{9\, D_N^6 \, N} \right)\s
\\ 
&\geq& \dfrac{1}{36\, D_N^6 \, (1 + \exp((3 + D_N) \, (\norm{\truth}_{\infty} + \epsilon^\star)))^4}\s
\\
&\geq& \dfrac{1}{36\, D_N^6 \, (2 \, \exp((3 + D_N) \, (\norm{\truth}_{\infty} + \epsilon^\star)))^4}\s
\\
\= \dfrac{1}{576\, D_N^6 \, \exp((3 + D_N) \, (\norm{\truth}_{\infty} + \epsilon^\star))^4}.
\ee
In conclusion, 
for all $N \geq N_0$ and all $\bx \in \mbH$,  
\beno
\sup\limits_{\nat\, \in\, \mB_{\infty}(\truth, \, \epsilon^\star)} \,
\norm{\cvec}_{\infty}
&\leq& 3 \, D_N^3
\ee
and 
\beno
\inf\limits_{\nat\, \in\, \mB_{\infty}(\truth, \, \epsilon^\star)} \, 
\dfrac{\vvar}{\cvec^{\top}\, \cvec}
\,\geq\, \dfrac{1}{576\; D^6_N \, \exp((3 + D_N) \, (\norm{\truth}_{\infty} + \epsilon^\star))^4}.
\ee
\qed

\begin{lemma}
\label{lem:Iscale}
Consider Models 2 and 3 and assume that $\truth \in \mR^p$ satisfies condition \eqref{condition.theta0} in Section \ref{sec:corollaries} and either condition \sone or condition \stwo is satisfied.
Then there exists an integer $N_0 \geq 3$,
independent of $N$ and $p$,
such that,
for all $N > N_0$, 
\beno
\mbP(\bX \in \mbH)
&\geq& 1 - \dfrac{2}{\max\{N,\, p\}^2},
\ee
where
\beno
\mbH
&\coloneqq& \left\{\bx \in \mbX:\; \dsum_{i=1}^{N}\, \norm{\bI_i(\bx)}_{\infty}\; \geq\; \dfrac{N}{2\, (1 + \exp((3 + D_N) \, (\norm{\truth}_{\infty} + \epsilon^\star)))^2}
\right\}.
\ee
\end{lemma}

\s

\llproof \ref{lem:Iscale}. 
Define
\beno
A(\bX)
&\coloneqq& \dsum_{i=1}^{N}\, \norm{\bI_{i}(\bX)}_{\infty}.
\ee
We prove Lemma \ref{lem:Iscale} as follows.
First,
we show that there exists an integer $N_0 \geq 3$,
independent of $N$ and $p$,
such that,
for all $N > N_0$,
\beno
\dfrac12\; \mbE\, A(\bX)
&\geq& \dfrac{N}{2\, (1 + \exp((3 + D_N) \, (\norm{\truth}_{\infty} + \epsilon^\star)))^2}.
\ee
Second,
we prove that
\beno
\mbP\left(\left|A(\bX) - \mbE\, A(\bX)\right|
\;<\; \dfrac{1}{2}\; \mbE\, A(\bX)\right) 
&\geq& 1 - \dfrac{2}{\max\{N,\, p\}^2},
\ee
which implies that the events
\beno
A(\bX) 
&\in& \left(\dfrac{1}{2}\; \mbE\, A(\bX),\; \dfrac{3}{2}\; \mbE\, A(\bX)\right)
\ee
and
\beno
A(\bX)
&>& \dfrac{1}{2}\; \mbE\, A(\bX)
&\geq& \dfrac{N}{2\, (1 + \exp((3 + D_N) \, (\norm{\truth}_{\infty} + \epsilon^\star)))^2}
\ee
occur with probability at least $1 - 2\, / \max\{N,\, p\}^2$,
proving the desired result.
\hide{
\beno
&& \mbP\left(\left|\dsum_{i=1}^{N}\, \norm{\bI_{i}(\bX)}_{\infty} - \mbE\, \dsum_{i=1}^{N}\, \norm{\bI_{i}(\bX)}_{\infty} \right|
\;<\; \dfrac{1}{2}\; \mbE\, \dsum_{i=1}^{N}\, \norm{\bI_{i}(\bX)}_{\infty}\right)\s
\\
&\geq& 1 - \dfrac{2}{\max\{N,\, p\}^2}
\ee
and
\beno
\dfrac12\; \mbE\, \dsum_{i=1}^N\, \norm{\bI_{i}(\bX)}_{\infty}
&\geq& \dfrac{N}{2\, (1 + \exp((3 + D_N) \, (\norm{\truth}_{\infty} + \epsilon^\star)))^2},
\ee
which establishes the desired result:
\beno
\mbP(\bX \in \mbH) 
&=& \mbP\left( \dsum_{i=1}^{N}\, \norm{\bI_{i}(\bX)}_{\infty} 
\;\geq\; \dfrac{N}{2\, (1 + \exp((3 + D_N) \, (\norm{\truth}_{\infty} + \epsilon^\star)))^2} \right) \s\\
&\geq& 1 - \dfrac{2}{\max\{N,\, p\}^2}.
\ee
}

\s

{\bf Bounding $(1 / 2)\, \mbE \sum_{i=1}^N\, \norm{\bI_{i}(\bX)}_{\infty}$ from below.} 
We first observe that,
for all $\bx \in \mbX$,  
\be
\label{eq:first_bound}
\dsum_{i=1}^N\, \norm{\bI_{i}(\bx)}_{\infty}
&\geq& \dsum_{i=1}^N \, I_{i,a_i}(\bx)
& \mbox{for all } a_i \in \mN \setminus \{i\} & (i \in \mN).  
\ee
Consider any $\nat \in \mB_\infty(\truth,\, \epsilon^\star)$ and note that each node $i \in \mN$ belongs to at least subpopulation $\mA_k$ ($k \in \{1, \dots, K\}$).
Since $\min_{1 \leq k \leq K}\, |\mA_k| \geq 3$,
there exists,
for any given node $i \in \mN$,
at least one other node $a_i \in \mN \setminus \{i\}$ such that $\mN_i \,\cap\, \mN_{a_i} \neq \emptyset$.
In addition,
there exists a node $b \in \mN_i\, \cap\, \mN_{a_i}$ so that 
\beno
\mbP_{\nat}(I_{i,a_i}(\bX) = 1) 
&\geq& \mbP_{\nat}(X_{i,b} \, X_{a_i,b} = 1),
\ee 
because the event $X_{i,b} \, X_{a_i,b} = 1$ implies the event $I_{i,a_i}(\bX) = 1$.
Thus, 
\beno
\mbE\, \dsum_{i=1}^{N} \, \norm{\bI_{i}(\bX)}_{\infty}
\gte \mbE \, \dsum_{i=1}^{N} \, I_{i,a_i}(\bX)
&=& \dsum_{i=1}^{N}\, \mbP_{\nat}(I_{i,a_i}(\bX) = 1)\s
\\
&& &\geq&  \dsum_{i=1}^{N} \, \mbP_{\nat}(X_{i,b} \, X_{a_i,b} = 1).  
\ee
Using Lemma \ref{lem:bound_P},
we obtain,
for all $\{i, j\} \subset \mN$ satisfying $\mN_i \cap \mN_j \neq \emptyset$ and all $\bx_{-\{i,j\}} \in \{0, 1\}^{\binom{N}{2}-1}$, 
the bounds 
\beno
\mbP_{\nat}(X_{i,j} = 0 \mid \bX_{-\{i,j\}} = \bx_{-\{i,j\}})
&\geq& \dfrac{1}{1 + \exp((3 + D_N) \, \norm{\nat}_{\infty})}\s\s
\\
\mbP_{\nat}(X_{i,j} = 1 \mid \bX_{-\{i,j\}} = \bx_{-\{i,j\}})
&\geq& \dfrac{1}{1 + \exp((3 + D_N) \, \norm{\nat}_{\infty})}.
\ee
As a result, 
\beno
&& \mbP_{\nat}(X_{i,b}\; X_{a_i,b} = 1)
\;=\; \mbP_{\nat}(X_{i,b} = 1 \,|\, X_{a_i,b} = 1) \, \mbP_{\nat}(X_{a_i,b} = 1)\s
\\
&\geq& \left(\min\limits_{\{i,\, j\}\, \subset\, \mN}\;\, \min\limits_{\bx_{-\{i,j\}}\, \in\, \{0, 1\}^{\binom{N}{2}-1}}\; 
\mbP_{\nat}(X_{i,j} = 1 \mid \bX_{-\{i,j\}} = \bx_{-\{i,j\}})\right)^2\s
\\
&\geq& \left( \dfrac{1}{1 + \exp((3 + D_N) \, \norm{\truth}_{\infty})}\right)^2\s
\\
&\geq& \dfrac{1}{(1 + \exp((3 + D_N) \, (\norm{\truth}_{\infty} + \epsilon^\star))^2},
\ee
because $\norm{\nat}_{\infty} \leq \norm{\truth}_{\infty} + \epsilon^\star$ for all $\nat \in \mB_\infty(\truth,\, \epsilon^\star)$. 
Using condition \eqref{condition.theta0} in Section \ref{sec:corollaries} with $\vartheta \in [0,\, 1/2 - \alpha)$,
we obtain
\beno
(1 + \exp((3 + D_N)\, (\norm{\truth}_{\infty} + \epsilon^\star)))^2
&\leq& (2\, \exp((3 + D_N)\, (\norm{\truth}_{\infty} + \epsilon^\star)))^2\s\s
\\
\hide{
&=& 4\, \exp(2\, (3 + D_N)\, (\norm{\truth}_{\infty} + \epsilon^\star))\s\s
\\
}
&\leq& 4\, \exp\left(2\, (3 + D_N)\, \dfrac{L + \vartheta\, \log N}{14\, (3 + D_N)}\right)\s\s
\\
\hide{
&=& 4\, \exp\left(\dfrac{L + \vartheta\, \log N}{7}\right)\s\s
\\
}
&=& 4\, \exp(L / 7)\, N^{\vartheta / 7},
\ee
which implies that
\beno
\mbP_{\nat}(X_{i,b}\; X_{a_i,b} = 1)
&\geq& \dfrac{1}{(1 + \exp((3 + D_N) \, (\norm{\truth}_{\infty} + \epsilon^\star))^2}\s
\\
&\geq& \dfrac{1}{4\, \exp(L / 7)\, N^{\vartheta / 7}}.
\ee
We have thus demonstrated that
\beno
\dfrac12\; \mbE\, \dsum_{i=1}^{N} \, \norm{\bI_{i}(\bX)}_{\infty}
&\geq& \dfrac12\; \dsum_{i=1}^{N} \, \mbP_{\nat}(X_{i,b} \, X_{a_i,b} = 1)\s
\\
&\geq& \dfrac{N}{2\, (1 + \exp((3 + D_N) \, (\norm{\truth}_{\infty} + \epsilon^\star)))^2}\s
\\
&\geq& \dfrac{N}{8\, \exp(L / 7)\, N^{\vartheta / 7}}.
\ee

\s

{\bf Concentrating $\sum_{i=1}^{N} \norm{\bI_{i}(\bX)}_{\infty}$.}
Applying Lemma \ref{theorem.mle.lemma} with 
\beno
s(\bX)
&\coloneqq&  \dsum_{i=1}^{N} \, \norm{\bI_{i}(\bX)}_{\infty}
&&\mbox{and}&& 
\bmu(\truth) 
&\coloneqq& \mbE \, \dsum_{i=1}^{N} \, \norm{\bI_{i}(\bX)}_{\infty}
\ee
shows that,
for all $t > 0$,  
\beno
\mbP\left(\left|\dsum_{i=1}^{N} \norm{\bI_{i}(\bX)}_{\infty} 
- \mbE \dsum_{i=1}^{N} \norm{\bI_{i}(\bX)}_{\infty}\right| < t \right)
\geq 1 - 2 \, \exp\left( - \dfrac{2 \, t^2}{\mnorm{\mD_N(\truth)}_2^2 \; \Psi_N^2} \right). 
\ee
Choosing
\beno
t 
&=& \sqrt{\log \max\{N,\, p\}} \, \mnorm{\mD_N(\truth)}_2 \, \Psi_N
&>& 0
\ee
gives
\beno
\mbP\left(\left|\dsum_{i=1}^{N} \norm{\bI_{i}(\bX)}_{\infty} - \mbE \dsum_{i=1}^{N} \norm{\bI_{i}(\bX)}_{\infty}\right| < \sqrt{\log \max\{N,\, p\}}\, \mnorm{\mD_N(\truth)}_2 \, \Psi_N\right)\s\s
\\
\geq\; 1 - \dfrac{2}{\max\{N,\, p\}^2}. 
\ee
We will demonstrate below that there exists an integer $N_0 \geq 3$,
independent of $N$ and $p$,
such that,
for all $N > N_0$,
\be
\label{lem_key1}
\dfrac12\; \mbE \, \dsum_{i=1}^{N}\, \norm{\bI_{i}(\bX)}_\infty
&\geq& \dfrac{N}{8\, \exp(L / 7)\, N^{\vartheta / 7}}\s\s
\\
&\geq& \sqrt{\log \max\{N,\, p\}}\, \mnorm{\mD_N(\truth)}_2\, \Psi_N,
\ee
which implies that
\beno
&& \mbP\left(\left|\dsum_{i=1}^{N} \norm{\bI_{i}(\bX)}_{\infty} - \mbE \dsum_{i=1}^{N} \norm{\bI_{i}(\bX)}_{\infty}\right| < \dfrac12\; \mbE \, \dsum_{i=1}^{N}\, \norm{\bI_{i}(\bX)}_\infty\right)\s\s
\\
&\geq& \mbP\left(\left|\dsum_{i=1}^{N} \norm{\bI_{i}(\bX)}_{\infty} - \mbE \dsum_{i=1}^{N} \norm{\bI_{i}(\bX)}_{\infty}\right| < \sqrt{\log \max\{N,\, p\}}\, \mnorm{\mD_N(\truth)}_2 \, \Psi_N\right)\s\s
\\
\gte 1 - \dfrac{2}{\max\{N,\, p\}^2}. 
\ee
\hide{
by bounding the terms $\sqrt{\log \max\{N,\, p\}}$,\,
$\Psi_N$,\,
and $\mnorm{\mD_N(\truth)}_2$ one by one.
}

\s

{\bf Bounding $\sqrt{\log \max\{N,\, p\}}$ from above.}
Since $p = N + 1$ under Models 2 and 3,
we obtain,
for all $N \geq 2$,
\beno
\sqrt{\log \max\{N,\, p\}} 
\,=\, \sqrt{\log(N + 1)} 
\,\leq\, \sqrt{\log 2\, N} 
\,\leq\, \sqrt{2\, \log N} 
\,\leq\, 2\, \sqrt{\log N}.
\ee 

\s

{\bf Bounding $\Psi_N$ from above.}
Lemma \ref{lem:brokerage_ham_bound_2} shows that,
for each pair of nodes $\{a, b\} \subset \mN$ with $\mN_a\, \cap\, \mN_b\, =\, \emptyset$,  
\beno
\max\limits_{(\bx,\, \bx^\prime)\, \in\, \mbX\, \times\, \mbX:\;\, x_{v,w}\, =\, x_{v,w}^\prime,\; \{v, w\}\, \neq\, \{a, b\}}\;
\left|\dsum_{i=1}^{N}\, I_{i,a_i}(\bx) - \dsum_{i=1}^{N} \, I_{i,a_i}(\bx^\prime)\right|
\;=\; 0,
\ee
and,
for each pair of nodes $\{a, b\} \subset \mN$ with $\mN_a\, \cap\, \mN_b\, \neq\, \emptyset$,
\beno
\max\limits_{(\bx,\, \bx^\prime)\, \in\, \mbX\, \times\, \mbX:\;\, x_{v,w}\, =\, x_{v,w}^\prime,\; \{v, w\}\, \neq\, \{a, b\}}\;
\left|\dsum_{i=1}^{N} \, I_{i,a_i}(\bx) - \dsum_{i=1}^{N} \, I_{i,a_i}(\bx^\prime)\right|
\;\leq\; D_N.
\ee
Using Lemma \ref{lem:DN},
the number of pairs $\{a, b\} \subset \mN$ with $\mN_a\, \cap\, \mN_b\, \neq\, \emptyset$ is bounded above by $N\, D_N^2$,
so that
\beno 
\Psi_N 
&\leq& \sqrt{N\, D_N^2} 
\= D_N\, \sqrt{N}. 
\ee

{\bf Bounding $\mnorm{\mD_N(\truth)}_2$ from above.}
We have shown that
\beno
\sqrt{\log \max\{N,\, p\}}\, \Psi_N\, \mnorm{\mD_N(\truth)}_2 
\lte 2\, \sqrt{\log N}\, D_N\, \sqrt{N}\, \mnorm{\mD_N(\truth)}_2\s\s
\\
\= 2\, D_N\, \sqrt{N \log N}\, \mnorm{\mD_N(\truth)}_2.
\ee
To bound $\mnorm{\mD_N(\truth)}_2$ from above,
we distinguish two scenarios:
\ben
\item {\bf Bounding $\mnorm{\mD_N(\truth)}_2$ from above when the subpopulations do not intersect ($\omega_1 = \omega_{2} = 0$):}\,
If condition \sone is satisfied,
Lemma \ref{prop:D_bound} implies that

\vspace{-.2cm}

\beno
\mnorm{\mD_N(\truth)}_2
&\leq& 1 + 4\, D_N^2
&\leq& 5\, D_N^2,
\ee
because $D_N \geq 1$ under Models 2 and 3.
As a result,
\beno
\sqrt{\log \max\{N,\, p\}}\, \mnorm{\mD_N(\truth)}_2\, \Psi_N
&\leq& 10\; D_N^3\, \sqrt{N \log N},
\ee
Since condition \sone ensures that $D_N = O(\log N)$,
there exist constants $C > 0$ and $N_0 \geq 3$,
independent of $N$ and $p$,
such that,
for all $N \geq N_0$,
\beno
\sqrt{\log \max\{N,\, p\}}\, \mnorm{\mD_N(\truth)}_2\, \Psi_N
\lte 10\; D_N^3\, \sqrt{N \log N} \s\s
\\
&\leq& C\, (\log N)^{3/2}\, \sqrt{N}\s
\\
&<& \dfrac{N}{8\, \exp(L / 7)\, N^{\vartheta / 7}}\s
\\
&\leq& \dfrac12\; \mbE\, \dsum_{i=1}^{N}\, \norm{\bI_{i}(\bX)}_\infty,
\ee
because $\sqrt{N} < N\, /\, N^{\vartheta / 7}$ owing to $\vartheta < 1/2 - \alpha < 1/2$. 

\s

\item {\bf Bounding $\mnorm{\mD_N(\truth)}_2$ from above when the subpopulations intersect ($\omega_1 > 0$, $\omega_2 \geq 0$):}
If condition \stwo is satisfied,
Lemma \ref{prop:D_bound} implies there exist constants 
$C_1,\, C_2,\, C_3 > 0$ and $N_1 \geq 3$, 
independent of $N$ and $p$,
such that,
for all $N > N_1$, 
\beno
\mnorm{\mD_N(\truth)}_{2}
&\leq& 1 + 4 \, D_N^2 + \omega_1 \, C_1  \, \exp(C_2 \, D_N^3)
&\leq& 
\hide{
3\; \omega_1 \, C_1 \, \exp(C_2 \, D_N^3) \s\\
\= 
}
C_3\, \exp(C_2 \, D_N^3).
\ee
Thus,
there exists constants $C_4 > 0$ and $N_2 \geq 3$,
independent of $N$ and $p$, 
such that,
for all $N > N_2$,
\beno
\sqrt{\log \max\{N,\, p\}}\, \mnorm{\mD_N(\truth)}_2\, \Psi_N
\lte 2\, \sqrt{\log N}\, D_N\, \sqrt{N}\; C_3\, \exp(C_2 \, D_N^3)\s\s
\\
\hide{
\= 2\; C_3\, \sqrt{N \log N}\, D_N\, \exp(C_2 \, D_N^3)\s\s
\\
\lte 2\; C_3\, \sqrt{N \log N}\, \exp(C_2 \, D_N^3 + \log D_N).
\lte 2\; C_3\, \sqrt{N \log N}\, \exp(C_2 \, D_N^3 + \log D_N).
}
\lte \exp(C_4\,  D_N^3)\, \sqrt{N \log N}.
\ee
We want to prove that
\beno
\sqrt{\log \max\{N,\, p\}}\, \mnorm{\mD_N(\truth)}_2\, \Psi_N
\lte \exp(C_4\,  D_N^3)\, \sqrt{N \log N}\s
\\
&\leq& \dfrac{N}{8\, \exp(L / 7)\, N^{\vartheta / 7}}.
\ee
Note that there exists an integer $N_3 \geq 3$, 
independent of $N$ and $p$, 
such that,
for all $N > N_3$, 
the condition 
\beno
\exp(C_4\, D_N^3)\, \sqrt{N \log N}
&\leq& \dfrac{N}{8\, \exp(L / 7)\, N^{\vartheta / 7}}
\ee
is satisfied,
by invoking the assumption that $\vartheta = 0$ and concluding that
\beno
\exp(C_4\, D_N^3)
&\leq& \dfrac{1}{8\, \exp(L / 7)}\; \sqrt{\dfrac{N}{\log N}},
\ee
provided 
\beno
D_N &=& o((\log(N\, /\, \log N))^{1/3}),
\ee
which is condition \stwo.
As a result,
for all $N > \max\{N_1,\, N_2,\, N_3\}\, \geq\, 3$,
\beno
\sqrt{\log \max\{N,\, p\}}\, \mnorm{\mD_N(\truth)}_2\, \Psi_N
&\leq& \dfrac{N}{8\, \exp(L / 7)\, N^{\vartheta / 7}}\s
\\
\lte \dfrac12\; \mbE\, \dsum_{i=1}^{N}\, \norm{\bI_{i}(\bX)}_\infty.
\ee
\een

{\bf Conclusion.}
We have demonstrated that there exists an integer $N_0 \geq 3$,
independent of $N$ and $p$, 
such that,
for all $N > N_0$,
\beno
\dfrac12\; \mbE\, \dsum_{i=1}^{N} \, \norm{\bI_{i}(\bX)}_{\infty}
&\geq& \dfrac{N}{2\, (1 + \exp((3 + D_N) \, (\norm{\truth}_{\infty} + \epsilon^\star)))^2}\s
\\
&\geq& \dfrac{N}{8\, \exp(L / 7)\, N^{\vartheta / 7}}\s\s
\\
\gte \sqrt{\log \max\{N,\, p\}}\, \mnorm{\mD_N(\truth)}_2\, \Psi_N
\ee
and
\beno
\mbP\left(\left|\dsum_{i=1}^{N} \norm{\bI_{i}(\bX)}_{\infty} - \mbE \dsum_{i=1}^{N} \norm{\bI_{i}(\bX)}_{\infty}\right| < \dfrac12\; \mbE \, \dsum_{i=1}^{N}\, \norm{\bI_{i}(\bX)}_\infty\right)\s\s
\\
\geq\; \mbP\left(\left|\dsum_{i=1}^{N} \norm{\bI_{i}(\bX)}_{\infty} - \mbE \dsum_{i=1}^{N} \norm{\bI_{i}(\bX)}_{\infty}\right| < \sqrt{\log \max\{N,\, p\}}\, \mnorm{\mD_N(\truth)}_2 \, \Psi_N\right)\s\s
\\
\geq\; 1 - \dfrac{2}{\max\{N,\, p\}^2}.
\ee
Combining these two results proves the desired result,
as explained at the beginning of the proof of Lemma \ref{lem:Iscale}.
\qed
\hide{
\beno
\mbP(\bX \in \mbH)
&=& \mbP\left( \dsum_{i=1}^{N}\, \norm{\bI_{i}(\bX)}_{\infty}
\;\geq\; \dfrac{N}{2\, (1 + \exp((3 + D_N) \, (\norm{\truth}_{\infty} + \epsilon^\star)))^2} \right) \s\\
&\geq& 1 - \dfrac{2}{\max\{N,\, p\}^2}.
\ee
}

\s

\subsection{Bounding \texorpdfstring{$\mnorm{\mD_N(\truth)}_2$}{D}}
\label{sup-sec:mrf_conc}

To bound the spectral norm $\mnorm{\mD_N(\truth)}_2$ of the coupling matrix $\mD_N(\truth)$,
we first review undirected graphical models encoding the conditional independence properties of generalized $\beta$-models with dependent edges in Appendices \ref{sup-sec:gm_mrf} and \ref{sup-sec:computing}.
We then bound\break 
$\mnorm{\mD_N(\truth)}_2$ by using these conditional independence properties in Appendix \ref{D}.
Auxiliary results can be found in Appendix \ref{sup-sec:D_aux_proofs}.

\subsubsection{Undirected graphical models of random graphs}
\label{sup-sec:gm_mrf}

Let $\mG(\mV, \mE)$ be an undirected graph with set of vertices $\mV$ and set of edges 
\beno 
\mE &\subseteq& \big\{\{v,w\} \,:\, v \in \mV, \, w \in \mV \setminus \{v\} \big\}.
\ee 
An undirected graphical model of a random graph \citep{LaRiSa17} is a family of probability measures 
$\{\mbP_{\nat}, \, \nat \in \Nat\}$ dominated by a $\sigma$-finite measure $\nu$,
with factorization and conditional independence properties \citep{Da79} of the form
\be
\label{ug.rg}
f_{\nat}\left(\bx\right)
&\coloneqq& \dfrac{\dd \mbP_{\nat}}{\dd \nu}(\bx)
&\propto&
\dprod_{\mathcal{C}\, \in\, \mathfrak{C}} g_{\mathcal{C}}(\bx_{\mathcal{C}};\, \nat),
&& \bx \in \mbX,
\ee
where $\mathfrak{C}$ is the set of all maximal complete subsets of the conditional independence graph $\mG(\mV, \mE)$ 
with set of vertices $\mV = \{X_1, \dots, X_{M}\}$ and set of edges 
$\mE \subset \{\{v,w\} \,:\, v \in \mV, \, w \in \mV \setminus \{v\}\}$.
The functions $g_{\mathcal{C}}: \mbX \times \Nat \mapsto \mbR^+ \cup \{0\}$ are non-negative functions defined on the maximal complete subsets $\mathcal{C} \in \mathfrak{C}$ of the conditional independence graph $\mG$.
A complete subset of the conditional independence graph $\mG$ is a subset of vertices such that each pair of vertices is connected by an edge,
and a complete subset is maximal complete if no vertices can be added without losing the property of completeness.

The probability density functions introduced in Section \ref{sec:models} are of the form
\be
\label{rg}
f_{\nat}\left(\bx\right)
&\propto& \dprod_{i\, <\, j}^N \varphi_{i,j}(x_{i,j}, \, \bx_{\mS_{i,j}};\, \nat),
&& \bx \in \mbX,
\ee
where $\mS_{i,j} \subset \{\{v,w\} : v \in \mN, \, w \in \mN \setminus \{w\}\} \setminus \{i,j\}$ 
for all $\{i,j\} \subset \mN$.
Probability density functions of the form \eqref{rg} can be represented as probability density functions of the form \eqref{ug.rg} by grouping the functions $\varphi_{i,j}$ in accordance with the maximal complete subsets of conditional independence graph $\mG$.
The conditional independence graph $\mG$ depends on the model:
e.g.,
the conditional independence graph of Model 1 has no edges,
because all edge variables are independent.
By contrast,
the conditional independence graph of Models 2 and 3 shown in Figure \ref{fig:cond_ind_graph} has edges,
which indicate the absence of conditional independence among edge variables due to brokerage in overlapping subpopulations.

\hide{
\begin{figure}[t]
\centering
\begin{tikzpicture}
\node[align = flush center, minimum size = .85cm] (center) [] {};
\node[node, fill = myred] (12) [right = 0cm of center] {$X_{1,2}$};
\node[minimum size = .85cm] (pivot1) [right = .5cm of 12] {};
\node[node, fill = myred] (13) [above = .01cm of pivot1] {$X_{1,3}$};
\node[node, fill = myred] (23) [below = .01cm of pivot1] {$X_{2,3}$};
\node[node, fill = myred] (34) [right = .5cm of pivot1] {$X_{3,4}$};
\node[node, fill = myorange] (14) [above = 1cm of 34] {$X_{1,4}$};
\node[node, fill = myorange] (24) [below = 1cm of 34] {$X_{2,4}$};
\node[node, fill = myred] (45) [right = .6cm of 34] {$X_{4,5}$};
\node[node, fill = myorange] (35) [above = 1cm of 45] {$X_{3,5}$};
\node[node, fill = myred] (56) [right = .6cm of 45] {$X_{5,6}$};
\node[node, fill = myorange] (46) [above = 1cm of 56] {$X_{4,6}$};
\node[node, fill = myorange] (47) [below = 1cm of 45] {$X_{4,7}$};
\node[node, fill = myred] (57) [below = 1cm of 56] {$X_{5,7}$};
\node[minimum size = .85cm] (pivot2) [right = .6cm of 57] {};
\node[node, fill = myred] (67) [right = .6 cm of 56] {$X_{6,7}$};

\node[node, fill = mygrey] (25) [below right = 3.5cm of 12] {$X_{2,5}$};
\node[node, fill = mygrey] (26) [right = .3cm of 25] {$X_{2,6}$};
\node[node, fill = mygrey] (17) [left = .3cm of 25] {$X_{1,7}$};
\node[node, fill = mygrey] (16) [left = .3cm of 17] {$X_{1,6}$};
\node[node, fill = mygrey] (15) [left = .3cm of 16] {$X_{1,5}$};
\node[node, fill = mygrey] (27) [right = .3cm of 26] {$X_{2,7}$};
\node[node, fill = mygrey] (36) [right = .3cm of 27] {$X_{3,6}$};
\node[node, fill = mygrey] (37) [right = .3cm of 36] {$X_{3,7}$};

\path[draw, thick]
(12) edge (13)
(12) edge (23)
(13) edge (23)
(13) edge (34)
(23) edge (34)
(13) edge (14)
(34) edge (14)
(24) edge (34)
(24) edge (23)
(34) edge (35)
(45) edge (34)
(45) edge (35)
(45) edge (56)
(56) edge (46)
(45) edge (46)
(45) edge (47)
(47) edge (57)
(45) edge (57)
(57) edge (56)
(57) edge (67)
(67) edge (56);

\end{tikzpicture}
\caption{\label{fig:cond_ind_graph}
The conditional independence graph of Models 2 and 3 with population of nodes $\mN = \{1, \dots, 7\}$,
consisting of overlapping subpopulations $\mA_1 = \{1, 2, 3\}$,
$\mA_2 = \{3, 4\}$,
$\mA_3 = \{4, 5\}$,
and $\mA_4 = \{5, 6, 7\}$.
If nodes $i$ and $j$ belong to the same subpopulation,
edge variable $X_{i,j}$ is colored red.
If nodes $i$ and $j$ do not belong to the same subpopulation,
edge variable $X_{i,j}$ is colored orange if the subpopulations of $i$ and $j$ overlap and is colored gray otherwise.
}
\end{figure}
}

\begin{figure}
\centering 
\includegraphics[width = .6\linewidth, keepaspectratio]{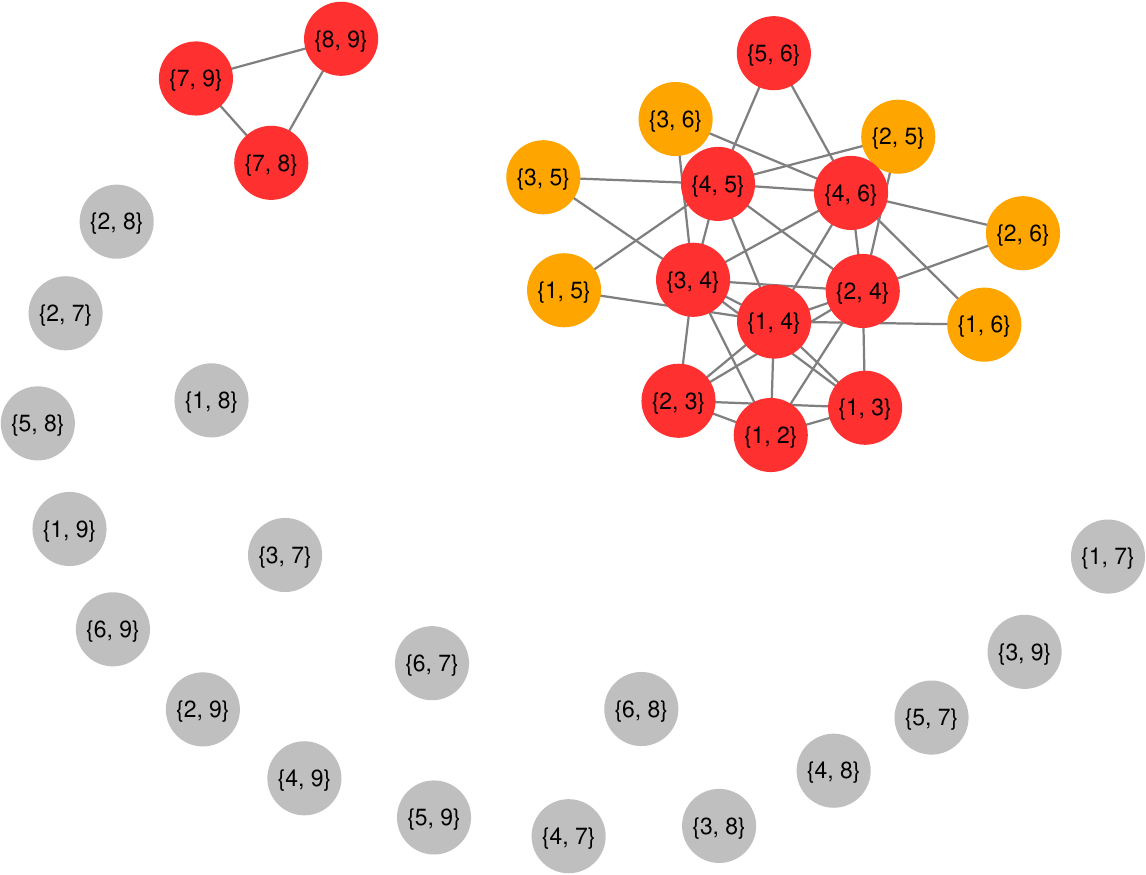}
\caption{\label{fig:cond_ind_graph}
The conditional independence graph of Models 2 and 3 with population of nodes $\mN \coloneqq \{1, \dots, 9\}$,
consisting of overlapping subpopulations $\mA_1 \coloneqq \{1,2,3,4\}$,\, 
$\mA_2 \coloneqq \{4, 5, 6\}$,\, 
and $\mA_3 \coloneqq \{7,8,9\}$. 
Edge variables $X_{i,j}$ are represented by circles with labels $\{i, j\}$.
If nodes $i$ and $j$ share a subpopulation,
$X_{i,j}$ is colored red. 
If nodes $i$ and $j$ do not share a subpopulation but belong to overlapping subpopulations,
$X_{i,j}$ is colored orange. 
Otherwise,
$X_{i,j}$ is colored gray.
}
\end{figure}

To distinguish the random graph (representing data structure) from the conditional independence graph $\mG$ (representing conditional independence structure, i.e., model structure),
we call elements of $\mV$ vertices rather than nodes,
and elements of $\mE$ edges rather than edge variables.

\subsubsection{Conditional independence properties}
\label{sup-sec:computing}

We prove selected conditional independence properties that help establish consistency results and convergence rates for generalized $\beta$-models with dependent edges.

By Equations \eqref{eq:factorization} and \eqref{parameterization},
for each $\{i,j\} \subset \mN$,  
\beno
\varphi_{i,j}(x_{i,j},\, \bx_{\mS_{i,j}};\, \nat)
\,\coloneqq\, a_{i,j}(x_{i,j})\, \exp\left((\theta_i + \theta_j)\, x_{i,j} + \theta_{N+1}\, b_{i,j}(x_{i,j},\, \bx_{\mS_{i,j}})\right),
\ee
where 
\beno
b_{i,j}(x_{i,j}, \bx_{\mS_{i,j}})
&\coloneqq&
\begin{cases}
0 & \mbox{if } \mN_i \cap \mN_j = \emptyset\s
\\
x_{i,j}\, \one\left(\dsum_{h\, \in\, \mN_i\, \cap\, \mN_j} x_{i,h}\, x_{j,h} \geq 1\right)
& \mbox{if } \mN_i \cap \mN_j \neq \emptyset.
\end{cases}
\ee

\s

\begin{definition}
\label{def:intersection}
{\bf Neighborhood intersection property.}
{\em
Consider a random graph model with a probability density function parameterized by \eqref{eq:factorization} and \eqref{parameterization}.
If $\mS_{i,j} = \left\{\{a,b\} \subset \mN \,:\, (a,b) \in \{i, j\} \times \{\mN_i\, \cap\, \mN_j\} \right\}$ for all pairs of nodes $\{i, j\} \subset \mN$,
then the random graph is said to satisfy the neighborhood intersection property.}
\end{definition}

\s 

By construction,
generalized $\beta$-models with dependent edges satisfy the neighborhood intersection property, 
which implies conditional independence properties,
including---but not limited to---the conditional independence properties established in Proposition \ref{prop:cond_ind} below.  
Define 
\beno
\dyads &\coloneqq& \big\{\{a,b\} \,:\, a \in \mN,\; b \in \mN \setminus \{a\} \big\}.
\ee
We will utilize $\dyads$ as the index set of all possible edge variables, 
which will be useful in constructing statements which exclude certain edge variables. 

\s 

\begin{proposition}
\label{prop:cond_ind}
A random graph with overlapping subpopulations $\mA_k$ of sizes $|\mA_k| \geq 3$ ($k = 1, \dots, K$) satisfying the neighborhood intersection property possesses the following conditional independence properties:
\bi
\item[1.] For all pairs of nodes $\{i,j\} \subset \mN$ such that $\mN_i\, \cap\, \mN_j = \emptyset$:
\beno
X_{i,j} &\orth& \bX_{\dyads \setminus\, \{i,j\}}. 
\ee
\item[2.] For all pairs of nodes $\{i,j\} \subset \mN$ such that $\mN_i\, \cap\, \mN_j \neq \emptyset$ and there exists $k \in \{1, \ldots, K\}$ such that $\{i,j\} \subset \mA_k$:  
\beno
X_{i,j} &\orth& \bX_{\dyads \setminus\, (\{i,j\} \,\cup\, \mathfrak{N}_{i,j})} \mid \bX_{\mathfrak{N}_{i,j}},
\ee
where $\mathfrak{N}_{i,j} = \mathcal{J}_{i}^{(j)} \cup \mathcal{J}_{j}^{(i)}$, 
using the definition 
\beno
\mathcal{J}_{t}^{(v)}
&\coloneqq&
\left[ \bigcup\limits_{b \in \mN_t \setminus \{v\}} \, \{v,b\}  \right]
\;\cup\;
\left[ \bigcup\limits_{b \in \mN_t} \, \big\{\{a, b\} : (a,b) \in \{v,b\} \times \mN_v \cap \mN_b \big\}  \right],
\ee
and with the property that 
\beno
\mathfrak{N}_{i,j}
&\subseteq& \big\{ \{a,b\} \,:\, a \in \mN_i \cup \mN_j, \, b \in \mN_i \cup \mN_j \setminus \{a\} \big\}.
\ee

\item[3.] For all pairs of nodes $\{i,j\} \subset \mN$ such that $\mN_i\, \cap\, \mN_j \neq \emptyset$ 
and there exists no $k \in \{1, \ldots, K\}$ such that $\{i,j\} \subset \mA_k$:
\beno
X_{i,j} &\orth& \bX_{\dyads \setminus\, (\{i,j\} \,\cup\, \mS_{i,j})}
\mid \bX_{\mS_{i,j}}.
\ee
where $\mS_{i,j}
= \left\{\{a,b\} \subset \mN \,:\, (a,b) \in \{i, j\} \times \{\mN_i\, \cap\, \mN_j\} \right\}$. 
\ei
\end{proposition}

\pproof \ref{prop:cond_ind}.
In the following,
we use the characterizations of conditional independence due to \citet{Da79},
which relate factorization properties of probability density functions to conditional independence properties.
Using these characterizations of conditional independence,
we establish the conditional independence properties of Proposition \ref{prop:cond_ind} by showing that,
for each pair of nodes $\{i,j\} \subset \mN$, 
there exists a subset of edge indices  
$\mathfrak{N}_{i,j} \subseteq \dyads \setminus  \{i,j\}$ 
and non-negative functions $g$ and $h$ 
such that the probability density function can be written as
\beno
f_{\nat}(\bx)
&\propto& g(x_{i,j},\, \bx_{\mathfrak{N}_{i,j}}) \; h(\bx_{\dyads \setminus \{i,j\}}), 
\ee
where $\dyads \coloneqq \{\{a,b\} \,:\, a \in \mN, \, b \in \mN \setminus \{a\}\}$,  
implying that 
\beno
X_{i,j} &\orth& \bX_{\dyads \setminus (\{i,j\} \,\cup\; \mathfrak{N}_{i,j})} \mid \bX_{\mathfrak{N}_{i,j}}.  
\ee
Proposition \ref{prop:cond_ind} assumes that the neighborhood intersection property is satisfied,
allowing us to write 
\beno
f_{\nat}(\bx) 
&\propto& \dprod_{i\, <\, j}^N \, \varphi_{i,j}(x_{i,j}, \, \bx_{\mS_{i,j}}), 
\ee
where 
\beno
\mS_{i,j} 
\= \left\{ \{v,w\} : (v,w) \in \{i,j\} \times \mN_i \cap \mN_j \right\}, 
\ee
recalling the definition of the node neighborhood sets $\mN_i$ ($i \in \mN$):
\beno
\mN_i 
\= \{h \in \mN \setminus \{i\} \,:\, \mbox{ exists } k \in \{1, \ldots, K\} \mbox{ such that } \{i,h\} \subset \mA_k \}. 
\ee

\s

{\bf Condition 1:} 
Consider any pair of nodes $\{i,j\} \subset \mN$ with $\mN_i\, \cap\, \mN_j = \emptyset$,
that is,
nodes $i$ and $j$ neither belong to a common subpopulations nor belong to distinct subpopulations that overlap.
Since $\{i, j\} \times \mN_i \cap \mN_j = \emptyset$,
\beno
\varphi_{i,j}(x_{i,j}, \, \bx_{\mS_{i,j}}) 
&\equiv& \varphi_{i,j}(x_{i,j}).
\ee
It remains to check whether any $\varphi_{a,b}$ with $\{a, b\} \neq \{i,j\}$ 
can be a function of $x_{i,j}$,
which can happen in one of two different ways: 
\bi
\item $i \in \{a, b\}$ and $j \in \mN_a \cap \mN_b$, 
in which case, 
by the definition of the node neighborhood sets $\mN_v$ ($v \in \mN$), 
it must be that $j \in \mN_i$; or \s
\item $j \in \{a,b\}$ and $i \in \mN_a \cap \mN_b$,
in which case,
similarly, 
$i \in \mN_j$. 
\ei
We prove that there exists no $\varphi_{a,b}$ with $\{a, b\} \neq \{i,j\}$
which is a function of $x_{i,j}$ by contradiction. 
Note $i \in \mN_j$ implies the existence of a $k \in \{1, \ldots, K\}$ 
such that $\{i,j\} \subset \mA_k$.
By assumption, 
$|\mA_k| \geq 3$ for all $k \in \{1, \ldots, K\}$, 
implying there exists at least one other node $h \in \mA_k \setminus \{i,j\}$. 
Thus,
if $\{i,j,h\} \subseteq \mA_k$,
then both $h \in \mN_i$ and $h \in \mN_j$ must hold,
violating the assumption that $\mN_i \cap \mN_j = \emptyset$. 
The case when $j \in \{a,b\}$ and $i \in \mN_a \cap \mN_b$ 
is proved similarly. 
Therefore,
there cannot exist a pair of nodes $\{a, b\} \neq \{i,j\}$ such that 
$\varphi_{a,b}$ 
is a function of $x_{i,j}$. 
As a consequence,
taking
\beno
g(x_{i,j}) 
\= \varphi_{i,j}(x_{i,j})
\ee
and
\beno
h(\bx_{\dyads \setminus \{i,j\}})
\= \dprod_{a<b:\; \{a,b\}\, \neq\, \{i,j\}}^N \, \varphi_{a,b}(x_{a,b},\, \bx_{\mS_{a,b}}) 
\ee
shows that $f_{\nat}(\bx)$ can be written as 
\beno
f_{\nat}(\bx)
&\propto& g(x_{i,j})\; h(\bx_{\dyads \setminus \{i,j\}}),
\ee
which implies 
$X_{i,j} \orth \bX_{\dyads \setminus \{i,j\}}$, 
i.e.,
$X_{i,j}$ is independent of all others edges.

\s\s

{\bf Condition 2:} 
Consider any pair of nodes $\{i,j\} \subset \mN$ with $\mN_i \cap \mN_j \neq \emptyset$ 
and such that there exists a $k \in \{1, \ldots, K\}$ such that $\{i,j\} \subset \mA_k$. 
By definition, 
$\varphi_{i,j}(x_{i,j}, \bx_{\mS_{i,j}})$ is a function of $x_{i,j}$.
Recall the key condition for the neighborhood intersection assumption,
which was that $\mS_{a,b}$ satisfies  
\beno
\mS_{a,b}
\= \left\{ \{v,w\} : (v,w) \in \{a,b\} \times \mN_a \,\cap\, \mN_b \right\},
&& \{a,b\} \subset \mN. 
\ee
For any $\varphi_{a,b}(x_{a,b}, \bx_{\mS_{a,b}})$
with $\{a, b\} \neq \{i,j\}$ to be a function of $x_{i,j}$,
one of the following must hold: 
\ben
\item $i \in \{a, b\}$ and $j \in \mN_a \cap \mN_b$,
in which case,
by the definition of the node neighborhood sets $\mN_v$ ($v \in \mN$),
$j \in \mN_i$; or \s
\item $j \in \{a,b\}$ and $i \in \mN_a \cap \mN_b$,
in which case,
similarly,
$i \in \mN_j$.
\een
Consider the first case:  
$i \in \{a, b\}$ and $j \in \mN_a \cap \mN_b$,
and without loss, 
take $a = i$. 
The condition for this case implies that $j \in \mN_i \cap \mN_b$. 
By assumption,
$\{i,j\} \subset \mA_k$ for some $k \in \{1, \ldots, K\}$, 
which implies $j \in \mN_i$.  
If $j \in \mN_b$,  
then $b \in \mN_j$,
implying $\varphi_{i,b}$ is a function of $x_{i,j}$ 
for all $\{i,b\} \subset \mN$ 
with $b \in \mN_j$. 
Applying the same argument to the second case  
where $j \in \{a,b\}$ and $i \in \mN_a \cap \mN_b$ 
reveals that $\varphi_{j,b}$ is a function of $x_{i,j}$ 
for all $\{j,b\}$ with $b \in \mN_i$.

Summarily, 
$\varphi_{a,b}$ is a function of $x_{i,j}$ if it is in the following list: 
\bi
\item $\varphi_{i,j}(x_{i,j}, \bx_{\mS_{i,j}})$, 
where $\{v,w\} \in \mS_{i,j}$ if $(v,w) \in \{i,j\} \times \mN_i \,\cap\,\mN_j$. \s

\item $\varphi_{i,b}(x_{i,b}, \bx_{\mS_{i,b}})$ ($b \in \mN_j$), 
where $\{v,w\} \in \mS_{i,b}$ if $(v,w) \in \{i,b\} \times \mN_i \,\cap\,\mN_b$. 

\vspace{.2cm}

\item $\varphi_{j,b}(x_{j,b}, \bx_{\mS_{j,b}})$ ($b \in \mN_i$),
where $\{v,w\} \in \mS_{j,b}$ if $(v,w) \in \{j,b\} \times \mN_j \,\cap\,\mN_b$. 

\vspace{.1cm}

\ei
This collection of functions is a function of all edge variables $X_{a,b}$ with indices $\{a,b\}$ 
in $\mathcal{J}_{i}^{(j)} \cup \mathcal{J}_{j}^{(i)} \cup \{i,j\}$,
where 
\beno
\mathcal{J}_{t}^{(v)} 
&\coloneqq& 
\left[ \bigcup\limits_{b \in \mN_t \setminus \{v\}} \, \{v,b\}  \right]
\;\cup\;  
\left[ \bigcup\limits_{b \in \mN_t} \, \big\{\{a, b\} : (a,b) \in \{v,b\} \times \mN_v \cap \mN_b \big\}  \right].  
\ee
Thus, 
there exist non-negative functions $g$ and $h$ such that 
the probability density function can be written as follows:
\beno
f_{\nat}(\bx)
&\propto& g(x_{i,j}, \bx_{\mathfrak{N}_{i,j}}) 
\; h(\bx_{\dyads \setminus \{i,j\}}), 
\ee
where $\mathfrak{N}_{i,j} = \mathcal{J}_{i}^{(j)} \cup \mathcal{J}_{j}^{(i)}$, 
implying that 
\beno
X_{i,j} &\orth& \bX_{\dyads \setminus (\{i,j\} \,\cup\, \mathfrak{N}_{i,j})} 
\mid \bX_{\mathfrak{N}_{i,j}}.  
\ee
As $\{i,j\} \subset \mN_i \cup \mN_j$, 
\beno
\mathfrak{N}_{i,j} 
&\subseteq& \big\{ \{a,b\} \,:\, a \in \mN_i \cup \mN_j, \, b \in \mN_i \cup \mN_j \setminus \{a\} \big\}. 
\ee

\s

{\bf Condition 3:} 
Consider any pair of nodes $\{i,j\} \subset \mN$ 
with $\mN_i \cap \mN_j \neq \emptyset$ 
and such that there exists no $k \in \{1, \ldots, K\}$ such that $\{i,j\} \subset \mA_k$.
It is clear that $\varphi_{i,j}$ is a function of $x_{i,j}$. 
For any $\varphi_{a,b}(x_{a,b}, \bx_{\mS_{a,b}})$ with $\{a, b\} \neq \{i,j\}$ to be a function of $x_{i,j}$,
one of the following must hold:
\bi
\item $i \in \{a, b\}$ and $j \in \mN_a \cap\mN_b$,
in which case,
by the definition of the node neighborhood sets $\mN_v$ ($v \in \mN$),
it must be that $j \in \mN_i$; or 

\vspace{.1cm}

\item $j \in \{a,b\}$ and $i \in \mN_a \cap \mN_b$,
in which case,
similarly,
$i \in \mN_j$.
\ei
In both conditions, 
$i \in \mN_j$ and $j \in \mN_i$,
which implies that  
$\{i,j\} \subset \mA_k$ for some $k \in \{1, \ldots, K\}$,
violating the assumption that no such $k$ exists. 
Thus, 
$\varphi_{a,b}(x_{a,b}, \bx_{\mS_{a,b}})$ 
is a function of $x_{i,j}$ if and only if $\{a,b\} = \{i,j\}$.  
As a result,
there exist non-negative functions $g$ and $h$ 
such that 

\vspace{-.2cm}

\beno
f_{\nat}(\bx)
&\propto& g(x_{i,j}, \, \bx_{\mS_{i,j}}) \; h(\bx_{\dyads \setminus \{i,j\}}) 
\ee
which implies  
$X_{i,j} \,\orth\, \bX_{\dyads \setminus (\{i,j\} \,\cup\, \mS_{i,j})}
\mid \bX_{\mS_{i,j}}$.\qed

\s

\begin{lemma}
\label{lem:DN}
Consider Models 2 and 3.
Then $\max_{t \in \mN} \, |\mN_t| \leq D_N$ and 
\beno
\max\limits_{t \in \mN} \,
\left|\left\{a \in \mN \setminus \{t\} \,:\, \mN_a \cap \mN_t \neq \emptyset \right\} \right|
&\leq& D_N^2,
\ee
where $D_N \coloneqq \max_{\{i,j\} \subset \mN} \, |\mathfrak{N}_{i,j}|$. 
\end{lemma}

\s

\llproof \ref{lem:DN}.
By Proposition \ref{prop:cond_ind},
for any $\{i,j\} \subset \mN$ satisfying $\mN_i\, \cap\, \mN_j \neq \emptyset$ 
and for which there exists $k \in \{1, \ldots, K\}$ such that 
$\{i,j\} \subset \mA_k$,
we have $\mathfrak{N}_{i,j} = \mathcal{J}_{i}^{(j)} \cup \mathcal{J}_{j}^{(i)}$,
where,
for all $\{t,v\} \subset \mN$, 
\beno
\mathcal{J}_{t}^{(v)}
&\coloneqq&
\left[ \bigcup\limits_{b \in \mN_t \setminus \{v\}} \, \{v,b\}  \right]
\;\cup\;
\left[ \bigcup\limits_{b \in \mN_t} \, \big\{\{a, b\} : (a,b) \in \{v,b\} \times \mN_v \cap \mN_b \big\}  \right].
\ee
Then, 
for each $t \in \mN$,
there exists $v \in \mN \setminus \{t\}$ 
and 
$k \in \{1, \ldots, K\}$ 
such that $\{t,v\} \subset \mA_k$,
due to the assumption that $|\mA_k| \geq 3$ for all $k \in \{1, \ldots, K\}$ under Models 2 and 3, 
implying that 
\beno
|\mathfrak{N}_{v,t}|
\,=\, 
|\mathcal{J}_{t}^{(v)} \cup \mathcal{J}_{v}^{(t)}|
\,\geq\, \left| \left[ \bigcup\limits_{b \in \mN_t \setminus \{v\}} \,
\{v,b\} \right] 
\cup 
\left[ \bigcup\limits_{b \in \mN_v \setminus \{t\}} \,
\{t,b\} \right]
\right|
\,\geq\, |\mN_t|. 
\ee
Thus $D_N \coloneqq \max_{\{i,j\} \subset \mN} \, |\mathfrak{N}_{i,j}| \geq |\mN_t|$ for all $t \in \mN$. 
Next,
for all $t \in \mN$, 
\beno
\left|\left\{a \in \mN \setminus \{t\} \,:\, \mN_a \cap \mN_t \neq \emptyset \right\} \right|
\,\leq\,
\left| \bigcup\limits_{r \in \mN_t} \, \mN_r \right|
\,\leq\, |\mN_t| \, \left(\max\limits_{r \in \mN} |\mN_r| \right)
\,\leq\, D_N^2, 
\ee
using the above-proven fact that $\max_{t \in \mN} \, |\mN_t| \leq D_N$.  
\qed

\s

\subsubsection{Bounding the spectral norm of the coupling matrix}
\label{D}

We bound the spectral norm $\mnorm{\mD_N(\truth)}_2$ of the coupling matrix $\mD_N(\truth)$.
Throughout,
we adopt the notation used in Section \ref{sec:stat_inf} of the manuscript and denote the number of edge variables by $M = \binom{N}{2}$ and edge variables by $X_1, \dots, X_M$. 

\vspace{-.1cm}

\begin{lemma}
\label{prop:D_bound}
Consider Models 2 and 3.
Assume that Assumption A is satisfied and that $\truth \in \Nat = \mR^p$ satisfies condition \eqref{condition.theta0} in Section \ref{sec:corollaries}.
%
%
\ben
\item If the subpopulations do not intersect ($\omega_1 = \omega_2 = 0$) and 
$\truth \in \mR^{N+1}$ satisfies condition \eqref{condition.theta0} with $\vartheta \in [0,\, 1/2 - \alpha)$,
then 

\vspace{-.2cm}

\beno
\mnorm{\mD_N(\truth)}_2
&\leq& 1 + 4 \, D_N^2.
\ee

\item If the subpopulations do intersect ($\omega_1 > 0$) and  
$\truth \in \mR^{N+1}$ satisfies condition \eqref{condition.theta0} with $\vartheta = 0$,
then there exists finite constants $C_1 > 0$ and $C_2 > 0$, 
independent of $N$ and $p$, 
such that 

\vspace{-.2cm}

\beno
\mnorm{\mD_N(\truth)}_{2} 
&\leq& 1 + 4 \, D_N^2 + \omega_1 \, C_1  \, \exp(C_2 \, D_N^3). 
\ee
\een 
\end{lemma}


\s

\llproof \ref{prop:D_bound}.
We adapt the coupling approach of \citet[][pp.~759--760]{BeMa94} 
from the literature on Gibbs measures and Markov random fields to coupling conditional distributions of subgraphs of random graphs.
Let $i \in \mV$ be any vertex of the conditional independence graph $\mG$,
corresponding to edge variable $X_i$,  
and consider any $\bx_{1:i-1} \in \{0, 1\}^{i-1}$.
Define
\beno
\mbP_{i,\bx_{1:i-1},0}(\bX_{i+1:M} = \bm{a})
&\coloneqq& \mbP(\bX_{i+1:M} = \bm{a} \mid \bX_{1:i-1} = \bx_{1:i-1}, X_i = 0)
\ee
and
\beno
\mbP_{i,\bx_{1:i-1},1}(\bX_{i+1:M} = \bm{a})
&\coloneqq& \mbP(\bX_{i+1:M} = \bm{a} \mid \bX_{1:i-1} = \bx_{1:i-1}, X_i = 1),
\ee
where $\bX_{1:i-1} = (X_1, \dots, X_{i-1})$,
$\bX_{i+1:M} = (X_{i+1}, \dots, X_{M})$,
and $\bm{a} \in \{0, 1\}^{M-i}$.

\s

\noindent We divide the proof into three parts:
\bi
\item[I.] Coupling conditional distributions of subgraphs.\s
\item[II.] Bounding the elements of the coupling matrix $\mD_N(\truth)$.\s
\item[III.] Bounding the spectral norm $\mnorm{\mD_N(\truth)}_2$ of the coupling matrix $\mD_N(\truth)$.
\ei

\s

\noindent {\bf I. Coupling conditional distributions of subgraphs.}
Given any vertex $i \in \mV$ of the conditional independence graph $\mG$ and any $\bx_{1:i-1} \in \{0, 1\}^{i-1}$,
we construct a coupling $(\bX^\star, \bX^{\star\star})$ 
of the conditional probability distributions $\mbP_{i, \bx_{1:i-1}, 0}$ and $\mbP_{i, \bx_{1:i-1}, 1}$.
Some background on coupling can be found in \citet{Li02}.

It will be convenient to assume 
that the coupling $(\bX^\star, \bX^{\star\star})$ 
takes on values in the set $\{0, 1\}^{M} \times \{0, 1\}^{M}$ rather than the set $\{0, 1\}^{M-i} \times \{0, 1\}^{M-i}$,
where  we set 
$(\bX_{1:i-1}^\star,\, X_i^\star) = (\bx_{1:i-1},\, 0)$ and $(\bX_{1:i-1}^{\star\star},\, X_i^{\star\star}) = (\bx_{1:i-1},\, 1)$
with probability $1$. 
As a consequence,
the random vectors $\bX^{\star} \in \{0, 1\}^M$ and $\bX^{\star\star} \in \{0, 1\}^M$ have the same dimension as random vector $\bX \in \{0, 1\}^M$.
We construct a coupling of the conditional probability distributions $\mbP_{i,\bx_{1:i-1},0}$ and $\mbP_{i,\bx_{1:i-1},1}$ as follows:
\ben
\item Initialize the subset of vertices $\mathfrak{V} = \{1, \ldots, i\}$.\s
\item Check whether there exists a vertex $j \in \mV\, \setminus\, \mathfrak{V}$ connected to a vertex $v \in \mathfrak{V}$ in the conditional independence graph $\mG$ such that the coupling disagrees at vertex $v \in \mathfrak{V}$,
in the sense that $X^{\star}_v \neq X_v^{\star\star}$.\s
\ben
\item If such a vertex $j$ exists,
pick the smallest such vertex,
and let $(X_j^{\star},\, X_j^{\star\star})$ be distributed according to an optimal coupling of
$\mbP( X_j = \,\cdot\, \mid \bX_{\mathfrak{V}} = \bx^{\star}_{\mathfrak{V}})$ and
$\mbP( X_j = \,\cdot\, \mid \bX_{\mathfrak{V}} = \bx^{\star\star}_{\mathfrak{V}})$. \s
\item If no such vertex $j$ exists,
select the smallest $j \in \mV\, \setminus\, \mathfrak{V}$ and let $(X_j^\star, \, X_j^{\star\star})$ be distributed according to an optimal coupling of $\mbP( X_j = \,\cdot\, \mid \bX_{\mathfrak{V}} = \bx^{\star}_{\mathfrak{V}})$ and $\mbP( X_j = \,\cdot\, \mid \bX_{\mathfrak{V}} = \bx^{\star\star}_{\mathfrak{V}})$.
In this case, 
an optimal coupling will ensure $X_j^{\star} = X_j^{\star\star}$ with probability $1$,
as conditional independence properties and the equality of edge variables in the conditioning statement
in this case 
will imply  
\beno
\mbP( X_j = a \mid \bX_{\mathfrak{V}} = \bx^{\star}_{\mathfrak{V}})
\= \mbP( X_j = a \mid \bX_{\mathfrak{V}} = \bx^{\star\star}_{\mathfrak{V}}),
& a \in \{0,1\},  
\ee
resulting in a  total variation distance of $0$. 
\een
\item[] In both steps,
an optimal coupling exists \citep[][Theorem 5.2, p.\ 19]{Li02},
but it may not be unique.
Any optimal coupling will do. \s
\item Replace $\mathfrak{V}$ by $\mathfrak{V}\, \cup \{j\}$ and repeat Step 2 until $\mV \setminus \mathfrak{V} = \emptyset$.\vspace{.1cm}
\een
Denote the resulting coupling distribution by $\mbQ_{i,\bx_{1:i-1}}^{}$.
Lemma \ref{lem:prove_coupling} verifies that the above algorithm constructs a valid coupling of the conditional distributions $\mbP_{i,\bx_{1:i-1},0}$ and $\mbP_{i,\bx_{1:i-1},1}$,
in the sense that the marginal distributions of $\bX^{\star}$ and $\bX^{\star\star}$ 
are $\mbP_{i,\bx_{1:i-1},0}$ and $\mbP_{i,\bx_{1:i-1},1}$,
respectively. 

For any two distinct vertices $i \in \mV$ and $j \in \{i + 1, \dots, M\}$ of the conditional independence graph $\mG$,
define the event $i \centernot{\longleftrightarrow} j$ to be the event that there exists 
a {\it path of disagreement} between $i$ and $j$ in $\mG$.  
A path of disagreement $i \centernot{\longleftrightarrow} j$ between vertices $X_i$ and $X_j$ is a sequence of two or more distinct vertices $(X_i, \dots, X_j)$ in the conditional independence graph $\mG$ starting at vertex $X_i$ and ending at vertex $X_j$,
such that 
\bi
\item each subsequent pair of vertices $(X_{v}, X_{w})$ in the sequence is connected by an edge in the conditional independence graph $\mG$,
which indicates the absence of conditional independence of vertices $X_{v}$ and $X_{w}$;\vspace{.1cm}
\item the coupling $(\bX^\star_{i+1:M}, \bX^{\star\star}_{i+1:M}) \in \{0, 1\}^{M-i} \times \{0, 1\}^{M-i}$ with joint probability mass function $\mbQ_{\truth,i,\bx_{1:i-1}}$ 
disagrees at each vertex $X_v$ in the sequence,
in the sense that $X_{v}^{\star} \neq X_{v}^{\star\star}$.
\ei

Theorem 1 of \citet[][p.\ 753]{BeMa94} shows that
\be
\label{eq:stochastic_ordering_probs}
\mbQ_{i,\bx_{1:i-1}}^{}(X_{j}^{\star} \neq X_{j}^{\star\star})
&=& \mbQ_{i,\bx_{1:i-1}}^{}(i \centernot{\longleftrightarrow} j) 
&\leq& \mbB_{\bm\pi}(i \centernot{\longleftrightarrow} j),
\ee
where $\mbB_{\bm\pi}$ is a Bernoulli product measure on $\{0, 1\}^{M}$ with probability vector $\bm\pi \in [0, 1]^{M}$.
The coordinates $\pi_v$ of $\bm\pi$ are given by 
\beno
\pi_v
&\coloneqq& 
\begin{cases}
0 & \mbox{ if } v \in \{1, \ldots, i-1\}\s
\\
1 & \mbox{ if } v = i\s
\\
\max\limits_{(\bx_{-v},\, \bx_{-v}^\prime)\, \in\, \{0, 1\}^{M-1} \times \{0, 1\}^{M-1}}\, \pi_{v,\, \bx_{-v},\, \bx_{-v}^\prime}
& \mbox{ if } v \in \{i+1, \ldots, M\},
\end{cases}
\ee
where
\beno
\pi_{v,\, \bx_{-v},\, \bx_{-v}^\prime}
&\coloneqq& \norm{\mbP(\,\cdot \mid \bX_{-v} = \bx_{-v}) - \mbP(\,\cdot \mid \bX_{-v}^\prime = \bx_{-v}^\prime)}_{\tv}. 
\ee
Observe that the total variation distance 
\beno
\norm{\mbP(\,\cdot \mid \bX_{-v} = \bx_{-v}) - \mbP(\,\cdot \mid \bX_{-v}^\prime = \bx_{-v}^\prime)}_{\tv}
\ee 
is equal to 
\beno 
\sup\limits_{x_{v}\, \in\, \{0,\, 1\}}\, \left| \mbP(X_v = x_v \mid \bX_{-v} = \bx_{-v}) - \mbP(X_v = x_v \mid \bX_{-v}^\prime = \bx_{-v}^\prime) \right|. 
\ee
The Bernoulli product measure $\mbB_{\bm\pi}$ assumes that independent Bernoulli experiments are carried out at vertices $v \in \{1, \dots, M\}$. 
The Bernoulli experiment at vertex $v \in \{i+1, \dots, M\}$ has two possible outcomes:
Either vertex $v$ is {\em open,}
in the sense that the event $\{X_v^{\star} \neq X_v^{\star\star}\}$ occurs and hence vertex $v$ allows a path of disagreement from $i$ to $j$ to pass through,
or vertex $v$ is {\em closed.}
A vertex $v$ is open with probability $\pi_v$,
and closed with probability $1-\pi_v$.
By construction, 
vertices $v \in \{1, \ldots, i-1\}$ are closed with probability $1$, 
and vertex $i$ is open with probability $1$. 

The coupling argument of \citet{BeMa94} is useful,
in that it translates the hard problem of bounding probabilities of events involving dependent random variables into the more convenient problem of bounding probabilities of events involving independent random variables.
Indeed,
we can bound the above-diagonal elements $\mD_{i,j}(\truth)$ of $\mD_N(\truth)$ by
\be
\label{boundingdij}
\mD_{i,j}(\truth)
\,\coloneqq\, \sup\limits_{\bx_{1:i-1} \in \{0,1\}^{i-1}} \,
\mbQ_{i,\bx_{1:i-1}}^{}(X_{j}^{\star} \neq X_{j}^{\star\star})
\,\leq\, \mbB_{\bm\pi}(i \centernot{\longleftrightarrow} j).
\ee
By the construction of $\mD_N(\truth)$,
the below-diagonal and diagonal elements of $\mD_N(\truth)$ are $0$ and $1$,
respectively.
We define $\pi^\star \in (0, 1)$ by
\beno
\pi^\star
&\coloneqq& \max\limits_{1 \leq v \leq M} \, \pi_v, 
\ee
and note that Lemma \ref{lem:pi_star_bound},
together with the assumption that $\truth \in \bTheta_N = \mR^p$ satisfies \eqref{condition.theta0}, 
implies that
\beno
\pi^{\star} &\leq& \dfrac{1}{1 + \exp(- L - \vartheta  \log N)} 
&<& 1,
\ee
where $L \in  [0, \infty)$ and $\vartheta \in [0, \infty)$ are the same constants as in \eqref{condition.theta0},
assumed to be independent of $N$ and $p$.  
Let
\beno
\pp 
&\coloneqq& \dfrac{1}{1 + \exp(-L - \vartheta \log N)}
\ee
and define the vector  
$\pvec \in [0, 1]^{M}$ by 
\beno
\xi_i 
&\coloneqq& 
\begin{cases}
0 & \mbox{ if } v \in \{1, \ldots, i-1\} \s\\
1 & \mbox{ if } v = i \s\\
\pp & \mbox{ if } v \in \{i+1, \ldots, M\}
\end{cases}.
\ee
Observe that the probabilities $\mbB_{\bm\pi}(i \centernot{\longleftrightarrow} j)$ of the events $\{i \centernot{\longleftrightarrow} j\}$ are non-decreasing in the coordinates of $\bm\pi$,
so that
\beno
\mbB_{\bm\pi}(i \centernot{\longleftrightarrow} j)
&\leq& \mbB_{\pvec}(i \centernot{\longleftrightarrow} j).
\ee
We bound the elements  $\mD_{i,j}(\truth)$ of $\mD_N(\truth)$ by bounding the probabilities $\mbB_{\pvec}(i \centernot{\longleftrightarrow} j)$ of the events $\{i \centernot{\longleftrightarrow} j\}$. 

\s\s

\noindent 
{\bf II. Bounding the elements of the coupling matrix $\mD_N(\truth)$.}
To bound the elements $\mD_{i,j}(\truth)$ of $\mD_N(\truth)$,
we bound the probabilities $\mbB_{\pvec}(i \dpath j)$ of the events $\{i \centernot{\longleftrightarrow} j\}$ using Assumption A.
To do so,
define
\beno
\mS_{\mG,i,k}
&\coloneqq& \left\{ v \in \mV \setminus \{i\} \,:\, d_{\mG}(i, v) = k \right\},
&& k = 1, \ldots, M-1,
\ee
where $d_{\mG}(i, v)$ is the graph distance (i.e., the length of the shortest path)
between vertices $i \in \mV$ and $v \in \mV$ in the conditional independence graph $\mG$. 
The set $\mS_{\mG,i,k} \subseteq \mV$ represents the subset of vertices in 
the conditional independence graph $\mG$ at graph distance $k$ from vertex $i$ in $\mG$.

We bound $\mbB_{\pvec}(i \dpath j)$ by placing restrictions on the subpopulation structure,
which determines which edges are present in $\mG$.  
To do so, 
define the {\it subpopulation graph} $\mG_{\mA}$ 
to be the graph with the set of subpopulations $\{\mA_1, \ldots, \mA_K\}$ as vertices  
and edges between vertices $\mA_r$ and $\mA_l$ if and only if $\mA_r \cap \mA_l \neq \emptyset$. 
In $\mG_{\mA}$, 
two vertices corresponding to subpopulations $\mA_r$ and $\mA_l$ are connected by an edge if and only if they overlap. 
Let $d_{\mG_{\mA}}(\mA_r, \mA_l)$ denote the graph distance (i.e., the length of the shortest path) 
between vertices $\mA_r$ and $\mA_l$ in $\mG_{\mA}$. 
Define,
for all $\mA_r \in \{\mA_1, \ldots, \mA_K\}$ and $k \in \{1, \ldots, K-1\}$,  
\beno
\mV_{\mA_r,k} 
&\coloneqq& \big\{\mA_l \in \{\mA_1, \ldots, \mA_K\} \setminus \{\mA_r\} \,:\, d_{\mG_{\mA}}(\mA_r, \, \mA_l) = k \big\}.  
\ee 
Let $g : \{1, 2, \ldots\} \mapsto [0, \, \infty)$ be such that,
for all $K \in \{1, 2, \ldots\}$, 
\beno
|\mV_{\mA_r,k}|
&\leq& g(k),
& k \in \{1, \ldots, K-1\}, 
&& \mbox{for all } \mA_r \in \{\mA_1, \ldots, \mA_K\}.
\ee
In words, 
$g(k)$ bounds the number of subpopulations at graph distance $k$ from any given subpopulation in $\mG_{\mA}$
for all conceivable subpopulations and thus all conceivable subpopulation structures,
i.e.,
for all $\mG_{\mA}$ defined for subpopulations $\mA_1, \ldots, \mA_K$ at all values of $K \in \{1, 2, \ldots\}$. 

Models 2 and 3 satisfy Definition \ref{def:intersection}
and posses the neighborhood intersection property.  
By Proposition \ref{prop:cond_ind},
the dependence neighborhood of any edge variable $X_i$ between nodes $\{a,b\} \subset \mN$
is not larger than the subset of edge indices contained in the set 
$\mM_{a,b} \coloneqq \big\{ 
\{c,d\} \,:\, 
c \in \mN_a \cup \mN_b, \, d \in \mN_a \cup \mN_b \setminus \{c\} 
\big\}$,
i.e.,
the edge variables contained in the set $\mathfrak{N}_i$ 
will correspond to edge variables between pairs of nodes in $\mM_{a,b}$.  
We construct a graph covering $\mG^\star$ of the conditional independence graph $\mG$ as follows:
\ben
\item Initialize $\mG^\star$ with the same set of vertices and edges as $\mG$. \s 
\item For each vertex $X_i$ in $\mG$ corresponding to an edge variable between nodes $\{a,b\} \subset \mN$ 
with degree greater than $0$ in $\mG$,
add edges between $X_i$ and any other edge variables $X_j$ contained in the subgraph $\bX_{\mM_{a,b}}$ 
which are not already present in $\mG$. 
\een
The construction of $\mG^\star$ ensures that the dependence neighborhood of any given vertex $X_i$ in $\mG^\star$  
corresponding to the edge variable between pair of nodes
$\{a,b\} \subset \mN$ is either empty or is equal to the set of vertices corresponding to 
edge variables contained in the subgraph $\bX_{\mM_{a,b}}$. 
Moreover, 
the fact that $\mG \subseteq \mG^\star$ implies  
\beno
\mbB_{\pvec}(i \centernot{\longleftrightarrow} j \mbox{ in } \mG)
&\leq& \mbB_{\pvec}(i \centernot{\longleftrightarrow} j \mbox{ in } \mG^\star).
\ee
In words, 
the probability of the existence of a path of disagreement 
does not decrease through the addition of edges in the graph. 
Henceforth and for ease of presentation, 
the event $i \centernot{\longleftrightarrow} j$ will represent a path of disagreement in the graph covering $\mG^\star$ of $\mG$
and we will assume  $\mS_{i,k} \equiv \mS_{\mG^\star,i,k}$.   

\s
 
We bound each $|\mS_{i,k}|$ ($k \in \{1, 2, \ldots\}$) for arbitrary $i \in \mV$ 
with non-zero degree in $\mG^\star$:  
\bi

\item {\it Bounding $|\mS_{i,1}|$.} Let $X_i$ denote the edge variable between pair of nodes $\{a,b\} \subset \mN$. 
By definition, 
$\mS_{i,1}$ contains all vertices in $\mG^\star$ corresponding to edge variables $X_j$
which lie in the dependence neighborhood of edge variable $X_i$ in $\mG^\star$.  
By the construction of $\mG^\star$, 
the dependence neighborhood of edge variable $X_i$ in $\mG^\star$ 
is equal to the set of edge variables contained in the subgraph $\bX_{\mM_{a,b}}$,
the number of which is bounded above by $4 \, D_N^2$: 
\beno
|\mM_{a,b}|
\leq |\mN_a \cup \mN_b|^2
\leq (|\mN_a| + |\mN_b|)^2
\leq (2 \, \max\{|\mN_a|, \, |\mN_b|\})^2
\leq 4 \, D_N^2,
\ee
where by Lemma \ref{lem:DN}, 
$D_N \geq \max_{t \in \mN} \, |\mN_t|$. 
Hence,
$|\mS_{i,1}| \leq 4 \, D_N^2$. \s 

\item {\it Bounding $|\mS_{i,2}|$.} 
Consider any $j \in \mS_{i,2}$ and let $X_{i}$ denote the edge variable 
between pair of nodes $\{a,b\} \subset \mN$. 
The shortest path between edge variables $X_i$ and $X_j$ in $\mG^\star$ is of length $2$,
implying the following:
\begin{enumerate}[leftmargin=.75in]
\item[(F.1)]  $X_j$ is not in the dependence neighborhood of $X_i$ in $\mG^\star$.
\item[(F.2)] In $\mG^\star$, there is at least one edge variable $X_l$ in the dependence neighborhood of $X_i$
such that $X_j$ is likewise in the dependence neighborhood of $X_l$.
\een
By the construction of $\mG^\star$,
facts (F.1) and (F.2) imply there exist 
\bi
\item a pair of nodes $\{v,w\} \subseteq \mN_a \cup \mN_b$
such that $X_l$ is the edge variable between $\{v,w\}$,
and
\item a pair of nodes $\{c,d\} \subseteq \mN_v \cup \mN_w$ such that
$X_j$ is the edge variable between $\{c,d\}$ and 
$\{c,d\} \not\subseteq \mN_a \cup \mN_b$,
otherwise $X_j$ would be in the dependence neighborhood of $X_i$,
violating the assumption that $j \in \mS_{i,2}$. 
\ei

\begin{center}
\begin{tikzpicture}[main node/.style={circle,fill=blue!20,draw,minimum size=1.5cm,inner sep=0pt,align=center}]
 \node[main node] (A) at (-3,0) {$X_i$\\ $\{a,b\}$}; 
 \node[main node] (B) at (0,0) {$X_l$\\ $\{v,w\}$};
 \node[main node] (C) at (3,0) {$X_j$\\ $\{c,d\}$};
  
 \path[draw, thick]
 (A) edge node {} (B) 
 (B) edge node {} (C); 

\end{tikzpicture}
\end{center}

Recall the definition, 
for all $v \in \mN$,
\beno
\mN_v
&\coloneqq&
\{w \in \mN \setminus \{v\} \,:\, \mbox{exists } r \in \{1, \ldots, K\} \mbox{ such that } \{v,w\} \subset \mA_r\}.
\ee
As $\{c,d\} \subseteq \mN_v \cup \mN_w$, 
there must exist $r,t \in \{1, \ldots, K\}$ such that: 
\bi  
\item 
either $\{c,v\} \subset \mA_r$ or $\{c,w\} \subset \mA_r$, and 
\item 
either $\{d,v\} \subset \mA_t$ or $\{d,w\} \subset \mA_t$. 
\ei
Since $\{c,d\} \not\subseteq \mN_a \cup \mN_b$,
therefore $\{a,b\} \not\subseteq  \mA_r \cup \mA_t \subseteq \mN_c \cup \mN_d$.
Finally, 
$\{v,w\} \subseteq \mN_a \cup \mN_b$ 
implies that there exists $n,m \in \{1, \ldots, K\}$ such that 
\bi
\item 
either $\{v,a\} \subset \mA_n$ or $\{v,b\} \subset \mA_n$, and 
\item 
either $\{w,a\} \subset \mA_m$ or $\{w,b\} \subset \mA_m$. 
\ei
As a result, 
$(\mA_n \cup \mA_m) \,\cap\, (\mA_r \cup \mA_t) \neq \emptyset$,
implying either $v$ or $w$ belong to a subpopulation 
$\mA_z \not\subseteq \mN_a \cup \mN_b$ ($z \in \{1, \ldots, K\}$)
for which 
$d_{\mG_{\mA}}(\mA_z, \mA_y) = 1$ for some $y \in \{1, \ldots, K\}$
with $\mA_y \subseteq \mN_a \cup \mN_b$,
i.e.,
a subpopulation with graph distance at least $1$ in $\mG_{\mA}$ 
from all subpopulations represented in $\mN_a \cup \mN_b$ 
and equal to $1$ for at least one such subpopulation. 
The same holds for either $c$ or $d$. 
Thus, 
\beno
|\mS_{i,2}|
&\leq& 2 \, D_N^3 \, (g(1) + 1) \, g(1), 
\ee
which follows from the following argument: 
\bi
\item First, the number of subpopulations contained in $\mN_a \cup \mN_b$
is bounded above by $2 (g(1) + 1)$,
because $g(1)$ bounds the number of subpopulations which overlap with any other subpopulation,
so that $g(1) + 1$ bounds the number of subpopulations to which any node $a \in \mN$ 
or $b \in \mN$ may belong; \s 

\item Second, 
the number of subpopulations with graph distance $1$ in the subpopulation graph $\mG_{\mA}$ 
to any subpopulation represented in 
$\mN_a \cup \mN_b$
is bounded above by $2 \, (g(1) + 1) \, g(1)$; \s

\item Third,
note that either $c$ or $d$ must be in one of the subpopulations 
with graph distance $1$ in $\mG_{\mA}$ to at least one of the subpopulations
represented in $\mN_a \cup \mN_b$. 
Without loss,
let this be $c$.  
Then the total number of such nodes $c$ which are contained in one of the aforementioned subpopulations at graph distance $1$ 
is bounded above by $2 \, D_N \, (g(1) + 1) \, g(1)$, 
using the bound $|\mA_k| \leq |\mN_i| \leq D_N$ from Lemma \ref{lem:DN}
which holds for all   
$k \in \{1, \ldots, K\}$ and 
$i \in \mA_k$. 

\s

\item Finally,
we bound the number of possible $d$ that may be paired with $c$. 
Note $X_j$ has non-zero degree in $\mG^\star$. 
By Proposition \ref{prop:cond_ind},
$X_j$ is independent of all other edges if $\mN_c \cap \mN_d = \emptyset$.
Thus, 
we bound the number of edge variables $X_j$ between node $c \in \mN$ and nodes $d \in \mN \setminus \{c\}$ 
for which $\mN_c \cap \mN_d \neq \emptyset$
using 
Lemma \ref{lem:DN}: 
\beno
\left| \left\{d \in \mN \setminus \{c\} \,:\, \mN_c \cap \mN_d \neq \emptyset \right\} \right|
&\leq& D_N^2.
\ee
Hence, 
the number of such $d$ (for a given $c$) 
numbers no more than $D_N^2$,
the total of which is bounded above by $2 \, D_N^3 \,  (g(1) + 1) \, g(1)$.

\ei

\hide{
\begin{figure}[t]
\centering
\includegraphics[width = .75 \linewidth, keepaspectratio]{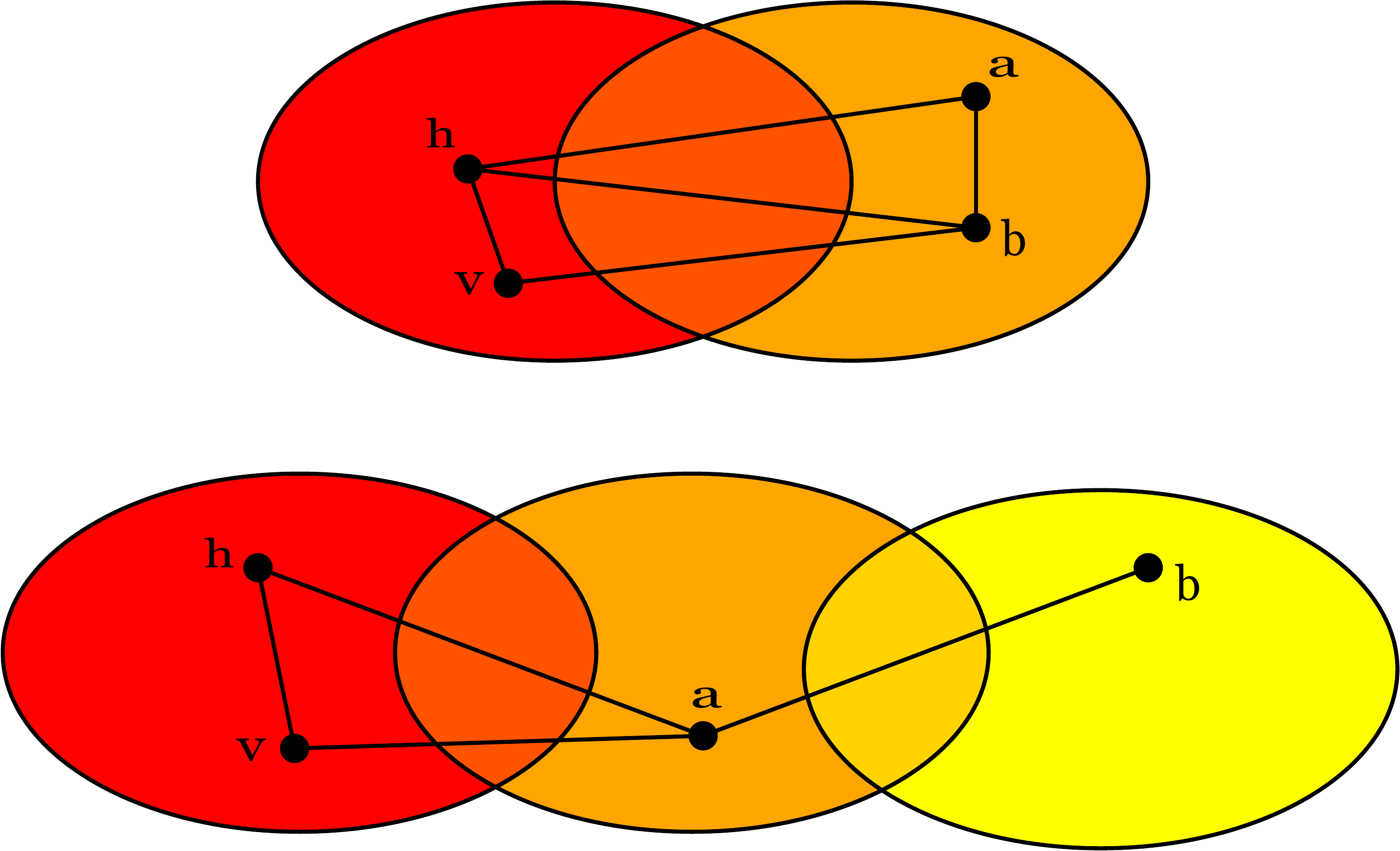}
\caption{\label{fig:lem11} 
Demonstrating how vertices $X_i$ and $X_j$ of the cover $\mG^\star$ of conditional independence graph $\mG$ can have graph distance $d_{\mG^\star}(i,\, j) = 2$ in $\mG^\star$.
Suppose that the indicator $X_j$ of an edge between nodes $v$ and $h$ depends on the indicator $X_i$ of an edge between nodes $a$ and $b$ through the indicator $X_l$ of an edge between nodes $a$ and $h$. 
By construction, 
$v$ either belongs to a subpopulation that intersects with one of the subpopulations of $b$,
or $a$ act as a broker for $v$ and $b$. 
}
\end{figure}
}

\item {\it Bounding $|\mS_{i,k}|$ for $k \in \{3, 4, \ldots\}$.} 
Consider any $k \in \{3, 4, \ldots\}$ and any $j \in \mS_{i,k}$.
Let $X_{i}$ be the edge variable between nodes $\{a,b\} \subset \mN$ 
and $X_j$ be the edge variable between pair of nodes $\{c,d\} \subset \mN$. 
For $j \in \mS_{i,k}$,
there must exist an $l \in \mS_{i,k-1}$ such that $j \in \mS_{l, 1}$. 
Let $X_l$ be the edge variable between nodes $\{v,w\} \subset \mN$. 
Leveraging arguments from the case bounding $|\mS_{i,1}|$ above, 
$\{c,d\} \subseteq \mN_v \cup \mN_w$, 
implying both  $c$ and $d$ belong to at least one subpopulation 
to which either $v$ or $w$ also belong. 
Assume the following: 
\begin{itemize}[leftmargin=.5in] 

\item[(A.1)] $h \in \mS_{i,k-1}$ if and only if 
$X_h$ 
is an edge variable between nodes $\{r,t\} \subset \mN$ 
and 
either $r$ or $t$ 
belongs to a subpopulation at graph distance 
$k-2$ in $\mG_{\mA}$ to a subpopulation represented in $\mN_a \cup \mN_b$
and not less than $k-2$ to all others in  $\mN_a \cup \mN_b$. 

\item[(A.2)] $h \in \mS_{i,k-2}$ if and only if 
$X_h$ is an edge variable between nodes $\{r,t\} \subset \mN$ 
and 
either $r$ or $t$ belongs to a subpopulation at graph distance $k-3$ in $\mG_{\mA}$ 
to a subpopulation represented in $\mN_a \cup \mN_b$
and not less than $k-3$ to all others in  $\mN_a \cup \mN_b$.  

\ei
We have shown above that (A.1) and (A.2) are satisfied when $k=3$:  
\bi
\item $h \in \mS_{i,1}$ corresponds to edge variables $X_h$ between nodes $\{r,t\} \subset \mN_a \cup \mN_b$,
in which case both $r$ and $t$ belong to a subpopulation represented in $\mN_a \cup \mN_b$
and therefore have graph distance $0$ in $\mG_{\mA}$ to a subpopulation represented in $\mN_a \cup \mN_b$.

\item $h \in \mS_{i,2}$ corresponds to edge variables $X_h$ between nodes $\{r,t\} \subset \mN$
for which either $r$ or $t$ belongs to a subpopulation which is not represented in $\mN_a \cup \mN_b$,
but which is at graph distance $1$ in $\mG_{\mA}$ to a subpopulation represented in $\mN_a \cup \mN_b$, 
and at graph distance no less than $1$ to all others in $\mN_a \cup \mN_b$. 
\ei
Assumptions (A.1) and (A.2) require 
that neither $v$ nor $w$ belong to a subpopulation at graph distance less than or equal to $k-2$ in $\mG_{\mA}$
from any subpopulation represented in $\mN_a \cup \mN_b$. 
Leveraging the argument used in the case bounding $|\mS_{i,2}|$: 
For $\{c,d\} \subseteq \mN_v \cup \mN_w$,
either $c$ or $d$ must belong to a subpopulation 
$\mA_r$ ($r \in \{1, \ldots, K\}$) 
jointly with either $v$ or $w$
which is at graph distance at least $k-1$ in in $\mG_{\mA}$ from any subpopulation represented in $\mN_i \cup \mN_j$. 
Thus, 
there must exist  some $q \in \{1, \ldots, K\}$ such that $\mA_r \cap \mA_q \neq \emptyset$
and for which both of the following are satisfied: 
\bi
\item $d_{\mG_{\mA}}(\mA_q, \mA_y) \geq k-2$ for all $y \in \{1, \ldots, K\}$ satisfying $\mA_y \subseteq \mN_a \cup \mN_b$;  
\item $d_{\mG_{\mA}}(\mA_q, \mA_y) = k-2$ for at least one $y \in \{1, \ldots, K\}$ satisfying $\mA_y \subseteq \mN_a \cup \mN_b$.
\ei 
Hence, 
there exists at least one $y \in \{1, \ldots, K\}$ satisfying $\mA_y \subseteq \mN_a \cup \mN_b$
for which $d_{\mG_{\mA}}(\mA_r, \mA_y) = k-1$.
Repeating the counting argument from the case bounding $|\mS_{i,2}|$ above shows that  
the number of such pairs $\{a,b\} \subset \mN$ is bounded above by 
\beno
|\mS_{i,k}|
&\leq& 2 \, D_N^3 \, (g(1) + 1) \, g(k-1), 
\ee
and the result is proved by induction.

\ei

From the above, 
we have the bounds $|\mS_{i,1}| \leq 4 \, D_N^2$
and 
\be
\label{eq:Si_bound}
|\mS_{i,k}|
&\leq& 2 \, D_N^3 \, (g(1) + 1) \, g(k-1),
&& k \in \{2, 3, \ldots\}. 
\ee
We proceed with bounding $\mD_{i,j}(\truth)$ under Assumption A.
Define the function $g: \{1, 2, \ldots\} \mapsto [0, \infty)$ by
\be
\label{eq:g_def}
g(k)
\= \omega_1 + \dfrac{\omega_2}{2 \, D^3_N} \, \log \, k, 
&& k \in \{1, 2, \dots\},
\ee
where $\omega_1 \geq 0$ and $\omega_2 \in [0, \omega_1]$ 
are independent of $N$ and $p$ by Assumption A, 
and $\omega_2 \in [0, \omega_1]$ additionally satisfies
\beno
\omega_2
&<& \dfrac{1}{(\omega_1 + 1) \, |\log(1 - \pp)|}. 
\ee
Using \eqref{eq:g_def} and \eqref{eq:Si_bound}, 
we obtain the bounds $|\mS_{i,1}| \leq 4 \, D_N^2$ 
and 
\beno
|\mS_{i,k}|
&\leq& (\omega_1 + 1) \, (2 \, D_N^3 \, \omega_1 + \omega_2 \, \log(k-1)),
&& k \in \{2, 3, \ldots\}. 
\ee
By construction, 
\beno
\mbB_{\pvec}(v \mbox{ is open in } \mG^\star)
&\leq& \pp 
&<& 1,
&& v \in \{i+1, \ldots, M\}.
\ee
For there to be a path of disagreement $i \dpath j$ in $\mG^\star$ 
between a vertex $i \in \{1, \ldots, M\}$ and $j \in \mS_{i,k}$,
there must be at least one open vertex in each of the sets $\mS_{i,1}, \ldots, \mS_{i,k-1}$ and $j$ must be open 
(placing no restrictions on the connectedness of vertices within sets or between two sequential sets);
note that $i$ is open with probability $1$.
The probability that there exists at least one open vertex $v \in \mS_{i,1}$
is bounded above by 
\beno
1 - (1 - \pp)^{|\mS_{i,1}|}
&\leq& 1 - (1 - \pp)^{4 \,D_N^2}
&\leq& 1,  
\ee
and for $k \in \{2, 3, \ldots\}$, 
the same is bounded above by
\beno
1 - (1 - \pp)^{|\mS_{i,k}|}
&\leq& 1 - (1 - \pp)^{(\omega_1 + 1) \, (2 \, D_N^3 \, \omega_1 + \omega_2 \, \log(k-1))} \s\\
\= 1 - (1 - \pp)^{C_1 D_N^3 + C_2 \log(k-1)},
\ee
defining $C_1 \coloneqq 2 \, \omega_1 \, (\omega_1 + 1) \in [0, \infty)$ 
and $C_2 \coloneqq \omega_2  \, (\omega_1 + 1) \in [0, \infty)$. 
Since the events that vertices are open are independent under the Bernoulli product measure $\mbB_{\pvec}$,
we obtain
\beno
\mbB_{\pvec}(i \dpath j)
&\leq& \pp \, 
\left[1 - (1 - \pp)^{4 D_N^2} \right]
\dprod\limits_{l=2}^{k-1} \,
\left[ 1 - (1 - \pp)^{C_1 D_N^3 + C_2 \log(l-1)} \right] \s\\
\lte \pp \, 
\left[1 - (1 - \pp)^{4 D_N^2} \right] \, 
\left[ 1 - (1 - \pp)^{C_1 D_N^3 + C_2 \log(k-2)} \right]^{k-2} \s\\
\lte \left[ 1 - (1 - \pp)^{C_1 D_N^3 + C_2 \log(k-2)} \right]^{k-2},  
\ee
by the monotonicity of logarithms. 
We then bound 
\beno
1 - (1 - \pp)^{C_1 D_N^3 + C_2 \log(k-2)} 
\lte \exp\left(- (1 - \pp)^{C_1 D_N^3 + C_2 \log(k-2)}\right),
\ee
using the inequality $1 - z \leq \exp(-z)$ ($z \in (0, 1)$).
We proceed by writing  
\beno 
&& \exp\left(- (1 - \pp)^{C_1 D_N^3 + C_2 \log(k-2)}\right) \s\s\\
\= \exp\left( - \exp\left(\left[C_1 \, D_N^3 + C_2 \log(k-2)\right] \, \log(1 - \pp)\right) \right) \s\s\\ 
\= \exp\left( - \exp\left(-\left[C_1 \, D_N^3 + C_2 \log(k-2)\right] \, |\log(1 - \pp)|\right) \right) \s\s\\
\= \exp\left( - \exp(-C_1 \, D_N^3 \, |\log(1 - \pp)|) \; (k-2)^{-C_2 \, |\log(1-\pp)|} \right), 
\ee
where in the above we used the fact that $\log(1 -\pp) < 0$. 
Define 
\beno
A &\coloneqq& \exp(- C_1  \, D_N^3 \, |\log(1 - \pp)|)
&\in& (0, 1), 
\ee
noting that $D_N \geq 1$ under Models 2 and 3. 
Then
\beno
\left[ 1 - (1 - \pp)^{C_1 D_N^3 + C_2 \log(k-2)} \right]^{k-2}
&\leq& \left[ \exp\left( -A \, (k-2)^{-C_2 \, |\log(1-\pp)|} \right) \right]^{k-2} \s\s\\
\=  \exp\left( -A \, (k-2)^{1-C_2 \, |\log(1-\pp)|} \right), 
\ee
demonstrating the bound (for $i \in \{1, \ldots, M\}$ and $j \in \mS_{i,k}$) 
\beno
\mbB_{\bm\pi}(i \dpath j)
&\leq& \mbB_{\pvec}(i \dpath j)
&\leq& \exp\left( -A \, (k-2)^{1-C_2 \, |\log(1-\pp)|} \right). 
\ee
We hence obtain (for $i \in \{1, \ldots, M\}$ and $j \in \mS_{i,k}$)
\be
\label{eq:bound_cond1}
\mD_{i,j}(\truth)
\lte \exp\left( -A \, (k-2)^{1-C_2 \, |\log(1-\pp)|} \right). 
\ee

\s\s

\noindent {\bf III. Bounding the spectral norm $\mnorm{\mD_N(\truth)}_2$ of the coupling matrix $\mD_N(\truth)$.}
To bound the spectral norm $\mnorm{\mD_N(\truth)}_2$ of the coupling matrix $\mD_N(\truth)$,
we first use H\"older's inequality to obtain
\beno
\mnorm{\mD_N(\truth)}_2
&\leq& \sqrt{\mnorm{\mD_N(\truth)}_1 \; \mnorm{\mD_N(\truth)}_{\infty}}.
\ee
Next, 
we form a symmetric $M \times M$ matrix $\mT$ by defining
\beno
\mT
&\coloneqq& \mD_N(\truth) + \mD_N(\truth)^{\top} - \mbox{diag}(\mD_N(\truth)),
\ee
where $\mD_N(\truth)^{\top}$ is the $M \times M$ transpose of $\mD_N(\truth)$ 
and $\mbox{diag}(\mD_N(\truth))$ is the $M \times M$ diagonal matrix with elements 
$\mD_{1,1}, \dots, \mD_{M,M}$ on the main diagonal. 
By the construction of $\mT$,
the elements $\mT_{i,j}$ of $\mT$ are given by
\beno
\mT_{i,j} \=
\begin{cases}
\mD_{i,j}(\truth) & \mbox{ if } j > i\s
\\
\mD_{i,i}(\truth) & \mbox{ if } i = j\s
\\
\mD_{j,i}(\truth) & \mbox{ if } j < i
\end{cases},
\ee
where $\mD_{i,i}(\truth) = 1$ ($i = 1, \ldots, M$) 
by the definition of $\mD_{N}(\truth)$. 
Using the fact that $\mT_{i,j} = \max(\mD_{i,j}(\truth),\, \mD_{j,i}(\truth))$ ($i, j = 1, \dots, M$),
we obtain
\beno
\mnorm{\mD_N(\truth)}_1
\= \max\limits_{1 \leq j \leq M} \, \dsum_{i = 1}^{M} \, |\mD_{i,j}(\truth)|
\lte \max\limits_{1 \leq j \leq M} \, \dsum_{i = 1}^{M} \, |\mT_{i,j}|
\= \mnorm{\mT}_1
\ee
and
\beno
\mnorm{\mD_N(\truth)}_{\infty}
\= \max\limits_{1 \leq i \leq M} \, \dsum_{j = 1}^{M} \, |\mD_{i,j}(\truth)|
\lte \max\limits_{1 \leq i \leq M} \, \dsum_{j = 1}^{M} \, |\mT_{i,j}|
\= \mnorm{\mT}_{\infty}.
\ee
In addition,
we know that $\mT_{i,j} = \mT_{j,i}$ ($i, j = 1, \dots, M$),
which implies that
\beno
\mnorm{\mT}_1
\= \mnorm{\mT^{\top}}_{\infty}
\= \mnorm{\mT}_{\infty}.
\ee
As a consequence,
we obtain
\beno
\mnorm{\mD_N(\truth)}_2
&\leq& \sqrt{ \mnorm{\mD_N(\truth)}_1 \; \mnorm{\mD_N(\truth)}_{\infty}}
&\leq& \sqrt{ \mnorm{\mT}_1 \; \mnorm{\mT}_{\infty}}
\= \mnorm{\mT}_{\infty},
\ee
where $\mnorm{\mT}_{\infty}$ can be bounded above by using \eqref{boundingdij}:
\beno
\mnorm{\mT}_{\infty}
\= \max\limits_{1 \leq i \leq M} \, \dsum_{j = 1}^{M} \, |\mT_{i,j}|
&\leq& 1 + \max\limits_{1 \leq i \leq M} \, \dsum_{j=1:\, j \neq i}^{M} \mbB_{\pvec}(i \centernot{\longleftrightarrow} j).
\ee
Hence using \eqref{eq:bound_cond1} along with Assumption A,
\beno
\mnorm{\mD_N(\truth)}_2 &\leq& 1 + \max\limits_{1 \leq i \leq M} \, \dsum_{j=1:\, j \neq i}^{M} \mbB_{\pvec}(i \centernot{\longleftrightarrow} j)
\ee
and,
for all $i \in \{1, \ldots, M\}$, 
\beno
&& \dsum_{j=1:\, j \neq i}^{M} \mbB_{\pvec}(i \centernot{\longleftrightarrow} j) \s\\
&\leq& |\mS_{i,1}| + \dsum_{k=2}^{\infty} \, |\mS_{i,k}| \, \exp\left( -A \, (k-2)^{1-C_2 \, |\log(1-\pp)|} \right) \s\\
&\leq& 4 \, D_N^2 
+\dsum_{k=2}^{\infty} \, 
(C_1 \, D_N^3 + C_2 \, \log(k-1))
\exp\left( -A \, (k-2)^{1-C_2 \, |\log(1-\pp)|} \right),
\ee
using the above bounds on the number of vertices $|\mS_{i,k}|$ which are at graph distance in $k$ 
to any given vertex $i \in \mV$ in $\mG^\star$. 
We focus on bounding the infinite series 
\beno
&& \dsum_{k=2}^{\infty} \,
(C_1 \, D_N^3 + C_2 \, \log(k-1))
\exp\left( -A \, (k-2)^{1-C_2 \, |\log(1-\pp)|} \right) \s\\
\= C_1 \, D_N^3  \dsum_{k=2}^{\infty} \exp\left( -A \, (k-2)^{1-C_2 \, |\log(1-\pp)|} \right) \s\\
&+& C_2 \dsum_{k=2}^{\infty} \, \log(k-1) \, \exp\left( -A \, (k-2)^{1-C_2 \, |\log(1-\pp)|} \right). 
\ee
For the first series, 
\beno
\dsum_{k=2}^{\infty} \exp\left( -A \, (k-2)^{1-C_2 \, |\log(1-\pp)|} \right)
\,=\, 1 + \dsum_{k=1}^{\infty} \exp\left(- A \, k^{1-C_2 \, |\log(1-\pp)|} \right), 
\ee
noting that $\exp\left( -A \, (k-2)^{1-C_2 \, |\log(1-\pp)|} \right) = 1$ when $k = 2$. 
By a Taylor expansion of $\exp(z)$, 
we can establish the inequality $\exp(z) > z^u \,/\, u!$ for any $z > 0$ and any $u \in \{1, 2, \dots\}$,
which in turn establishes the inequality $\exp(-z) < u! \,/\, z^u$.
Using this inequality, 
\be
\label{eq:exp_inq}
\exp\left(- A \, k^{1-C_2 \, |\log(1-\pp)|} \right)
&<& \dfrac{u!}{A^u \, k^{u \, (1 - C_2 \, |\log(1-\pp)|)}}. 
\ee
Assume that 
$1 - C_2 \, |\log(1-\pp)| > 0$,
which is satisfied when 
\beno
\omega_2 &<& \dfrac{1}{(\omega_1 + 1) \, |\log(1 - \pp)|}, 
\ee
recalling $C_2 \coloneqq \omega_2 \, (\omega_1 + 1) > 0$. 
Taking $u = \lceil 2 \,/\, (1 - C_2 \, |\log(1-\pp)|) \rceil > 0$,
\beno
\dsum_{k=1}^{\infty} \exp\left(- A \, k^{1-C_2 \, |\log(1-\pp)|} \right)
&\leq&  \dfrac{u!}{A^u} \, \dsum_{k=1}^{\infty} \, \dfrac{1}{k^2}
\= \dfrac{u!}{A^{u}} \, \left(\dfrac{\pi^2}{6} \right).  
\ee
Thus, 
the first infinite series is bounded above by 
\be
\label{eq:series_b1}
C_1 \, D_N^3 \, \left( 1 + \dfrac{u!}{A^{u}} \, \left(\dfrac{\pi^2}{6} \right) \right). 
\ee
For the second infinite series, 
we write 
\beno
&& C_2 \, \dsum_{k=2}^{\infty} \, \log(k-1) \, \exp\left( -A \, (k-2)^{1-C_2 \, |\log(1-\pp)|} \right) \s\\ 
\= C_2 \, \dsum_{k=3}^{\infty} \, \log(k-1) \, \exp\left( -A \, (k-2)^{1-C_2 \, |\log(1-\pp)|} \right) \s\\
\= C_2 \, \dsum_{k=1}^{\infty} \, \log(k+1) \, \exp\left(-A \, k^{1-C_2 \, |\log(1-\pp)|} \right). 
\ee
We employ \eqref{eq:exp_inq} once more to show that 
\beno
\exp\left(-A \, k^{1-C_2 \, |\log(1-\pp)|} \right)
&\leq& \dfrac{u!}{A^u \, k^3},  
\ee
taking $u = \lceil 3 \,/\, (1-C_2 \, |\log(1-\pp)|) \rceil > 0$ this time. 
Thus, 
\beno
\dsum_{k=1}^{\infty} \, \log(k+1) \, \exp\left(-A \, k^{1-C_2 \, |\log(1-\pp)|} \right) 
\;\;\leq\;\; \dfrac{u!}{A^u} \, \dsum_{k=1}^{\infty} \, \dfrac{\log(k+1)}{k^3} \s\\
\;\;\leq\;\; \dfrac{u!}{A^u} \, \dsum_{k=1}^{\infty} \, \dfrac{k}{k^3} 
\;\;=\;\;  \dfrac{u!}{A^u} \, \dsum_{k=1}^{\infty} \, \dfrac{1}{k^2} 
\;\;=\;\; \dfrac{u!}{A^{u}} \, \left(\dfrac{\pi^2}{6} \right), 
\ee
using the inequality $\log(z+1) \leq z$ for $z \in  (0, \infty)$. 
As a result, 
the second infinite series is bounded above by 
\be
\label{eq:series_b2}
C_2 \, \dfrac{u!}{A^{u}} \, \left(\dfrac{\pi^2}{6} \right). 
\ee
Combining \eqref{eq:series_b1}, 
\eqref{eq:series_b2},
and the bound $2 \geq \pi^2 \,/\, 6$,   
\beno
\mnorm{\mD_N(\truth)}_{2}
&\leq& 1 + 4 \, D_N^2 + 
C_1 \, D_N^3 \, \left(1 + 2 \, \dfrac{u!}{A^{u}} \right) 
+ 2 \, C_2 \, \dfrac{u!}{A^{u}}\s\\ 
\lte 1 + 4 \, D_N^2 + \max\{C_1, C_2\} \, A^{-u} \left((A^{u} + 2 \, u!) \, D_N^3  + 2 \, u! \right) \s\\
\lte 1 + 4 \, D_N^2 + \max\{C_1, C_2\} \, A^{-u} \left( (1 + 2 \, u!) \, D_N^3 + 1 + 2 \, u! \right)\s\\ 
\lte 1 + 4 \, D_N^2 + \max\{C_1, C_2\} \, A^{-u} \, (1 + 2 \, u!)
\left(D_N^3 + 1\right), 
\ee
recalling that 
$A \coloneqq \exp(- C_1  \, D_N^3 \, |\log(1 - \pp)|)$,
which implies $A^{u} \in (0, 1)$. 
Next,
using the definitions of $C_1$ and  $C_2$,
and the assumption that $\omega_2 \leq \omega_1$, 
\beno
\max\{C_1, \, C_2\}
\= \max\{2 \, \omega_1 \, (\omega_1 + 1), \, \omega_2 \, (\omega_1+1)\}
&\leq& 2 \, \omega_1 \, (\omega_1 + 1).  
\ee
Then, 
there exist finite constants 
\beno
C_3 
&\coloneqq& 2 \, (\omega_1 + 1) \, (1 + 2 \, u!) 
&>& 0
\ee
and $C_4 \coloneqq C_1 \, |\log(1- \pp)| > 0$,
independent of $N$ and $p$, 
such that 
\beno
\mnorm{\mD_N(\truth)}_{2}
&\leq& 1 +  4 \,D_N^2  + \omega_1 \, C_3 \, \exp(C_4\, D_N^3) \, (1 + D_N^3). 
\ee
We complete the proof by noticing two key facts:
\bi
\item If $\omega_1 = 0$, 
then $\omega_1 \, (1 + D_N^3) \leq 2 \, \omega_1 \, D_N^3$ for all $D_N^3$. \s 
\item Since $D_N \geq 1$ under Models 2 and 3, 
$\omega_1 \, (1 + D_N^3) \leq 2 \, \omega_1 \, D_N^3$. 
\ei
Thus, 
there exists $C_5 \coloneqq 2 \, C_3 > 0$, 
independent of $N$ and $p$, 
such that  
\beno
\mnorm{\mD_N(\truth)}_{2}
&\leq& 1 + 4 \, D_N^2 + \omega_1 \, C_5  \, D_N^3 \, \exp(C_4 \, D_N^3) \s\\
&\leq& 1 + 4 \, D_N^2 + \omega_1 \, C_5 \, \exp(C_4 \, D_N^3 + 3 \, \log D_N) \s\\
&\leq& 1 + 4 \, D_N^2 + \omega_1 \, C_5 \, \exp(C_6 \, D_N^3),  
\ee
taking $C_6 \coloneqq (3 + C_4) > 0$ and since $\log D_N \leq D_N^3$.
Recall that  
\beno
\pp 
&\coloneqq& \dfrac{1}{1 + \exp(-L - \vartheta \, \log N)}
\ee
is bounded away from $1$ when $\vartheta = 0$. 

We then have the following cases:   
\bi
\item If subpopulations overlap,
i.e., 
$\omega_1 > 0$ and $\omega_2 \in [0, \omega_1]$,
provided 
\beno
\omega_2 &<& \dfrac{1}{(1 + \omega_1) \, |\log(1 - \pp)|},
\ee
then 
\beno
\mnorm{\mD_N(\truth)}_{2}
&\leq& 1 + 4 \,D_N^2 + \omega_1  \, C_5 \, \exp(C_6 \, D_N^3), 
\ee
provided $\vartheta = 0$ to ensure the constants are independent of $N$ and $p$.  

\s

\item If subpopulations do not overlap, 
i.e.,
$\omega_1 = \omega_2 = 0$, 
then 
\beno
\mnorm{\mD_N(\truth)}_{2}
&\leq& 1 + 4 \, D_N^2.
\ee
Since $D_N$ does not depend on $\pp$, 
we allow $\vartheta > 0$. 
\qed
\ei

\hide{

\begin{lemma}
\label{lem:bound_P} 
Consider Models 1, 2, and 3 with $\nat \in \mbR^p$ and $\alpha \in [0, \, 1/2)$. 
Then there exist functions $L_{k} : \mbR^p \mapsto (0, 1)$ and $U_{k} : \mbR^p \mapsto (0, 1)$ 
($k = 0, 1$) such that, 
for all $\{i,j\} \subset \mN$ and $\bx_{-\{i,j\}} \in \{0, 1\}^{M-1}$,
\beno
0 &<& L_k(\nat) 
&\leq&  \mbP_{\nat}(X_{i,j} = k \mid \bX_{-\{i,j\}} = \bx_{-\{i,j\}}) 
&\leq& U_k(\nat) &<& 1. 
\ee 
The functions $L_k(\nat)$ and $U_k(\nat)$ ($k = 0,1$) are given by 
\beno
L_1(\nat)
&\coloneqq& \begin{cases}
\dfrac{1}{1 + \exp((3 + D_N) \, \norm{\nat}_{\infty})} & \hspace{.12in} \mbox{if } \mN_i \cap \mN_j \neq \emptyset \s\s\\
\dfrac{N^{-\alpha}}{1 + \exp(3 \, \norm{\nat}_{\infty})} & \hspace{.12in} \mbox{if } \mN_i \cap \mN_j = \emptyset
\end{cases} 

\s\s\\

U_1(\nat)
&\coloneqq& \begin{cases}
\dfrac{1}{1 + \exp(-(3 + D_N) \, \norm{\nat}_{\infty})} & \hspace{.01in} \mbox{if } \mN_i \cap \mN_j \neq \emptyset \s\s\\
\dfrac{1}{1 + \exp(-3 \, \norm{\nat}_{\infty}) \, N^{\alpha}} & \hspace{.01in} \mbox{if } \mN_i \cap \mN_j = \emptyset
\end{cases}

\s\s\\

L_0(\nat)
&\coloneqq& \begin{cases}
\dfrac{1}{1 + \exp((3+D_N) \, \norm{\nat}_{\infty})} & \hspace{.12in} \mbox{if } \mN_i \cap \mN_j \neq \emptyset \s\s\\
\dfrac{1}{1 + \exp(3 \, \norm{\nat}_{\infty}) \, N^{-\alpha}} & \hspace{.12in} \mbox{if } \mN_i \cap \mN_j = \emptyset
\end{cases}

\s\s\\

U_0(\nat)
&\coloneqq& \begin{cases}
\dfrac{1}{1 + \exp(-(3+D_N) \, \norm{\nat}_{\infty})} & \hspace{0in} \mbox{if } \mN_i \cap \mN_j \neq \emptyset \s\s\\ 
\dfrac{1}{1 + \exp(-3 \, \norm{\nat}_{\infty}) \, N^{-\alpha}} & \hspace{0in} \mbox{if } \mN_i \cap \mN_j = \emptyset.
\end{cases},
\ee
\end{lemma}

\s 

\llproof \ref{lem:bound_P}. 
Consider any pair of nodes $\{i,j\} \subset \mN$ and any $\bx_{-\{i,j\}} \in \{0, 1\}^{\binom{N}{2}-1}$.
We can express the full conditional probability  
\beno
\mbP_{\nat}(X_{i,j} = x_{i,j} \mid \bX_{-\{i,j\}} = \bx_{-\{i,j\}})
\ee
two different ways depending whether $\mN_i$ and $\mN_j$ are disjoint. 

First, 
if $\mN_i \cap \mN_j = \emptyset$, 
\beno
\label{eq:eq_sparsity} 
&& \mbP_{\nat}(X_{i,j} = x_{i,j} \mid \bX_{-\{i,j\}} = \bx_{-\{i,j\}}) \s\s\\
&=& \dfrac{\exp(\langle \nat, \, s(\bx_{-\{i,j\}},\, x_{i,j}) \rangle) \, N^{-\alpha \, x_{i,j}}}{\exp(\langle \nat, \, s(\bx_{-\{i,j\}},\, x_{i,j}=0) \rangle) + \exp(\langle \nat, \, s(\bx_{-\{i,j\}},\, x_{i,j}=1) \rangle)\; N^{-\alpha}}\s\s
\\
\= \dfrac{1}{g(0;\, \bx_{-\{i,j\}},\, x_{i,j},\, \nat)\; N^{\alpha \, x_{i,j}} + g(1;\, \bx_{-\{i,j\}},\, x_{i,j},\, \nat)\; N^{-\alpha \, (1-x_{i,j})}}, 
\ee
defining,
for $y \in \{0, 1\}$,  
\beno
g(y;\, \bx_{-\{i,j\}},\, x_{i,j},\, \nat)
&\coloneqq& \exp(\langle \nat, \, s(\bx_{-\{i,j\}},\, y) - s(\bx_{-\{i,j\}},\, x_{i,j})\rangle). 
\ee
Note $g(x_{i,j}; \, \bx_{-\{i,j\}},\, x_{i,j},\, \nat) = 1$
for all $\bx_{-\{i,j\}} \in \{0,1\}^{\binom{N}{2}-1}$ and all $\nat \in \mbR^p$.  

Second, 
if $\mN_i \cap \mN_j \neq \emptyset$, 
\beno
\label{eq:eq_no_sparsity}
&& \mbP_{\nat}(X_{i,j} = x_{i,j} \mid \bX_{-\{i,j\}} = \bx_{-\{i,j\}}) \s\s\\
&=& \dfrac{\exp(\langle \nat, \, s(\bx_{-\{i,j\}},\, x_{i,j}) \rangle)}{\exp(\langle \nat, \, s(\bx_{-\{i,j\}},\, x_{i,j}=0) \rangle) + \exp(\langle \nat, \, s(\bx_{-\{i,j\}},\, x_{i,j}=1) \rangle)}\s\s
\\
\= \dfrac{1}{1 + g(1-x_{i,j};\, \bx_{-\{i,j\}},\, x_{i,j},\, \nat)}. 
\ee

Next, 
observe that 
\beno
&& \max\limits_{\bx_{-\{i,j\}}\, \in\, \{0, 1\}^{M-1}} 
\, \left| s_l(\bx_{-\{i,j\}},\, x_{i,j}=0) - s_l(\bx_{-\{i,j\}},\, x_{i,j}=1)  \right|\s\s
\\
&\leq& \begin{cases}
0 & \mbox{ if } l \in \{1, \dots, N\} \setminus \{i,j\}\s
\\
1 & \mbox{ if } l \in \{i,j\}\s 
\\
1+D_N & \mbox{ if } l = N + 1 \mbox{ and } \mN_i \cap \mN_j \neq \emptyset\s
\\ 0 & \mbox{ if } l = N + 1 \mbox{ and } \mN_i \cap \mN_j = \emptyset
\end{cases}.
\ee
The bound on $s_{N+1}$ follows from Lemma \ref{lem:brokerage_ham_bound},
whereas the conditions follow from Proposition \ref{prop:cond_ind},
which implies $X_{i,j}$ is independent of all other edges in the graph when $\mN_i \cap \mN_j = \emptyset$. 
The bound on $s_{l}(\bx)$ ($l \in \mN$) follows because $s_l(\bx) = \sum_{h \in \mN \setminus \{l\}} x_{l,h}$ is a 
function of $x_{i,j}$ if and only if $l \in \{i,j\}$. 
As a result, 
the triangle inequality and the above bounds reveal 
\beno 
&& \left|\left\langle \nat, \, s(\bx_{-\{i,j\}},\, x_{i,j}=1)\right\rangle - \left\langle\nat, \, s(\bx_{-\{i,j\}},\, x_{i,j}=0)\right\rangle \right| \s\s\\  
&\leq& 
\begin{cases}
\left(3 + D_N  \right) \, \norm{\nat}_{\infty} & \mbox{if } \mN_i \cap \mN_j \neq \emptyset \s\\ 
3 \, \norm{\nat}_{\infty} & \mbox{if } \mN_i \cap \mN_j = \emptyset 
\end{cases},
\ee
implying,
for $\{i,j\} \subseteq \mN$ with $\mN_i \cap \mN_j \neq \emptyset$, 
\beno
\exp(-(3+D_N) \, \norm{\nat}_{\infty})
\,\leq\, g(1-x_{i,j};\, \bx_{-\{i,j\}},\, x_{i,j},\, \nat) 
\,\leq\, \exp((3+D_N) \, \norm{\nat}_{\infty}),
\ee
and for $\{i,j\} \subseteq \mN$ with $\mN_i \cap \mN_j = \emptyset$, 
\beno
\exp(- 3 \, \norm{\nat}_{\infty})
\,\leq\, g(1-x_{i,j};\, \bx_{-\{i,j\}},\, x_{i,j},\, \nat)
\,\leq\, \exp(3 \, \norm{\nat}_{\infty}). 
\ee
As a result, 
for all $k \in \{0, 1\}$, 
\beno
0 &<& L_k(\nat) 
&\leq& \mbP(X_{i,j} = k \,|\, \bX_{-\{i,j\}} =  \bx_{-\{i,j\}}) 
&\leq& U_k(\nat)
&<& 1,
\ee
where 
\beno
L_1(\nat)
&\coloneqq& \begin{cases}
\dfrac{1}{1 + \exp((3 + D_N) \, \norm{\nat}_{\infty})} & \mbox{if } \mN_i \cap \mN_j \neq \emptyset \s\s\\
\dfrac{N^{-\alpha}}{1 + \exp(3 \, \norm{\nat}_{\infty})} & \mbox{if } \mN_i \cap \mN_j = \emptyset  
\end{cases}
\ee
and 
\beno
U_1(\nat)
&\coloneqq& \begin{cases}
\dfrac{1}{1 + \exp(-(3 + D_N) \, \norm{\nat}_{\infty})} & \mbox{if } \mN_i \cap \mN_j \neq \emptyset \s\s\\
\dfrac{1}{1 + \exp(-3 \, \norm{\nat}_{\infty}) \, N^{\alpha}} & \mbox{if } \mN_i \cap \mN_j = \emptyset
\end{cases}.
\ee
We obtain $L_0(\nat)$ and $U_0(\nat)$ 
by noting 
\beno
\mbP(X_{i,j} = 0 \,|\, \bX_{-\{i,j\}} =  \bx_{-\{i,j\}})
\= 1 - \mbP(X_{i,j} = 1 \,|\, \bX_{-\{i,j\}} =  \bx_{-\{i,j\}}),
\ee
implying 
\beno
1 - U_1(\nat)
&\leq& \mbP(X_{i,j} = 0 \,|\, \bX_{-\{i,j\}} =  \bx_{-\{i,j\}})
&\leq& 1 - L_1(\nat),
\ee
which allows us to obtain 
\beno
L_0(\nat)
&\coloneqq& \begin{cases}
\dfrac{1}{1 + \exp((3+D_N) \, \norm{\nat}_{\infty})} & \mbox{if } \mN_i \cap \mN_j \neq \emptyset \s\s\\
\dfrac{1}{1 + \exp( 3\, \norm{\nat}_{\infty}) \, N^{-\alpha}} & \mbox{if } \mN_i \cap \mN_j = \emptyset 
\end{cases}
\ee
and 
\beno
U_0(\nat) 
&\coloneqq& \begin{cases}
\dfrac{1}{1 + \exp(-(3+D_N) \, \norm{\nat}_{\infty})} & \mbox{if } \mN_i \cap \mN_j \neq \emptyset \s\s\\
\dfrac{1}{1 + \exp(-3 \, \norm{\nat}_{\infty}) \, N^{-\alpha}} & \mbox{if } \mN_i \cap \mN_j = \emptyset
\end{cases}. 
\ee
\qed 

\s\s 

\begin{lemma}
\label{lem:brokerage_ham_bound}
Consider 
\beno
s_{N+1}(\bx) 
\= \dsum_{i < j}^N\, x_{i,j} \, I_{i,j}(\bx),
\ee
where
\beno
I_{i,j}(\bx) 
\= 
\begin{cases}
0 & \mbox{if } \mN_i \cap \mN_j = \emptyset\s
\\
\one\left( \, \dsum_{h \in \mN_i \cap \mN_j} \, x_{i,h} \, x_{j,h} \,\geq\, 1\right)
& \mbox{if } \mN_i \cap \mN_j \neq \emptyset.
\end{cases}
\ee
Then,
for all $\{i,j\} \subset \mN$,  
\beno
\max\limits_{(\bx,\,\bx^\prime)\, \in\, \mbX \times \mbX:\; x_{v,w} = x_{v,w}^\prime,\; \{v,\, w\} \neq \{i,\, j\}}
|s_{N+1}(\bx) - s_{N+1}(\bx^\prime)|
&\leq& 1+D_N,
\ee
where $D_N \coloneqq \max_{\{i,j\} \subset \mN} \, |\mathfrak{N}_{i,j}|$. 
\end{lemma}

\s\s 

\llproof \ref{lem:brokerage_ham_bound}. 
Consider any pair of nodes $\{i,j\} \subset \mN$. 
The number of $x_{a,b}\, I_{a,b}(\bx)$ ($\{a,b\} \neq \{i,j\}$)
which are a function of $x_{i,j}$ 
includes all $\{a,b\} \neq \{i,j\}$ satisfying one of the following: 
\bi
\item $\{a,b\} = \{i,b\}$ and $j \in \mN_i \cap \mN_b \subseteq \mN_i$. 
\item $\{a,b\} = \{j,b\}$ and $i \in \mN_j \cap \mN_b \subseteq \mN_j$. 
\ei

From the proof of Proposition \ref{prop:cond_ind},
the number of $x_{a,b}\, I_{a,b}(\bx)$ ($\{a,b\} \neq \{i,j\}$)
which are a function of $x_{i,j}$ is bounded above by 
$|\{i,j\} \times \mN_i \cup \mN_j|$. 
Proposition 2 establishes that,
in the case when $\mN_i \cap \mN_j \neq \emptyset$, 
$|\mathfrak{N}_{i,j}|$ is at least as big as $|\{i,j\} \times \mN_i \cap \mN_j|$. 
Defining $D_N \coloneqq \max_{\{i,j\} \subset \mN} |\mathfrak{N}_{i,j}|$, 
the number of summands which are a function of $x_{i,j}$ is bounded above by $1 + D_N$,
now counting the case when $\{a,b\} = \{i,j\}$.  
Consider any $(\bx,\, \bx^\prime) \in \mbX \times \mbX$ such that $x_{v,w} = x_{v,w}^\prime$ for all $\{v,w\} \neq \{i,j\}$.  
Then, 
by the triangle inequality,
\beno
|s_{N+1}(\bx) - s_{N+1}(\bx^\prime)|
&\leq& \dsum_{\{a,b\} \subset \mN} \, |x_{a,b} \, I_{a,b}(\bx) - x_{a,b}^\prime \, I_{a,b}(\bx^\prime)|
&\leq& 1 + D_N,
\ee
using $x_{a,b} \, I_{a,b}(\bx) \leq 1$ for all $\{a,b\} \subset \mN$ and all $\bx \in \mbX$.
\qed

}

\subsubsection{Auxiliary results}
\label{sup-sec:D_aux_proofs}

We prove Lemmas \ref{lem:bound_P}--\ref{lem:pi_star_bound},
which establish auxiliary results utilized in the proof of Lemma \ref{prop:D_bound}.

\begin{lemma}
\label{lem:bound_P}
Consider Models 1, 2, and 3 with $\nat \in \mbR^p$ and $\alpha \in [0, \, 1/2)$.
Then there exist functions $L_{k} : \mbR^p \mapsto (0, 1)$ and $U_{k} : \mbR^p \mapsto (0, 1)$
($k = 0, 1$) such that,
for all $\{i,j\} \subset \mN$ and $\bx_{-\{i,j\}} \in \{0, 1\}^{M-1}$,
\beno
0 &<& L_k(\nat)
&\leq&  \mbP_{\nat}(X_{i,j} = k \mid \bX_{-\{i,j\}} = \bx_{-\{i,j\}})
&\leq& U_k(\nat) &<& 1.
\ee
The functions $L_k(\nat)$ and $U_k(\nat)$ ($k = 0,1$) are given by
\beno
L_1(\nat)
&\coloneqq& \begin{cases}
\dfrac{1}{1 + \exp((3 + D_N) \, \norm{\nat}_{\infty})} & \hspace{.12in} \mbox{if } \mN_i \cap \mN_j \neq \emptyset \s\s\\
\dfrac{N^{-\alpha}}{1 + \exp(3 \, \norm{\nat}_{\infty})} & \hspace{.12in} \mbox{if } \mN_i \cap \mN_j = \emptyset
\end{cases}

\s\s\\

U_1(\nat)
&\coloneqq& \begin{cases}
\dfrac{1}{1 + \exp(-(3 + D_N) \, \norm{\nat}_{\infty})} & \hspace{.01in} \mbox{if } \mN_i \cap \mN_j \neq \emptyset \s\s\\
\dfrac{1}{1 + \exp(-3 \, \norm{\nat}_{\infty}) \, N^{\alpha}} & \hspace{.01in} \mbox{if } \mN_i \cap \mN_j = \emptyset
\end{cases}

\s\s\\

L_0(\nat)
&\coloneqq& \begin{cases}
\dfrac{1}{1 + \exp((3+D_N) \, \norm{\nat}_{\infty})} & \hspace{.12in} \mbox{if } \mN_i \cap \mN_j \neq \emptyset \s\s\\
\dfrac{1}{1 + \exp(3 \, \norm{\nat}_{\infty}) \, N^{-\alpha}} & \hspace{.12in} \mbox{if } \mN_i \cap \mN_j = \emptyset
\end{cases}

\s\s\\

U_0(\nat)
&\coloneqq& \begin{cases}
\dfrac{1}{1 + \exp(-(3+D_N) \, \norm{\nat}_{\infty})} & \hspace{0in} \mbox{if } \mN_i \cap \mN_j \neq \emptyset \s\s\\
\dfrac{1}{1 + \exp(-3 \, \norm{\nat}_{\infty}) \, N^{-\alpha}} & \hspace{0in} \mbox{if } \mN_i \cap \mN_j = \emptyset.
\end{cases},
\ee
\end{lemma}

\llproof \ref{lem:bound_P}.
Consider any pair of nodes $\{i,j\} \subset \mN$ and any $\bx_{-\{i,j\}} \in \{0, 1\}^{\binom{N}{2}-1}$.
We can express the full conditional probability
\beno
\mbP_{\nat}(X_{i,j} = x_{i,j} \mid \bX_{-\{i,j\}} = \bx_{-\{i,j\}})
\ee
two different ways depending whether $\mN_i$ and $\mN_j$ are disjoint.

First,
if $\mN_i \cap \mN_j = \emptyset$,
\beno
\label{eq:eq_sparsity}
&& \mbP_{\nat}(X_{i,j} = x_{i,j} \mid \bX_{-\{i,j\}} = \bx_{-\{i,j\}}) \s\s\\
&=& \dfrac{\exp(\langle \nat, \, s(\bx_{-\{i,j\}},\, x_{i,j}) \rangle) \, N^{-\alpha \, x_{i,j}}}{\exp(\langle \nat, \, s(\bx_{-\{i,j\}},\, x_{i,j}=0) \rangle) + \exp(\langle \nat, \, s(\bx_{-\{i,j\}},\, x_{i,j}=1) \rangle)\; N^{-\alpha}}\s\s
\\
\= \dfrac{1}{g(0;\, \bx_{-\{i,j\}},\, x_{i,j},\, \nat)\; N^{\alpha \, x_{i,j}} + g(1;\, \bx_{-\{i,j\}},\, x_{i,j},\, \nat)\; N^{-\alpha \, (1-x_{i,j})}},
\ee
defining,
for $y \in \{0, 1\}$,
\beno
g(y;\, \bx_{-\{i,j\}},\, x_{i,j},\, \nat)
&\coloneqq& \exp(\langle \nat, \, s(\bx_{-\{i,j\}},\, y) - s(\bx_{-\{i,j\}},\, x_{i,j})\rangle).
\ee
Note $g(x_{i,j}; \, \bx_{-\{i,j\}},\, x_{i,j},\, \nat) = 1$
for all $\bx_{-\{i,j\}} \in \{0,1\}^{\binom{N}{2}-1}$ and all $\nat \in \mbR^p$.

Second,
if $\mN_i \cap \mN_j \neq \emptyset$,
\beno
\label{eq:eq_no_sparsity}
&& \mbP_{\nat}(X_{i,j} = x_{i,j} \mid \bX_{-\{i,j\}} = \bx_{-\{i,j\}}) \s\s\\
&=& \dfrac{\exp(\langle \nat, \, s(\bx_{-\{i,j\}},\, x_{i,j}) \rangle)}{\exp(\langle \nat, \, s(\bx_{-\{i,j\}},\, x_{i,j}=0) \rangle) + \exp(\langle \nat, \, s(\bx_{-\{i,j\}},\, x_{i,j}=1) \rangle)}\s\s
\\
\= \dfrac{1}{1 + g(1-x_{i,j};\, \bx_{-\{i,j\}},\, x_{i,j},\, \nat)}.
\ee

Next,
observe that
\beno
&& \max\limits_{\bx_{-\{i,j\}}\, \in\, \{0, 1\}^{M-1}}
\, \left| s_l(\bx_{-\{i,j\}},\, x_{i,j}=0) - s_l(\bx_{-\{i,j\}},\, x_{i,j}=1)  \right|\s\s
\\
&\leq& \begin{cases}
0 & \mbox{ if } l \in \{1, \dots, N\} \setminus \{i,j\}\s
\\
1 & \mbox{ if } l \in \{i,j\}\s
\\
1+D_N & \mbox{ if } l = N + 1 \mbox{ and } \mN_i \cap \mN_j \neq \emptyset\s
\\ 0 & \mbox{ if } l = N + 1 \mbox{ and } \mN_i \cap \mN_j = \emptyset
\end{cases}.
\ee
The bound on $s_{N+1}$ follows from Lemma \ref{lem:brokerage_ham_bound},
whereas the conditions follow from Proposition \ref{prop:cond_ind}:
$s_{N+1}$ is a function of only dependent edge variables 
and $X_{i,j}$ is independent of all other edges in the graph when $\mN_i \cap \mN_j = \emptyset$. 
The bound on $s_{l}(\bx)$ ($l \in \mN$) follows because $s_l(\bx) = \sum_{h \in \mN \setminus \{l\}} x_{l,h}$ is a
function of $x_{i,j}$ if and only if $l \in \{i,j\}$.
As a result,
the triangle inequality shows that 
\beno
&& \left|\left\langle \nat, \, s(\bx_{-\{i,j\}},\, x_{i,j}=1)\right\rangle - \left\langle\nat, \, s(\bx_{-\{i,j\}},\, x_{i,j}=0)\right\rangle \right| \s\s\\
&\leq&
\begin{cases}
\left(3 + D_N  \right) \, \norm{\nat}_{\infty} & \mbox{if } \mN_i \cap \mN_j \neq \emptyset \s\\
3 \, \norm{\nat}_{\infty} & \mbox{if } \mN_i \cap \mN_j = \emptyset
\end{cases},
\ee
implying,
for $\{i,j\} \subseteq \mN$ with $\mN_i \cap \mN_j \neq \emptyset$,
\beno
\exp(-(3+D_N) \, \norm{\nat}_{\infty})
\,\leq\, g(1-x_{i,j};\, \bx_{-\{i,j\}},\, x_{i,j},\, \nat)
\,\leq\, \exp((3+D_N) \, \norm{\nat}_{\infty}),
\ee
and for $\{i,j\} \subseteq \mN$ with $\mN_i \cap \mN_j = \emptyset$,
\beno
\exp(- 3 \, \norm{\nat}_{\infty})
\,\leq\, g(1-x_{i,j};\, \bx_{-\{i,j\}},\, x_{i,j},\, \nat)
\,\leq\, \exp(3 \, \norm{\nat}_{\infty}).
\ee
As a result,
for all $k \in \{0, 1\}$,
\beno
0 &<& L_k(\nat)
&\leq& \mbP(X_{i,j} = k \,|\, \bX_{-\{i,j\}} =  \bx_{-\{i,j\}})
&\leq& U_k(\nat)
&<& 1,
\ee
where
\beno
L_1(\nat)
&\coloneqq& \begin{cases}
\dfrac{1}{1 + \exp((3 + D_N) \, \norm{\nat}_{\infty})} & \mbox{if } \mN_i \cap \mN_j \neq \emptyset \s\s\\
\dfrac{N^{-\alpha}}{1 + \exp(3 \, \norm{\nat}_{\infty})} & \mbox{if } \mN_i \cap \mN_j = \emptyset
\end{cases}
\ee
and
\beno
U_1(\nat)
&\coloneqq& \begin{cases}
\dfrac{1}{1 + \exp(-(3 + D_N) \, \norm{\nat}_{\infty})} & \mbox{if } \mN_i \cap \mN_j \neq \emptyset \s\s\\
\dfrac{1}{1 + \exp(-3 \, \norm{\nat}_{\infty}) \, N^{\alpha}} & \mbox{if } \mN_i \cap \mN_j = \emptyset
\end{cases}.
\ee
We obtain $L_0(\nat)$ and $U_0(\nat)$
by noting
\beno
\mbP(X_{i,j} = 0 \,|\, \bX_{-\{i,j\}} =  \bx_{-\{i,j\}})
\= 1 - \mbP(X_{i,j} = 1 \,|\, \bX_{-\{i,j\}} =  \bx_{-\{i,j\}}),
\ee
implying
\beno
1 - U_1(\nat)
&\leq& \mbP(X_{i,j} = 0 \,|\, \bX_{-\{i,j\}} =  \bx_{-\{i,j\}})
&\leq& 1 - L_1(\nat),
\ee
which allows us to obtain
\beno
L_0(\nat)
&\coloneqq& \begin{cases}
\dfrac{1}{1 + \exp((3+D_N) \, \norm{\nat}_{\infty})} & \mbox{if } \mN_i \cap \mN_j \neq \emptyset \s\s\\
\dfrac{1}{1 + \exp( 3\, \norm{\nat}_{\infty}) \, N^{-\alpha}} & \mbox{if } \mN_i \cap \mN_j = \emptyset
\end{cases}
\ee
and
\beno
U_0(\nat)
&\coloneqq& \begin{cases}
\dfrac{1}{1 + \exp(-(3+D_N) \, \norm{\nat}_{\infty})} & \mbox{if } \mN_i \cap \mN_j \neq \emptyset \s\s\\
\dfrac{1}{1 + \exp(-3 \, \norm{\nat}_{\infty}) \, N^{-\alpha}} & \mbox{if } \mN_i \cap \mN_j = \emptyset
\end{cases}.
\ee
\qed

\s

\begin{lemma}
\label{lem:brokerage_ham_bound}
Consider
\beno
s_{N+1}(\bx)
\= \dsum_{i < j}^N\, x_{i,j} \, I_{i,j}(\bx),
\ee
where
\beno
I_{i,j}(\bx)
\=
\begin{cases}
0 & \mbox{if } \mN_i \cap \mN_j = \emptyset\s
\\
\one\left( \, \dsum_{h \in \mN_i \cap \mN_j} \, x_{i,h} \, x_{j,h} \,\geq\, 1\right)
& \mbox{if } \mN_i \cap \mN_j \neq \emptyset.
\end{cases}
\ee
Then,
for all $\{i,j\} \subset \mN$,
\beno
\max\limits_{(\bx,\,\bx^\prime)\, \in\, \mbX \times \mbX:\; x_{v,w} = x_{v,w}^\prime,\; \{v,\, w\} \neq \{i,\, j\}}
|s_{N+1}(\bx) - s_{N+1}(\bx^\prime)|
&\leq& 1+D_N,
\ee
where $D_N \coloneqq \max_{\{i,j\} \subset \mN} \, |\mathfrak{N}_{i,j}|$.
\end{lemma}

\s\s

\llproof \ref{lem:brokerage_ham_bound}.
Consider any pair of nodes $\{i,j\} \subset \mN$.
The number of $x_{a,b}\, I_{a,b}(\bx)$ ($\{a,b\} \neq \{i,j\}$)
which are a function of $x_{i,j}$
includes 
\bi
\item $\{a,b\} = \{i,b\}$ ($b \in \mN \setminus \{i,j\}$) satisfying $j \in \mN_i \cap \mN_b \subseteq \mN_i$, and \s 
\item $\{a,b\} = \{j,b\}$ ($b \in \mN \setminus \{i,j\}$) satisfying $i \in \mN_j \cap \mN_b \subseteq \mN_j$.
\ei
As a result, 
the number of summands $x_{a,b}\, I_{a,b}(\bx)$ ($\{a,b\} \neq \{i,j\}$)
which can change value due to changing the value of $x_{i,j}$ 
is bounded above by 
$|(\{i\} \times \mN_i) \cup (\{j\} \times \mN_j)|$.
Proposition \ref{prop:cond_ind} establishes that, 
for any  $k \in  \{1, \ldots, K\}$ and $\{i,j\} \subset \mA_k$,
\beno
\big\{\{a,b\} \,:\, (a,b) \in \{i\} \times \mN_i \mbox{ or } (a,b) \in \{j\} \times \mN_j \big\}
&\subseteq& \mathfrak{N}_{i,j}.
\ee
Hence  
$|(\{i\} \times \mN_i) \cup (\{j\} \times \mN_j)| \leq D_N$,
noting $D_N \coloneqq \max_{\{i,j\} \subset \mN} \, |\mathfrak{N}_{i,j}|$. 
As a result,
the number of total summands $x_{a,b} \, I_{a,b}(\bx)$ which are a function of $x_{i,j}$ is bounded above by $1 + D_N$,
now counting the case when $\{a,b\} = \{i,j\}$.
Consider any $(\bx,\, \bx^\prime) \in \mbX \times \mbX$ such that $x_{v,w} = x_{v,w}^\prime$ for all $\{v,w\} \neq \{i,j\}$.
Then,
by the triangle inequality,
\beno
|s_{N+1}(\bx) - s_{N+1}(\bx^\prime)|
&\leq& \dsum_{\{a,b\} \subset \mN} \, |x_{a,b} \, I_{a,b}(\bx) - x_{a,b}^\prime \, I_{a,b}(\bx^\prime)|
&\leq& 1 + D_N,
\ee
using $x_{a,b} \, I_{a,b}(\bx) \in \{0, 1\}$ for all $\{a,b\} \subset \mN$ and all $\bx \in \mbX$.
\qed

\s\s

\begin{lemma}
\label{lem:brokerage_ham_bound_2}
Consider the function 
\beno
I(\bx)
\= \dsum_{i < j}^N\, I_{i,j}(\bx),
\ee
where
\beno
I_{i,j}(\bx)
\=
\begin{cases}
0 & \mbox{if } \mN_i \cap \mN_j = \emptyset\s
\\
\one\left( \, \dsum_{h \in \mN_i \cap \mN_j} \, x_{i,h} \, x_{j,h} \,\geq\, 1\right)
& \mbox{if } \mN_i \cap \mN_j \neq \emptyset.
\end{cases}
\ee
Then,
for all $\{i,j\} \subset \mN$,
\beno
\max\limits_{(\bx,\,\bx^\prime)\, \in\, \mbX \times \mbX:\; x_{v,w} = x_{v,w}^\prime,\; \{v,\, w\} \neq \{i,\, j\}}
|I(\bx) - I(\bx^\prime)|
&\leq& \begin{cases}
D_N & \mN_i \,\cap\, \mN_j \neq \emptyset \s\\
0 & \mN_i \,\cap\, \mN_j = \emptyset 
\end{cases},
\ee
where $D_N \coloneqq \max_{\{i,j\} \subset \mN} \, |\mathfrak{N}_{i,j}|$.
\end{lemma}

\s\s

\llproof \ref{lem:brokerage_ham_bound_2}.
Consider any pair of nodes $\{i,j\} \subset \mN$.
Note that $I_{i,j}(\bx)$ is not a function of $x_{i,j}$.
Additionally, 
note that each $I_{i,j}(\bx)$ is only a function of edge variables for which 
$\mN_i \cap \mN_j \neq \emptyset$,
and is constant in all edge variables $x_{a,b}$ with $\mN_a \cap \mN_b = \emptyset$.   
The number of $I_{a,b}(\bx)$ ($\{a,b\} \subset \mN$) 
which are a function of $x_{i,j}$
includes 
\bi
\item $\{a,b\} = \{i,b\}$ ($b \in \mN \setminus \{i,j\}$) satisfying $j \in \mN_i \cap \mN_b \subseteq \mN_i$, and \s 
\item $\{a,b\} = \{j,b\}$ ($b \in \mN \setminus \{i,j\}$) satisfying $i \in \mN_j \cap \mN_b \subseteq \mN_j$.
\ei
As a result, 
the number of summands $I_{a,b}(\bx)$ ($\{a,b\} \subset \mN$) 
which can change value due to changing the value of $x_{i,j}$ 
is bounded above by 
\beno
|(\{i\} \times \mN_i) \cup (\{j\} \times \mN_j)|.
\ee
Proposition \ref{prop:cond_ind} establishes that, 
for any  $k \in  \{1, \ldots, K\}$ and $\{i,j\} \subset \mA_k$,
\beno
\big\{\{a,b\} \,:\, (a,b) \in \{i\} \times \mN_i \mbox{ or } (a,b) \in \{j\} \times \mN_j \big\}
&\subseteq& \mathfrak{N}_{i,j}.
\ee
Hence  
$|(\{i\} \times \mN_i) \cup (\{j\} \times \mN_j)| \leq D_N$,
noting $D_N \coloneqq \max_{\{i,j\} \subset \mN} \, |\mathfrak{N}_{i,j}|$. 
As a result,
the number of total summands $I_{a,b}(\bx)$ which are a function of $x_{i,j}$ is bounded above by $D_N$,
Consider any $(\bx,\, \bx^\prime) \in \mbX \times \mbX$ such that $x_{v,w} = x_{v,w}^\prime$ for all $\{v,w\} \neq \{i,j\}$.
Then,
by the triangle inequality,
\beno
|I(\bx) - I(\bx^\prime)|
&\leq& \dsum_{\{a,b\} \subset \mN} \, |I_{a,b}(\bx) -I_{a,b}(\bx^\prime)|
&\leq& D_N,
\ee
using $I_{a,b}(\bx) \in \{0, 1\}$ for all $\{a,b\} \subset \mN$ and all $\bx \in \mbX$.
\qed

\s

\begin{lemma}
\label{lem:prove_coupling}
Choose any $i \in \{1, \dots, M\}$ and any $\bx_{1:i-1} \in \{0, 1\}^{i-1}$.
Then the coupling of the conditional distributions
\beno
\mbP(\,\cdot \mid \bX_{1:i-1} = \bx_{1:i-1},\; X_i = 0)
&\mbox{and}&
\mbP(\,\cdot \mid \bX_{1:i-1} = \bx_{1:i-1},\; X_i = 1)
\ee
of $\bX_{(i+1):M}$ 
constructed in Lemma \ref{prop:D_bound} is a valid coupling.
\end{lemma}

\s

\llproof \ref{lem:prove_coupling}.
Denote the coupling distribution generated by the algorithm in Lemma \ref{prop:D_bound} by $\mbQ_{i,\bx_{1:i-1}}^{}$ and let $v_1, \dots, v_{M-i}$ be the vertices added to the set $\mathfrak{V}$ at iteration $1, \dots, M-i$ of the algorithm.
To reduce the notational burden, 
define
\beno
&& q(\bx^\star_{a:b},\; \bx^{\star\star}_{a:b} \mid \bx^\star_{c:d},\; \bx^{\star\star}_{c:d})\s\s
\\
&\coloneqq& \mbQ_{i,\bx_{1:i-1}}^{}(\bX^\star_{a:b} = \bx^\star_{a:b},\; \bX^{\star\star}_{a:b} = \bx^{\star\star}_{a:b} \mid \bX^\star_{c:d} = \bx^\star_{c:d},\; \bX^{\star\star}_{c:d} = \bx^{\star\star}_{c:d}),
\ee
where $a, b, c, d \in \{1, \dots, M\}$ are distinct and $\{a, \dots, b\}\, \cap\, \{c, \dots, d\} = \emptyset$.
By construction,
\beno
q(\bx^\star_{i+1:M},\; \bx^{\star\star}_{i+1:M})
\= q(x_{v_1}^\star,\; x_{v_1}^{\star\star})\; \dprod_{l=2}^{M-i} q(x_{v_l}^\star,\; x_{v_l}^{\star\star} \mid \bx_{v_1, \ldots, v_{l-1}}^{\star},\; \bx_{v_1, \ldots, v_{l-1}}^{\star\star}).
\ee
Observe that
\beno
\dsum_{x_{v_{M-i}}^\star\, \in\, \{0, 1\}} q(x_{v_{M-i}}^\star,\; x_{v_{M-i}}^{\star\star} \mid \bx_{v_1, \ldots, v_{M-i-1}}^{\star},\; \bx_{v_1, \ldots, v_{M-i-1}}^{\star\star})\s\s
\\
=\; \mbP(X_{v_{M-i}} = x_{v_{M-i}}^{\star\star} \mid \bX_{1:i-1} = \bx_{1:i-1},\; X_i = 1,\; \bX_{v_1,\ldots,v_{M-i-1}} = \bx_{v_1, \ldots, v_{M-i-1}}^{\star\star})
\ee
and
\beno
\dsum_{x_{v_{M-i}}^{\star\star}\, \in\, \{0, 1\}} q(x_{v_{M-i}}^\star,\; x_{v_{M-i}}^{\star\star} \mid \bx_{v_1, \ldots, v_{M-i-1}}^{\star},\; \bx_{v_1, \ldots, v_{M-i-1}}^{\star\star})\s\s
\\
=\; \mbP(X_{v_{M-i}} = x_{v_{M-i}}^{\star} \mid \bX_{1:i-1} = \bx_{1:i-1},\; X_i = 0,\; \bX_{v_1,\ldots,v_{M-i-1}} = \bx_{v_1, \ldots, v_{M-i-1}}^{\star}),
\ee
owing to the fact that $(X_{v_{M-i}}^\star,\, X_{v_{M-i}}^{\star\star})$ is distributed according to the optimal coupling of the conditional distributions
\beno
\mbP(X_{v_{M-i}} = \,\cdot \mid \bX_{1:i-1} = \bx_{1:i-1},\; X_i = 0,\; \bX_{v_1,\ldots,v_{M-i-1}} = \bx_{v_1, \ldots, v_{M-i-1}}^{\star})
\ee
and
\beno
\mbP(X_{v_{M-i}} = \,\cdot \mid \bX_{1:i-1} = \bx_{1:i-1},\; X_i = 1,\; \bX_{v_1,\ldots,v_{M-i-1}} = \bx_{v_1, \ldots, v_{M-i-1}}^{\star\star}).
\ee
\hide{
Consider any $l \in \{1, \ldots, M-i-1\}$ and suppose that 
Then
\beno
\dsum_{x_{v_l}^\star\, \in\, \{0, 1\}} q(x_{v_l}^\star,\; x_{v_l}^{\star\star} \mid \bx_{v_1, \ldots, v_{l-1}}^{\star},\; \bx_{v_1, \ldots, v_{l-1}}^{\star\star})\s\s
\\
=\; \mbP(X_{v_l} = x_{v_l}^{\star\star} \mid \bX_{1:i-1} = \bx_{1:i-1},\; X_i = 1,\; \bX_{v_1,\ldots,v_{l-1}} = \bx_{v_1, \ldots, v_{l-1}}^{\star\star})
\ee
and
\beno
\dsum_{x_{v_l}^{\star\star}\, \in\, \{0, 1\}} q(x_{v_l}^\star,\; x_{v_l}^{\star\star} \mid \bx_{v_1, \ldots, v_{l-1}}^{\star},\; \bx_{v_1, \ldots, v_{l-1}}^{\star\star})\s\s
\\
=\; \mbP(X_{v_l} = x_{v_l}^\star \mid \bX_{1:i-1} = \bx_{1:i-1},\; X_i = 0,\; \bX_{v_1,\ldots,v_{l-1}} = \bx_{v_1, \ldots, v_{l-1}}^{\star}).
\ee
}
We can repeat the same argument to show that 
\beno
&& \dsum_{x_{v_1}^\star \in \{0, 1\}} \cdots \dsum_{x_{v_{M-i}}^\star \in \{0, 1\}} q(\bx_{v_1, \ldots, v_{M-i}}^\star,\; \bx^{\star\star}_{i+1:M})\s
\\
\hide{
\= \dsum_{x_{v_1}^\star \in \{0, 1\}}\, q(x_{v_1}^\star, x_{v_1}^{\star\star})
\\
&\times& \dprod_{l=2}^{M-i} \dsum_{x_{v_l}^\star \in \{0, 1\}} q(x_{v_l}^\star,\; x_{v_l}^{\star\star} \mid \bx_{v_1, \ldots, v_{l-1}}^{\star},\; \bx^{\star\star}_{v_1, \ldots, v_{l-1}})\s
\\
\= \mbP(X_{v_1} = x_{v_1}^{\star\star} \mid \bX_{1:i-1} = \bx_{1:i-1}, X_i = 1)
\\
&\times& \dprod_{l=2}^{M-i} \mbP(X_{v_{l}} = x_{v_l}^{\star\star} \mid \bX_{1:i-1} = \bx_{1:i-1}, X_i = 1, \bX_{v_1,\ldots,v_{l-1}} = \bX_{v_1, \ldots, v_{l-1}}^{\star\star})
\\
}
\= \mbP(\bX_{i+1:M} = \bx_{1+i:M}^{\star\star} \mid \bX_{1:i-1} = \bx_{1:i-1},\; X_i = 1)
\ee
and
\beno
&& \dsum_{x_{v_1}^{\star\star}\, \in\, \{0, 1\}}\, \cdots\, \dsum_{x_{v_{M-i}^{\star\star}}\, \in\, \{0, 1\}} q(\bx_{i+1:M}^\star,\; \bx_{v_1, \ldots, v_{M-i}}^{\star\star})\s
\\
\= \mbP(\bX_{i+1:M} = \bx_{1+i:M}^{\star} \mid \bX_{1:i-1} = \bx_{1:i-1},\; X_i = 0),
\ee
so the coupling is indeed a valid coupling of the conditional distributions
\beno
\mbP(\bX_{i+1:M} = \,\cdot \mid \bX_{1:i-1} = \bx_{1:i-1},\; X_i = 0)
\ee
and
\beno
\mbP(\bX_{i+1:M} = \,\cdot \mid \bX_{1:i-1} = \bx_{1:i-1},\; X_i = 1). 
\ee
\qed

\begin{lemma}
\label{lem:pi_star_bound}
Consider Models 1, 2, and 3,
any $v \in \{1, \dots, M\}$, 
and any $(\bx_{-v},\, \bx_{-v}^\prime) \in \{0, 1\}^{M-1} \times \{0, 1\}^{M-1}$. 
Define
\beno
\pi_{v,\, \bx_{-v},\, \bx_{-v}^\prime}
&\coloneqq& \norm{\mbP(\,\cdot \mid \bX_{-v} = \bx_{-v}) -  \mbP(\,\cdot \mid \bX_{-v} = \bx_{-v}^\prime)}_{\tv}
\ee
and
\beno
\pi^\star
&\coloneqq& \max\limits_{1\, \leq\, v\, \leq\, M}\;\, \max\limits_{(\bx_{-v},\, \bx_{-v}^\prime)\; \in\; \{0, 1\}^{M-1} \times \{0, 1\}^{M-1}} \; \pi_{v,\, \bx_{-v},\, \bx_{-v}^\prime},
\ee
and define 
$D_N \coloneqq \max_{\{i,j\} \subset \mN}\, |\mathfrak{N}_{i,j}|$.
Then  
\beno
\pi^\star
\lte
\begin{cases}
0 & \mbox{ under Model 1} \s
\\
\dfrac{1}{1 + \exp(-(3+ D_N) \, \norm{\truth}_{\infty})} 
& \mbox{ under Models 2 and 3}.
\end{cases}
\ee
\end{lemma}

\s

\llproof \ref{lem:pi_star_bound}.
Under Model 1, 
edge variables $X_v$ are independent, 
which implies that $\pi_{v,\, \bx_{-v},\, \bx_{-v}^\prime} = 0$ for all $v \in \{1, \ldots, M\}$ and all $(\bx_{-v},\, \bx_{-v}^\prime) \in \{0, 1\}^{M - 1} \times \{0, 1\}^{M - 1}$,
which in turn implies that $\pi^\star = 0$.  
To bound $\pi^\star$ under Models 2 and 3,
we distinguish two cases: 
\ben
\item[(a)] If edge variable $X_v$ corresponds to a pair of nodes with non-intersecting node neighborhoods,
i.e.,
a pair $\{i,j\} \subset \mN$ with $\mN_i \cap \mN_j = \emptyset$, 
then $X_v$ is independent of all other edge variables by Proposition \ref{prop:cond_ind}.
As a result,
$\pi_{v,\, \bx_{-v},\, \bx_{-v}^\prime} = 0$ 
for all $(\bx_{-v},\, \bx_{-v}^\prime) \in \{0, 1\}^{M - 1} \times \{0, 1\}^{M - 1}$. \s
\item[(b)] If edge variable $X_v$ corresponds to a pair of nodes with intersecting node neighborhoods,
i.e.,
a pair $\{i,j\} \subset \mN$ with $\mN_i \cap \mN_j \neq \emptyset$,
then $X_v$ is not independent of all other edges, 
implying $\pi_{v,\, \bx_{-v},\, \bx_{-v}^\prime} > 0$ for some or all $(\bx_{-v},\, \bx_{-v}^\prime) \in \{0, 1\}^{M - 1} \times \{0, 1\}^{M - 1}$.\s
\een
We focus henceforth on case (b).
Consider any $v \in \{1, \ldots, M\}$ such that 
$\pi_{v,\, \bx_{-v},\, \bx_{-v}^\prime} > 0$ 
for some 
$(\bx_{-v},\, \bx_{-v}^\prime) \in \{0, 1\}^{M - 1} \times \{0, 1\}^{M - 1}$ and define
\beno
a_0 \= \mbP(X_v = 0 \mid \bX_{-v} = \bx_{-v})
&\mbox{and}&
a_1 \= \mbP(X_v = 1 \mid \bX_{-v} = \bx_{-v}) \s\s \\
b_0 \= \mbP(X_v = 0 \mid \bX_{-v} = \bx_{-v}^\prime)
&\mbox{and}&
b_1 \= \mbP(X_v = 1 \mid \bX_{-v} = \bx_{-v}^\prime).
\ee
Then 
\beno
\pi_{v,\, \bx_{-v},\, \bx_{-v}^\prime} 
= \dfrac{1}{2} \left( |(1 - a_1) - (1 - b_1) | + |a_1 - b_1| \right)
= |a_1 - b_1|
\leq \max\{a_1,\, b_1\}.
\ee
By symmetry,
\beno
\pi_{v,\, \bx_{-v},\, \bx_{-v}^\prime}
\lte \max\{a_0,\, b_0\},
\ee
which implies that
\beno
\pi_{v,\, \bx_{-v},\, \bx_{-v}^\prime}
\lte \min\left\{\max\{a_0,\, b_0\},\; \max\{a_1,\, b_1\}\right\}. 
\ee 
Lemma \ref{lem:bound_P} shows that,
under Models 2 and 3,
\beno
\mbP(X_{v} = 0 \mid \bX_{-v} = \bx_{-v})
\lte \dfrac{1}{1 + \exp(- (3+D_N) \, \norm{\truth}_{\infty})} 
\ee
and
\beno
\mbP(X_{v} = 1 \mid \bX_{-v} = \bx_{-v})
\lte \dfrac{1}{1 + \exp(-(3 + D_N) \, \norm{\truth}_{\infty})}.
\ee
We therefore conclude that,
under Models 2 and 3,
\beno
\pi^\star
\lte \min\left\{\max\{a_0,\, b_0\},\; \max\{a_1,\, b_1\}\right\}  \s \\ 
&\leq& \dfrac{1}{1 + \exp(-(3 + D_N) \, \norm{\truth}_{\infty})}. 
\ee
\qed

\s

\section{Proof of Proposition 1}
\label{proof.degrees}

We prove Proposition \ref{p:expected.degrees} stated in Section \ref{model3} of the manuscript.

\pproof \ref{p:expected.degrees}.
The expected degree of any node $i \in \mN$ under Model 3 with $\alpha \in (0,\, 1]$ is given by  
\beno
\mbE_{\nat}\left(\dsum_{j \neq i }^N \, X_{i,j}\right)
\= \dsum_{j\, \in\, \mathfrak{A}_{i,1}}\, \mbE_{\nat}\, X_{i,j} + \dsum_{j\, \in\, \mathfrak{A}_{i,2}}\, \mbE_{\nat}\, X_{i,j}\s\s
\\
\lte |\mathfrak{A}_{i,1}|\, \max\limits_{j\, \in\, \mathfrak{A}_{i,1}}\, \mbE_{\nat}\, X_{i,j} + |\mathfrak{A}_{i,2}|\, \max\limits_{j\, \in\, \mathfrak{A}_{i,2}}\, \mbE_{\nat}\, X_{i,j},
\ee 
where
\bi
\item $\mathfrak{A}_{i,1} = \{j \in \mN \setminus \{i\}:\, \mN_i\, \cap\, \mN_j \neq  \emptyset\}$;\s
\item $\mathfrak{A}_{i,2} = \{j \in  \mN \setminus \{i\}:\, \mN_i\, \cap\, \mN_j = \emptyset\}$.\s
\ei
We bound the expectations of edges $\mbE_{\nat}\, X_{i,j}$ by using the bound  
\beno
\mbE_{\nat}\, X_{i,j}
\= \mbP_{\nat}(X_{i,j} = 1)\s\s
\\
&\leq& \max\limits_{\bx_{-\{i,\, j\}}\, \in\, \{0, 1\}^{\binom{N}{2} - 1}}\; \mbP_{\nat}(X_{i,j} = 1 \mid \bX_{-\{i,\, j\}} = \bx_{-\{i,\, j\}}). 
\ee
For any $j \in \mathfrak{A}_{i,1}$,\,
$\mbP_{\nat}(X_{i,j} = 1 \,|\, \bX_{-\{i,\, j\}} = \bx_{-\{i,\, j\}}) \leq 1 \leq \exp(3 \, \norm{\nat}_\infty)$ for all $\bx_{-\{i,\, j\}} \in \{0, 1\}^{\binom{N}{2} - 1}$.
In addition,
for any $j \in \mathfrak{A}_{i,2}$,\,
Lemma \ref{lem:bound_P} in Appendix \ref{D} shows that
\beno
\mbP_{\nat}(X_{i,j} = 1 \mid \bX_{-\{i,\, j\}} = \bx_{-\{i,\, j\}})
\leq \dfrac{1}{1 + \exp(- 3\, \norm{\nat}_\infty) \, N^{\alpha}}
< \dfrac{\exp(3 \, \norm{\nat}_\infty)}{N^{\alpha}}, 
\ee
for all $\bx_{-\{i,\, j\}}\, \in\, \{0, 1\}^{\binom{N}{2} - 1}$.
Hence, 
\beno
\mbE_{\nat}\left(\dsum_{j \neq i }^N \, X_{i,j}\right)
\lte \exp(3 \, \norm{\nat}_{\infty}) \, (|\mathfrak{A}_{i,1}| + |\mathfrak{A}_{i,2}| \, N^{-\alpha}). 
\ee

\s

{\bf Bounding $|\mathfrak{A}_{i,1}|$.}
To bound $|\mathfrak{A}_{i,1}|$,
we distinguish two cases:
\bi
\item $\mN_i \cap \mN_j \neq \emptyset$ and $j \in  \mN_i$,
which implies that there exists $k \in \{1, \dots, K\}$ such that $\{i,\, j\} \subset \mA_k$,
in which case $j\, \in\, \mN_i$.\s
\item $\mN_i \cap \mN_j \neq \emptyset$ and $j \not\in \mN_i$,
in which case there exists $h\, \in\, \mN_i \cap \mN_j$,
which implies that $h \in \mN_i$ and $h \in \mN_j$,
which further implies $j \in \mN_h$.  
\ei
The number of nodes $j\, \in\, \mN$ satisfying the first case is bounded above by 
$|\mN_i| \leq \max_{1 \leq r \leq N} |\mN_r| \leq \left(\max_{1 \leq r \leq N} |\mN_r|\right)^2$,
since $|\mN_r| \in \{0, 1, \ldots, N-1\}$, 
and the number of $j\, \in\, \mN$ satisfying the second case is bounded above by 
\beno
\left| \bigcup\limits_{h \in \mN_i} \, \mN_h \right| 
&\leq&
|\mN_i| \; \max\limits_{h \in \mN_i} |\mN_h|
&\leq& \left(\max\limits_{1 \leq h \leq N} |\mN_h|\right)^2.
\ee
In conclusion,
\beno
|\mathfrak{A}_{i,1}| 
\lte 2\; \left(\max\limits_{1 \leq h \leq N} |\mN_h|\right)^2.
\ee

\s

{\bf Bounding $|\mathfrak{A}_{i,2}|$.}
For each node $i \in \mN$,
there are at most $N - 1 < N$ other nodes $j \in \mN \setminus \mN_i$,
hence $|\mathfrak{A}_{i,2}| \leq N \leq 2 \, N$.

\s 

{\bf Conclusion.}
By collecting terms,
for all nodes $i \in \mN$,
\beno
\mbE_{\nat}\left(\dsum_{j \neq i }^N \, X_{i,j}\right)
&\leq& 2\; \exp(3\, \norm{\nat}_{\infty}) \; \left(\left(\max\limits_{1 \leq h \leq N} |\mN_h|\right)^2 + N^{1-\alpha}\right).
\ee
\qed

\vspace{.1cm}

\section{Simulation results}
\label{sec:simulationsstudy}

\begin{figure}[t]
\centering
\includegraphics[width = .975 \linewidth, keepaspectratio]{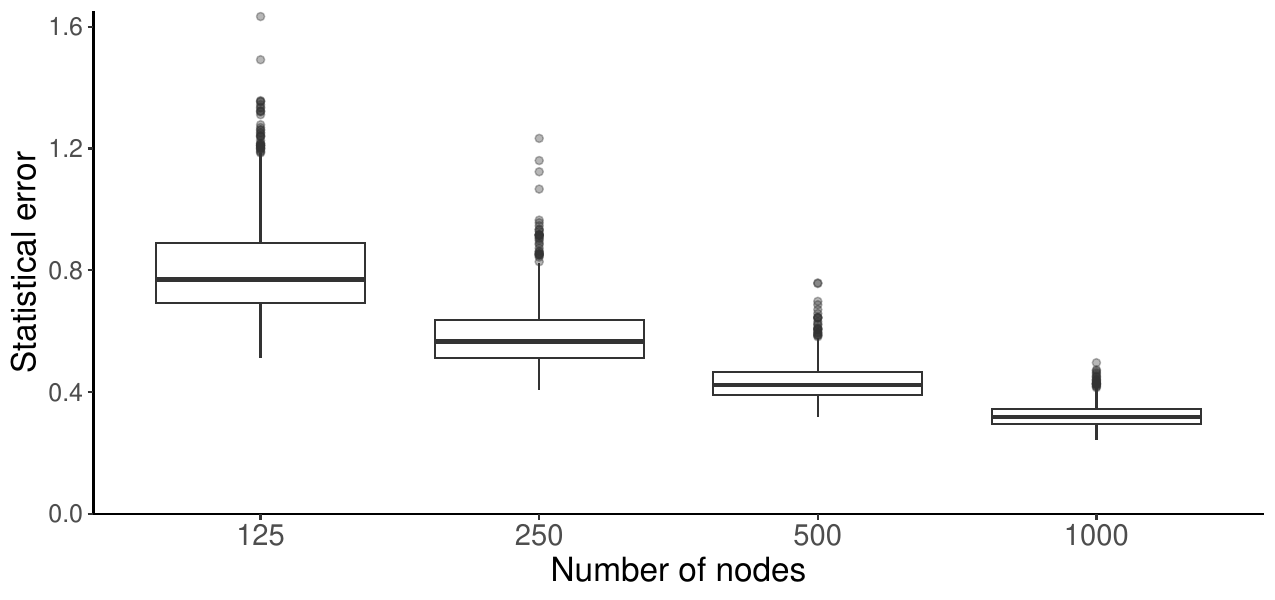}
\caption{\label{fig:sim_res_1}
The statistical error $\norm{\mple - \truth}_{\infty}$ of maximum pseudo-likelihood estimator $\mple$ as an estimator of $\truth \in \mR^{N+1}$ plotted against the number of nodes $N$.
}
\end{figure}


We study the performance of maximum pseudo-likelihood estimators by considering populations with $N =$ 125, 250, 500, and 1,000 nodes.
We focus on maximum pseudo-likelihood estimators,
because computing maximum likelihood and Monte Carlo maximum likelihood estimators is too time-consuming when $N$ is large (e.g., when $N = 500$ and $N =$ 1,000).
For each value of $N$,
we generate 1,000 populations with overlapping subpopulations as follows:
\bi
\item The number of subpopulations $K$ is $N \,/\, 25$.
\item Each node $i \in \mN$ belongs to $1 + Y_i$ subpopulations,
where\break 
$Y_i\; \iid\; \mbox{Binomial}(K - 1,\; 1\, /\, K)$ ($i = 1, \dots, N$).
\item For node $i \in \{1, \dots, N\}$,
the $1 + Y_i$ subpopulation memberships are sampled from the Multinomial$(p_1^{(i)}, \dots, p_K^{(i)})$ distribution with
\beno
p_k^{(i)} &=& 
\begin{cases}
\dfrac1K & \mbox{ if } i = 1\s
\\
\dfrac{1}{K-1}\; \left(1 - \dfrac{N_k^{(i-1)}}{N_1^{(i-1)} + \ldots + N_K^{(i-1)}}\right) & \mbox{ if } i \in \{2, \dots, N\},
\end{cases}
\ee
where $N_k^{(i-1)}$ is the number of nodes in $\{1, \dots, i-1\}$ that belong to subpopulation $\mA_k$ ($k = 1, \dots, K$)
at the current time. 
\ei
We consider Model 2 with degree parameters $\theta_1^\star, \ldots, \theta_{N}^\star$ drawn from\break 
Uniform($-1.25,\, -.75$) and brokerage parameter $\theta_{N+1}^\star = .25$.
For each population size $N \in \{125, 250, 500, 1000\}$,
we generate a graph from Model 2 and compute the maximum pseudo-likelihood estimator from the generated graph.
For each value of $N$,
the gradient ascent algorithm used to compute the maximum pseudo-likelihood estimator converged for at least 95\% of the simulated data sets,
and the following simulation results are based on the simulated data sets for which the gradient ascent algorithm converged.

\begin{figure}[t]
\centering
\includegraphics[width = \linewidth, keepaspectratio]{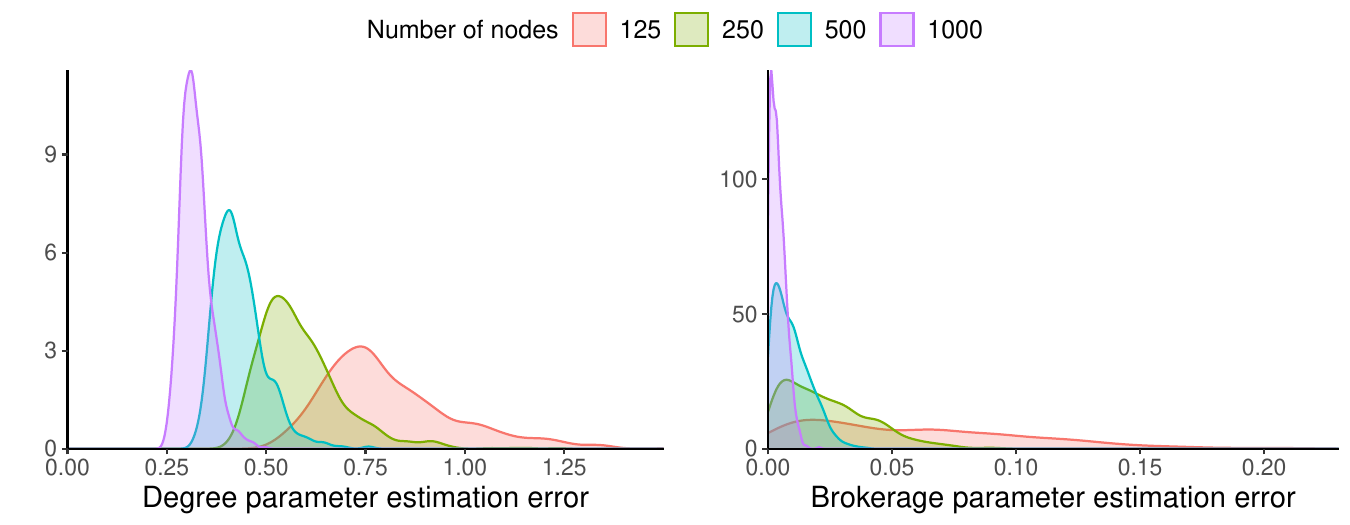}
\caption{\label{fig:den_plot}
The maximum deviation $\max_{1 \leq i \leq N} |\widetilde\theta_i - \theta_i^\star|$ of the maximum pseudo-likelihood estimators $\widetilde\theta_i$ from the data-generating degree parameters $\theta_i^\star$ ($i = 1, \dots, N$) (left) and the deviation $|\widetilde\theta_{N+1} - \theta_{N+1}^\star|$ of the maximum pseudo-likelihood estimator $\widetilde\theta_{N+1}$ from the data-generating brokerage parameter $\theta_{N+1}^\star$ (right). 
}
\end{figure}

Figure \ref{fig:sim_res_1} demonstrates that the statistical error $\norm{\mple - \truth}_{\infty}$ of $\mple$ as an estimator of the data-generating parameter vector $\truth \in \mR^{N+1}$ decreases as the number of nodes $N$ increases.
Figure \ref{fig:den_plot} decomposes the statistical error of $\mple$ into the statistical error of the degree parameter estimators $\widetilde\theta_1, \dots, \widetilde\theta_N$ and the statistical error of the brokerage parameter estimator $\widetilde\theta_{N+1}$.
Figure \ref{fig:den_plot} reveals that the brokerage parameter is estimated with greater accuracy than the degree parameters,
which makes sense as the degree parameters are greater in absolute value than the brokerage parameter
and there are $N$ estimated degree parameters $\widetilde\theta_1, \ldots, \widetilde\theta_N$, 
compared with the single estimated brokerage parameter $\widetilde\theta_{N+1}$. 

\label{last.page}

\end{document}